\pdfoutput=1
\RequirePackage[T1]{fontenc}
\documentclass[12pt,twoside]{amsart}
\calclayout

\usepackage[utf8]{inputenc}
\usepackage{amsmath}
\usepackage{amssymb}
\usepackage{mathtools}
\allowdisplaybreaks[1]
\usepackage{amsthm}
\usepackage[svgnames]{xcolor}
\usepackage[colorlinks,citecolor=DarkGreen,linkcolor=FireBrick,urlcolor=FireBrick,linktocpage,breaklinks=true]{hyperref}
\usepackage{cite}
\usepackage[pdftex]{graphicx}
\usepackage{tikz}
\usepackage{newtxtext}
\usepackage{courier}

\usepackage{enumitem}

\usetikzlibrary{arrows.meta,calc,decorations.markings}

\theoremstyle{plain}
    
    \newtheorem{theorem}{Theorem}[section]
    \newtheorem{proposition}[theorem]{Proposition}
    \newtheorem{corollary}[theorem]{Corollary}
    \newtheorem{lemma}[theorem]{Lemma}

\theoremstyle{definition}
    \newtheorem{definition}[theorem]{Definition}

\theoremstyle{remark}
	\newtheorem{remark}[theorem]{Remark}%

\makeatletter
\newcommand{\cornermark@sym}{$\lrcorner$}
\let\oldrem\remark
\let\oldendrem\endremark
\def\remark{\pushQED{\qed}\let\qedpic@sym\cornermark@sym\global\qedpic@lower=\z@\oldrem}
\def\endremark{\popQED\oldendrem}
\let\olddef\definition
\let\oldenddef\enddefinition
\def\definition{\pushQED{\qed}\let\qedpic@sym\cornermark@sym\global\qedpic@lower=\z@\olddef}
\def\enddefinition{\popQED\oldenddef}
\makeatother

\newtheoremstyle{restated}
    {\topsep}{\topsep} 
    {\itshape}         
    {}                 
    {\bfseries}        
    {.}                
    { }                
    {\thmname{#1} \ref{#3} {\normalfont(Restated)}}
    
\theoremstyle{restated}
    \newtheorem{restate-theorem}{Theorem}
    \newtheorem{restate-proposition}{Proposition}
    \newtheorem{restate-corollary}{Corollary}
    
\numberwithin{equation}{section}
\numberwithin{table}{section}

\newcommand{\FF}{\mathbb{F}}

\newcommand{\id}{\mathit{id}}

\newcommand{\bR}{\mathbb{R}}

\DeclareMathOperator{\op}{\mathrm{op}}

\DeclareMathOperator{\Ob}{Ob}

\DeclareMathOperator{\Hom}{Hom}

\DeclareMathOperator{\End}{End}

\DeclareMathOperator{\Ima}{Im}
\DeclareMathOperator{\Ker}{Ker}

\DeclareMathOperator{\Fix}{Fix}

\DeclareMathOperator{\Mat}{Mat}

\DeclareMathOperator{\ev}{ev}
\DeclareMathOperator{\coev}{coev}
\DeclareMathOperator{\Tr}{Tr}
\DeclareMathOperator{\Inv}{Inv}

\newcommand{\sP}{\mathsf{P}}
\newcommand{\sT}{\mathsf{T}}

\newcommand{\bZ}{\mathbb{Z}}

\newcommand{\bC}{\mathbb{C}}

\newcommand{\cC}{\mathcal{C}}
\newcommand{\cD}{\mathcal{D}}

\newcommand{\cE}{\mathcal{E}}

\newcommand{\cT}{\mathcal{T}}

\newcommand{\Irr}{\mathrm{Irr}}

\providecommand{\one}{\mathbf{1}}

\providecommand{\cH}{\mathfrak{H}}
\providecommand{\Mat}{\operatorname{Mat}}

\newcommand{\fV}{\mathfrak{V}}
\newcommand{\fR}{\mathfrak{R}}
\newcommand{\fX}{\mathfrak{M}}
\newcommand{\fF}{\mathfrak{F}}

\definecolor{mshade}{RGB}{228,233,243}
\tikzset{
  mb/.style={baseline={([yshift=-0.7ex]current bounding box.center)},
    execute at end picture={\qedpicrecord},
    x=8.5mm, y=7mm, line cap=round, line join=round,
    every node/.append style={font=\scriptsize}},
  mbar/.style={line width=0.8pt},
  mbartop/.style={mbar, red, -{Stealth[length=2mm,width=1.8mm]}},
  mbarbot/.style={mbar, red, {Stealth[length=2mm,width=1.8mm]}-},
  mstr/.style={line width=0.6pt},
  mover/.style={line width=0.6pt, preaction={draw, white, line width=4pt,
    line cap=butt, shorten <=1.4pt, shorten >=1.4pt}},
  mcpn/.style={draw, line width=0.6pt, fill=white, rectangle,
    inner sep=1.4pt, minimum height=4.2mm, minimum width=5.2mm,
    font=\scriptsize},
  mcpnc/.style={mcpn, fill=black!15},
}
\newcommand{\mbars}[3]{\draw[mbarbot] (#1,0) -- (#2,0); \draw[mbartop] (#1,#3) -- (#2,#3);}
\newcommand{\blab}[2]{\node[anchor=base,inner sep=0pt,yshift=-8.8pt] at (#1,0) {$#2$};}
\newcommand{\tlab}[3]{\node[anchor=base,inner sep=0pt,yshift=4pt] at (#1,#2) {$#3$};}
\makeatletter
\newdimen\qedpic@depth
\newdimen\qedpic@tmp
\everydisplay\expandafter{\the\everydisplay\global\qedpic@depth=\z@}
\newcommand{\qedpicrecord}{%
  \qedpic@tmp=\pgf@picmaxy \advance\qedpic@tmp-\pgf@picminy
  \divide\qedpic@tmp 2 \advance\qedpic@tmp-0.7ex
  \ifdim\qedpic@tmp>\qedpic@depth \global\qedpic@depth=\qedpic@tmp\fi}
\newdimen\qedpic@lower
\let\qedpic@sym\qedsymbol
\renewcommand{\qedsymbol}{\raisebox{-\qedpic@lower}[\height][\z@]{\qedpic@sym}}
\AddToHook{env/proof/begin}{\global\qedpic@lower=\z@}
\newcommand{\qedpic}[1][\qedpic@depth]{\global\qedpic@lower=#1\relax\qedhere}
\makeatother
\newcommand{\sstr}[5][mstr]{\draw[#1] (#2,#3) .. controls (#2,{(#3+#5)/2}) and (#4,{(#3+#5)/2}) .. (#4,#5);}
\newcommand{\mone}[3]{%
\begin{tikzpicture}[mb]
\mbars{-0.55}{0.55}{2}
\draw[mstr] (0,0) -- (0,2);
\node[mcpn] at (0,1) {$#1$};
\blab{0}{#2}\tlab{0}{2}{#3}
\end{tikzpicture}}

\newcommand{\tubemark}[2]{\draw[mstr] ({#1-0.11},{#2-0.11}) -- ({#1+0.11},{#2-0.02}); \draw[mstr] ({#1-0.11},{#2+0.02}) -- ({#1+0.11},{#2+0.11});}
\newcommand{\tubeelt}[4]{%
\begin{tikzpicture}[mb]
\draw[mstr] (-0.25,0) -- (-0.25,1.5);
\draw[mstr] (0.25,1.5) -- (0.25,3);
\draw[mstr,rounded corners=4pt] (1.3,0.75) -- (0.25,0.75) -- (0.25,1.5);
\draw[mstr,rounded corners=4pt] (-0.25,1.5) -- (-0.25,2.25) -- (-1.3,2.25);
\tubemark{-0.25}{0.22}\tubemark{0.25}{2.78}
\node[mcpn,minimum width=8mm] at (0,1.5) {$#1$};
\node[anchor=west,inner sep=1.5pt] at (1.3,0.75) {$#3$};
\node[anchor=east,inner sep=1.5pt] at (-1.3,2.25) {$#4$};
\blab{-0.25}{#2}\tlab{0.25}{3}{#2}
\end{tikzpicture}}

\newcommand{\Ubox}[1]{\langle #1\rangle_{\sP}}
\newcommand{\ubox}[4]{\draw[densely dashed,rounded corners=2pt,line width=0.4pt] (#1,#2) rectangle (#3,#4);}
\tikzset{medgeline/.style={draw=none}}
\newcommand{\medge}[3]{\path[medgeline] (#1,#2) -- (#1,#3);}
\newcommand{\legend}[3]{\fill (#1,#2) circle (1.1pt);
  \node[anchor=east,inner sep=1.5pt] at (#1,#2) {$#3$};}
\newcommand{\picZg}{%
\begin{tikzpicture}[mb]
\mbars{-1}{0.6}{2}\medge{-1}{0}{2}
\draw[mstr] (0,0) -- (0,1);
\draw[mstr] (0.15,1) -- (0.15,2);
\draw[mstr] (-0.15,1) .. controls (-0.15,1.45) and (-0.6,1.45) .. (-1,1.45);
\node[mcpn] at (0,1) {$g$};
\legend{-1}{1.45}{E}
\blab{0}{Z}\tlab{0.15}{2}{\sP Z}
\end{tikzpicture}}
\newcommand{\picXg}{%
\begin{tikzpicture}[mb]
\mbars{-1}{0.6}{2}\medge{-1}{0}{2}
\draw[mstr] (0,0) -- (0,1);
\draw[mstr] (0.15,1) -- (0.15,2);
\draw[mstr] (-0.15,1) .. controls (-0.15,1.45) and (-0.6,1.45) .. (-1,1.45);
\node[mcpn] at (0,1) {$g$};
\legend{-1}{1.45}{E}
\blab{0}{X}\tlab{0.15}{2}{\sP X}
\end{tikzpicture}}

\newcommand{\mob}{\mathrm{M}}
\newcommand{\bdy}{\mathrm{L}}
\newcommand{\Cas}{\Xi}
\newcommand{\obstruction}{O}
\newcommand{\obs}{o}

\begin{document}
\title{Möbius Fusion Algebra}

\author{Ippo Orii}
\email{ippo.orii@ipmu.jp}
\author{Yuji Tachikawa}
\email{yuji.tachikawa@ipmu.jp}
\address{Kavli Institute for the Physics and Mathematics of the Universe \textsc{(wpi)},
  The University of Tokyo,
  5-1-5 Kashiwanoha, Kashiwa, Chiba, 277-8583,
  Japan
}

\begin{abstract}
  Given a pseudo-unitary modular tensor category $\cC$ 
  with an anti-monoidal involution $\sP$,
  we define a generally non-commutative algebra from this data,
  which we call the M\"obius fusion algebra $\fX$ of $(\cC,\sP)$,
  in terms of ribbon graphs 
  drawn on the M\"obius band, where the orientation-reversal
  of the M\"obius band acts on the graphs by $\sP$.
  We show that this algebra is semisimple,
  and give a formula for the dimensions of its simple modules,
  in terms of the modular $S$-matrix of $\cC$ and 
  a certain set of quantities extracted from the action of $\sP$.
  This formula is expected to play the role of 
  a non-orientable analogue of the Verlinde formula,
  in the context of Reshetikhin--Turaev-type topological
  quantum field theories on non-orientable manifolds.
\end{abstract}

\date{November 2026}

\maketitle
\setcounter{tocdepth}{2}
\enlargethispage{\baselineskip}
\tableofcontents

\newpage

\section{Introduction}
\label{sec:overview}
\addtocontents{toc}{\protect\setcounter{tocdepth}{1}}

\subsection{Motivation from physics}

We start from a physics motivation;
those who do not trust physics may skip this subsection
and go directly to the next one,
from which everything is mathematically well-defined.
Let $\cC$ be a modular tensor category.  Its fusion coefficients
$N_{AB}^C=\dim\Hom(A\otimes B,C)$ are given in terms of the modular
$S$-matrix by the Verlinde formula \cite{Verlinde:1988sn,Moore:1988uz,Moore:1988qv}
\begin{equation}
N_{AB}^C=\sum_{X\in\Irr(\cC)}
\frac{S_{A,X}\,S_{B,X}\,\overline{S_{C,X}}}{S_{\one,X}} .
\label{eq:intro-verlinde}
\end{equation}
In the context of $(2+1)$-dimensional topological quantum
field theory (TQFT) of Reshetikhin and Turaev
\cite{ReshetikhinTuraev1991,Turaev1994}, a mathematical construction of
the Chern--Simons theories of \cite{Witten:1988hf},
this equation rewrites the dimension of the state space
on $S^2$ with three points marked by $A$, $B$ and $C^\vee$
in terms of the modular data of $\cC$ alone.
If we are presented with the right hand side alone
without any prior knowledge,
it would be very difficult to understand that 
it would always turn out to be a non-negative integer.

More recently, TQFTs have also played important roles in condensed matter physics,
where $\cC$ describes the anyons of a $(2+1)$-dimensional topological phase of matter \cite{Kitaev_2006}.
In this context, topological phases may also have additional symmetries. Internal symmetries which do not act on the spacetime 
are well understood mathematically, in terms of
$G$-crossed braided extensions \cite{Etingof:2009yvg}, homotopy quantum
field theories \cite{Turaev:2014yha} and twisted Verlinde formulas
\cite{Deshpande2018},
and there are corresponding physics interpretations,
 see e.g.~\cite{Barkeshli_2019}.

Symmetry which reverses the orientation of the spacetime
is also extensively studied in physics
\cite{Barkeshli:2016mew,Wang_2017,TachikawaYonekura2017ptep},
but its mathematical formulation is much less developed.
It is expected that such an orientation-reversing
symmetry acts on $\cC$ by an \emph{anti-monoidal involution} $\sP$,
an autoequivalence  which reverses the order
of tensor products and inverts the braiding,
such that there is a natural isomorphism $\beta\colon\id_\cC\Rightarrow\sP^2$.
From this data of $(\cC,\sP)$, together with a further coherence
condition on $\beta$,
it is believed in the physics community that 
we can construct a topological phase, or a TQFT,
defined not only on orientable but also on
non-orientable spacetimes \cite{Barkeshli:2016mew}.  
Such a TQFT has
not been constructed from $(\cC,\sP)$ in general,
although partial results exist. 
For example, state sum descriptions have been
proposed at the physics level of rigor when $\cC$ is a Drinfeld center
\cite{Barkeshli:2016mew,Bhardwaj2017}, and unoriented twisted
Dijkgraaf--Witten theories are constructed mathematically in
\cite{Young2020}.  

One interesting observation by physicists \cite{Barkeshli:2016mew} 
is that, assuming that such a TQFT defined on non-orientable manifolds exists, 
its space of states on the real
projective plane $\mathbb{RP}^2$ with one point marked by $X\in \cC$
should have dimension
\begin{equation}\label{eq:intro-MX}
M_X=\sum_{Y\in\Irr(\cC),\ \sP Y\cong Y^\vee}S_{X,Y}\,\eta_Y,
\end{equation}
via a cutting-and-pasting argument of the putative TQFT.
Here, $\eta_Y$ is a root of unity determined by how 
$\sP$ acts on the one-dimensional space $\Hom(\one,Y\otimes \sP Y)$.
Furthermore, when $\beta$ satisfies a coherence condition,
which is usually assumed in the physics literature,
$\eta_Y$ reduces to  signs satisfying $\eta_{Y^\vee}=\eta_Y$.

The analogy with the ordinary Verlinde formula \eqref{eq:intro-verlinde}
is obvious:
both formulas, \eqref{eq:intro-verlinde} and \eqref{eq:intro-MX}, 
rewrite the dimensions of the respective state spaces
in terms of the modular $S$-matrix of $\cC$,
with an additional piece of information coming from $\sP$.

For this reason we call \eqref{eq:intro-MX} the
\emph{non-orientable Verlinde formula}, as $\mathbb{RP}^2$ 
with a marked point
is the basic building block of non-orientable surfaces,
just as $S^2$ with three marked points is
the basic building block of orientable surfaces.
Analogously with the right hand side of \eqref{eq:intro-verlinde},
the right hand side of \eqref{eq:intro-MX} is
not at all manifestly a non-negative integer.
That said, if the putative non-orientable TQFT existed,
the left hand side should clearly be a non-negative integer,
as the dimension of a vector space.
Indeed, 
in many explicit computations performed by physicists,
the right hand side has always turned out to be an integer.

This raises an obvious question.
Let us \emph{define} the quantities $M_X$ via \eqref{eq:intro-MX},
given the data of $(\cC,\sP)$.
Then physics gives us a conjecture that 
this $M_X$ is always a non-negative integer.\footnote{
Some partial results toward this direction exist in the physics literature.
In \cite{QiJianWang2019},
the system was folded along the reflection axis into a bilayer theory
$\cC\boxtimes\overline{\cC}$, in which the reflection becomes the exchange of the
two layers and the axis becomes a gapped boundary,
and it was shown that $M_X$ appears as a condensation multiplicity
\cite[(42)]{QiJianWang2019},  giving a physical argument for its integrality.
As another more mathematical approach, for pointed $\cC$, the
integrality of $M_X$ was proved in \cite{Orii_2025}.  }

In this paper we prove that $M_X$ thus defined by \eqref{eq:intro-MX} 
from $(\cC,\sP)$ is always a non-negative integer,
but without fully constructing the putative non-orientable TQFT.
Instead, 
from $(\cC,\sP)$, we construct a finite-dimensional semisimple algebra $\fX$,
which we call the \emph{M\"obius fusion algebra}, and show that the $M_X$ are the
dimensions of its irreducible modules.  
Here, the vector space $\fX$ should be thought of as the space of states on the Klein bottle,
just as the fusion algebra $\fF$ of $\cC$ is the space of states on the torus.

\subsection{M\"obius fusion algebra and its properties}

\newcommand{\skbcyl}[1][1]{%
  \draw[black!45,line width=0.4pt] (-#1,0) -- (-#1,3) (#1,0) -- (#1,3);
  \draw[black!45,line width=0.4pt] (0,3) ellipse ({#1} and 0.3);
  \draw[black!45,line width=0.4pt] (-#1,0) arc (180:360:{#1} and 0.3);
  \draw[black!45,line width=0.4pt,densely dashed] (#1,0) arc (0:180:{#1} and 0.3);
  \draw[mbarbot] (-#1,0) -- (#1,0);
  \draw[mbartop] (-#1,3) -- (#1,3);}
\newcommand{\stcyl}[1][1]{%
  \draw[black!45,line width=0.4pt] (-#1,0) -- (-#1,3) (#1,0) -- (#1,3);
  \draw[black!45,line width=0.4pt] (0,3) ellipse ({#1} and 0.3);
  \draw[black!45,line width=0.4pt] (-#1,0) arc (180:360:{#1} and 0.3);
  \draw[black!45,line width=0.4pt,densely dashed] (#1,0) arc (0:180:{#1} and 0.3);
  \draw[mbartop] (-#1,0) -- (#1,0);
  \draw[mbartop] (-#1,3) -- (#1,3);}
\newcommand{\skbbot}[2]{\node[anchor=base,inner sep=0pt] at (#1,-0.78) {$#2$};}
\newcommand{\skbtop}[2]{\node[anchor=base,inner sep=0pt] at (#1,3.45) {$#2$};}

Let $\cC$ be a pseudo-unitary modular tensor category,
and $\Irr(\cC)$ be a set of representatives of the isomorphism classes of simple objects in $\cC$.  The fusion
algebra $\fF$ of $\cC$ is the algebra of ribbon graphs in the solid
torus; we construct its analogue for the solid Klein bottle.

Let us start by recalling the relation of the solid torus to the fusion algebra $\fF$ of $\cC$.
The solid torus is obtained from the cylinder $D^2\times[0,1]$ by gluing
its top disk to its bottom disk by the identity map.  
We consider ribbon graphs in it in the
sense of \cite{Reshetikhin:1990pr,Turaev1994}.  A strand colored by $X$
along the core closes up with a coupon colored by some $f\in\End(X)$;
we write this element as \[
[X,f]=\;\begin{tikzpicture}[mb]
\stcyl
\draw[mstr] (0,0) -- (0,3);
\node[mcpn] at (0,1.5) {$f$};
\skbbot{0}{X}\skbtop{0}{X}
\end{tikzpicture}
\;=\;\begin{tikzpicture}[mb]
\draw[mbartop] (-1,0) -- (1,0); \draw[mbartop] (-1,3) -- (1,3);
\draw[mstr] (0,0) -- (0,3);
\node[mcpn] at (0,1.5) {$f$};
\blab{0}{X}\tlab{0}{3}{X}
\end{tikzpicture}.
\]
Here, in the middle picture, we drew the surrounding solid torus,
but in the right picture we projected it to the central plane,
to reduce clutter. 
In the introduction we use both types of pictures,
but in the main text we will only use the projected pictures.

A coupon colored by $a\in\End(X)$, slid once around the core, comes back
unchanged, so that \[
[X,f\circ a]=\;
\begin{tikzpicture}[mb]
\draw[mbartop] (-0.75,0) -- (0.75,0); \draw[mbartop] (-0.75,3) -- (0.75,3);
\draw[mstr] (0,0) -- (0,3);
\node[mcpn] at (0,0.9) {$a$};
\node[mcpn] at (0,2.1) {$f$};
\blab{0}{X}\tlab{0}{3}{X}
\end{tikzpicture}
\;=\;
\begin{tikzpicture}[mb]
\draw[mbartop] (-0.75,0) -- (0.75,0); \draw[mbartop] (-0.75,3) -- (0.75,3);
\draw[mstr] (0,0) -- (0,3);
\node[mcpn] at (0,0.9) {$f$};
\node[mcpn] at (0,2.1) {$a$};
\blab{0}{X}\tlab{0}{3}{X}
\end{tikzpicture}
\;=[X,a\circ f].
\]
The product places two such
strands side by side, so that \[
[X,f]\cdot[Y,g]=\;
\begin{tikzpicture}[mb]
\stcyl
\draw[mstr] (-0.5,-0.2) -- (-0.5,2.8);
\draw[mstr] (0.5,0.2) -- (0.5,3.2);
\node[mcpn] at (-0.5,0.55) {$f$};
\node[mcpn] at (0.5,0.75) {$g$};
\skbbot{-0.5}{X}\skbbot{0.5}{Y}\skbtop{-0.5}{X}\skbtop{0.5}{Y}
\end{tikzpicture}
\;=
[X\otimes Y,f\otimes g].
\]
More precisely, these elements take values in \[
\fF\cong\bigl(\bigoplus_{X\in\Ob(\cC)}\End(X)\bigr)/\fR_0,
\] where 
we write an element $f\in\End(X)$ of the summand labeled by $X$ as the
pair $(X,f)$,
the subspace $\fR_0$ is spanned by
\begin{equation}
(X\oplus Y,f\oplus g)-(X,f)-(Y,g),\qquad \text{and}\qquad  (X,a\circ f)-(X,f\circ a)
\label{eq:intro-relations}
\end{equation}
for $a\in\End(X)$,
and the class of $(X,f)$ in the quotient is $[X,f]$.
It is not difficult to check that $\fF$ thus defined
is isomorphic to the fusion algebra of $\cC$ defined as $K_0(\cC)\otimes_\bZ\bC$,
and the class $[X]$ for $X\in\Ob(\cC)$ in $\fF$ in our notation 
is then $[X]=[X,\id_X]$.

Let us next study 
the solid Klein bottle, which is obtained instead by gluing the two ends of $D^2\times [0,1]$ with a reflection of $D^2$.  
A strand passing through this
gluing is affected by the reflection. Let us denote this effect on the objects of $\cC$ by $\sP$.
As the reflection reverses the order of the
points along the cut, the natural isomorphisms
\[ U_{X,Y}\colon\sP Y\otimes\sP X\to\sP(X\otimes Y)\] 
reverse the order of the tensor product.
Together with 
$U_0\colon\one\to\sP\one$, this makes $(\sP,U)$ a monoidal functor
$\cC\to\cC^{\otimes\op}$, where $\cC^{\otimes\op}$ is $\cC$ with the opposite tensor product $X\otimes^{\op}Y:=Y\otimes X$. 
Similarly, as the orientation-reversal turns a crossing into its
mirror image, we should have 
\[\sP(c_{X,Y})\circ U_{X,Y}=U_{Y,X}\circ c_{\sP X,\sP Y}^{-1},\] 
meaning that  this functor is braided for the braiding $c_{X,Y}^{-1}$ of
$\cC^{\otimes\op}$.
As the square of the reflection is the identity,
we should also have  a monoidal natural isomorphism 
\[\beta\colon\id_\cC\Rightarrow\sP^2.\]
We summarize these conditions by saying that $(\sP,\beta)$ is an \emph{anti-monoidal involution} of $\cC$.

With these preparations, let us consider strands placed in a solid Klein bottle.
A strand colored by $X$ along the core now closes up only with a coupon
colored by some $f\colon X\to\sP X$; this is
\[
[X,f]=\;\begin{tikzpicture}[mb]
\skbcyl
\draw[mstr] (0,0) -- (0,3);
\node[mcpn] at (0,1.5) {$f$};
\skbbot{0}{X}\skbtop{0}{\sP X}
\end{tikzpicture}
\;=\;\begin{tikzpicture}[mb]
\mbars{-1}{1}{3}
\draw[mstr] (0,0) -- (0,3);
\node[mcpn] at (0,1.5) {$f$};
\blab{0}{X}\tlab{0}{3}{\sP X}
\end{tikzpicture}.
\]
A coupon colored by $a\in\End(X)$, slid once around the core, comes back
as $\sP(a)$, so that
\[
[X,f\circ a]=\;
\begin{tikzpicture}[mb]
\mbars{-0.75}{0.75}{3}
\draw[mstr] (0,0) -- (0,3);
\node[mcpn] at (0,0.9) {$a$};
\node[mcpn] at (0,2.1) {$f$};
\blab{0}{X}\tlab{0}{3}{\sP X}
\end{tikzpicture}
\;=\;
\begin{tikzpicture}[mb]
\mbars{-0.75}{0.75}{3}
\draw[mstr] (0,0) -- (0,3);
\node[mcpn] at (0,0.9) {$f$};
\node[mcpn] at (0,2.1) {$\sP(a)$};
\blab{0}{X}\tlab{0}{3}{\sP X}
\end{tikzpicture}
\;=[X,\sP(a)\circ f].
\]
In the product the reflection reverses the order of the two strands,
which therefore cross:
\begin{equation}\label{eq:skb-product}
[X,f]\cdot[Y,g]=\;
\begin{tikzpicture}[mb]
\skbcyl
\draw[mstr] (0.5,0.2) -- (0.5,1.2);
\sstr{0.5}{1.2}{-0.5}{3.2}
\draw[mstr] (-0.5,-0.2) -- (-0.5,1.2);
\sstr[mover]{-0.5}{1.2}{0.5}{2.8}
\node[mcpn] at (-0.5,0.55) {$f$};
\node[mcpn] at (0.5,0.75) {$g$};
\skbbot{-0.5}{X}\skbbot{0.5}{Y}\skbtop{-0.5}{\sP Y}\skbtop{0.5}{\sP X}
\end{tikzpicture}
\;=\bigl[X\otimes Y,\ U_{X,Y}\circ c_{\sP X,\sP Y}\circ(f\otimes g)\bigr].
\end{equation}
Here, we put the strand $X$ in front of the strand $Y$,
as shown in the picture above.
This is the \emph{M\"obius fusion algebra}
\[
\fX:=\Bigl(\bigoplus_{X\in\Ob(\cC)}\Hom(X,\sP X)\Bigr)\Big/\fR ,
\]
where $\fR$ is spanned by the same vectors as in \eqref{eq:intro-relations} with $a\circ f$
replaced by $\sP(a)\circ f$, with the product \eqref{eq:skb-product}.
The classes
$\mob(X):=[X,\xi_X]$ for
$X\in\Fix(\sP):=\{X\in\Irr(\cC)\mid X\cong\sP X\}$, where
$\xi_X\colon X\to\sP X$ is an isomorphism with
$\sP(\xi_X)\circ\xi_X=\beta_X$, turn out to form a basis of $\fX$.
We work with the projections of these ribbon graphs to the M\"obius
band, which we call \emph{M\"obius diagrams}.
The equation \eqref{eq:skb-product} shows that the
product is in general not commutative: exchanging the factors changes
the crossing of the two strands.  This is unlike $\fF$, where the strands do not cross. 
Closing up the strands at the cut
gives a trace $\varepsilon\colon\fX\to\bC$, 
which makes $\fX$ a finite-dimensional $C^*$-algebra.

\begin{remark}
If a corresponding TQFT exists, $\fX$
should be its space of states on the Klein bottle, just as $\fF$ is the
space of states on the torus.
From this point on, physics interpretations, such as these, are enclosed in remarks.
\end{remark}

We next consider the \emph{boundary loop}
\[\bdy(Y):=[Y\otimes\sP Y,\,U_{Y,\sP Y}\circ(\beta_Y\otimes\id_{\sP Y})]
\;=\;
\begin{tikzpicture}[mb]
\mbars{-0.75}{1.65}{3}
\draw[mstr] (0,0) -- (0,3); \draw[mstr] (0.9,0) -- (0.9,3);
\node[mcpn] at (0,1.5) {$\beta_Y$};
\blab{0}{Y}\blab{0.9}{\sP Y}\tlab{0}{3}{\sP^2Y}\tlab{0.9}{3}{\sP Y}
\end{tikzpicture}\;,
\] which is
a strand colored by $Y$ running twice around the core, parallel to the
boundary of the M\"obius band. 
Unlike a strand running once, it closes up by itself: after one round its color
is $\sP Y$, and after two it is $\sP^2Y$, which $\beta_Y$ identifies with
$Y$.  A neighborhood of this loop is a solid torus in the solid Klein
bottle, so ribbon graphs in the solid torus also give elements of $\fX$; this
explains why $\bdy$ lies in the center of $\fX$ and extends to an algebra
homomorphism $\bdy\colon\fF\to\fX$. 
We put $p_X:=\bdy(e_X)$ for
the primitive idempotents $e_X$ of $\fF$, labeled by
$X\in\Irr(\cC)$.

Now, consider a linear map $R$ on $\Hom(\one,Y\otimes\sP Y)$
given by 
\[
\phi \mapsto  (\beta_Y^{-1}\otimes\id_{\sP Y})  U_{Y, \sP Y}^{-1} \sP( \phi) U_0.
\] 
When $Y$ is simple and $\sP Y\cong Y^\vee$, 
$\Hom(\one,Y\otimes\sP Y)$ is one-dimensional, and $R$ acts on it by a scalar,
which we denote by $\eta_Y$.
When $\beta$ is \emph{coherent}, that
is, $\beta_{\sP X}=\sP(\beta_X)$ for all $X$, 
we have
$\eta_Y=\eta_{Y^\vee}=\pm1$, 
and $\eta_Y$ can be considered as an analogue of the Frobenius--Schur
indicator. 
Finally we define the numbers $M_X$ by the formula \eqref{eq:intro-MX},
which we recall here for convenience:
\[M_X:=\sum_{Y\in\Irr(\cC),\ \sP Y\cong Y^\vee} S_{X,Y}\,\eta_Y.\]
Our first main result is the following theorem:

\begin{theorem}\label{thm:intro}
\begin{enumerate}[label=(\alph*)]
\item The center of $\fX$ is $\bdy(\fF)$.
\item For each $X\in\Irr(\cC)$ with $p_X\neq 0$,
the algebra $p_X\fX$ is simple, 
\[p_X\fX\cong\End_\bC(V_X),\]
where $V_X$ is an irreducible $\fX$-module; we put $V_X:=0$ if $p_X=0$.
Furthermore,  every irreducible $\fX$-module arises in this manner for a unique $X$.
\item We have $\dim V_X=M_X$, with $M_X$ given by the formula above.
In
  particular the numbers $M_X$ are non-negative
  integers, and $\fX\cong\bigoplus_{X\in\Irr(\cC)}\Mat_{M_X}(\bC)$.
\end{enumerate}
\end{theorem}

\begin{remark}
If the TQFT exists, $V_X$ should be its space
of states on $\mathbb{RP}^2$ with a marked point colored by $X$, whose
dimension was predicted to be $M_X$ in \cite{Barkeshli:2016mew} 
and $\fX\cong\bigoplus_X\End_\bC(V_X)$ is what cutting the Klein bottle into two
M\"obius bands would give.
\end{remark}
A rough outline of the proof is as follows.
As in any semisimple algebra,
the center of $\fX$ is the image of an
averaging map $\Cas\colon u\mapsto \sum_i b^i u b_i$
where $b^i$ and $b_i$ are dual bases of $\fX$ with respect to the trace $\varepsilon$.
Using the fact that $\mob(X)$ for $X\in\Fix(\sP)$ form a basis of $\fX$, 
the averaging map takes the form
$\sum_{X\in\Fix(\sP)}\mob(X)^\star u\mob(X)$,
where $\star$ is the involution of the $C^*$-algebra $\fX$.  A
computation with M\"obius diagrams shows that $\Cas(u)$ is a linear
combination of boundary loops $\bdy(Y)$, which gives (a).
Then (b) follows easily.  For (c), we use the same formula for
the unit $u=\mob(\one)$ which gives
\[
\sum_{X\in\Fix(\sP)}\mob(X)^\star \mob(X)
=\sum_{Y\in\Irr(\cC),\ \sP Y\cong Y^\vee}\eta_Y\,\bdy(Y).
\]
Multiplying by $p_X$ and applying $\varepsilon$ gives $(\dim V_X)^2$ on
the left and $M_X^2$ on the right, since $\varepsilon(p_X)=S_{\one,X}M_X$;
the positivity of $\varepsilon(p_X)$ then gives $\dim V_X=M_X$.  Thus the
integrality comes from counting the dimension of a block in two ways.

Our second main result concerns the twists $\theta_X$
which we explain here for coherent $\beta$ for brevity:

\begin{theorem}
\label{thm:constancy-intro}
Assume that $\beta$ is coherent.
For each $X\in\Irr(\cC)$ with $M_X\neq 0$, 
its $\theta_X$ is given by the formula
\[
\theta_X=\lambda\sigma
\] which is independent of $X$, where \[
\lambda=\frac1D\sum_{Y\in\Irr(\cC)}d_Y^2\,\theta_Y,\qquad
\sigma=\frac1D\sum_{Y\in\Irr(\cC),\ \sP Y\cong Y^\vee}\eta_Y\,d_Y\,\theta_Y,\qquad
D=(\sum_Yd_Y^2)^{1/2}.
\]
\end{theorem}

We will actually prove a more general result, Theorem~\ref{thm:theta},
which is applicable to the non-coherent case as well.
This result will be shown using the following intermediate result,
which is of independent interest.

For $E\in\Irr(\cC)$ let
$\cH_E$ be the quotient of
$\bigoplus_{Z\in\Ob(\cC)}\Hom(Z,E\otimes\sP Z)$ by the analogues
of the two relations spanning $\fR$ we gave above.
The class $[Z,g]$ of $g\colon Z\to E\otimes\sP Z$ is
the M\"obius diagram
\[\picZg\]
in which a strand colored by $E$ leaves the M\"obius band through
its boundary.
We have $\cH_\one=\fX$,
and $\cH_E$ is naturally an $\fX$-bimodule, by stacking M\"obius diagrams
from above and from below.
We then have the following structure theorem for $\cH_E$:
\begin{theorem}
\label{thm:bimodule-intro}
For each $E\in\Irr(\cC)$, we have an isomorphism of $\fX$-bimodules
\[\cH_E\cong\bigoplus_{X,X'\in\Irr(\cC)}(V_X\otimes V_{X'}^*)^{\oplus N_{EX'}^X},\]
where $N_{EX'}^X$ is the fusion coefficient of $\cC$.
\end{theorem}

Now, when $\beta$ is coherent, we can show that carrying the strand $E$ once around the band
acts as the identity on $\cH_E$, i.e.~\[
\begin{tikzpicture}[mb]
\mbars{-1.3}{2.0}{2.6}\medge{-1.3}{0}{2.6}
\draw[mstr] (-0.7,0) .. controls (-0.7,0.45) and (-1,0.6) .. (-1.3,0.6);
\draw[mstr] (0.3,0) -- (0.3,0.8);
\draw[mstr] (0.15,0.8) -- (0.15,2.6);
\draw[mstr] (0.45,0.8) .. controls (0.45,1.3) and (0.95,1.2) .. (0.95,1.7) -- (0.95,2.6);
\draw[mstr] (1.6,0) -- (1.6,2.6);
\node[mcpn] at (0.3,0.8) {$g$};
\node[mcpn] at (0.15,1.8) {$\beta_E$};
\legend{-1.3}{0.6}{E}
\blab{-0.7}{E}\blab{0.3}{Z}\blab{1.6}{\sP E}
\tlab{0.15}{2.6}{\sP^2E}\tlab{0.95}{2.6}{\sP Z}\tlab{1.6}{2.6}{\sP E}
\end{tikzpicture}
\;=\; \picZg\;.
\]

\begin{remark}
If the expected TQFT exists, $\cH_E$ should be its space of states on the Klein bottle with a marked point colored by $E$.
Then this is a consequence of the geometric fact that moving the marked point 
around the two-sided curve which cuts the Klein bottle into two M\"obius bands is isotopic to the identity.
\end{remark}

Now, under the isomorphism given in Theorem~\ref{thm:bimodule-intro}, 
this operation acts on the summand $(V_X\otimes V_{X'}^*)^{\oplus N_{EX'}^X}$ 
by $\theta_X\theta_{X'}^{-1}$.
Since for any two $X,X'$ we have some $E$ with $N_{EX'}^X\neq 0$, 
$\theta_X$ is constant on the $X$ with $M_X\neq0$.
Its value $\lambda\sigma$ in Theorem~\ref{thm:constancy-intro} then follows by a further
computation. 

\subsection{Organization of the paper}
The paper is organized as follows.  Section~\ref{sec:setting} fixes
conventions and the pictures of the action of $\sP$.
In Section~\ref{sec:mobius} we define $\fX$, show that it is a
$C^*$-algebra, and construct the boundary loops.
Our first and second main results are then
proved in Sections~\ref{sec:integrality} and \ref{sec:bimodules}, respectively.\footnote{
Two crucial propositions, Propositions~\ref{prop:two-steps} and~\ref{prop:W-factor} there, are proved by long
sequences of pictures; animated versions of these two proofs are available at
\url{https://ippo0127.github.io/movies/mobius_prop4_23.html} and \url{https://ippo0127.github.io/movies/mobius_prop5_6.html}.}
Finally in Section~\ref{sec:examples}, we study three classes of examples
to illustrate our results:
(i) pointed categories, where the non-zero
$M_X$ are powers of $2$, 
(ii) the categories $\cC=\cE\boxtimes\overline{\cE}$
with $\sP$ exchanging the factors, where $\fX$ of $\cC$ turns out to be the fusion algebra of $\cE$, and 
(iii) two non-pointed categories of rank $6$.
We also have Appendix~\ref{app:skeletal} where we explain how to find an anti-monoidal involution, 
or to show that none exists, from the $F$- and $R$-symbols of a skeletal model.

\begin{remark}
There will be two sequels to this paper.
First, in a paper in preparation
by the first author and Sakie Suzuki,
Reshetikhin--Turaev-type invariants of non-orientable three-manifolds 
will be constructed, building on the results of this paper.
Second, in a forthcoming companion paper intended for physicists,
various related topics will be discussed,
including 
physical interpretations, 
the action of the mapping class group of the Klein bottle, 
and the generalizations to fermionic cases (i.e.~super-modular categories). 
\end{remark}

\subsection*{Acknowledgments}
IO thanks Yasuyuki Kawahigashi, Yuji Terashima and Sakie Suzuki for helpful
discussions, and Sakie Suzuki for comments on the draft.
IO is supported by the Forefront Physics and Mathematics Program to Drive
Transformation (FoPM), a World-leading Innovative Graduate Study (WINGS)
Program, and the JSR Fellowship at the University of Tokyo.
YT is supported in part by WPI Initiative, MEXT, Japan at Kavli IPMU, the
University of Tokyo and by JSPS KAKENHI Grant-in-Aid (Kiban-C),
No.~24K06883.
AI assistants based on large language models were used in preparing this
paper, for drafting and checking proofs, drawing diagrams, numerical
computations and editing the text.  All results were checked by the
authors, who take full responsibility for the content.

\addtocontents{toc}{\protect\setcounter{tocdepth}{2}}

\section{Setting}
\label{sec:setting}

\subsection{Conventions}
\label{subsec:conventions}

In the following,  $\cC$ is a pseudo-unitary\footnote{%
Recall that a fusion category is \emph{pseudo-unitary} if its global
dimension equals its Frobenius--Perron dimension; it then has a unique
spherical structure for which the dimensions of all simple objects are
positive, and they coincide with the Frobenius--Perron dimensions
\cite[Proposition~8.23]{ENO2005}.
We always take this spherical structure, so that our assumption is
equivalent to: $\cC$ is a spherical braided fusion category in which the quantum
dimension $d_X$ of every simple object $X$ is positive.
Unitary modular tensor categories are the main examples, but no
unitary structure on $\cC$ will be used (Remark~\ref{rem:no-unitarity}).
} braided fusion category over $\bC$,
which we can assume small and strict without loss of generality. 
We write $\Irr(\cC)$ for a set of
representatives of the isomorphism classes of simple objects, $d_X>0$ for the quantum
dimension and
\[
D:=\Bigl(\sum_{X\in\Irr(\cC)}d_X^2\Bigr)^{1/2}.
\]
Without a subscript, $\Hom$ and $\End$ are taken in $\cC$, or in the
category under discussion; for a vector space $V$, $\End_\bC(V)$ is the
algebra of linear maps $V\to V$.
We write $g\circ f$, or simply $gf$, for the composition of morphisms.
Starting from Section~\ref{sec:integrality},
we will also assume that $\cC$ is modular,
as will be explicitly stated there.

We use $c_{X,Y}: X\otimes Y\to Y\otimes X$ for the braiding
and $\theta_X$ for the twist. Then the balancing axiom reads
\[
\theta_{X\otimes Y}=c_{Y,X}\circ c_{X,Y}\circ(\theta_X\otimes\theta_Y).
\]
We use the graphical calculus of $\cC$
\cite{Turaev1994}: diagrams are read from bottom to top,
a morphism $f$ is drawn as a coupon containing $f$,
and the braiding $c_{X,Y}$ is drawn with the $X$-strand passing over the
$Y$-strand.

Direct sums are
treated with the canonical distributivity isomorphisms, which we
suppress.
We use the left duals of $\cC$, $\ev_V\colon V^\vee\otimes V\to\one$ and
$\coev_V\colon\one\to V\otimes V^\vee$, and write
$\ev'_V\colon V\otimes V^\vee\to\one$ and $\coev'_V\colon\one\to V^\vee\otimes V$
for the right evaluation and coevaluation of the spherical structure.
In pictures a cap or a cup is drawn as follows:
\[
\ev_V=
\begin{tikzpicture}[mb]
\draw[mstr] (0,0) -- (0,0.6) .. controls (0,1.2) and (1,1.2) .. (1,0.6) -- (1,0);
\blab{0}{V^\vee}\blab{1}{V}
\end{tikzpicture},\qquad
\coev_V=
\begin{tikzpicture}[mb]
\draw[mstr] (0,1.3) -- (0,0.7) .. controls (0,0.1) and (1,0.1) .. (1,0.7) -- (1,1.3);
\tlab{0}{1.3}{V}\tlab{1}{1.3}{V^\vee}
\end{tikzpicture},\qquad
\ev'_V=
\begin{tikzpicture}[mb]
\draw[mstr] (0,0) -- (0,0.6) .. controls (0,1.2) and (1,1.2) .. (1,0.6) -- (1,0);
\blab{0}{V}\blab{1}{V^\vee}
\end{tikzpicture},\qquad
\coev'_V=
\begin{tikzpicture}[mb]
\draw[mstr] (0,1.3) -- (0,0.7) .. controls (0,0.1) and (1,0.1) .. (1,0.7) -- (1,1.3);
\tlab{0}{1.3}{V^\vee}\tlab{1}{1.3}{V}
\end{tikzpicture}.
\]
For $f\colon V_1\to V_2$ the
transpose $f^\vee\colon V_2^\vee\to V_1^\vee$ is such that
\begin{equation}\label{eq:transpose}
\begin{tikzpicture}[mb]
\draw[mstr] (0,0) -- (0,1.5) .. controls (0,2.1) and (1,2.1) .. (1,1.5) -- (1,0);
\node[mcpn] at (0,0.75) {$f^\vee$};
\node[anchor=east,inner sep=1pt] at (0,1.45) {$V_1^\vee$};
\blab{0}{V_2^\vee}\blab{1}{V_1}
\end{tikzpicture}
=
\begin{tikzpicture}[mb]
\draw[mstr] (0,0) -- (0,1.5) .. controls (0,2.1) and (1,2.1) .. (1,1.5) -- (1,0);
\node[mcpn] at (1,0.75) {$f$};
\node[anchor=west,inner sep=1pt] at (1,1.45) {$V_2$};
\blab{0}{V_2^\vee}\blab{1}{V_1}
\end{tikzpicture},
\qquad
\begin{tikzpicture}[mb]
\draw[mstr] (0,2.1) -- (0,0.6) .. controls (0,0) and (1,0) .. (1,0.6) -- (1,2.1);
\node[mcpn] at (0,1.35) {$f$};
\node[anchor=east,inner sep=1pt] at (0,0.65) {$V_1$};
\tlab{0}{2.1}{V_2}\tlab{1}{2.1}{V_1^\vee}
\end{tikzpicture}
=
\begin{tikzpicture}[mb]
\draw[mstr] (0,2.1) -- (0,0.6) .. controls (0,0) and (1,0) .. (1,0.6) -- (1,2.1);
\node[mcpn] at (1,1.35) {$f^\vee$};
\node[anchor=west,inner sep=1pt] at (1,0.65) {$V_2^\vee$};
\tlab{0}{2.1}{V_2}\tlab{1}{2.1}{V_1^\vee}
\end{tikzpicture}.
\end{equation}
For $V\otimes Y$ we use the left dual $Y^\vee\otimes V^\vee$, whose cap and cup are
nested,
\[
\ev_{V\otimes Y}=
\begin{tikzpicture}[mb]
\draw[mstr] (0.6,0) -- (0.6,0.5) .. controls (0.6,1.0) and (1.2,1.0) .. (1.2,0.5) -- (1.2,0);
\draw[mstr] (0,0) -- (0,0.5) .. controls (0,1.6) and (1.8,1.6) .. (1.8,0.5) -- (1.8,0);
\blab{0}{Y^\vee}\blab{0.6}{V^\vee}\blab{1.2}{V}\blab{1.8}{Y}
\end{tikzpicture},\qquad
\coev_{V\otimes Y}=
\begin{tikzpicture}[mb]
\draw[mstr] (0.6,1.5) -- (0.6,1.0) .. controls (0.6,0.5) and (1.2,0.5) .. (1.2,1.0) -- (1.2,1.5);
\draw[mstr] (0,1.5) -- (0,1.0) .. controls (0,-0.1) and (1.8,-0.1) .. (1.8,1.0) -- (1.8,1.5);
\tlab{0}{1.5}{V}\tlab{0.6}{1.5}{Y}\tlab{1.2}{1.5}{Y^\vee}\tlab{1.8}{1.5}{V^\vee}
\end{tikzpicture},
\]
and $\one^\vee=\one$; then $(f\otimes h)^\vee=h^\vee\otimes f^\vee$.
Two left
duals of an object are related by a unique isomorphism compatible with
evaluation and coevaluation \cite[Proposition~2.10.5]{EGNO}, and everything below
is compatible with these isomorphisms, so we suppress them.
In particular, for $X\in\Irr(\cC)$,
we always choose its dual $X^\vee$ from among $\Irr(\cC)$,
i.e.~we write $X^\vee\in\Irr(\cC)$.

A closed loop labeled $V$ is $\ev_V\coev'_V=d_V$, and the partial trace
$\operatorname{ptr}_V(g):=(\ev_V\otimes\id_W)(\id_{V^\vee}\otimes g)
(\coev'_V\otimes\id_W)$ over $V$ of $g\in\End(V\otimes W)$ closes the strand $V$ on the left,
and for $W=\one$ it is the trace $\operatorname{tr}(g)$ of $g\in\End(V)$
($\Tr$ denotes the trace of a linear map):
\[
\begin{tikzpicture}[mb]
\draw[mstr] (0,0.7) ellipse (0.45 and 0.5);
\node[anchor=west,inner sep=1.5pt] at (0.45,0.7) {$V$};
\end{tikzpicture}
=d_V,\qquad
\operatorname{ptr}_V(g)=
\begin{tikzpicture}[mb]
\draw[mstr] (0,0.9) -- (0,0.45) .. controls (0,-0.05) and (-0.7,-0.05) .. (-0.7,0.45) -- (-0.7,1.95)
  .. controls (-0.7,2.45) and (0,2.45) .. (0,1.95) -- (0,1.5);
\draw[mstr] (0.8,0) -- (0.8,2.4);
\node[mcpn,minimum width=13mm] at (0.4,1.2) {$g$};
\node[anchor=west,inner sep=1pt] at (0,0.67) {$V$};
\blab{0.8}{W}\tlab{0.8}{2.4}{W}
\end{tikzpicture}.
\]

To lighten the notation, we also use the following convention:
\begin{definition}\label{def:repeated}
Two shaded coupons carrying the same label $e$, one in $\Hom(A,B)$ and one
in $\Hom(B,A)$, stand for an element
$e=\sum_if_i\otimes g^i\in\Hom(A,B)\otimes\Hom(B,A)$: a diagram containing
them denotes the sum over $i$ of the diagrams with $f_i$ and $g^i$ in the
two coupons, for example
\[
\begin{tikzpicture}[mb]
\draw[mstr] (0,0) -- (0,2.0);
\node[mcpnc] at (0,1.0) {$e$};
\blab{0}{A}\tlab{0}{2.0}{B}
\draw[mstr] (1.2,0) -- (1.2,2.0);
\node[mcpnc] at (1.2,1.0) {$e$};
\blab{1.2}{B}\tlab{1.2}{2.0}{A}
\end{tikzpicture}
\;=\;\sum_i\;
\begin{tikzpicture}[mb]
\draw[mstr] (0,0) -- (0,2.0);
\node[mcpn] at (0,1.0) {$f_i$};
\blab{0}{A}\tlab{0}{2.0}{B}
\draw[mstr] (1.2,0) -- (1.2,2.0);
\node[mcpn] at (1.2,1.0) {$g^i$};
\blab{1.2}{B}\tlab{1.2}{2.0}{A}
\end{tikzpicture}\,.
\]
Shaded pairs with different labels are summed independently.
Unless stated otherwise, $e$ is the canonical element
\[
e_{A,B}=\sum_if_i\otimes g^i\in\Hom(A,B)\otimes\Hom(B,A),\qquad
\operatorname{tr}(g^if_j)=\delta_{ij},
\]
of the trace pairing $(f,g)\mapsto\operatorname{tr}(g\circ f)$, which is
non-degenerate since $\cC$ is semisimple with non-zero dimensions; $e_{A,B}$
does not depend on the choice of the basis $f_i$.
Since $\operatorname{tr}(g\circ f)=\operatorname{tr}(f\circ g)$, the element
$e_{B,A}$ is $e_{A,B}$ with its two factors exchanged, so the two coupons of
the pair play symmetric roles.
\end{definition}

\begin{proposition}\label{prop:canonical-element}
Shaded coupons carrying the canonical elements have the following properties.
\begin{enumerate}[label=(E\arabic*)]
\item\label{item:E1} \emph{Sliding.}  A morphism attached to one end of $e$ may be
moved to the other: for $h\colon B\to B'$,
\[
\begin{tikzpicture}[mb]
\draw[mstr] (0,0) -- (0,2.8);
\node[mcpnc] at (0,0.8) {$e$}; \node[mcpn] at (0,2.0) {$h$};
\node[anchor=west,inner sep=1pt] at (0,1.4) {$B$};
\blab{0}{A}\tlab{0}{2.8}{B'}
\draw[mstr] (1.4,0) -- (1.4,2.8);
\node[mcpnc] at (1.4,0.8) {$e$};
\blab{1.4}{B}\tlab{1.4}{2.8}{A}
\end{tikzpicture}
\;=\;
\begin{tikzpicture}[mb]
\draw[mstr] (0,0) -- (0,2.8);
\node[mcpnc] at (0,2.0) {$e$};
\blab{0}{A}\tlab{0}{2.8}{B'}
\draw[mstr] (1.4,0) -- (1.4,2.8);
\node[mcpn] at (1.4,0.8) {$h$}; \node[mcpnc] at (1.4,2.0) {$e$};
\node[anchor=west,inner sep=1pt] at (1.4,1.4) {$B'$};
\blab{1.4}{B}\tlab{1.4}{2.8}{A}
\end{tikzpicture},
\]
and similarly at the ends $A$.
\item\label{item:E2} \emph{Closing.}  For a linear map $\Phi$ of $\Hom(A,B)$, write $\Phi(e)$ for
$\Phi$ applied to the end of $e$ in $\Hom(A,B)$.
Then
\[
\begin{tikzpicture}[mb]
\draw[mstr] (0,0.2) -- (0,2.4) .. controls (0,2.9) and (0.8,2.9) .. (0.8,2.4) -- (0.8,0.2)
  .. controls (0.8,-0.3) and (0,-0.3) .. (0,0.2);
\node[mcpnc] at (0,0.9) {$\Phi(e)$}; \node[mcpnc] at (0,1.9) {$e$};
\node[anchor=west,inner sep=1pt] at (0,0.4) {$A$};
\node[anchor=west,inner sep=1pt] at (0,1.4) {$B$};
\end{tikzpicture}
=\Tr\Phi,
\qquad\text{and, if $A$ is simple,}\qquad
\begin{tikzpicture}[mb]
\draw[mstr] (0,0) -- (0,2.8);
\node[mcpnc] at (0,0.8) {$\Phi(e)$}; \node[mcpnc] at (0,2.0) {$e$};
\node[anchor=west,inner sep=1pt] at (0,1.4) {$B$};
\blab{0}{A}\tlab{0}{2.8}{A}
\end{tikzpicture}
=\frac{\Tr\Phi}{d_A}\;
\begin{tikzpicture}[mb]
\draw[mstr] (0,0) -- (0,2.8);
\blab{0}{A}\tlab{0}{2.8}{A}
\end{tikzpicture}.
\]
\item\label{item:E3} \emph{Resolution of the identity.}  We have
\[
\begin{tikzpicture}[mb]
\draw[mstr] (0,0) -- (0,2.8);
\blab{0}{Z}\tlab{0}{2.8}{Z}
\end{tikzpicture}
=\sum_{A\in\Irr(\cC)}d_A\;
\begin{tikzpicture}[mb]
\draw[mstr] (0,0) -- (0,0.6);
\draw[mstr] (0,0.9) -- (0,1.9);
\draw[mstr] (0,2.2) -- (0,2.8);
\node[mcpnc] at (0,0.75) {$e$};
\node[mcpnc] at (0,2.05) {$e$};
\node[anchor=west,inner sep=1pt] at (0,1.4) {$A$};
\blab{0}{Z}\tlab{0}{2.8}{Z}
\end{tikzpicture}
\]
for every $Z\in\cC$.
\item\label{item:E4} \emph{Bending.}  Bending a leg of both ends of $e$ in the same way
gives the canonical element of the resulting $\Hom$-spaces.
For example, if $e$ is the canonical element of
$\Hom(X,A\otimes B)\otimes\Hom(A\otimes B,X)$ and $e'$ that of
$\Hom(A^\vee\otimes X,B)\otimes\Hom(B,A^\vee\otimes X)$, then
\[
\begin{tikzpicture}[mb]
\draw[mstr] (1,0) -- (1,0.9);
\node[mcpnc, minimum width=8mm] at (1,1.2) {$e$};
\draw[mstr] (0.65,1.5) -- (0.65,1.7) .. controls (0.65,2.3) and (0,2.3) .. (0,1.7) -- (0,0);
\draw[mstr] (1.35,1.5) -- (1.35,2.6);
\node[anchor=west,inner sep=1pt] at (0.65,1.75) {$A$};
\blab{0}{A^\vee}\blab{1}{X}\tlab{1.35}{2.6}{B}
\end{tikzpicture}
\qquad
\begin{tikzpicture}[mb]
\draw[mstr] (0,2.6) -- (0,0.7) .. controls (0,0.1) and (0.65,0.1) .. (0.65,0.7) -- (0.65,0.9);
\draw[mstr] (1.35,0) -- (1.35,0.9);
\node[mcpnc, minimum width=8mm] at (1,1.2) {$e$};
\draw[mstr] (1,1.5) -- (1,2.6);
\node[anchor=west,inner sep=1pt] at (0.65,0.65) {$A$};
\blab{1.35}{B}\tlab{0}{2.6}{A^\vee}\tlab{1}{2.6}{X}
\end{tikzpicture}
\quad=\quad
\begin{tikzpicture}[mb]
\draw[mstr] (0,0) -- (0,0.9); \draw[mstr] (0.8,0) -- (0.8,0.9);
\node[mcpnc, minimum width=11mm] at (0.4,1.2) {$e'$};
\draw[mstr] (0.4,1.5) -- (0.4,2.6);
\blab{0}{A^\vee}\blab{0.8}{X}\tlab{0.4}{2.6}{B}
\end{tikzpicture}
\qquad
\begin{tikzpicture}[mb]
\draw[mstr] (0.4,0) -- (0.4,0.9);
\node[mcpnc, minimum width=11mm] at (0.4,1.2) {$e'$};
\draw[mstr] (0,1.5) -- (0,2.6); \draw[mstr] (0.8,1.5) -- (0.8,2.6);
\blab{0.4}{B}\tlab{0}{2.6}{A^\vee}\tlab{0.8}{2.6}{X}
\end{tikzpicture}.
\]
\end{enumerate}
\end{proposition}

\begin{proof}
Write $e=\sum_if_i\otimes g^i$ as in Definition~\ref{def:repeated}.
Then
$x=\sum_i\operatorname{tr}(g^ix)\,f_i$ for $x\in\Hom(A,B)$ and
$y=\sum_i\operatorname{tr}(yf_i)\,g^i$ for $y\in\Hom(B,A)$.

\medskip\noindent \ref{item:E1}
An element of $\Hom(A,B')\otimes\Hom(B,A)$ is determined by its pairings
with $y\otimes x$ for $y\in\Hom(B',A)$ and $x\in\Hom(A,B)$, under
$(f\otimes g,\,y\otimes x)\mapsto\operatorname{tr}(yf)\operatorname{tr}(gx)$.
By the expansions above, both sides pair to $\operatorname{tr}(yhx)$.
The ends $A$ are treated in the same way.

\medskip\noindent \ref{item:E2}
The closed diagram is $\sum_i\operatorname{tr}(g^i\,\Phi(f_i))$, which is
$\Tr\Phi$ by the first expansion.
If $A$ is simple, $\sum_ig^i\,\Phi(f_i)\in\End(A)$ is a multiple of $\id_A$,
and taking the trace determines the multiple.

\medskip\noindent \ref{item:E3}
Write $Z$ as a direct sum with inclusions $\iota_i\colon A\to Z$ and
projections $\pi_i\colon Z\to A$ of the summands isomorphic to $A$, so that
$\pi_j\iota_i=\delta_{ij}\id_A$.
Then $\operatorname{tr}(\pi_j\iota_i)=d_A\delta_{ij}$, so we may take
$f_i=\iota_i$ and $g^i=d_A^{-1}\pi_i$, and $\sum_if_ig^i=d_A^{-1}\sum_i\iota_i\pi_i$.
Summing the projections over $A\in\Irr(\cC)$ gives $\id_Z$, since $\cC$ is semisimple.

\medskip\noindent \ref{item:E4}
Bending a leg is a linear isomorphism of $\Hom$-spaces, and bending the
two ends in the same way does not change the trace of their composite, by
sphericality.
So the bent $f_i$ and $g^i$ are again dual bases.
\end{proof}

\subsection{Anti-monoidal involution}
\label{subsec:P}

\begin{definition}\label{def:anti-monoidal}
An \emph{anti-monoidal involution} on $\cC$ consists of a
$\bC$-linear autoequivalence $\sP\colon\cC\to\cC$ with natural isomorphisms
\[
U_{X,Y}\colon \sP Y\otimes\sP X\xrightarrow{\ \cong\ }\sP(X\otimes Y),
\qquad U_0\colon\one\xrightarrow{\ \cong\ }\sP\one,
\]
satisfying
\begin{align}
&U_{X,\,Y\otimes Z}\circ(U_{Y,Z}\otimes\id_{\sP X})
 =U_{X\otimes Y,\,Z}\circ(\id_{\sP Z}\otimes U_{X,Y}),
 \tag{A1}\label{eq:U-coherence}\\
&U_{\one,X}\circ(\id_{\sP X}\otimes U_0)=\id_{\sP X}
 =U_{X,\one}\circ(U_0\otimes\id_{\sP X}),
 \tag{A2}\label{eq:U-unit}\\
&\sP(c_{X,Y})\circ U_{X,Y}=U_{Y,X}\circ c_{\sP X,\sP Y}^{-1},
 \tag{A3}\label{eq:U-braiding}
\end{align}
together with a monoidal natural isomorphism 
$\beta\colon\id_\cC\Rightarrow\sP^2$, that is,
\begin{align}
&\beta_{X\otimes Y}
 =\sP(U_{X,Y})\circ U_{\sP Y,\sP X}\circ(\beta_X\otimes\beta_Y),
 \qquad \beta_\one=\sP(U_0)\circ U_0.
 \tag{A4}\label{eq:beta-monoidal}
\end{align} 
\end{definition}

We always assume that a datum
$(\sP,U,\beta)$ satisfying
\eqref{eq:U-coherence}--\eqref{eq:beta-monoidal} is fixed.

\begin{definition}\label{def:Pbox}
For objects $V_1,\dots,V_m$ let
$U_{V_1,\dots,V_m}\colon\sP V_m\otimes\cdots\otimes\sP V_1\to
\sP(V_1\otimes\cdots\otimes V_m)$ be the iterated anti-monoidal
structure; by \eqref{eq:U-coherence} it does not depend on the
bracketing, and it is $\id_{\sP V_1}$ for $m=1$ and $U_0$ for $m=0$.
We use it to draw images under $\sP$ of morphisms of $\cC$.
Let $\Gamma$ be a diagram of $\cC$
for a morphism $V_1\otimes\dots\otimes V_m\to W_1\otimes\dots\otimes W_n$.
We
draw the morphism
\[
\Ubox{\Gamma}:=U^{-1}_{W_1,\dots,W_n}\circ\sP(\Gamma)\circ U_{V_1,\dots,V_m}\colon
\sP V_m\otimes\dots\otimes\sP V_1\to\sP W_n\otimes\dots\otimes\sP W_1
\]
as the \emph{mirror image} of $\Gamma$ inside a dashed box, the \emph{$\sP$-box}:
\[
\Ubox{\Gamma}\;=\;
\begin{tikzpicture}[mb]
\draw[mstr] (0,0) -- (0,2.9); \draw[mstr] (2.4,0) -- (2.4,2.9);
\node[mcpn, minimum width=22mm] at (1.2,1.45) {$\Gamma$};
\ubox{-0.55}{0.55}{2.95}{2.35}
\node[anchor=west,inner sep=1pt] at (0,0.85) {$V_m$}; \node[anchor=east,inner sep=1pt] at (2.4,0.85) {$V_1$};
\node[anchor=west,inner sep=1pt] at (0,2.05) {$W_n$}; \node[anchor=east,inner sep=1pt] at (2.4,2.05) {$W_1$};
\blab{0}{\sP V_m}\blab{1.2}{\cdots}\blab{2.4}{\sP V_1}
\tlab{0}{2.9}{\sP W_n}\tlab{1.2}{2.9}{\cdots}\tlab{2.4}{2.9}{\sP W_1}
\end{tikzpicture}
\qedpic
\]
\end{definition}

The picture inside the box is the reflection of $\Gamma$ in a vertical
line, read in a mirror: a coupon labeled $g$ inside a box has its legs
in the reverse order.
The strands
outside the box carry the images under $\sP$ of the labels inside, in
the same order.
At an end with one strand $U$ is the identity, and at
an end without strands it is $U_0$.
For example, for a morphism
$f\colon V_1\to V_2$ and a
morphism $g\colon V_1\otimes V_2\otimes V_3\to W_1\otimes W_2$,
\[
\begin{tikzpicture}[mb]
\draw[mstr] (0,0) -- (0,2.2);
\node[mcpn] at (0,1.1) {$f$};
\ubox{-0.45}{0.55}{0.45}{1.65}
\blab{0}{\sP V_1}\tlab{0}{2.2}{\sP V_2}
\end{tikzpicture}
=\sP(f),\qquad
\begin{tikzpicture}[mb]
\draw[mstr] (0,0) -- (0,1.2); \draw[mstr] (1,0) -- (1,1.2); \draw[mstr] (2,0) -- (2,1.2);
\draw[mstr] (0.5,1.6) -- (0.5,2.6); \draw[mstr] (1.5,1.6) -- (1.5,2.6);
\node[mcpn, minimum width=20mm] at (1,1.4) {$g$};
\ubox{-0.45}{0.7}{2.45}{2.0}
\blab{0}{\sP V_3}\blab{1}{\sP V_2}\blab{2}{\sP V_1}\tlab{0.5}{2.6}{\sP W_2}\tlab{1.5}{2.6}{\sP W_1}
\end{tikzpicture}
=U^{-1}_{W_1,W_2}\,\sP(g)\,U_{V_1,V_2,V_3}.
\]
In the first example both ends have one strand, so no $U$ appears.

\begin{proposition}\label{prop:box-rules}
The boxes obey the following rules, each an equality of morphisms of
$\cC$.
\begin{enumerate}[label=(B\arabic*)]
\item\label{item:B1} \emph{Composition.}  Stacked boxes merge:
\[
\begin{tikzpicture}[mb]
\draw[mstr] (0,0) -- (0,3.4);
\node[mcpn] at (0,0.95) {$\Gamma_1$};
\node[mcpn] at (0,2.45) {$\Gamma_2$};
\ubox{-0.45}{0.45}{0.45}{1.45}
\ubox{-0.45}{1.95}{0.45}{2.95}
\end{tikzpicture}
=
\begin{tikzpicture}[mb]
\draw[mstr] (0,0) -- (0,3.4);
\node[mcpn] at (0,0.95) {$\Gamma_1$};
\node[mcpn] at (0,2.45) {$\Gamma_2$};
\ubox{-0.45}{0.45}{0.45}{2.95}
\end{tikzpicture}.
\]
\item\label{item:B2} \emph{Tensor product.}  Boxes side by side merge:
\[
\begin{tikzpicture}[mb]
\draw[mstr] (0,0) -- (0,2.4); \draw[mstr] (1.1,0) -- (1.1,2.4);
\node[mcpn] at (0,1.2) {$\Gamma_1$}; \node[mcpn] at (1.1,1.2) {$\Gamma_2$};
\ubox{-0.45}{0.6}{0.45}{1.8}
\ubox{0.65}{0.6}{1.55}{1.8}
\end{tikzpicture}
=
\begin{tikzpicture}[mb]
\draw[mstr] (0,0) -- (0,2.4); \draw[mstr] (1.1,0) -- (1.1,2.4);
\node[mcpn] at (0,1.2) {$\Gamma_1$}; \node[mcpn] at (1.1,1.2) {$\Gamma_2$};
\ubox{-0.45}{0.6}{1.55}{1.8}
\end{tikzpicture}.
\]
Since the box reverses the order of the strands, in formulas this
reads $\Ubox{\Gamma_1}\otimes\Ubox{\Gamma_2}=\Ubox{\Gamma_2\otimes \Gamma_1}$.
\item\label{item:B3} \emph{Identity.}  A box around an identity strand may be removed:
\[
\begin{tikzpicture}[mb]
\draw[mstr] (0,0) -- (0,2.4);
\ubox{-0.45}{0.6}{0.45}{1.8}
\end{tikzpicture}
=
\begin{tikzpicture}[mb]
\draw[mstr] (0,0) -- (0,2.4);
\end{tikzpicture}.
\]
\item\label{item:B4} \emph{Crossings.}  Crossings pass
through the walls of a box unchanged:
\[
\begin{tikzpicture}[mb]
\draw[mstr] (0,0) -- (0,0.5); \draw[mstr] (0.9,0) -- (0.9,0.5);
\draw[mstr] (0,1.7) -- (0,2.2); \draw[mstr] (0.9,1.7) -- (0.9,2.2);
\sstr{0}{0.5}{0.9}{1.7}
\sstr[mover]{0.9}{0.5}{0}{1.7}
\ubox{-0.4}{0.4}{1.3}{1.8}
\blab{0}{\sP Y}\blab{0.9}{\sP X}\tlab{0}{2.2}{\sP X}\tlab{0.9}{2.2}{\sP Y}
\end{tikzpicture}
=
\begin{tikzpicture}[mb]
\draw[mstr] (0,0) -- (0,0.5); \draw[mstr] (0.9,0) -- (0.9,0.5);
\draw[mstr] (0,1.7) -- (0,2.2); \draw[mstr] (0.9,1.7) -- (0.9,2.2);
\sstr{0}{0.5}{0.9}{1.7}
\sstr[mover]{0.9}{0.5}{0}{1.7}
\blab{0}{\sP Y}\blab{0.9}{\sP X}\tlab{0}{2.2}{\sP X}\tlab{0.9}{2.2}{\sP Y}
\end{tikzpicture},
\]
and in the same way for the other crossing.
\item\label{item:B5} \emph{Caps and cups, I.}
$\sP(V^\vee)$, with $\Ubox{\ev'_V}\colon\sP(V^\vee)\otimes\sP V\to\one$ and
$\Ubox{\coev'_V}\colon\one\to\sP V\otimes\sP(V^\vee)$, is a left dual of $\sP V$.
Following the convention of Section~\ref{subsec:conventions}, we use it as $(\sP V)^\vee$, also written $\sP V^\vee$, so that
$\Ubox{\ev'_V}=\ev_{\sP V}$ and $\Ubox{\coev'_V}=\coev_{\sP V}$, that is,
\[
\begin{tikzpicture}[mb]
\draw[mstr] (0,0) -- (0,1.3); \draw[mstr] (1,0) -- (1,1.3);
\draw[mstr] (0,1.3) .. controls (0,1.9) and (1,1.9) .. (1,1.3);
\ubox{-0.4}{0.75}{1.4}{1.95}
\node[anchor=west,inner sep=1pt] at (0,1.05) {$V^\vee$}; \node[anchor=east,inner sep=1pt] at (1,1.05) {$V$};
\blab{0}{\sP(V^\vee)}\blab{1}{\sP V}
\end{tikzpicture}
=
\begin{tikzpicture}[mb]
\draw[mstr] (0,0) -- (0,1.3); \draw[mstr] (1,0) -- (1,1.3);
\draw[mstr] (0,1.3) .. controls (0,1.9) and (1,1.9) .. (1,1.3);
\blab{0}{\sP(V^\vee)}\blab{1}{\sP V}
\end{tikzpicture},
\qquad
\begin{tikzpicture}[mb]
\draw[mstr] (0,1.95) -- (0,0.65); \draw[mstr] (1,1.95) -- (1,0.65);
\draw[mstr] (0,0.65) .. controls (0,0.05) and (1,0.05) .. (1,0.65);
\ubox{-0.4}{0}{1.4}{1.2}
\node[anchor=west,inner sep=1pt] at (0,0.9) {$V$}; \node[anchor=east,inner sep=1pt] at (1,0.9) {$V^\vee$};
\tlab{0}{1.95}{\sP V}\tlab{1}{1.95}{\sP(V^\vee)}
\end{tikzpicture}
=
\begin{tikzpicture}[mb]
\draw[mstr] (0,1.95) -- (0,0.65); \draw[mstr] (1,1.95) -- (1,0.65);
\draw[mstr] (0,0.65) .. controls (0,0.05) and (1,0.05) .. (1,0.65);
\tlab{0}{1.95}{\sP V}\tlab{1}{1.95}{\sP(V^\vee)}
\end{tikzpicture}.
\]
\item\label{item:B6} \emph{Caps and cups, II.}  With this identification,
$\Ubox{\ev_V}=\ev'_{\sP V}$ and $\Ubox{\coev_V}=\coev'_{\sP V}$, that is,
\[
\begin{tikzpicture}[mb]
\draw[mstr] (0,0) -- (0,1.3); \draw[mstr] (1,0) -- (1,1.3);
\draw[mstr] (0,1.3) .. controls (0,1.9) and (1,1.9) .. (1,1.3);
\ubox{-0.4}{0.75}{1.4}{1.95}
\node[anchor=west,inner sep=1pt] at (0,1.05) {$V$}; \node[anchor=east,inner sep=1pt] at (1,1.05) {$V^\vee$};
\blab{0}{\sP V}\blab{1}{\sP(V^\vee)}
\end{tikzpicture}
=
\begin{tikzpicture}[mb]
\draw[mstr] (0,0) -- (0,1.3); \draw[mstr] (1,0) -- (1,1.3);
\draw[mstr] (0,1.3) .. controls (0,1.9) and (1,1.9) .. (1,1.3);
\blab{0}{\sP V}\blab{1}{\sP(V^\vee)}
\end{tikzpicture},
\qquad
\begin{tikzpicture}[mb]
\draw[mstr] (0,1.95) -- (0,0.65); \draw[mstr] (1,1.95) -- (1,0.65);
\draw[mstr] (0,0.65) .. controls (0,0.05) and (1,0.05) .. (1,0.65);
\ubox{-0.4}{0}{1.4}{1.2}
\node[anchor=west,inner sep=1pt] at (0,0.9) {$V^\vee$}; \node[anchor=east,inner sep=1pt] at (1,0.9) {$V$};
\tlab{0}{1.95}{\sP(V^\vee)}\tlab{1}{1.95}{\sP V}
\end{tikzpicture}
=
\begin{tikzpicture}[mb]
\draw[mstr] (0,1.95) -- (0,0.65); \draw[mstr] (1,1.95) -- (1,0.65);
\draw[mstr] (0,0.65) .. controls (0,0.05) and (1,0.05) .. (1,0.65);
\tlab{0}{1.95}{\sP(V^\vee)}\tlab{1}{1.95}{\sP V}
\end{tikzpicture}.
\]
\item\label{item:B7} \emph{Double box.}  A box around a box puts $\beta$ on every strand:
\[
\begin{tikzpicture}[mb]
\draw[mstr] (0,0) -- (0,3.4); \draw[mstr] (1.9,0) -- (1.9,3.4);
\node at (0.95,0.3) {$\cdots$}; \node at (0.95,3.1) {$\cdots$};
\node[mcpn, minimum width=22mm] at (0.95,1.7) {$\Gamma$};
\ubox{-0.45}{1.15}{2.35}{2.25}
\ubox{-0.8}{0.55}{2.7}{2.85}
\blab{0}{\sP^2V_1}\blab{1.9}{\sP^2V_m}\tlab{0}{3.4}{\sP^2W_1}\tlab{1.9}{3.4}{\sP^2W_n}
\end{tikzpicture}
=
\begin{tikzpicture}[mb]
\draw[mstr] (0,0) -- (0,3.4); \draw[mstr] (1.9,0) -- (1.9,3.4);
\node[mcpn] at (0,0.65) {$\beta^{-1}_{V_1}$}; \node at (0.95,0.65) {$\cdots$}; \node[mcpn] at (1.9,0.65) {$\beta^{-1}_{V_m}$};
\node[mcpn, minimum width=22mm] at (0.95,1.7) {$\Gamma$};
\node[mcpn] at (0,2.75) {$\beta_{W_1}$}; \node at (0.95,2.75) {$\cdots$}; \node[mcpn] at (1.9,2.75) {$\beta_{W_n}$};
\blab{0}{\sP^2V_1}\blab{1.9}{\sP^2V_m}\tlab{0}{3.4}{\sP^2W_1}\tlab{1.9}{3.4}{\sP^2W_n}
\end{tikzpicture},
\]
that is,
$\Ubox{\Ubox{\Gamma}}=(\beta_{W_1}\otimes\cdots\otimes\beta_{W_n})\circ \Gamma\circ(\beta_{V_1}\otimes\cdots\otimes\beta_{V_m})^{-1}$,
with no $\beta$ at an end without strands.
\end{enumerate}
\end{proposition}

\begin{proof}
\ref{item:B1}
The middle factor $UU^{-1}$ cancels, and $\sP$ is a functor.

\medskip\noindent \ref{item:B2}
This is the naturality of $U$ and \eqref{eq:U-coherence}.
At an end without
strands one uses the unit axioms \eqref{eq:U-unit}.

\medskip\noindent \ref{item:B3}
This holds since $\sP(\id)=\id$.

\medskip\noindent \ref{item:B4}
The mirror image of $c_{X,Y}$ drawn inside the box is the
picture of $c^{-1}_{\sP X,\sP Y}$ with $\sP$ dropped from the labels, and
$U^{-1}_{Y,X}\sP(c_{X,Y})U_{X,Y}=c^{-1}_{\sP X,\sP Y}$ by \eqref{eq:U-braiding}.

\medskip\noindent \ref{item:B5}
Taking the boxes of the zigzag identities of $\ev'_V$ and $\coev'_V$, by
\ref{item:B1}--\ref{item:B3} we obtain the zigzag identities of a left dual.

\medskip\noindent \ref{item:B6}
By \ref{item:B1}--\ref{item:B3}, $\Ubox{\ev_V}$ and $\Ubox{\coev_V}$ satisfy
the zigzag identities which make $\sP V$ a left dual of $\sP(V^\vee)=(\sP V)^\vee$, as do
$\ev'_{\sP V}$ and $\coev'_{\sP V}$.
So $\Ubox{\ev_V}=\ev'_{\sP V}(g_V\otimes\id)$ and
$\Ubox{\coev_V}=(\id\otimes g_V^{-1})\coev'_{\sP V}$ for an automorphism $g_V$ of
$\sP V$, additive in $V$, and we may assume $V$ simple.
Then $g_V$ is a scalar, as $\sP V$ is simple.
The box of a scalar is the same scalar, since $\sP$ is linear.
Taking the boxes of the closed loops $\ev_V\coev'_V=d_V$ and $\ev'_V\coev_V=d_V$, by
\ref{item:B1} and \ref{item:B5} we get
\[
d_V=\ev'_{\sP V}(g_V\otimes\id)\coev_{\sP V}=g_V\,d_{\sP V},\qquad
d_V=\ev_{\sP V}(\id\otimes g_V^{-1})\coev'_{\sP V}=g_V^{-1}\,d_{\sP V}.
\]
Hence $g_V^2=1$, and $g_V=1$ since $d_V>0$ and $d_{\sP V}>0$.

\medskip\noindent \ref{item:B7}
This is
\eqref{eq:beta-monoidal}, applied repeatedly, and the naturality of
$\beta$.
\end{proof}

Since $\sP$ is $\bC$-linear, $\Ubox{\Gamma}$ depends linearly on each coupon of $\Gamma$; so the convention of Definition~\ref{def:repeated} and \ref{item:E1}--\ref{item:E4} apply also to shaded coupons inside boxes.

\begin{corollary}\label{cor:P-dim}
$d_{\sP X}=d_X$ for $X\in\cC$.
\end{corollary}

\begin{proof}
Take the box of the closed loop $\ev_X\coev'_X=d_X$ and use \ref{item:B1}, \ref{item:B5} and \ref{item:B6}:
$d_X=\Ubox{\ev_X}\Ubox{\coev'_X}=\ev'_{\sP X}\coev_{\sP X}=d_{\sP X}$.
\end{proof}

%
%

\section{The M\"obius fusion algebra}
\label{sec:mobius}

\subsection{Definition and basic properties}
\label{subsec:definition}

Let us first recall the fusion algebra of $\cC$.
\begin{definition}\label{def:fusion-algebra}
The \emph{fusion algebra} $\fF$ of $\cC$ is \[
\fF:=K_0(\cC)\otimes_\bZ\bC.
\]
We denote the class of an object $X$ by $[X]\in\fF$.
\end{definition}
The following is standard:
\begin{proposition}
  $\fF$ is a commutative associative unital algebra
  with a basis $[X]$ labeled by the simple objects $X\in\Irr(\cC)$, 
  with a product \[
  [X][Y]=\sum_{Z\in\Irr(\cC)}N_{XY}^Z[Z],\quad
  N_{XY}^Z:=\dim\Hom(X\otimes Y,Z).
  \]
\end{proposition}

A more explicit description of $\fF$ is as follows,
whose proof is straightforward and omitted:
\begin{proposition}
  \label{prop:fusion-algebra-explicit}
  Let 
  \[\fV_0:=\bigoplus_{X\in\Ob(\cC)}\Hom(X,X),\]
  and denote by $(X,f)$ the vector $f\in\Hom(X,X)$ in the summand labeled by $X$.
  Let $\fR_0\subset \fV_0$ be the subspace generated
by the vectors \[
  (X\oplus Y,f\oplus g)-(X,f)-(Y,g)
\] and \[
  (X,a\circ f)-(X,f\circ a),
\] where $f\in\Hom(X,X)$, $g\in\Hom(Y,Y)$ and $a\in\End(X)$.
  Then $\fF\cong\fV_0/\fR_0$ as algebras, where, writing $[X,f]$ for the
  class of $(X,f)$, we have $[X]=[X,\id_X]$ and
  $[X,f]\cdot[Y,g]=[X\otimes Y,f\otimes g]$; the isomorphism sends $[X,f]$ to
  $\sum_{Z\in\Irr(\cC)}\Tr(f_*)\,[Z]$, where $f_*$ is the linear map
  $\varphi\mapsto f\circ\varphi$ of $\Hom(Z,X)$, the action of $f$ on the
  multiplicity space of $Z$ in $X$; for $f=\id_X$ its trace is the
  multiplicity $\dim\Hom(Z,X)$.
\end{proposition}

We are going to define a M\"obius version of the fusion algebra, 
incorporating the anti-monoidal involution $\sP$ of $\cC$.
We define it by generalizing the description of $\fF$ 
in Proposition~\ref{prop:fusion-algebra-explicit} above.

\begin{definition}\label{def:X}
Let
\[
\fV:=\bigoplus_{X\in\Ob(\cC)}\Hom(X,\sP X),
\]
and write $(X,f)$ for the vector $f\in\Hom(X,\sP X)$ in the summand labeled by $X$.
Let $\fR\subset \fV$ be the subspace spanned by
\begin{enumerate}[label=(R\arabic*)]
\item (direct sums)
  $(X\oplus Y,\,f\oplus g)-(X,f)-(Y,g)$, where
  $f\oplus g:=\sP(\iota_X)\,f\,\pi_X+\sP(\iota_Y)\,g\,\pi_Y$ for the
  inclusions $\iota$ and projections $\pi$ of $X\oplus Y$;
\item (coinvariance) $(X,\sP(a)\circ f)-(X,f\circ a)$ for
  $a\in\End(X)$.
\end{enumerate}
The \emph{M\"obius fusion algebra} is the vector space $\fX:=\fV/\fR$, with the
product later defined in Proposition~\ref{prop:product},
which is not commutative in general.
The class of $(X,f)$ is
denoted $[X,f]$.
\end{definition}

\begin{definition}\label{def:bar}
We draw elements of $\fX$ in the graphical calculus of $\cC$
(Section~\ref{subsec:conventions}) between two horizontal arrows in the opposite directions:
\[
\mone{f}{X}{\sP X}\;:=\;[X,f].
\]
The top and the bottom arrows indicate that the top boundary $\sP X$ is glued back to the
bottom boundary $X$ through $\sP$, so that the picture describes a line
running once around the core circle of the M\"obius band.
We call the diagram we just introduced a \emph{M\"obius diagram}, 
and the line running around the core circle a \emph{core loop}.
\end{definition}

The two relations of Definition~\ref{def:X} read
\[
\begin{tikzpicture}[mb]
\mbars{-0.8}{0.8}{2}
\draw[mstr] (0,0) -- (0,2);
\node[mcpn] at (0,1) {$f\oplus g$};
\blab{0}{X\oplus Y}\tlab{0}{2}{\sP(X\oplus Y)}
\end{tikzpicture}
=\mone{f}{X}{\sP X}+\mone{g}{Y}{\sP Y},
\qquad
\begin{tikzpicture}[mb]
\mbars{-0.6}{0.6}{2.9}
\draw[mstr] (0,0) -- (0,2.9);
\node[mcpn] at (0,0.8) {$a$};
\node[mcpn] at (0,1.9) {$f$};
\blab{0}{X}\tlab{0}{2.9}{\sP X}
\end{tikzpicture}
=
\begin{tikzpicture}[mb]
\mbars{-0.6}{0.6}{2.9}
\draw[mstr] (0,0) -- (0,2.9);
\node[mcpn] at (0,0.75) {$f$};
\node[mcpn] at (0,1.85) {$a$};
\ubox{-0.45}{1.4}{0.45}{2.3}
\blab{0}{X}\tlab{0}{2.9}{\sP X}
\end{tikzpicture}
\]
where $a\in\End(X)$.

Next we extend the M\"obius diagrams so that they can include more than one strand at the top and bottom arrows.

\begin{definition}\label{def:bar-multi}
For objects $X_1,\dots,X_n$ and a morphism
$f\colon X_1\otimes\cdots\otimes X_n\to\sP X_n\otimes\cdots\otimes\sP X_1$
we put
\[
\begin{tikzpicture}[mb]
\mbars{-0.6}{3.6}{2}
\foreach \x in {0,1.2,3.0} \draw[mstr] (\x,0) -- (\x,2);
\node[mcpn, minimum width=32mm] at (1.5,1) {$f$};
\blab{0}{X_1}\blab{1.2}{X_2}\blab{2.1}{\cdots}\blab{3.0}{X_n}
\tlab{0}{2}{\sP X_n}\tlab{1.2}{2}{\sP X_{n-1}}\tlab{2.1}{2}{\cdots}\tlab{3.0}{2}{\sP X_1}
\end{tikzpicture}
\;:=\;\bigl[X_1\otimes\cdots\otimes X_n,\;U_{X_1,\dots,X_n}\circ f\bigr],
\]
where
$U_{X_1,\dots,X_n}\colon\sP X_n\otimes\cdots\otimes\sP X_1\to
\sP(X_1\otimes\cdots\otimes X_n)$ is the iterated anti-monoidal
structure of Definition~\ref{def:Pbox}.
\end{definition}

Thus the labels of the top arrow are those of the bottom arrow, with $\sP$
applied and the order reversed.
Strands may be grouped or split at
will, provided that the corresponding $U$'s are inserted or omitted.
For example,
\[
\bigl[X\otimes Y,\;U_{X,Y}\circ f\bigr]
\;=\;
\begin{tikzpicture}[mb]
\mbars{-0.5}{1.3}{2}
\draw[mstr] (0,0) -- (0,2); \draw[mstr] (0.8,0) -- (0.8,2);
\node[mcpn, minimum width=14.5mm] at (0.4,1) {$f$};
\blab{0}{X}\blab{0.8}{Y}\tlab{0}{2}{\sP Y}\tlab{0.8}{2}{\sP X}
\end{tikzpicture}
\;=\;
\begin{tikzpicture}[mb]
\mbars{-0.5}{1.3}{3.2}
\draw[mstr] (0.4,0) -- (0.4,0.3);
\draw[mstr] (0.4,0.3) .. controls (0.4,0.6) and (0,0.5) .. (0,0.85);
\draw[mstr] (0.4,0.3) .. controls (0.4,0.6) and (0.8,0.5) .. (0.8,0.85);
\draw[mstr] (0,1.15) .. controls (0,1.5) and (0.4,1.4) .. (0.4,1.7);
\draw[mstr] (0.8,1.15) .. controls (0.8,1.5) and (0.4,1.4) .. (0.4,1.7);
\draw[mstr] (0.4,1.7) -- (0.4,3.2);
\node[mcpn, minimum width=14.5mm] at (0.4,1.0) {$f$};
\node[mcpn] at (0.4,2.4) {$U_{X,Y}$};
\blab{0.4}{X\otimes Y}\tlab{0.4}{3.2}{\sP(X\otimes Y)}
\end{tikzpicture}
\]
for $f\colon X\otimes Y\to\sP Y\otimes\sP X$.

\begin{lemma}\label{lem:trace-property}
For $a\colon X\to Y$ and $h\colon Y\to\sP X$ one has
\[
[X,\,h\circ a]=[Y,\,\sP(a)\circ h],
\]
that is, in pictures, the morphism $a$ moves once around the core:
\[
\begin{tikzpicture}[mb]
\mbars{-0.6}{0.6}{2.4}
\draw[mstr] (0,0) -- (0,2.4);
\node[mcpn] at (0,0.75) {$a$};
\node[mcpn] at (0,1.65) {$h$};
\blab{0}{X}\tlab{0}{2.4}{\sP X}
\end{tikzpicture}
=
\begin{tikzpicture}[mb]
\mbars{-0.6}{0.6}{2.4}
\draw[mstr] (0,0) -- (0,2.4);
\node[mcpn] at (0,0.6) {$h$};
\node[mcpn] at (0,1.6) {$a$};
\ubox{-0.45}{1.15}{0.45}{2.05}
\blab{0}{Y}\tlab{0}{2.4}{\sP Y}
\end{tikzpicture}.
\]
In particular $[X,f]=[Y,\sP(u)\circ f\circ u^{-1}]$ for every
isomorphism $u\colon X\to Y$, that is,
\[
\begin{tikzpicture}[mb]
\mbars{-0.6}{0.6}{3.0}
\draw[mstr] (0,0) -- (0,3.0);
\node[mcpn] at (0,1.5) {$f$};
\blab{0}{X}\tlab{0}{3.0}{\sP X}
\end{tikzpicture}
=
\begin{tikzpicture}[mb]
\mbars{-0.6}{0.6}{3.0}
\draw[mstr] (0,0) -- (0,3.0);
\node[mcpn] at (0,0.55) {$u^{-1}$};
\node[mcpn] at (0,1.37) {$f$};
\node[mcpn] at (0,2.2) {$u$};
\ubox{-0.45}{1.75}{0.45}{2.65}
\blab{0}{Y}\tlab{0}{3.0}{\sP Y}
\end{tikzpicture}\,.
\]
\end{lemma}

\begin{proof}
Put $Z=X\oplus Y$ with inclusions $\iota$ and projections $\pi$.
By
relation (R1) (with the second summand carrying $0$),
$[Z,\sP(\iota_X)\varphi \pi_X]=[X,\varphi]$ for every
$\varphi\colon X\to\sP X$, and similarly for $Y$.
Put
$F:=\sP(\iota_X)\,h\,\pi_Y\colon Z\to\sP Z$ and
$b:=\iota_Y\,a\,\pi_X\in\End(Z)$.
Then $F b=\sP(\iota_X)(ha)\pi_X$ and
$\sP(b)F=\sP(\iota_Y)\bigl(\sP(a)h\bigr)\pi_Y$, so the claim follows from relation (R2).
The special case is $(a,h)=(u,fu^{-1})$.
\end{proof}

\begin{remark}\label{rem:box-trace}
When there is more than one strand in the diagram,
the lemma above results in a more general identity.
For example,
for $a\colon X_1\otimes X_2\to Y_1\otimes Y_2$ and
$h\colon Y_1\otimes Y_2\to\sP X_2\otimes\sP X_1$,
we have
\[
\begin{tikzpicture}[mb]
\mbars{-0.5}{1.8}{3.0}
\draw[mstr] (0,0) -- (0,3.0); \draw[mstr] (1.3,0) -- (1.3,3.0);
\node[mcpn, minimum width=16mm] at (0.65,0.8) {$a$};
\node[mcpn, minimum width=16mm] at (0.65,2.0) {$h$};
\node[anchor=east,inner sep=1pt] at (0,1.4) {$Y_1$};
\node[anchor=west,inner sep=1pt] at (1.3,1.4) {$Y_2$};
\blab{0}{X_1}\blab{1.3}{X_2}\tlab{0}{3.0}{\sP X_2}\tlab{1.3}{3.0}{\sP X_1}
\end{tikzpicture}
\;=\;
\begin{tikzpicture}[mb]
\mbars{-0.5}{1.8}{3.8}
\draw[mstr] (0,0) -- (0,3.8); \draw[mstr] (1.3,0) -- (1.3,3.8);
\node[mcpn, minimum width=16mm] at (0.65,0.6) {$h$};
\node[mcpn, minimum width=16mm] at (0.65,2.5) {$a$};
\ubox{-0.45}{1.6}{1.75}{3.35}
\node[anchor=east,inner sep=1pt] at (0,1.15) {$\sP X_2$};
\node[anchor=west,inner sep=1pt] at (1.3,1.15) {$\sP X_1$};
\node[anchor=west,inner sep=1pt] at (0,1.95) {$X_2$}; \node[anchor=east,inner sep=1pt] at (1.3,1.95) {$X_1$};
\node[anchor=west,inner sep=1pt] at (0,3.05) {$Y_2$}; \node[anchor=east,inner sep=1pt] at (1.3,3.05) {$Y_1$};
\blab{0}{Y_1}\blab{1.3}{Y_2}\tlab{0}{3.8}{\sP Y_2}\tlab{1.3}{3.8}{\sP Y_1}
\end{tikzpicture}\ .
\qedpic
\]
\end{remark}


\begin{lemma}\label{lem:crossing-move}
Let $f\colon X_1\otimes\cdots\otimes X_n\to\sP X_n\otimes\cdots\otimes\sP X_1$
and $1\le i<n$.
Then
\[
\begin{tikzpicture}[mb]
\mbars{-0.5}{5.1}{2.4}
\foreach \x in {0,1.7,2.9,4.6} \draw[mstr] (\x,0) -- (\x,2.4);
\node[mcpn, minimum width=45mm] at (2.3,1.2) {$f$};
\blab{0}{X_1}\blab{0.85}{\cdots}\blab{1.7}{X_i}\blab{2.9}{X_{i+1}}\blab{3.75}{\cdots}\blab{4.6}{X_n}
\tlab{0}{2.4}{\sP X_n}\tlab{0.85}{2.4}{\cdots}\tlab{1.7}{2.4}{\sP X_{i+1}}\tlab{2.9}{2.4}{\sP X_i}\tlab{3.75}{2.4}{\cdots}\tlab{4.6}{2.4}{\sP X_1}
\end{tikzpicture}
\;=\;
\begin{tikzpicture}[mb]
\mbars{-0.5}{5.1}{3.2}
\draw[mstr] (0,0) -- (0,3.2); \draw[mstr] (4.6,0) -- (4.6,3.2);
\sstr{1.7}{0}{2.9}{1.1}
\sstr[mover]{2.9}{0}{1.7}{1.1}
\sstr{1.7}{2.1}{2.9}{3.2}
\sstr[mover]{2.9}{2.1}{1.7}{3.2}
\draw[mstr] (1.7,1.1) -- (1.7,2.1); \draw[mstr] (2.9,1.1) -- (2.9,2.1);
\node[mcpn, minimum width=45mm] at (2.3,1.6) {$f$};
\blab{0}{X_1}\blab{0.85}{\cdots}\blab{1.7}{X_{i+1}}\blab{2.9}{X_i}\blab{3.75}{\cdots}\blab{4.6}{X_n}
\tlab{0}{3.2}{\sP X_n}\tlab{0.85}{3.2}{\cdots}\tlab{1.7}{3.2}{\sP X_i}\tlab{2.9}{3.2}{\sP X_{i+1}}\tlab{3.75}{3.2}{\cdots}\tlab{4.6}{3.2}{\sP X_1}
\end{tikzpicture}
\;.
\]
\end{lemma}

\begin{proof}
We draw the case $n=2$; for $n>2$ the other strands are spectators.
\[
\begin{tikzpicture}[mb]
\mbars{-0.5}{1.3}{2}
\draw[mstr] (0,0) -- (0,2); \draw[mstr] (0.8,0) -- (0.8,2);
\node[mcpn, minimum width=14.5mm] at (0.4,1) {$f$};
\blab{0}{X}\blab{0.8}{Y}\tlab{0}{2}{\sP Y}\tlab{0.8}{2}{\sP X}
\end{tikzpicture}
\;=\;
\begin{tikzpicture}[mb]
\mbars{-0.5}{1.3}{2.8}
\sstr{0.8}{0}{0}{0.8}
\sstr[mover]{0}{0}{0.8}{0.8}
\sstr{0}{0.8}{0.8}{1.6}
\sstr[mover]{0.8}{0.8}{0}{1.6}
\draw[mstr] (0,1.6) -- (0,2.8); \draw[mstr] (0.8,1.6) -- (0.8,2.8);
\node[mcpn, minimum width=14.5mm] at (0.4,2.2) {$f$};
\blab{0}{X}\blab{0.8}{Y}\tlab{0}{2.8}{\sP Y}\tlab{0.8}{2.8}{\sP X}
\end{tikzpicture}
\;\overset{\text{Lem.~\ref{lem:trace-property}}}{=}\;
\begin{tikzpicture}[mb]
\mbars{-0.5}{1.3}{3.4}
\sstr{0}{0}{0.8}{0.95}
\sstr[mover]{0.8}{0}{0}{0.95}
\draw[mstr] (0,0.95) -- (0,2.05); \draw[mstr] (0.8,0.95) -- (0.8,2.05);
\sstr{0}{2.05}{0.8}{2.95}
\sstr[mover]{0.8}{2.05}{0}{2.95}
\draw[mstr] (0,2.95) -- (0,3.4); \draw[mstr] (0.8,2.95) -- (0.8,3.4);
\node[mcpn, minimum width=14.5mm] at (0.4,1.4) {$f$};
\ubox{-0.35}{1.95}{1.15}{3.05}
\blab{0}{Y}\blab{0.8}{X}\tlab{0}{3.4}{\sP X}\tlab{0.8}{3.4}{\sP Y}
\end{tikzpicture}
\;\overset{\text{\ref{item:B4}}}{=}\;
\begin{tikzpicture}[mb]
\mbars{-0.5}{1.3}{2.8}
\sstr{0}{0}{0.8}{0.95}
\sstr[mover]{0.8}{0}{0}{0.95}
\sstr{0}{1.85}{0.8}{2.8}
\sstr[mover]{0.8}{1.85}{0}{2.8}
\draw[mstr] (0,0.95) -- (0,1.85); \draw[mstr] (0.8,0.95) -- (0.8,1.85);
\node[mcpn, minimum width=14.5mm] at (0.4,1.4) {$f$};
\blab{0}{Y}\blab{0.8}{X}\tlab{0}{2.8}{\sP X}\tlab{0.8}{2.8}{\sP Y}
\end{tikzpicture}\,.
\qedpic
\]
\end{proof}

\begin{remark}
  Here and in the following, when we have an equation
  with graphical calculus or M\"obius diagrams,
  we will often mark the equality sign with a reference to the proposition,
  the lemma, or the rule which justifies it.
  An unmarked equality usually involves only the basic properties of the 
  graphical calculus, which are general for any braided tensor category. 
\end{remark}


\begin{definition}\label{def:product}
For morphisms $f\colon X\to\sP X$ and $g\colon Y\to\sP Y$ put
\[
f\cdot g:=U_{X,Y}\circ c_{\sP X,\sP Y}\circ(f\otimes g)\colon
X\otimes Y\to\sP(X\otimes Y).
\]
In pictures, this is \[
f\cdot g=U_{X,Y}\circ
\begin{tikzpicture}[mb]
\draw[mstr] (0,0) -- (0,0.9); \draw[mstr] (0.9,0) -- (0.9,0.9);
\sstr{0.9}{0.9}{0}{2.2}
\sstr[mover]{0}{0.9}{0.9}{2.2}
\node[mcpn] at (0,0.5) {$f$}; \node[mcpn] at (0.9,0.5) {$g$};
\blab{0}{X}\blab{0.9}{Y}\tlab{0}{2.2}{\sP Y}\tlab{0.9}{2.2}{\sP X}
\end{tikzpicture}\ .
\] We will often omit the $U_{X,Y}$ in the picture.
\end{definition}

\begin{lemma}\label{lem:product-unit}
$U_0: \one\to \sP \one$ is the unit of this product.
\end{lemma}
\begin{proof}
We have
$U_0\cdot f=U_{\one,X}(f\otimes U_0)c_{\one,X}=f$ and
$f\cdot U_0=U_{X,\one}(U_0\otimes f)=f$ by the naturality of $c$,
$c_{\one,-}=c_{-,\one}=\id$ and \eqref{eq:U-unit}.
\end{proof}

\begin{lemma}\label{lem:product-assoc}
For morphisms $f\colon X\to\sP X$, $g\colon Y\to\sP Y$ and $h\colon Z\to\sP Z$,
\[
(f\cdot g)\cdot h=U_{X,Y,Z}\circ\tau\circ(f\otimes g\otimes h)=f\cdot(g\cdot h),
\]
where $\tau$ is the positive half twist of $\sP X\otimes\sP Y\otimes\sP Z$,
in which every strand passes in front of the strands to its right.
\end{lemma}

\begin{proof}
For a product of three factors, $U_{X,Y}$ slides up through the
crossing with $\sP Z$ (naturality of the braiding) and merges with
$U_{X\otimes Y,Z}$ into the iterated $U_{X,Y,Z}$ of
Definition~\ref{def:Pbox}.
Therefore, we have
\[
(f\cdot g)\cdot h=U_{X,Y,Z}\circ
\begin{tikzpicture}[mb]
\foreach \x in {0,0.8,1.6} \draw[mstr] (\x,0) -- (\x,0.75);
\node[mcpn] at (0,0.42) {$f$}; \node[mcpn] at (0.8,0.42) {$g$}; \node[mcpn] at (1.6,0.42) {$h$};
\draw[mstr] (1.6,0.75) -- (1.6,1.35);
\sstr{0.8}{0.75}{0}{1.35}
\sstr[mover]{0}{0.75}{0.8}{1.35}
\draw[mstr] (0,1.35) -- (0,1.95);
\sstr{1.6}{1.35}{0.8}{1.95}
\sstr[mover]{0.8}{1.35}{1.6}{1.95}
\draw[mstr] (1.6,1.95) -- (1.6,2.55);
\sstr{0.8}{1.95}{0}{2.55}
\sstr[mover]{0}{1.95}{0.8}{2.55}
\blab{0}{X}\blab{0.8}{Y}\blab{1.6}{Z}
\tlab{0}{2.55}{\sP Z}\tlab{0.8}{2.55}{\sP Y}\tlab{1.6}{2.55}{\sP X}
\end{tikzpicture}
=U_{X,Y,Z}\circ
\begin{tikzpicture}[mb]
\foreach \x in {0,0.8,1.6} \draw[mstr] (\x,0) -- (\x,0.75);
\node[mcpn] at (0,0.42) {$f$}; \node[mcpn] at (0.8,0.42) {$g$}; \node[mcpn] at (1.6,0.42) {$h$};
\draw[mstr] (0,0.75) -- (0,1.35);
\sstr{1.6}{0.75}{0.8}{1.35}
\sstr[mover]{0.8}{0.75}{1.6}{1.35}
\draw[mstr] (1.6,1.35) -- (1.6,1.95);
\sstr{0.8}{1.35}{0}{1.95}
\sstr[mover]{0}{1.35}{0.8}{1.95}
\draw[mstr] (0,1.95) -- (0,2.55);
\sstr{1.6}{1.95}{0.8}{2.55}
\sstr[mover]{0.8}{1.95}{1.6}{2.55}
\blab{0}{X}\blab{0.8}{Y}\blab{1.6}{Z}
\tlab{0}{2.55}{\sP Z}\tlab{0.8}{2.55}{\sP Y}\tlab{1.6}{2.55}{\sP X}
\end{tikzpicture}
=f\cdot(g\cdot h).
\qedpic
\]
\end{proof}

\begin{proposition}\label{prop:product}
The rule
\[
[X,f]\cdot[Y,g]:=\bigl[X\otimes Y,\;f\cdot g\bigr],
\] that is, \[
\mone{f}{X}{\sP X}\;\cdot\;\mone{g}{Y}{\sP Y}
\;:=\;
\begin{tikzpicture}[mb]
\mbars{-0.5}{1.4}{2.4}
\draw[mstr] (0,0) -- (0,0.9); \draw[mstr] (0.9,0) -- (0.9,0.9);
\sstr{0.9}{0.9}{0}{2.4}
\sstr[mover]{0}{0.9}{0.9}{2.4}
\node[mcpn] at (0,0.6) {$f$}; \node[mcpn] at (0.9,0.6) {$g$};
\blab{0}{X}\blab{0.9}{Y}\tlab{0}{2.4}{\sP Y}\tlab{0.9}{2.4}{\sP X}
\end{tikzpicture}\; ,
\]
gives a well-defined associative product on $\fX$
with unit
$[\one,U_0]$.
\end{proposition}

\begin{proof}
Relation (R1) is respected by the naturality of $U$ and $c$ with
respect to the distributivity isomorphisms.
Relation (R2) is respected in the first
factor: for $a\in\End(X)$,
\[
[X,\sP(a)f]\cdot[Y,g]
=
\begin{tikzpicture}[mb]
\mbars{-0.55}{1.45}{3.1}
\draw[mstr] (0,0) -- (0,1.65); \draw[mstr] (0.9,0) -- (0.9,1.65);
\node[mcpn] at (0,0.5) {$f$}; \node[mcpn] at (0.9,0.5) {$g$};
\node[mcpn] at (0,1.25) {$\sP a$};
\sstr{0.9}{1.65}{0}{3.1}
\sstr[mover]{0}{1.65}{0.9}{3.1}
\blab{0}{X}\blab{0.9}{Y}\tlab{0}{3.1}{\sP Y}\tlab{0.9}{3.1}{\sP X}
\end{tikzpicture}
=
\begin{tikzpicture}[mb]
\mbars{-0.55}{1.45}{3.1}
\draw[mstr] (0,0) -- (0,0.9); \draw[mstr] (0.9,0) -- (0.9,0.9);
\node[mcpn] at (0,0.5) {$f$}; \node[mcpn] at (0.9,0.5) {$g$};
\sstr{0.9}{0.9}{0}{2.1}
\sstr[mover]{0}{0.9}{0.9}{2.1}
\draw[mstr] (0,2.1) -- (0,3.1); \draw[mstr] (0.9,2.1) -- (0.9,3.1);
\node[mcpn] at (0.9,2.6) {$\sP a$};
\blab{0}{X}\blab{0.9}{Y}\tlab{0}{3.1}{\sP Y}\tlab{0.9}{3.1}{\sP X}
\end{tikzpicture}
\overset{\text{Lem.~\ref{lem:trace-property}}}{=}
\begin{tikzpicture}[mb]
\mbars{-0.55}{1.45}{3.1}
\draw[mstr] (0,0) -- (0,1.65); \draw[mstr] (0.9,0) -- (0.9,1.65);
\node[mcpn] at (0,0.5) {$a$};
\node[mcpn] at (0,1.25) {$f$}; \node[mcpn] at (0.9,1.25) {$g$};
\sstr{0.9}{1.65}{0}{3.1}
\sstr[mover]{0}{1.65}{0.9}{3.1}
\blab{0}{X}\blab{0.9}{Y}\tlab{0}{3.1}{\sP Y}\tlab{0.9}{3.1}{\sP X}
\end{tikzpicture}
=[X,fa]\cdot[Y,g].
\]
\noindent The second factor is treated in the same way.  
The unit property and associativity follow from 
Lemma~\ref{lem:product-unit} and Lemma~\ref{lem:product-assoc}.
\end{proof}


To describe the basis of $\fX$,
we introduce the following notion.

\begin{definition}\label{def:equivariant}
A \emph{$\sP$-structure} on $X\in\cC$ is an isomorphism
$\xi\colon X\to\sP X$ with
\begin{equation}\label{eq:equivariant}
\sP(\xi)\circ\xi=\beta_X ,
\end{equation}
that is,
\[
\begin{tikzpicture}[mb]
\draw[mstr] (0,0) -- (0,3.2);
\node[mcpn] at (0,0.7) {$\xi$};
\node[mcpn] at (0,2.15) {$\xi$};
\ubox{-0.45}{1.6}{0.45}{2.7}
\node[anchor=west,inner sep=1.5pt] at (0,1.25) {$\sP X$};
\blab{0}{X}\tlab{0}{3.2}{\sP^2X}
\end{tikzpicture}
=
\begin{tikzpicture}[mb]
\draw[mstr] (0,0) -- (0,3.2);
\node[mcpn] at (0,1.6) {$\beta_X$};
\blab{0}{X}\tlab{0}{3.2}{\sP^2X}
\end{tikzpicture}\,.
\]
By \eqref{eq:beta-monoidal}, $U_0\colon\one\to\sP\one$ is a $\sP$-structure on $\one$.
\end{definition}

\begin{remark}\label{rem:P-structure-box}
More generally, for
$\hat b\colon A\otimes B\to\sP B\otimes\sP A$ the morphism
$U_{A,B}\hat b$ is a $\sP$-structure on $A\otimes B$ if and only if
\[
\begin{tikzpicture}[mb]
\draw[mstr] (0,0) -- (0,3.2); \draw[mstr] (0.9,0) -- (0.9,3.2);
\node[mcpn, minimum width=11mm] at (0.45,0.7) {$\hat b$};
\node[mcpn, minimum width=11mm] at (0.45,2.15) {$\hat b$};
\ubox{-0.45}{1.6}{1.35}{2.7}
\node[anchor=east,inner sep=1.5pt] at (0,1.25) {$\sP B$};
\node[anchor=west,inner sep=1.5pt] at (0.9,1.25) {$\sP A$};
\blab{0}{A}\blab{0.9}{B}\tlab{0}{3.2}{\sP^2A}\tlab{0.9}{3.2}{\sP^2B}
\end{tikzpicture}
=
\begin{tikzpicture}[mb]
\draw[mstr] (0,0) -- (0,3.2); \draw[mstr] (0.9,0) -- (0.9,3.2);
\node[mcpn] at (0,1.6) {$\beta_A$}; \node[mcpn] at (0.9,1.6) {$\beta_B$};
\blab{0}{A}\blab{0.9}{B}\tlab{0}{3.2}{\sP^2A}\tlab{0.9}{3.2}{\sP^2B}
\end{tikzpicture}\,.
\]
Indeed,
$\sP(U_{A,B}\hat b)\,U_{A,B}\hat b=\sP(U_{A,B})U_{\sP B,\sP A}\,\Ubox{\hat b}\,\hat b$,
and $\sP(U_{A,B})U_{\sP B,\sP A}(\beta_A\otimes\beta_B)=\beta_{A\otimes B}$
by \eqref{eq:beta-monoidal}.
When $\Ubox{\hat b}\hat b$ is invertible, so is $\hat b$, since $A\otimes B$ and
$\sP B\otimes\sP A$ have the same length.
\end{remark}

\begin{remark}\label{rem:equivariantization}
When $\beta$ is coherent, i.e.\ satisfies \eqref{eq:beta-coherence} of
Definition~\ref{def:coherent} below,
the pairs $(X,\xi)$ of an object $X$ and a $\sP$-structure $\xi$
are the objects of the
equivariantization $\cC^{\bZ/2}$ \cite[\S4.15]{EGNO} of $\cC$ with respect
to the $\bZ/2$-action given by the monoidal functor $(\sP,J)$,
$J_{X,Y}:=U_{X,Y}\,c_{\sP X,\sP Y}\colon\sP X\otimes\sP Y\to\sP(X\otimes Y)$,
and the isomorphism $\beta\colon\id\Rightarrow\sP^2$.
The morphisms
$(X_1,\xi_1)\to(X_2,\xi_2)$ of $\cC^{\bZ/2}$ are the fixed points of
$R_{\xi_1,\xi_2}$ of Definition~\ref{def:R}, 
the tensor product is the product of Lemma~\ref{lem:J},
and the dual is the natural dual
of Lemma~\ref{lem:natural-dual}.
Many statements of this section 
follow, in the coherent case, from standard properties of the
equivariantization. 
We chose to prove them directly, so that the paper could be understood
without knowledge of equivariantization, and to avoid the coherence assumption.
\end{remark}

\begin{definition}\label{def:InvP}
We set
\[
\Fix(\sP):=\{X\in\Irr(\cC)\mid X\cong\sP X\}.
\qedhere
\]
\end{definition}
\begin{lemma}\label{lem:xiT}
For $X\in\Fix(\sP)$ there is a $\sP$-structure $\xi_X$ on $X$, unique
up to sign.
If $\xi\colon Z\to\sP Z$ is a $\sP$-structure on an object $Z$ and
$u\colon Z\to X$ is an isomorphism, then $\sP(u)\,\xi\,u^{-1}=\pm\xi_X$.
\end{lemma}

\begin{proof}
Choose any isomorphism $\varphi\colon X\to\sP X$.
Since $X$ is simple,
$\sP(\varphi)\varphi=\lambda\beta_X$ for some $\lambda\in\bC^\times$, and
$\lambda^{-1/2}\varphi$ satisfies \eqref{eq:equivariant}.
Another solution is $c\,\xi_X$
with $c\in\bC^\times$, and $\sP(c\xi_X)\,c\xi_X=c^2\beta_X$ forces
$c=\pm1$.
Finally $\sP(u)\xi u^{-1}\colon X\to\sP X$ satisfies \eqref{eq:equivariant}
by the naturality of $\beta$, so it is $\pm\xi_X$ by the above.
\end{proof}

\begin{definition}\label{def:M}
For each $X\in\Fix(\sP)$ we fix, once and for all, one of the two
$\sP$-structures $\xi_X$ of Lemma~\ref{lem:xiT}, with $\xi_\one:=U_0$, and put
\[
\mob(X):=[X,\xi_X]\in\fX .
\qedhere
\]
\end{definition}

\begin{definition}\label{def:TrX}
Let $X\in\Fix(\sP)$ and $f\colon Z\to\sP Z$.
We define $\operatorname{Tr}_X(Z,f)\in\bC$ by
\[
\operatorname{Tr}_X(Z,f)\;:=\;
\begin{tikzpicture}[mb]
\draw[mstr] (0,0) -- (0,4.75) .. controls (0,5.25) and (0.9,5.25) .. (0.9,4.75)
  -- (0.9,0) .. controls (0.9,-0.5) and (0,-0.5) .. (0,0);
\node[mcpnc] at (0,0.6) {$e$};
\node[mcpn] at (0,1.8) {$f$};
\node[mcpnc] at (0,2.95) {$e$};
\ubox{-0.45}{2.50}{0.45}{3.40}
\node[mcpn] at (0,4.0) {$\xi_X^{-1}$};
\node[anchor=west,inner sep=1pt] at (0,1.2) {$Z$};
\node[anchor=west,inner sep=1pt] at (0,4.52) {$X$};
\end{tikzpicture}\,,
\]
where $e$ is the canonical element of $\Hom(X,Z)\otimes\Hom(Z,X)$
(Section~\ref{subsec:conventions}).
We call it the \emph{closing} of $(Z,f)$ at $X$. 
\end{definition}

\begin{lemma}\label{lem:TrX}
$(Z,f)\mapsto\operatorname{Tr}_X(Z,f)$ vanishes on $\fR$.
Hence it defines a linear map $\operatorname{Tr}_X\colon\fX\to\bC$.
\end{lemma}

\begin{proof}
The
closing is linear on $\fV$.
For relation (R1), the canonical element of $\Hom(X,Z\oplus W)\otimes\Hom(Z\oplus W,X)$
is the sum of those for $Z$ and for $W$, composed with the inclusions and
projections; this gives the sum of the closings of $(Z,f)$ and $(W,g)$.
For relation
(R2), let $a\in\End(Z)$.
Then
\[
\begin{tikzpicture}[mb]
\draw[mstr] (0,0) -- (0,4.5) -- (0,5.65) .. controls (0,6.15) and (0.9,6.15) .. (0.9,5.65)
  -- (0.9,0) .. controls (0.9,-0.5) and (0,-0.5) .. (0,0);
\node[mcpnc] at (0,0.6) {$e$};
\node[mcpn] at (0,1.55) {$f$};
\node[mcpn] at (0,2.65) {$a$};
\ubox{-0.45}{2.20}{0.45}{3.10}
\node[mcpnc] at (0,3.75) {$e$};
\ubox{-0.45}{3.30}{0.45}{4.20}
\node[mcpn] at (0,4.85) {$\xi_X^{-1}$};
\node[anchor=west,inner sep=1pt] at (0,5.4) {$X$};
\end{tikzpicture}
\overset{\text{\ref{item:B1}}}{=}
\begin{tikzpicture}[mb]
\draw[mstr] (0,0) -- (0,4.5) -- (0,5.65) .. controls (0,6.15) and (0.9,6.15) .. (0.9,5.65)
  -- (0.9,0) .. controls (0.9,-0.5) and (0,-0.5) .. (0,0);
\node[mcpnc] at (0,0.6) {$e$};
\node[mcpn] at (0,1.55) {$f$};
\node[mcpn] at (0,2.65) {$a$};
\node[mcpnc] at (0,3.6) {$e$};
\ubox{-0.45}{2.20}{0.45}{4.05}
\node[mcpn] at (0,4.85) {$\xi_X^{-1}$};
\node[anchor=west,inner sep=1pt] at (0,5.4) {$X$};
\end{tikzpicture}
\overset{\text{\ref{item:E1}}}{=}
\begin{tikzpicture}[mb]
\draw[mstr] (0,0) -- (0,4.5) -- (0,5.65) .. controls (0,6.15) and (0.9,6.15) .. (0.9,5.65)
  -- (0.9,0) .. controls (0.9,-0.5) and (0,-0.5) .. (0,0);
\node[mcpnc] at (0,0.6) {$e$};
\node[mcpn] at (0,1.55) {$a$};
\node[mcpn] at (0,2.5) {$f$};
\node[mcpnc] at (0,3.6) {$e$};
\ubox{-0.45}{3.15}{0.45}{4.05}
\node[mcpn] at (0,4.85) {$\xi_X^{-1}$};
\node[anchor=west,inner sep=1pt] at (0,5.4) {$X$};
\end{tikzpicture}\,,
\]
and the two ends are the closings of $(Z,\sP(a)f)$ and $(Z,fa)$.
\end{proof}

\begin{proposition}\label{prop:basis}
The elements $\mob(X)$, $X\in\Fix(\sP)$, form a basis of $\fX$, and the
coefficient of $\mob(X)$ in $u\in\fX$ is $\operatorname{Tr}_X(u)$.
In particular $\dim\fX=|\Fix(\sP)|$.
\end{proposition}

\begin{proof}
\emph{Spanning.}  Let $(Z,f)\in\fV$.
As $\cC$ is semisimple, there are $Z_i\in\Irr(\cC)$ and morphisms
$\iota_i\colon Z_i\to Z$, $\pi_i\colon Z\to Z_i$ with
$\pi_i\iota_j=\delta_{ij}\id_{Z_i}$ and $\sum_i\iota_i\pi_i=\id_Z$.
With $f_{ij}:=\sP(\pi_i)\,f\,\iota_j\colon Z_j\to\sP Z_i$,
\[
\mone{f}{Z}{\sP Z}
=\sum_{i,j}
\begin{tikzpicture}[mb]
\mbars{-0.6}{0.6}{3}
\draw[mstr] (0,0) -- (0,3);
\node[mcpn] at (0,0.6) {$\pi_j$};
\node[mcpn] at (0,1.5) {$f_{ij}$};
\node[mcpn] at (0,2.4) {$\sP\iota_i$};
\blab{0}{Z}\tlab{0}{3}{\sP Z}
\end{tikzpicture}
\overset{\text{Lem.~\ref{lem:trace-property}}}{=}\sum_{i,j}
\begin{tikzpicture}[mb]
\mbars{-0.6}{0.6}{3}
\draw[mstr] (0,0) -- (0,3);
\node[mcpn] at (0,0.6) {$f_{ij}$};
\node[mcpn] at (0,1.5) {$\sP\iota_i$};
\node[mcpn] at (0,2.4) {$\sP\pi_j$};
\blab{0}{Z_j}\tlab{0}{3}{\sP Z_j}
\end{tikzpicture}
=\sum_i\mone{f_{ii}}{Z_i}{\sP Z_i}.
\]
Finally, $f_{ii}=0$ unless $Z_i\in\Fix(\sP)$, and in that case
$[Z_i,f_{ii}]\in\bC \mob(Z_i)$, since $\Hom(Z_i,\sP Z_i)=\bC\xi_{Z_i}$.

\noindent \emph{Linear independence.}  We have
$\operatorname{Tr}_X(\mob(X'))=\delta_{X,X'}$, since $\Hom(X,X')=0$ for $X'\neq X$,
and the canonical element of $\End(X)\otimes\End(X)$ is $d_X^{-1}\id_X\otimes\id_X$.
Hence the $\mob(X)$ are linearly independent, and $\operatorname{Tr}_X(u)$ is
the coefficient of $\mob(X)$ in $u$.
\end{proof}

\begin{proposition}\label{prop:product-formula}
For $X,Y\in\Fix(\sP)$,
\[
\mob(X)\,\mob(Y)=\sum_{Z\in\Fix(\sP)}\;
\begin{tikzpicture}[mb]
\draw[mstr] (0.45,0.2) -- (0.45,0.7);
\draw[mstr] (0,0.7) -- (0,2.0); \draw[mstr] (0.9,0.7) -- (0.9,2.0);
\sstr{0.9}{2.0}{0}{3.2}
\sstr[mover]{0}{2.0}{0.9}{3.2}
\draw[mstr] (0,3.2) -- (0,3.8); \draw[mstr] (0.9,3.2) -- (0.9,3.8);
\draw[mstr] (0.45,3.8) -- (0.45,5.5) .. controls (0.45,6.1) and (1.8,6.1) .. (1.8,5.5)
  -- (1.8,0.2) .. controls (1.8,-0.4) and (0.45,-0.4) .. (0.45,0.2);
\node[mcpnc, minimum width=11mm] at (0.45,0.7) {$e$};
\node[mcpn] at (0,1.55) {$\xi_X$}; \node[mcpn] at (0.9,1.55) {$\xi_Y$};
\node[mcpnc, minimum width=11mm] at (0.45,3.8) {$e$};
\ubox{-0.4}{3.4}{1.3}{4.2}
\node[mcpn] at (0.45,4.8) {$\xi_Z^{-1}$};
\node[anchor=west,inner sep=1pt] at (0.45,5.3) {$Z$};
\end{tikzpicture}
\;\mob(Z),
\]
where the shaded pair is the canonical element of
$\Hom(Z,X\otimes Y)\otimes\Hom(X\otimes Y,Z)$ (Definition~\ref{def:repeated}).
\end{proposition}

\begin{proof}
By Proposition~\ref{prop:basis} the coefficient of $\mob(Z)$ is the closing
$\operatorname{Tr}_Z(X\otimes Y,\xi_X\cdot\xi_Y)$ of
Definition~\ref{def:TrX}.
In its picture, the factor $U_{X,Y}$ of $\xi_X\cdot\xi_Y$
(Definition~\ref{def:product}) and $\sP$ of the upper end of $e$ combine
into the box of Definition~\ref{def:Pbox}.
\end{proof}

In the following corollary we drop the strictness assumption and work in
a skeletal model of $\cC$; this is harmless, since
Proposition~\ref{prop:product-formula} involves no associativity
constraint.

\begin{corollary}\label{cor:skeletal-product}
Suppose that $\cC$ is skeletal.
Then $\sP$ permutes the simple objects, and
$\beta_a\in\End(a)=\bC$ for simple $a$ with $\sP^2a=a$.
For simple $a,b,c$ define the linear maps
\begin{alignat*}{2}
R^{ab}_c&\colon\Hom(c,a\otimes b)\to\Hom(c,b\otimes a),&\qquad v&\mapsto c_{a,b}\,v,\\
u^{ab}_c&\colon\Hom(c,a\otimes b)\to\Hom(\sP c,\sP b\otimes\sP a),&\qquad v&\mapsto\Ubox{v}.
\end{alignat*}
Then $\xi_X\in\End(X)=\bC$ is a square root of $\beta_X$ for
$X\in\Fix(\sP)$, and for $X,Y\in\Fix(\sP)$
\[
\mob(X)\,\mob(Y)=\sum_{Z\in\Fix(\sP)}
\frac{\xi_X\,\xi_Y}{\xi_Z}\,\Tr\bigl((u^{XY}_Z)^{-1}R^{XY}_Z\bigr)\,\mob(Z),
\]
where $\Tr$ is the trace of an endomorphism of $\Hom(Z,X\otimes Y)$.
\end{corollary}

\begin{proof}
For $X\in\Fix(\sP)$ we have $\sP X=X$, and $\sP(\xi_X)=\xi_X$ since $\sP$
is linear, so \eqref{eq:equivariant} reads $\xi_X^2=\beta_X$.
For $q\in\Hom(X\otimes Y,Z)$ and $w\in\Hom(Z,Y\otimes X)$,
\[
\Ubox{q}\,w=\Ubox{q}\,\Ubox{(u^{XY}_Z)^{-1}w}=\Ubox{q\,(u^{XY}_Z)^{-1}w}
=q\,(u^{XY}_Z)^{-1}w
\]
by \ref{item:B1} and \ref{item:B3}, as $q\,(u^{XY}_Z)^{-1}w\in\End(Z)=\bC$.
So, after the scalars $\xi_X,\xi_Y,\xi_Z^{-1}$ are taken out, the picture of
Proposition~\ref{prop:product-formula} is the closed diagram of $e$ and
$\Phi(e)$ with $\Phi=(u^{XY}_Z)^{-1}R^{XY}_Z$, which is $\Tr\Phi$ by
\ref{item:E2}.
\end{proof}

\subsection{Semisimplicity}
\label{subsec:positivity}
Here we prove that $\fX$ is semisimple,
by giving $\fX$  the structure of a $C^*$-algebra.
We start with a few lemmas.

\begin{lemma}\label{lem:natural-dual}
Let $\xi$ be a $\sP$-structure on $Z$.
\begin{enumerate}[label=(\alph*)]
\item There is a unique $\sP$-structure $\xi^\natural\colon Z^\vee\to\sP(Z^\vee)$ with
\begin{equation}\label{eq:natural-dual}
\begin{tikzpicture}[mb]
\draw[mstr] (0,1.75) -- (0,0.6) .. controls (0,0.1) and (1,0.1) .. (1,0.6) -- (1,1.75);
\node[mcpn] at (0,1.35) {$\xi$}; \node[mcpn] at (1,1.35) {$\xi^\natural$};
\sstr{1}{1.75}{0}{2.65}
\sstr[mover]{0}{1.75}{1}{2.65}
\draw[mstr] (0,2.65) -- (0,3.0); \draw[mstr] (1,2.65) -- (1,3.0);
\node[anchor=east,inner sep=1pt] at (0,0.65) {$Z$}; \node[anchor=west,inner sep=1pt] at (1,0.65) {$Z^\vee$};
\tlab{0}{3.0}{\sP(Z^\vee)}\tlab{1}{3.0}{\sP Z}
\end{tikzpicture}
=
\begin{tikzpicture}[mb]
\draw[mstr] (0,3.0) -- (0,0.6) .. controls (0,0.1) and (1,0.1) .. (1,0.6) -- (1,3.0);
\tlab{0}{3.0}{\sP(Z^\vee)}\tlab{1}{3.0}{\sP Z}
\end{tikzpicture}\,,
\qquad
\begin{tikzpicture}[mb]
\draw[mstr] (0,0) -- (0,0.95); \draw[mstr] (1,0) -- (1,0.95);
\node[mcpn] at (0,0.5) {$\xi^\natural$}; \node[mcpn] at (1,0.5) {$\xi$};
\sstr{1}{0.95}{0}{1.85}
\sstr[mover]{0}{0.95}{1}{1.85}
\draw[mstr] (0,1.85) -- (0,2.65) .. controls (0,3.2) and (1,3.2) .. (1,2.65) -- (1,1.85);
\node[anchor=east,inner sep=1pt] at (0,2.4) {$\sP Z$}; \node[anchor=west,inner sep=1pt] at (1,2.4) {$\sP(Z^\vee)$};
\blab{0}{Z^\vee}\blab{1}{Z}
\end{tikzpicture}
=
\begin{tikzpicture}[mb]
\draw[mstr] (0,0) -- (0,2.65) .. controls (0,3.2) and (1,3.2) .. (1,2.65) -- (1,0);
\blab{0}{Z^\vee}\blab{1}{Z}
\end{tikzpicture}\,,
\end{equation}
  which we call  the \emph{natural dual} of $\xi$.
\item $(\xi^\natural)^\natural=\xi$, where $Z^{\vee\vee}=Z$, that is:
\[
\begin{tikzpicture}[mb]
\draw[mstr] (0,1.75) -- (0,0.6) .. controls (0,0.1) and (1,0.1) .. (1,0.6) -- (1,1.75);
\node[mcpn] at (0,1.35) {$\xi^\natural$}; \node[mcpn] at (1,1.35) {$\xi$};
\sstr{1}{1.75}{0}{2.65}
\sstr[mover]{0}{1.75}{1}{2.65}
\draw[mstr] (0,2.65) -- (0,3.0); \draw[mstr] (1,2.65) -- (1,3.0);
\node[anchor=east,inner sep=1pt] at (0,0.65) {$Z^\vee$}; \node[anchor=west,inner sep=1pt] at (1,0.65) {$Z$};
\tlab{0}{3.0}{\sP Z}\tlab{1}{3.0}{\sP(Z^\vee)}
\end{tikzpicture}
=
\begin{tikzpicture}[mb]
\draw[mstr] (0,3.0) -- (0,0.6) .. controls (0,0.1) and (1,0.1) .. (1,0.6) -- (1,3.0);
\tlab{0}{3.0}{\sP Z}\tlab{1}{3.0}{\sP(Z^\vee)}
\end{tikzpicture}\,,
\qquad
\begin{tikzpicture}[mb]
\draw[mstr] (0,0) -- (0,0.95); \draw[mstr] (1,0) -- (1,0.95);
\node[mcpn] at (0,0.5) {$\xi$}; \node[mcpn] at (1,0.5) {$\xi^\natural$};
\sstr{1}{0.95}{0}{1.85}
\sstr[mover]{0}{0.95}{1}{1.85}
\draw[mstr] (0,1.85) -- (0,2.65) .. controls (0,3.2) and (1,3.2) .. (1,2.65) -- (1,1.85);
\node[anchor=east,inner sep=1pt] at (0,2.4) {$\sP(Z^\vee)$}; \node[anchor=west,inner sep=1pt] at (1,2.4) {$\sP Z$};
\blab{0}{Z}\blab{1}{Z^\vee}
\end{tikzpicture}
=
\begin{tikzpicture}[mb]
\draw[mstr] (0,0) -- (0,2.65) .. controls (0,3.2) and (1,3.2) .. (1,2.65) -- (1,0);
\blab{0}{Z}\blab{1}{Z^\vee}
\end{tikzpicture}\,.
\]
\item For $X\in\Fix(\sP)$ we have $X^\vee\in\Fix(\sP)$, and there is a sign $\nu_X$ with
  $\xi_X^\natural=\nu_X\,\xi_{X^\vee}$.  Moreover $\nu_{X^\vee}=\nu_X$.
\end{enumerate}
\end{lemma}

\begin{proof}
\medskip\noindent (a) 

The first identity of \eqref{eq:natural-dual} clearly determines $\xi^\natural$ uniquely as a morphism.
We then have
\[
\begin{tikzpicture}[mb]
\draw[mstr] (0,3.3) -- (0,0.6) .. controls (0,0.1) and (1.2,0.1) .. (1.2,0.6) -- (1.2,3.3);
\node[mcpn] at (0,1.75) {$\beta_Z$};
\node[mcpn] at (1.2,1.1) {$\xi^\natural$};
\node[mcpn] at (1.2,2.4) {$\xi^\natural$};
\ubox{0.75}{1.95}{1.65}{2.85}
\tlab{0}{3.3}{\sP^2Z}\tlab{1.2}{3.3}{\sP^2(Z^\vee)}
\end{tikzpicture}
\overset{\eqref{eq:equivariant}}{=}
\begin{tikzpicture}[mb]
\draw[mstr] (0,3.3) -- (0,0.6) .. controls (0,0.1) and (1.2,0.1) .. (1.2,0.6) -- (1.2,3.3);
\node[mcpn] at (0,1.1) {$\xi$};
\node[mcpn] at (0,2.4) {$\xi$};
\ubox{-0.45}{1.95}{0.45}{2.85}
\node[mcpn] at (1.2,1.1) {$\xi^\natural$};
\node[mcpn] at (1.2,2.4) {$\xi^\natural$};
\ubox{0.75}{1.95}{1.65}{2.85}
\tlab{0}{3.3}{\sP^2Z}\tlab{1.2}{3.3}{\sP^2(Z^\vee)}
\end{tikzpicture}
\overset{\text{\ref{item:B2}}}{\underset{\text{\ref{item:B4}}}{=}}
\begin{tikzpicture}[mb]
\draw[mstr] (0,1.4) -- (0,0.6) .. controls (0,0.1) and (1.2,0.1) .. (1.2,0.6) -- (1.2,1.4);
\node[mcpn] at (0,1.0) {$\xi$}; \node[mcpn] at (1.2,1.0) {$\xi^\natural$};
\sstr{1.2}{1.4}{0}{2.3}
\sstr[mover]{0}{1.4}{1.2}{2.3}
\draw[mstr] (0,2.3) -- (0,3.2); \draw[mstr] (1.2,2.3) -- (1.2,3.2);
\node[mcpn] at (0,2.8) {$\xi^\natural$}; \node[mcpn] at (1.2,2.8) {$\xi$};
\sstr{0}{3.2}{1.2}{4.1}
\sstr[mover]{1.2}{3.2}{0}{4.1}
\draw[mstr] (0,4.1) -- (0,4.5); \draw[mstr] (1.2,4.1) -- (1.2,4.5);
\ubox{-0.5}{2.4}{1.7}{4.25}
\tlab{0}{4.5}{\sP^2Z}\tlab{1.2}{4.5}{\sP^2(Z^\vee)}
\end{tikzpicture}
\overset{\eqref{eq:natural-dual}}{\underset{\text{\ref{item:B6}}}{=}}
\begin{tikzpicture}[mb]
\draw[mstr] (0,2.4) -- (0,0.75) .. controls (0,0.25) and (1.2,0.25) .. (1.2,0.75) -- (1.2,2.4);
\ubox{-0.4}{0.15}{1.6}{1.35}
\node[anchor=west,inner sep=1pt] at (0,0.85) {$Z^\vee$}; \node[anchor=east,inner sep=1pt] at (1.2,0.85) {$Z$};
\node[mcpn] at (0,2.0) {$\xi^\natural$}; \node[mcpn] at (1.2,2.0) {$\xi$};
\sstr{0}{2.4}{1.2}{3.3}
\sstr[mover]{1.2}{2.4}{0}{3.3}
\draw[mstr] (0,3.3) -- (0,3.7); \draw[mstr] (1.2,3.3) -- (1.2,3.7);
\ubox{-0.5}{1.6}{1.7}{3.45}
\tlab{0}{3.7}{\sP^2Z}\tlab{1.2}{3.7}{\sP^2(Z^\vee)}
\end{tikzpicture}
\overset{\eqref{eq:natural-dual}}{\underset{\text{\ref{item:B1},\ref{item:B6}}}{=}}
\begin{tikzpicture}[mb]
\draw[mstr] (0,2.8) -- (0,0.75) .. controls (0,0.25) and (1.2,0.25) .. (1.2,0.75) -- (1.2,2.8);
\ubox{-0.4}{0.15}{1.6}{1.35}
\ubox{-0.75}{0.0}{1.95}{1.95}
\node[anchor=west,inner sep=1pt] at (0,0.85) {$Z$}; \node[anchor=east,inner sep=1pt] at (1.2,0.85) {$Z^\vee$};
\tlab{0}{2.8}{\sP^2Z}\tlab{1.2}{2.8}{\sP^2(Z^\vee)}
\end{tikzpicture}
\overset{\text{\ref{item:B7}}}{=}
\begin{tikzpicture}[mb]
\draw[mstr] (0,2.8) -- (0,0.6) .. controls (0,0.1) and (1.2,0.1) .. (1.2,0.6) -- (1.2,2.8);
\node[mcpn] at (0,1.5) {$\beta_Z$}; \node[mcpn] at (1.2,1.5) {$\beta_{Z^\vee}$};
\tlab{0}{2.8}{\sP^2Z}\tlab{1.2}{2.8}{\sP^2(Z^\vee)}
\end{tikzpicture}.
\]
Canceling $\beta_Z$, we find 
$\sP(\xi^\natural)\xi^\natural=\beta_{Z^\vee}$.
In particular $\xi^\natural$ has a left inverse and $\sP(\xi^\natural)$ a right inverse, so
$\xi^\natural$ is invertible, as $\sP$ is an equivalence.
Therefore $\xi^\natural$ is a $\sP$-structure.

Let us now check the second identity. 
Inserting the zigzag identity on the leg $Z$ of the left hand side of the second identity of \eqref{eq:natural-dual}, we get
\[
\begin{tikzpicture}[mb]
\draw[mstr] (0,0) -- (0,0.95); \draw[mstr] (1,0) -- (1,0.95);
\node[mcpn] at (0,0.5) {$\xi^\natural$}; \node[mcpn] at (1,0.5) {$\xi$};
\sstr{1}{0.95}{0}{1.85}
\sstr[mover]{0}{0.95}{1}{1.85}
\draw[mstr] (0,1.85) -- (0,2.65) .. controls (0,3.2) and (1,3.2) .. (1,2.65) -- (1,1.85);
\blab{0}{Z^\vee}\blab{1}{Z}
\end{tikzpicture}
\;\overset{\text{zigzag}}{=}\;
\begin{tikzpicture}[mb]
\draw[mstr] (0,0) -- (0,2.4);
\draw[mstr] (1,2.4) -- (1,1.1) .. controls (1,0.6) and (2,0.6) .. (2,1.1) -- (2,4.3)
  .. controls (2,4.85) and (3,4.85) .. (3,4.3) -- (3,0);
\node[mcpn] at (0,1.75) {$\xi^\natural$}; \node[mcpn] at (1,1.75) {$\xi$};
\sstr{1}{2.4}{0}{3.2}
\sstr[mover]{0}{2.4}{1}{3.2}
\draw[mstr] (0,3.2) -- (0,3.6) .. controls (0,4.15) and (1,4.15) .. (1,3.6) -- (1,3.2);
\node[anchor=west,inner sep=1pt] at (2,1.6) {$Z^\vee$};
\blab{0}{Z^\vee}\blab{3}{Z}
\end{tikzpicture}
\;\overset{\eqref{eq:natural-dual}}{=}\;
\begin{tikzpicture}[mb]
\draw[mstr] (0,0) -- (0,2.6);
\draw[mstr] (1,1.5) -- (1,1.1) .. controls (1,0.6) and (2,0.6) .. (2,1.1) -- (2,1.5);
\sstr{1}{1.5}{2}{2.3}
\sstr[mover]{2}{1.5}{1}{2.3}
\draw[mstr] (1,2.3) -- (1,2.6);
\draw[mstr] (2,2.3) -- (2,4.5) .. controls (2,5.05) and (3,5.05) .. (3,4.5) -- (3,0);
\node[mcpn] at (0,1.0) {$\xi^\natural$};
\sstr{1}{2.6}{0}{3.4}
\sstr[mover]{0}{2.6}{1}{3.4}
\draw[mstr] (0,3.4) -- (0,3.8) .. controls (0,4.35) and (1,4.35) .. (1,3.8) -- (1,3.4);
\node[mcpn] at (2,3.5) {$(\xi^\natural)^{-1}$};
\node[anchor=west,inner sep=1pt] at (2,1.15) {$\sP Z$};
\node[anchor=east,inner sep=1pt] at (0,1.95) {$\sP(Z^\vee)$};
\blab{0}{Z^\vee}\blab{3}{Z}
\end{tikzpicture}
\;\overset{\text{isotopy}}{=}\;
\begin{tikzpicture}[mb]
\draw[mstr] (0,0) -- (0,3.5) .. controls (0,4.1) and (1,4.1) .. (1,3.5) -- (1,0);
\blab{0}{Z^\vee}\blab{1}{Z}
\end{tikzpicture}\,.
\]

\medskip\noindent (b) With $Z^{\vee\vee}=Z$, the cup $\coev'_Z$ is $\coev_{Z^\vee}$, so the
first identity is the first one of \eqref{eq:natural-dual} for $\xi^\natural$ with $\xi$ in place of
$(\xi^\natural)^\natural$, and by the uniqueness in (a) it is equivalent to
$(\xi^\natural)^\natural=\xi$.
For the ribbon structure, $\coev'_Z=(\id_{Z^\vee}\otimes\theta_Z)\,c_{Z,Z^\vee}\coev_Z$, so
\[
\begin{tikzpicture}[mb]
\draw[mstr] (0,1.75) -- (0,0.6) .. controls (0,0.1) and (1,0.1) .. (1,0.6) -- (1,1.75);
\node[mcpn] at (0,1.35) {$\xi^\natural$}; \node[mcpn] at (1,1.35) {$\xi$};
\sstr{1}{1.75}{0}{2.65}
\sstr[mover]{0}{1.75}{1}{2.65}
\draw[mstr] (0,2.65) -- (0,3.0); \draw[mstr] (1,2.65) -- (1,3.0);
\node[anchor=east,inner sep=1pt] at (0,0.65) {$Z^\vee$}; \node[anchor=west,inner sep=1pt] at (1,0.65) {$Z$};
\tlab{0}{3.0}{\sP Z}\tlab{1}{3.0}{\sP(Z^\vee)}
\end{tikzpicture}
\;\overset{\text{ribbon}}{=}\;
\begin{tikzpicture}[mb]
\draw[mstr] (0,1.0) -- (0,0.6) .. controls (0,0.1) and (1,0.1) .. (1,0.6) -- (1,1.0);
\sstr{1}{1.0}{0}{1.8}
\sstr[mover]{0}{1.0}{1}{1.8}
\draw[mstr] (0,1.8) -- (0,3.6); \draw[mstr] (1,1.8) -- (1,3.6);
\node[mcpn] at (1,2.25) {$\theta_Z$};
\node[mcpn] at (0,3.1) {$\xi^\natural$}; \node[mcpn] at (1,3.1) {$\xi$};
\sstr{1}{3.6}{0}{4.5}
\sstr[mover]{0}{3.6}{1}{4.5}
\draw[mstr] (0,4.5) -- (0,4.8); \draw[mstr] (1,4.5) -- (1,4.8);
\node[anchor=east,inner sep=1pt] at (0,0.65) {$Z$}; \node[anchor=west,inner sep=1pt] at (1,0.65) {$Z^\vee$};
\tlab{0}{4.8}{\sP Z}\tlab{1}{4.8}{\sP(Z^\vee)}
\end{tikzpicture}
\;=\;
\begin{tikzpicture}[mb]
\draw[mstr] (0,1.6) -- (0,0.6) .. controls (0,0.1) and (1,0.1) .. (1,0.6) -- (1,1.6);
\node[mcpn] at (0,1.2) {$\xi$}; \node[mcpn] at (1,1.2) {$\xi^\natural$};
\sstr{1}{1.6}{0}{2.4}
\sstr[mover]{0}{1.6}{1}{2.4}
\draw[mstr] (0,2.4) -- (0,3.6); \draw[mstr] (1,2.4) -- (1,3.6);
\node[mcpn] at (1,3.0) {$\theta_{\sP Z}$};
\sstr{1}{3.6}{0}{4.5}
\sstr[mover]{0}{3.6}{1}{4.5}
\draw[mstr] (0,4.5) -- (0,4.8); \draw[mstr] (1,4.5) -- (1,4.8);
\node[anchor=east,inner sep=1pt] at (0,0.65) {$Z$}; \node[anchor=west,inner sep=1pt] at (1,0.65) {$Z^\vee$};
\tlab{0}{4.8}{\sP Z}\tlab{1}{4.8}{\sP(Z^\vee)}
\end{tikzpicture}
\;\overset{\eqref{eq:natural-dual}}{=}\;
\begin{tikzpicture}[mb]
\draw[mstr] (0,1.9) -- (0,0.6) .. controls (0,0.1) and (1,0.1) .. (1,0.6) -- (1,1.9);
\node[mcpn] at (1,1.3) {$\theta_{\sP Z}$};
\sstr{1}{1.9}{0}{2.8}
\sstr[mover]{0}{1.9}{1}{2.8}
\draw[mstr] (0,2.8) -- (0,3.1); \draw[mstr] (1,2.8) -- (1,3.1);
\node[anchor=east,inner sep=1pt] at (0,0.65) {$\sP(Z^\vee)$}; \node[anchor=west,inner sep=1pt] at (1,0.65) {$\sP Z$};
\tlab{0}{3.1}{\sP Z}\tlab{1}{3.1}{\sP(Z^\vee)}
\end{tikzpicture}
\;\overset{\text{ribbon}}{=}\;
\begin{tikzpicture}[mb]
\draw[mstr] (0,3.0) -- (0,0.6) .. controls (0,0.1) and (1,0.1) .. (1,0.6) -- (1,3.0);
\tlab{0}{3.0}{\sP Z}\tlab{1}{3.0}{\sP(Z^\vee)}
\end{tikzpicture}\,,
\]
where in the second step we used the naturality of $\theta$.
The second identity is similar.

\medskip\noindent (c) By \ref{item:B5}, $\sP(X^\vee)\cong(\sP X)^\vee\cong X^\vee$, so $X^\vee\in\Fix(\sP)$.
Now, both $\xi^\natural_X$ and $\xi_{X^\vee}$ are $\sP$-structures on
$X^\vee$, so $\xi_X^\natural=\lambda\,\xi_{X^\vee}$ with $\lambda=\pm1$ by Lemma~\ref{lem:xiT}.
Put $\nu_X:=\lambda$.
Since $(-\xi)\cdot(-\xi^\natural)=\xi\cdot\xi^\natural$, we have
$(-\xi)^\natural=-\xi^\natural$ by the uniqueness in (a).
Hence, by (b),
\[
\xi_X=(\xi_X^\natural)^\natural=(\nu_X\,\xi_{X^\vee})^\natural
=\nu_X\,\xi_{X^\vee}^\natural=\nu_X\nu_{X^\vee}\,\xi_X ,
\]
that is, $\nu_{X^\vee}=\nu_X$.
\end{proof}

\begin{lemma}\label{lem:J}
For $\sP$-structures $\xi$ and $\xi'$ on $X$ and $Y$, the product
$\xi\cdot\xi'$ (Definition~\ref{def:product}) is a $\sP$-structure on
$X\otimes Y$.
\end{lemma}

\begin{proof}
By Remark~\ref{rem:P-structure-box} it suffices to show that the
picture of $\xi\cdot\xi'$, followed by its box, is
$\beta_X\otimes\beta_Y$:
\[
\begin{tikzpicture}[mb]
\draw[mstr] (0,0) -- (0,1.0); \draw[mstr] (1.2,0) -- (1.2,1.0);
\node[mcpn] at (0,0.55) {$\xi$}; \node[mcpn] at (1.2,0.55) {$\xi'$};
\sstr{1.2}{1.0}{0}{1.9}
\sstr[mover]{0}{1.0}{1.2}{1.9}
\draw[mstr] (0,1.9) -- (0,2.85); \draw[mstr] (1.2,1.9) -- (1.2,2.85);
\node[mcpn] at (0,2.45) {$\xi'$}; \node[mcpn] at (1.2,2.45) {$\xi$};
\sstr{0}{2.85}{1.2}{3.75}
\sstr[mover]{1.2}{2.85}{0}{3.75}
\draw[mstr] (0,3.75) -- (0,4.4); \draw[mstr] (1.2,3.75) -- (1.2,4.4);
\ubox{-0.5}{2.05}{1.7}{3.95}
\blab{0}{X}\blab{1.2}{Y}\tlab{0}{4.4}{\sP^2X}\tlab{1.2}{4.4}{\sP^2Y}
\end{tikzpicture}
\overset{\text{\ref{item:B1},\ref{item:B2}}}{\underset{\text{\ref{item:B4}}}{=}}
\begin{tikzpicture}[mb]
\draw[mstr] (0,0) -- (0,1.0); \draw[mstr] (1.2,0) -- (1.2,1.0);
\node[mcpn] at (0,0.55) {$\xi$}; \node[mcpn] at (1.2,0.55) {$\xi'$};
\sstr{1.2}{1.0}{0}{1.9}
\sstr[mover]{0}{1.0}{1.2}{1.9}
\draw[mstr] (0,1.9) -- (0,3.3); \draw[mstr] (1.2,1.9) -- (1.2,3.3);
\node[mcpn] at (0,2.45) {$\xi'$}; \node[mcpn] at (1.2,2.45) {$\xi$};
\ubox{-0.45}{2.0}{0.45}{2.9}
\ubox{0.75}{2.0}{1.65}{2.9}
\sstr{0}{3.3}{1.2}{4.1}
\sstr[mover]{1.2}{3.3}{0}{4.1}
\draw[mstr] (0,4.1) -- (0,4.4); \draw[mstr] (1.2,4.1) -- (1.2,4.4);
\blab{0}{X}\blab{1.2}{Y}\tlab{0}{4.4}{\sP^2X}\tlab{1.2}{4.4}{\sP^2Y}
\end{tikzpicture}
=
\begin{tikzpicture}[mb]
\draw[mstr] (0,0) -- (0,3.2); \draw[mstr] (1.2,0) -- (1.2,3.2);
\node[mcpn] at (0,0.6) {$\xi$}; \node[mcpn] at (1.2,0.6) {$\xi'$};
\node[mcpn] at (0,1.9) {$\xi$}; \node[mcpn] at (1.2,1.9) {$\xi'$};
\ubox{-0.45}{1.45}{0.45}{2.35}
\ubox{0.75}{1.45}{1.65}{2.35}
\blab{0}{X}\blab{1.2}{Y}\tlab{0}{3.2}{\sP^2X}\tlab{1.2}{3.2}{\sP^2Y}
\end{tikzpicture}
\overset{\eqref{eq:equivariant}}{=}
\begin{tikzpicture}[mb]
\draw[mstr] (0,0) -- (0,2.4); \draw[mstr] (1.2,0) -- (1.2,2.4);
\node[mcpn] at (0,1.2) {$\beta_X$}; \node[mcpn] at (1.2,1.2) {$\beta_Y$};
\blab{0}{X}\blab{1.2}{Y}\tlab{0}{2.4}{\sP^2X}\tlab{1.2}{2.4}{\sP^2Y}
\end{tikzpicture}\,.
\qedpic
\]
\end{proof}

\begin{lemma}\label{lem:natural-dual-product}
For $\sP$-structures $\xi$, $\xi'$ on $Z$, $Z'$, we have
$(\xi\cdot\xi')^\natural=\xi'^\natural\cdot\xi^\natural$, where
$(Z\otimes Z')^\vee=Z'^\vee\otimes Z^\vee$.
\end{lemma}

\begin{proof}
The first picture below is
\[
\bigl((\xi\cdot\xi')\cdot(\xi'^\natural\cdot\xi^\natural)\bigr)\circ\coev_{Z\otimes Z'},\qquad
\coev_{Z\otimes Z'}=(\id_Z\otimes\coev_{Z'}\otimes\id_{Z^\vee})\coev_Z ,
\]
with the products of $\sP$-structures drawn as after Definition~\ref{def:product}.
We have
\[
\begin{tikzpicture}[mb,x=6.8mm]
\draw[mstr] (0,1.85) -- (0,0.95) .. controls (0,0.1) and (3.9,0.1) .. (3.9,0.95) -- (3.9,1.85);
\draw[mstr] (1.45,1.85) -- (1.45,1.1) .. controls (1.45,0.6) and (2.7,0.6) .. (2.7,1.1) -- (2.7,1.85);
\node[mcpn] at (0,1.45) {$\xi$}; \node[mcpn] at (1.45,1.45) {$\xi'$};
\node[mcpn] at (2.7,1.45) {$\xi'^\natural$}; \node[mcpn] at (3.9,1.45) {$\xi^\natural$};
\sstr{1.45}{1.85}{0}{2.75}
\sstr[mover]{0}{1.85}{1.45}{2.75}
\sstr{3.9}{1.85}{2.7}{2.75}
\sstr[mover]{2.7}{1.85}{3.9}{2.75}
\sstr{2.7}{2.75}{0}{4.3}
\sstr{3.9}{2.75}{1.45}{4.3}
\sstr[mover]{0}{2.75}{2.7}{4.3}
\sstr[mover]{1.45}{2.75}{3.9}{4.3}
\foreach \x in {0,1.45,2.7,3.9} \draw[mstr] (\x,4.3) -- (\x,5.0);
\tlab{0}{5.0}{\sP(Z^\vee)}\tlab{1.45}{5.0}{\sP(Z'^\vee)}\tlab{2.7}{5.0}{\sP Z'}\tlab{3.9}{5.0}{\sP Z}
\end{tikzpicture}
\overset{\text{isotopy}}{=}
\begin{tikzpicture}[mb,x=6.8mm]
\draw[mstr] (0,1.85) -- (0,0.95) .. controls (0,0.1) and (3.9,0.1) .. (3.9,0.95) -- (3.9,1.85);
\node[mcpn] at (0,1.45) {$\xi$}; \node[mcpn] at (3.9,1.45) {$\xi^\natural$};
\sstr{3.9}{1.85}{0}{2.9}
\sstr[mover]{0}{1.85}{3.9}{2.9}
\draw[mstr] (0,2.9) -- (0,5.0); \draw[mstr] (3.9,2.9) -- (3.9,5.0);
\draw[mstr] (1.45,4.25) -- (1.45,3.55) .. controls (1.45,3.05) and (2.7,3.05) .. (2.7,3.55) -- (2.7,4.25);
\node[mcpn] at (1.45,3.85) {$\xi'$}; \node[mcpn] at (2.7,3.85) {$\xi'^\natural$};
\sstr{2.7}{4.25}{1.45}{5.0}
\sstr[mover]{1.45}{4.25}{2.7}{5.0}
\tlab{0}{5.0}{\sP(Z^\vee)}\tlab{1.45}{5.0}{\sP(Z'^\vee)}\tlab{2.7}{5.0}{\sP Z'}\tlab{3.9}{5.0}{\sP Z}
\end{tikzpicture}
\overset{\eqref{eq:natural-dual}}{=}
\begin{tikzpicture}[mb,x=6.8mm]
\draw[mstr] (0,1.85) -- (0,0.95) .. controls (0,0.1) and (3.9,0.1) .. (3.9,0.95) -- (3.9,1.85);
\node[mcpn] at (0,1.45) {$\xi$}; \node[mcpn] at (3.9,1.45) {$\xi^\natural$};
\sstr{3.9}{1.85}{0}{2.9}
\sstr[mover]{0}{1.85}{3.9}{2.9}
\draw[mstr] (0,2.9) -- (0,5.0); \draw[mstr] (3.9,2.9) -- (3.9,5.0);
\draw[mstr] (1.45,5.0) -- (1.45,3.9) .. controls (1.45,3.4) and (2.7,3.4) .. (2.7,3.9) -- (2.7,5.0);
\tlab{0}{5.0}{\sP(Z^\vee)}\tlab{1.45}{5.0}{\sP(Z'^\vee)}\tlab{2.7}{5.0}{\sP Z'}\tlab{3.9}{5.0}{\sP Z}
\end{tikzpicture}
\overset{\eqref{eq:natural-dual}}{=}
\begin{tikzpicture}[mb,x=6.8mm]
\draw[mstr] (0,5.0) -- (0,0.95) .. controls (0,0.1) and (3.9,0.1) .. (3.9,0.95) -- (3.9,5.0);
\draw[mstr] (1.45,5.0) -- (1.45,3.9) .. controls (1.45,3.4) and (2.7,3.4) .. (2.7,3.9) -- (2.7,5.0);
\tlab{0}{5.0}{\sP(Z^\vee)}\tlab{1.45}{5.0}{\sP(Z'^\vee)}\tlab{2.7}{5.0}{\sP Z'}\tlab{3.9}{5.0}{\sP Z}
\end{tikzpicture}.
\]
The last picture is the cup of $\sP(Z\otimes Z')$, so
$\xi'^\natural\cdot\xi^\natural$ satisfies the defining identity of
$(\xi\cdot\xi')^\natural$, and the uniqueness in Lemma~\ref{lem:natural-dual}(a) gives the claim.
\end{proof}

\begin{definition}\label{def:R}
For isomorphisms $\varphi_1\colon X_1\to\sP X_1$ and $\varphi_2\colon X_2\to\sP X_2$ and a
morphism $\Phi\colon X_1\to X_2$ put
\[
R_{\varphi_1,\varphi_2}(\Phi):=\varphi_2^{-1}\circ\sP(\Phi)\circ\varphi_1
=
\begin{tikzpicture}[mb]
\draw[mstr] (0,0) -- (0,5.6);
\node[mcpn] at (0,0.6) {$\varphi_1$};
\ubox{-0.55}{1.7}{0.8}{3.9}
\node[mcpn] at (0,2.8) {$\Phi$};
\node[mcpn] at (0,4.9) {$\varphi_2^{-1}$};
\node[anchor=west,inner sep=1pt] at (0,1.35) {$\sP X_1$};
\node[anchor=west,inner sep=1pt] at (0,2.05) {$X_1$};
\node[anchor=west,inner sep=1pt] at (0,3.55) {$X_2$};
\node[anchor=west,inner sep=1pt] at (0,4.25) {$\sP X_2$};
\blab{0}{X_1}\tlab{0}{5.6}{X_2}
\end{tikzpicture}
\colon X_1\to X_2 .
\]
The fixed points of $R_{\varphi_1,\varphi_2}$ are the morphisms $\Phi$ with
$\sP(\Phi)\,\varphi_1=\varphi_2\,\Phi$.
\end{definition}

\begin{lemma}\label{lem:R}
Let $\xi_i$ be $\sP$-structures on $X_i$ $(i=1,2,3)$ and $\xi'_i$
$\sP$-structures on $Y_i$ $(i=1,2)$.
\begin{enumerate}[label=(\alph*)]
\item $R_{\xi_1,\xi_2}$ is a linear involution of $\Hom(X_1,X_2)$, and
  $R_{\xi_1,\xi_1}(\id)=\id$.
\item $R_{\xi_1,\xi_3}(\Psi\Phi)=R_{\xi_2,\xi_3}(\Psi)\,R_{\xi_1,\xi_2}(\Phi)$
  for $\Phi\colon X_1\to X_2$ and $\Psi\colon X_2\to X_3$.
\item $R_{\xi_1\cdot\xi'_1,\,\xi_2\cdot\xi'_2}(\Phi\otimes\Phi')
  =R_{\xi_1,\xi_2}(\Phi)\otimes R_{\xi'_1,\xi'_2}(\Phi')$
  for $\Phi\colon X_1\to X_2$ and $\Phi'\colon Y_1\to Y_2$.
\end{enumerate}
\end{lemma}

\begin{proof}
(a) By Definition~\ref{def:R}, applied twice,
\[
R_{\xi_1,\xi_2}\bigl(R_{\xi_1,\xi_2}(\Phi)\bigr)
=
\begin{tikzpicture}[mb]
\draw[mstr] (0,0) -- (0,5.2);
\node[mcpn] at (0,0.5) {$\xi_1$};
\node[mcpn] at (0,1.5) {$\xi_1$};
\node[mcpn] at (0,2.55) {$\Phi$};
\ubox{-0.45}{2.1}{0.45}{3.0}
\node[mcpn] at (0,3.6) {$\xi_2^{-1}$};
\ubox{-0.8}{1.0}{0.8}{4.1}
\node[mcpn] at (0,4.65) {$\xi_2^{-1}$};
\blab{0}{X_1}\tlab{0}{5.2}{X_2}
\end{tikzpicture}
\overset{\text{\ref{item:B1}}}{=}
\begin{tikzpicture}[mb]
\draw[mstr] (0,0) -- (0,5.8);
\node[mcpn] at (0,0.5) {$\xi_1$};
\node[mcpn] at (0,1.5) {$\xi_1$};
\ubox{-0.45}{1.05}{0.45}{1.95}
\node[mcpn] at (0,2.9) {$\Phi$};
\ubox{-0.45}{2.45}{0.45}{3.35}
\ubox{-0.8}{2.2}{0.8}{3.6}
\node[mcpn] at (0,4.3) {$\xi_2^{-1}$};
\ubox{-0.55}{3.85}{0.55}{4.75}
\node[mcpn] at (0,5.25) {$\xi_2^{-1}$};
\blab{0}{X_1}\tlab{0}{5.8}{X_2}
\end{tikzpicture}
\overset{\eqref{eq:equivariant}}{=}
\begin{tikzpicture}[mb]
\draw[mstr] (0,0) -- (0,3.7);
\node[mcpn] at (0,0.5) {$\beta_{X_1}$};
\node[mcpn] at (0,1.85) {$\Phi$};
\ubox{-0.45}{1.4}{0.45}{2.3}
\ubox{-0.8}{1.15}{0.8}{2.55}
\node[mcpn] at (0,3.15) {$\beta_{X_2}^{-1}$};
\blab{0}{X_1}\tlab{0}{3.7}{X_2}
\end{tikzpicture}
\overset{\text{\ref{item:B7}}}{=}
\begin{tikzpicture}[mb]
\draw[mstr] (0,0) -- (0,2.2);
\node[mcpn] at (0,1.1) {$\Phi$};
\blab{0}{X_1}\tlab{0}{2.2}{X_2}
\end{tikzpicture}\ .
\]

\medskip\noindent (b) is clear, and (c) follows from
$(\xi_2\cdot\xi'_2)^{-1}\sP(\Phi\otimes\Phi')(\xi_1\cdot\xi'_1)
=(\xi_2\otimes\xi'_2)^{-1}\,c^{-1}\,U^{-1}\sP(\Phi\otimes\Phi')U\,c\,(\xi_1\otimes\xi'_1)$,
since $U^{-1}\sP(\Phi\otimes\Phi')U=\sP\Phi'\otimes\sP\Phi$ (naturality of
$U$) and $c^{-1}(\sP\Phi'\otimes\sP\Phi)c=\sP\Phi\otimes\sP\Phi'$ (naturality
of $c$).
\end{proof}

\begin{lemma}\label{lem:natural-dual-R}
For $\sP$-structures $\xi_1$, $\xi_2$ on $Z_1$, $Z_2$ and $f\colon Z_1\to Z_2$, we have
$R_{\xi_1,\xi_2}(f)^\vee=R_{\xi_2^\natural,\xi_1^\natural}(f^\vee)$.
\end{lemma}

\begin{proof}
By Lemma~\ref{lem:R}(b), (c) and Lemma~\ref{lem:product-assoc}, $R$ acts on each piece of a picture,
with the products of the $\sP$-structures on the intermediate strands.
So
\[
R_{\xi_2^\natural,\xi_1^\natural}\Biggl(
\begin{tikzpicture}[mb]
\draw[mstr] (0,0) -- (0,3.6) .. controls (0,4.15) and (1.3,4.15) .. (1.3,3.6)
  -- (1.3,1.5) .. controls (1.3,0.95) and (2.6,0.95) .. (2.6,1.5) -- (2.6,4.6);
\node[mcpn] at (1.3,2.55) {$f$};
\node[anchor=east,inner sep=1pt] at (1.3,1.8) {$Z_1$}; \node[anchor=east,inner sep=1pt] at (1.3,3.3) {$Z_2$};
\blab{0}{Z_2^\vee}\tlab{2.6}{4.6}{Z_1^\vee}
\end{tikzpicture}
\Biggr)
=
\begin{tikzpicture}[mb]
\draw[mstr] (0,0) -- (0,3.6) .. controls (0,4.15) and (1.3,4.15) .. (1.3,3.6)
  -- (1.3,1.5) .. controls (1.3,0.95) and (2.6,0.95) .. (2.6,1.5) -- (2.6,4.6);
\node[mcpn] at (1.3,2.55) {$R_{\xi_1,\xi_2}(f)$};
\node[anchor=east,inner sep=1pt] at (1.3,1.8) {$Z_1$}; \node[anchor=east,inner sep=1pt] at (1.3,3.3) {$Z_2$};
\blab{0}{Z_2^\vee}\tlab{2.6}{4.6}{Z_1^\vee}
\end{tikzpicture},
\]
since $R$ fixes the strands by Lemma~\ref{lem:R}(a), fixes the cup of $Z_1$ by
\eqref{eq:natural-dual} and the cap of $Z_2$ by Lemma~\ref{lem:natural-dual}(a).
The left picture
is $f^\vee$, the right one $R_{\xi_1,\xi_2}(f)^\vee$.
\end{proof}

\begin{lemma}\label{lem:P-modular-data}
For $X\in\Irr(\cC)$ one has $\theta_{\sP X}=\theta_X^{-1}$; in
particular $\theta_X=\pm1$ for $X\in\Fix(\sP)$.  
\end{lemma}

\begin{proof}
The twist $\theta_X$ is the curl, the partial trace of $c_{X,X}$.
Its box $\sP(\theta_X)$ is drawn as the mirror image of the curl; by
\ref{item:B1}--\ref{item:B6} the box can be removed, the crossing passing through the
walls unchanged (\ref{item:B4}) and the cap and cup becoming those of $\sP X$
(\ref{item:B5}, \ref{item:B6}).
What remains is the mirror image of the curl on $\sP X$, which is the curl with the
opposite crossing, $\theta_{\sP X}^{-1}$.
For simple $X$, $\sP$ of the scalar $\theta_X$ is the same scalar, as $\sP$ is linear;
so $\theta_{\sP X}=\theta_X^{-1}$, and for $X\cong\sP X$ this reads $\theta_X^2=1$.
\end{proof}

We have prepared all the necessary tools to show that $\fX$ is a finite-dimensional $C^*$-algebra.  
\begin{definition}\label{def:eps}
Let $\varepsilon:=\operatorname{Tr}_\one\colon\fX\to\bC$ be the closing at
$\one$ of Definition~\ref{def:TrX}.
In pictures, we have 
\[
\varepsilon\Bigl(\mone{f}{X}{\sP X}\Bigr)=
\begin{tikzpicture}[mb]
\draw[mstr] (0,0.45) -- (0,3.35);
\node[mcpnc] at (0,0.45) {$e$};
\node[mcpn] at (0,1.75) {$f$};
\node[mcpnc] at (0,3.35) {$e$};
\ubox{-0.5}{2.85}{0.5}{3.85}
\node[anchor=west,inner sep=2pt] at (0,1.1) {$X$};
\node[anchor=west,inner sep=2pt] at (0,2.4) {$\sP X$};
\end{tikzpicture},
\]
where $e$ is the canonical element of $\Hom(\one,X)\otimes\Hom(X,\one)$.
By
Proposition~\ref{prop:basis}, $\varepsilon(u)$ is the coefficient of
$\mob(\one)$ in $u$.
\end{definition}

\begin{proposition}\label{prop:hilbert}
\begin{enumerate}[label=(\alph*)]
\item For a $\sP$-structure $(Z,\xi)$ and $X\in\Fix(\sP)$, the
  coefficient of $\mob(X)$ in $[Z,\xi]$ is the trace of the involution
  $R_{\xi_X,\xi}$ of Definition~\ref{def:R}; in particular it is an
  integer.  
\item For $X,X'\in\Fix(\sP)$ one has
  $\varepsilon\bigl(\mob(X)\mob(X')\bigr)=\nu_X\,\delta_{X',X^\vee}$,
  where $\nu_X$ is the sign of Lemma~\ref{lem:natural-dual}(c).
  Since $\nu_{X^\vee}=\nu_X$ (Lemma~\ref{lem:natural-dual}(c)), this
  implies $\varepsilon(uv)=\varepsilon(vu)$ for general $u,v\in \fX$.
\item The antilinear map $u\mapsto u^\star$ with
  $\mob(X)^\star=\nu_X\,\mob(X^\vee)$ is an anti-involution of $\fX$, and
  $[Z,\xi]^\star=[Z^\vee,\xi^\natural]$ for every $\sP$-structure
  $(Z,\xi)$.
\item $\langle u,v\rangle:=\varepsilon(u^\star v)$ is a positive definite
  Hermitian form, $\{\mob(X)\}_{X\in\Fix(\sP)}$ is an orthonormal basis,
  and $\langle au,v\rangle=\langle u,a^\star v\rangle$.
\item $(\fX,\star)$ is a finite-dimensional $C^*$-algebra; in
  particular $\fX$ is semisimple.
For each irreducible $\fX$-module
  $V$ there is $\mu_V>0$ with
  $\varepsilon=\sum_V\mu_V\operatorname{Tr}_V$, where $\operatorname{Tr}_V$
  denotes the trace of the action on $V$.
\end{enumerate}
\end{proposition}

\begin{proof}
(a) By Proposition~\ref{prop:basis}, the coefficient $c$ of $\mob(X)$ in
$[Z,\xi]$ is the closing $\operatorname{Tr}_X(Z,\xi)$ of Definition~\ref{def:TrX}.
Writing
$R:=R_{\xi_X,\xi}$ for this proof, and $R(e)$ for $R$ applied to the end
of $e$ in $\Hom(X,Z)$,
\[
c=
\begin{tikzpicture}[mb]
\draw[mstr] (0,0) -- (0,5.15) .. controls (0,5.65) and (0.9,5.65) .. (0.9,5.15)
  -- (0.9,0) .. controls (0.9,-0.5) and (0,-0.5) .. (0,0);
\node[mcpnc] at (0,0.55) {$e$};
\node[mcpn] at (0,1.8) {$\xi$};
\node[mcpnc] at (0,3.1) {$e$};
\ubox{-0.45}{2.65}{0.45}{3.55}
\node[mcpn] at (0,4.35) {$\xi_X^{-1}$};
\node[anchor=west,inner sep=1pt] at (0,1.17) {$Z$};
\node[anchor=west,inner sep=1pt] at (0,2.4) {$\sP Z$};
\node[anchor=west,inner sep=1pt] at (0,4.9) {$X$};
\end{tikzpicture}
\overset{\text{Lem.~\ref{lem:R}(a)}}{=}
\begin{tikzpicture}[mb]
\draw[mstr] (0,0) -- (0,5.55) .. controls (0,6.05) and (1.25,6.05) .. (1.25,5.55)
  -- (1.25,0) .. controls (1.25,-0.5) and (0,-0.5) .. (0,0);
\node[mcpn] at (0,0.55) {$\xi_X$};
\node[mcpnc, minimum width=12mm] at (0,2.0) {$R(e)$};
\ubox{-0.85}{1.55}{0.85}{2.45}
\node[mcpnc] at (0,3.5) {$e$};
\ubox{-0.45}{3.05}{0.45}{3.95}
\node[mcpn] at (0,4.75) {$\xi_X^{-1}$};
\node[anchor=west,inner sep=1pt] at (0,1.22) {$\sP X$};
\node[anchor=west,inner sep=1pt] at (0,2.75) {$\sP Z$};
\node[anchor=west,inner sep=1pt] at (0,5.3) {$X$};
\end{tikzpicture}
\overset{\text{\ref{item:B1}}}{=}
\begin{tikzpicture}[mb]
\draw[mstr] (0,0) -- (0,5.25) .. controls (0,5.75) and (1.25,5.75) .. (1.25,5.25)
  -- (1.25,0) .. controls (1.25,-0.5) and (0,-0.5) .. (0,0);
\node[mcpn] at (0,0.55) {$\xi_X$};
\node[mcpnc, minimum width=12mm] at (0,2.0) {$R(e)$};
\node[mcpnc] at (0,3.2) {$e$};
\ubox{-0.85}{1.55}{0.85}{3.65}
\node[mcpn] at (0,4.45) {$\xi_X^{-1}$};
\node[anchor=west,inner sep=1pt] at (0,5.0) {$X$};
\end{tikzpicture}
\overset{\text{\ref{item:E2}}}{\underset{\text{\ref{item:B3}}}{=}}\frac{\operatorname{Tr}R}{d_X}\;
\begin{tikzpicture}[mb]
\draw[mstr] (0,0) -- (0,2.9) .. controls (0,3.4) and (0.9,3.4) .. (0.9,2.9)
  -- (0.9,0) .. controls (0.9,-0.5) and (0,-0.5) .. (0,0);
\node[mcpn] at (0,0.7) {$\xi_X$};
\node[mcpn] at (0,2.2) {$\xi_X^{-1}$};
\node[anchor=west,inner sep=1pt] at (0,2.7) {$X$};
\end{tikzpicture}
=\operatorname{Tr}R\ .
\]
So $c$ is the trace of an involution which is an integer.  

\medskip\noindent (b) By (a), applied to the coefficient of $\mob(\one)$,
$\varepsilon(\mob(X)\mob(X'))$ is the trace of the involution $R_{U_0,\,\xi_X\cdot\xi_{X'}}$ of
$\Hom(\one,X\otimes X')$, which is zero unless $X'\cong X^\vee$.
For
$X'=X^\vee$  the space is spanned by $\coev_X$, and
\[
R_{U_0,\,\xi_X\cdot\xi_{X^\vee}}(\coev_X)
=\nu_X\,R_{U_0,\,\xi_X\cdot\xi_X^\natural}(\coev_X)=\nu_X\coev_X
\]
by Lemma~\ref{lem:natural-dual}(c) and \eqref{eq:natural-dual}; so $\varepsilon(\mob(X)\mob(X^\vee))=\nu_X$.

\medskip\noindent (c) Let $(Z,\xi)$ be a $\sP$-structure,
and write 
$[Z,\xi]=\sum_Xc_X\mob(X)$ and $[Z^\vee,\xi^\natural]=\sum_Xc'_X\mob(X^\vee)$.
Let us show that $c'_X=\nu_Xc_X$ for each $X\in\Fix(\sP)$.
Let $e$ and $R=R_{\xi_X,\xi}$ be as in (a).
Transposing both ends of $e$ gives the canonical element $e^\vee$ of
$\Hom(X^\vee,Z^\vee)\otimes\Hom(Z^\vee,X^\vee)$, as transposition preserves
traces.
So, by Proposition~\ref{prop:basis} applied to $[Z^\vee,\xi^\natural]$ at $X^\vee\in\Fix(\sP)$,
\[
c'_X=
\begin{tikzpicture}[mb]
\draw[mstr] (0,0) -- (0,4.15) .. controls (0,4.65) and (1.0,4.65) .. (1.0,4.15)
  -- (1.0,0) .. controls (1.0,-0.5) and (0,-0.5) .. (0,0);
\node[mcpnc] at (0,0.5) {$e^\vee$};
\node[mcpn] at (0,1.25) {$\xi^\natural$};
\node[mcpnc] at (0,2.15) {$e^\vee$};
\ubox{-0.5}{1.75}{0.5}{2.55}
\node[mcpn] at (0,3.05) {$\xi_{X^\vee}^{-1}$};
\node[anchor=west,inner sep=1pt] at (0,3.9) {$X^\vee$};
\end{tikzpicture}
\overset{\text{Lem.~\ref{lem:natural-dual}(c)}}{=}\nu_X
\begin{tikzpicture}[mb]
\draw[mstr] (0,0) -- (0,4.15) .. controls (0,4.65) and (1.0,4.65) .. (1.0,4.15)
  -- (1.0,0) .. controls (1.0,-0.5) and (0,-0.5) .. (0,0);
\node[mcpnc] at (0,0.5) {$e^\vee$};
\node[mcpn] at (0,1.25) {$\xi^\natural$};
\node[mcpnc] at (0,2.15) {$e^\vee$};
\ubox{-0.5}{1.75}{0.5}{2.55}
\node[mcpn] at (0,3.05) {$(\xi_X^\natural)^{-1}$};
\node[anchor=west,inner sep=1pt] at (0,3.9) {$X^\vee$};
\end{tikzpicture}
\overset{\text{Lem.~\ref{lem:natural-dual-R}}}{=}\nu_X
\begin{tikzpicture}[mb]
\draw[mstr] (0,0) -- (0,2.6) .. controls (0,3.1) and (1.0,3.1) .. (1.0,2.6)
  -- (1.0,0) .. controls (1.0,-0.5) and (0,-0.5) .. (0,0);
\node[mcpnc] at (0,0.65) {$e^\vee$};
\node[mcpnc] at (0,1.75) {$R(e)^\vee$};
\node[anchor=west,inner sep=1pt] at (0,2.33) {$X^\vee$};
\end{tikzpicture}
=\nu_X
\begin{tikzpicture}[mb]
\draw[mstr] (0,0) -- (0,2.6) .. controls (0,3.1) and (1.0,3.1) .. (1.0,2.6)
  -- (1.0,0) .. controls (1.0,-0.5) and (0,-0.5) .. (0,0);
\node[mcpnc] at (0,0.65) {$R(e)$};
\node[mcpnc] at (0,1.75) {$e$};
\node[anchor=west,inner sep=1pt] at (0,0.1) {$X$};
\end{tikzpicture}
\overset{\text{\ref{item:E2}}}{=}\nu_X\operatorname{Tr}R
\overset{\text{(a)}}{=}\nu_X\,c_X.
\]
Since $X\mapsto X^\vee$ permutes $\Fix(\sP)$
(Lemma~\ref{lem:natural-dual}(c)), the identity $c'_X=\nu_Xc_X$ for every $X$
says that
\[
[Z^\vee,\xi^\natural]=\sum_X\nu_X\,c_X\,\mob(X^\vee),
\qquad\text{while}\qquad
[Z,\xi]^\star=\sum_X\overline{c_X}\,\nu_X\,\mob(X^\vee)
\]
because $\star$ is antilinear.
The $c_X$ are integers by (a), so
$\overline{c_X}=c_X$, and $[Z^\vee,\xi^\natural]=[Z,\xi]^\star$.
For
$(Z,\xi)=(X\otimes X',\xi_X\cdot\xi_{X'})$, we have by Lemma~\ref{lem:natural-dual-product}
\[
(\mob(X)\mob(X'))^\star=[X'^\vee,\xi_{X'}^\natural]\,[X^\vee,\xi_X^\natural]
=\mob(X')^\star \mob(X)^\star,
\]
so $\star$ is anti-multiplicative.
It is involutive because
$\nu_X\nu_{X^\vee}=1$.

\medskip\noindent (d) By (b) and (c),
\[
\langle \mob(X),\mob(X')\rangle=\nu_X\,\varepsilon(\mob(X^\vee)\mob(X'))=\delta_{X,X'},
\qquad
\langle au,v\rangle=\varepsilon(u^\star a^\star v)=\langle u,a^\star v\rangle.
\]

\medskip\noindent (e) By (d), left multiplication realizes $\fX$ as a $\star$-closed
subalgebra of $\End_\bC(\fX)$, hence as a $C^*$-algebra.
Since
$\varepsilon$ is a faithful positive trace, its restriction to each
simple summand $\End_\bC(V)$ is $\mu_V\operatorname{Tr}_V$ with
$\mu_V=\varepsilon(q)>0$ for a minimal projection $q$.
\end{proof}

\begin{remark}\label{rem:no-unitarity}
The $C^*$-structure of Proposition~\ref{prop:hilbert} is not inherited
from a unitary structure on $\cC$, and none is assumed: the involution
$\star$ is defined by the natural duals of Lemma~\ref{lem:natural-dual},
and the positivity of $\langle\,,\rangle$ comes from the integrality (a) of
the coefficients and from the signs $\nu_X=\pm1$.
Of the pseudo-unitarity of $\cC$ we use only $d_X>0$. 
\end{remark}

\subsection{Algebra homomorphism from the ordinary fusion algebra}
\label{subsec:loops}

\begin{definition}\label{def:L}
For $Y\in\cC$ put
\[
b_Y:=U_{Y,\sP Y}\circ(\beta_Y\otimes\id_{\sP Y})\colon
Y\otimes\sP Y\to\sP(Y\otimes\sP Y),
\]
and $\bdy(Y):=[Y\otimes\sP Y,b_Y]\in\fX$, i.e.
\[
\bdy(Y)=
\begin{tikzpicture}[mb]
\mbars{-0.5}{1.4}{2}
\draw[mstr] (0,0) -- (0,2); \draw[mstr] (0.9,0) -- (0.9,2);
\node[mcpn] at (0,1) {$\beta_Y$};
\blab{0}{Y}\blab{0.9}{\sP Y}\tlab{0}{2}{\sP^2Y}\tlab{0.9}{2}{\sP Y}
\end{tikzpicture} \; .
\]
We call $\bdy(Y)$ a \emph{boundary loop}.
We extend $\bdy$
linearly to the fusion algebra $\fF$ of $\cC$.
\end{definition}

\begin{proposition}\label{prop:L}
\begin{enumerate}[label=(\alph*)]
\item $\bdy(X)\,\bdy(Y)=\bdy(X\otimes Y)$.
\item $\bdy(X\oplus Y)=\bdy(X)+\bdy(Y)$, $\bdy(\one)=\mob(\one)$, and $\bdy(X)=\bdy(X')$ if
  $X\cong X'$.
\item $\bdy(X)$ commutes with every element of $\fX$.
\end{enumerate}
Consequently $\bdy\colon\fF\to\fX$ is a unital algebra homomorphism
with values in the center $Z(\fX)$.
\end{proposition}

\begin{proof}
(a) We have
\begin{align*}
\bdy(X)\,\bdy(Y)
&\overset{\text{Prop.~\ref{prop:product}}}{=}
\begin{tikzpicture}[mb]
\mbars{-0.5}{3.2}{2.8}
\foreach \x in {0,0.9,1.8,2.7} \draw[mstr] (\x,0) -- (\x,0.95);
\sstr{1.8}{0.95}{0}{2.8}
\sstr{2.7}{0.95}{0.9}{2.8}
\sstr[mover]{0}{0.95}{1.8}{2.8}
\sstr[mover]{0.9}{0.95}{2.7}{2.8}
\node[mcpn] at (0,0.6) {$\beta_X$}; \node[mcpn] at (1.8,0.6) {$\beta_Y$};
\blab{0}{X}\blab{0.9}{\sP X}\blab{1.8}{Y}\blab{2.7}{\sP Y}
\tlab{0}{2.8}{\sP^2Y}\tlab{0.9}{2.8}{\sP Y}\tlab{1.8}{2.8}{\sP^2X}\tlab{2.7}{2.8}{\sP X}
\end{tikzpicture}
\overset{\text{Lem.~\ref{lem:crossing-move}}}{=}
\begin{tikzpicture}[mb]
\mbars{-0.5}{3.2}{2}
\foreach \x in {0,0.9,1.8,2.7} \draw[mstr] (\x,0) -- (\x,2);
\node[mcpn] at (0,1) {$\beta_X$}; \node[mcpn] at (0.9,1) {$\beta_Y$};
\blab{0}{X}\blab{0.9}{Y}\blab{1.8}{\sP Y}\blab{2.7}{\sP X}
\tlab{0}{2}{\sP^2X}\tlab{0.9}{2}{\sP^2Y}\tlab{1.8}{2}{\sP Y}\tlab{2.7}{2}{\sP X}
\end{tikzpicture}\\
&\overset{\text{Lem.~\ref{lem:trace-property}}}{=}\bigl[X\otimes Y\otimes\sP(X\otimes Y),\;U_{X\otimes Y,\sP(X\otimes Y)}
\bigl(\sP(U_{X,Y})\,U_{\sP Y,\sP X}(\beta_X\otimes\beta_Y)\otimes\id\bigr)\bigr]\\
&\overset{\eqref{eq:beta-monoidal}}{=}\bdy(X\otimes Y).
\end{align*}

\medskip\noindent (b) Decompose $(X\oplus Y)\otimes\sP(X\oplus Y)$ into
$X\otimes\sP X$, $Y\otimes\sP Y$ and the mixed part
$X\otimes\sP Y\oplus Y\otimes\sP X$.
By naturality, $b_{X\oplus Y}$ is
$b_X\oplus b_Y$ on the first two summands and exchanges the two mixed
summands.
By the spanning argument of Proposition~\ref{prop:basis},
the class of a morphism $Z\to\sP Z$ is the sum of the classes of its
diagonal blocks for any decomposition of $Z$, and here these blocks are
$b_X$, $b_Y$ and zero.
Hence $\bdy(X\oplus Y)=\bdy(X)+\bdy(Y)$.
Next,
\[
\bdy(\one)=[\one\otimes\sP\one,b_\one]=[\one,\sP(U_0)^{-1}b_\one U_0]
\]
by Lemma~\ref{lem:trace-property} for the isomorphism
$U_0\colon\one\to\sP\one$, and
$b_\one U_0=\beta_\one=\sP(U_0)U_0$ by \eqref{eq:U-unit} and
\eqref{eq:beta-monoidal}; so $\bdy(\one)=[\one,U_0]=\mob(\one)$.
Finally,
for an isomorphism $a\colon X\to X'$ we have
$b_{X'}(a\otimes\sP(a))=\sP(a\otimes\sP(a))\,b_X$ by naturality of
$\beta$ and of $U$, so $\bdy(X)=\bdy(X')$ by Lemma~\ref{lem:trace-property}.

\medskip\noindent (c) It suffices to take $u=[Y,f]$.
In pictures:
\begin{align*}
\bdy(X)\cdot[Y,f]&\overset{\text{Prop.~\ref{prop:product}}}{=}
\begin{tikzpicture}[mb]
\mbars{-0.5}{2.3}{2.8}
\foreach \x in {0,0.9,1.8} \draw[mstr] (\x,0) -- (\x,0.95);
\sstr{1.8}{0.95}{0}{2.8}
\sstr[mover]{0}{0.95}{0.9}{2.8}
\sstr[mover]{0.9}{0.95}{1.8}{2.8}
\node[mcpn] at (0,0.6) {$\beta_X$}; \node[mcpn] at (1.8,0.6) {$f$};
\blab{0}{X}\blab{0.9}{\sP X}\blab{1.8}{Y}
\tlab{0}{2.8}{\sP Y}\tlab{0.9}{2.8}{\sP^2X}\tlab{1.8}{2.8}{\sP X}
\end{tikzpicture}
\overset{\text{Lem.~\ref{lem:crossing-move}}}{=}
\begin{tikzpicture}[mb]
\mbars{-0.5}{2.3}{2}
\foreach \x in {0,0.9,1.8} \draw[mstr] (\x,0) -- (\x,2);
\node[mcpn] at (0,1) {$\beta_X$}; \node[mcpn] at (0.9,1) {$f$};
\blab{0}{X}\blab{0.9}{Y}\blab{1.8}{\sP X}
\tlab{0}{2}{\sP^2X}\tlab{0.9}{2}{\sP Y}\tlab{1.8}{2}{\sP X}
\end{tikzpicture}
\\
&\overset{\text{Lem.~\ref{lem:crossing-move}}}{=}
\begin{tikzpicture}[mb]
\mbars{-0.5}{2.3}{2.8}
\foreach \x in {0,0.9,1.8} \draw[mstr] (\x,0) -- (\x,0.95);
\sstr{0.9}{0.95}{0}{2.8}
\sstr{1.8}{0.95}{0.9}{2.8}
\sstr[mover]{0}{0.95}{1.8}{2.8}
\node[mcpn] at (0,0.6) {$f$}; \node[mcpn] at (0.9,0.6) {$\beta_X$};
\blab{0}{Y}\blab{0.9}{X}\blab{1.8}{\sP X}
\tlab{0}{2.8}{\sP^2X}\tlab{0.9}{2.8}{\sP X}\tlab{1.8}{2.8}{\sP Y}
\end{tikzpicture}
\overset{\text{Prop.~\ref{prop:product}}}{=}[Y,f]\cdot \bdy(X).
\end{align*}
For the consequences: $\bdy$ is invariant under isomorphism, unital and
additive by (b), multiplicative by (a), and central by (c).
\end{proof}

The fusion algebra $\fF$ has the antilinear involution $*$ with $[Y]^*:=[Y^\vee]$.

\begin{lemma}\label{lem:L-star}
The homomorphism $\bdy\colon\fF\to\fX$ preserves the involutions:
$\bdy(a^*)=\bdy(a)^\star$ for $a\in\fF$.
In particular $\bdy(Y)^\star=\bdy(Y^\vee)$ for $Y\in\cC$.
\end{lemma}

\begin{proof}
The Hermitian form on $\fF$ for which the $[Y]$, $Y\in\Irr(\cC)$, are orthonormal
satisfies $(ab,c)=(b,a^*c)$, as $N_{XY}^Z=N_{X^\vee Z}^Y$.
So, as in the proof of Proposition~\ref{prop:hilbert}(e), $(\fF,*)$ is a commutative
finite-dimensional $C^*$-algebra, and $\fF$ is spanned by idempotents $e$ with $e^*=e$.

By Proposition~\ref{prop:hilbert}(e), $\fX\cong\bigoplus_V\End_\bC(V)$; let $1_V$
be the unit of the summand $\End_\bC(V)$.
The anti-involution $\star$ permutes the minimal central idempotents
$1_V$, and $1_V^\star1_V\neq0$ since $\varepsilon(1_V^\star1_V)>0$; so
$1_V^\star=1_V$ for each $V$.
For an idempotent $e=e^*$ of $\fF$, $\bdy(e)$ is a central idempotent of $\fX$
by Proposition~\ref{prop:L}, hence a sum of some of the $1_V$, so
$\bdy(e)^\star=\bdy(e)=\bdy(e^*)$.
As both sides of the claim are antilinear in $a$, it holds for all $a\in\fF$.
\end{proof}

\section{Properties of the simple modules}
\label{sec:integrality}

\subsection{Generalities on modular tensor categories}
\label{subsec:statement}

From now on, we assume that $\cC$ is not just braided but modular.
In terms of the braiding $c_{X,Y}$,
the unnormalized and normalized $S$-matrices are
\[
s_{X,Y}:=\operatorname{tr}\bigl(c_{Y,X}\circ c_{X,Y}\bigr),\qquad
S_{X,Y}:=D^{-1}s_{X,Y}.
\]
$S$ is symmetric and unitary, $S_{\one,X}=d_X/D$, and
$S_{X^\vee,Y}=\overline{S_{X,Y}}$.\footnote{These hold in any modular
tensor category with real $d_X$: the twists are roots of unity by Vafa's
theorem \cite{Vafa1988,AndersonMoore1988} (see also \cite[Corollary~8.18.2]{EGNO}), so
$\overline{s_{X,Y}}=\operatorname{tr}\bigl((c_{Y,X}c_{X,Y})^{-1}\bigr)=s_{X^\vee,Y}$,
the last equality by duality; and $S^2$ is the map
$X\mapsto X^\vee$, so $S\,\overline{S}=\id$.}

Let $\fF$ be the fusion algebra of $\cC$.
By the
Verlinde formula its characters are
\[
\chi_Y([X]):=\frac{S_{X,Y}}{S_{\one,Y}}\qquad(Y\in\Irr(\cC)),
\]
and the corresponding primitive idempotents are
\begin{equation}\label{eq:fusion-idempotent}
e_X:=S_{\one,X}\sum_{A\in\Irr(\cC)}\overline{S_{A,X}}\,[A],
\qquad \chi_Y(e_X)=\delta_{X,Y}.
\end{equation}

\begin{definition}
We let $p_X=\bdy(e_X)$, which by Proposition~\ref{prop:L} and \eqref{eq:fusion-idempotent} are orthogonal central
idempotents of $\fX$ with $\sum_Xp_X=\mob(\one)$.
Then 
$\fX=\bigoplus_{X\in\Irr(\cC)}p_X\fX$ as an algebra.
\end{definition}

\begin{definition}\label{def:loop}
For
$A\in\Irr(\cC)$ let $\ell_A\in\End(\id_\cC)$ be the \emph{$A$-loop},
$\ell_{A,X}:=\operatorname{ptr}_A(c_{X,A}c_{A,X})$, the partial trace over $A$
of $c_{X,A}c_{A,X}\in\End(A\otimes X)$.
In the picture the $A$-strand passes
over $X$ from left to right, comes back under $X$, and is closed up on
the left:
\[
\ell_{A,X}
\;=\;
\begin{tikzpicture}[mb]
\draw[mstr] (0,0) -- (0,0.6); \draw[mstr] (0,0.78) -- (0,1.7);
\draw[mstr] (0,0.85) ++(0:0.42 and 0.16)
  arc[start angle=0, end angle=70, x radius=0.42, y radius=0.16];
\draw[mstr] (0,0.85) ++(110:0.42 and 0.16)
  arc[start angle=110, end angle=360, x radius=0.42, y radius=0.16];
\draw[mstr,-{Stealth[length=1.5mm,width=1.5mm]}] (0,0.85) ++(232:0.42 and 0.16)
  arc[start angle=232, end angle=256, x radius=0.42, y radius=0.16];
\node[anchor=west,inner sep=1.5pt] at (0.42,0.85) {$A$};
\blab{0}{X}\tlab{0}{1.7}{X}
\end{tikzpicture}
\]
More generally we put $\ell_{E,X}:=\operatorname{ptr}_E(c_{X,E}c_{E,X})$ for every
$E\in\cC$.
\end{definition}

The loop $\ell_{E,X}$ is natural in $X$ because the braiding is.
For simple $Y$ we have
$\ell_{A,Y}=\lambda\,\id_Y$ with
$\lambda\,d_Y=\operatorname{tr}(\ell_{A,Y})=\operatorname{tr}(c_{Y,A}c_{A,Y})=s_{A,Y}$,
and $S_{\one,Y}=d_Y/D$ gives
\[
\ell_{A,Y}=\chi_Y([A])\,\id_Y\qquad(Y\in\Irr(\cC)).
\]
Reversing the orientation of the loop replaces $A$ by $A^\vee$.

\begin{definition}\label{def:Omega}
Let $\Omega:=\sum_{V\in\Irr(\cC)}d_V\,[V]\in\fF$.
A closed strand labeled $\Omega$ in a diagram stands for the
sum over $V\in\Irr(\cC)$, weighted by $d_V$, of the diagram with this
strand labeled $V$.
\end{definition}

\begin{remark}\label{rem:Omega-mirror}
An $\Omega$-strand needs no orientation, since reversing the
orientation of a $V$-strand replaces $V$ by $V^\vee$, and
$d_{V^\vee}=d_V$.
So the $\Omega$-loop is drawn without an arrow:
\[
\ell_{\Omega,X}=\sum_{V\in\Irr(\cC)}d_V\,\ell_{V,X}=
\begin{tikzpicture}[mb]
\draw[mstr] (0,0) -- (0,0.95); \draw[mstr] (0,1.13) -- (0,2.3);
\draw[mstr] (0,1.2) ++(0:0.42 and 0.16)
  arc[start angle=0, end angle=70, x radius=0.42, y radius=0.16];
\draw[mstr] (0,1.2) ++(110:0.42 and 0.16)
  arc[start angle=110, end angle=360, x radius=0.42, y radius=0.16];
\node[anchor=west,inner sep=1.5pt] at (0.42,1.2) {$\Omega$};
\blab{0}{X}\tlab{0}{2.3}{X}
\end{tikzpicture}
\]
Moreover, since $V\mapsto\sP V$ permutes $\Irr(\cC)$ with $d_{\sP V}=d_V$
(Corollary~\ref{cor:P-dim}), a sum over $V\in\Irr(\cC)$, weighted by $d_V$, of
diagrams with a strand labeled $\sP V$ is the diagram with this strand
labeled $\Omega$; in particular an $\Omega$-strand passes through the walls of
a $\sP$-box unchanged.
\end{remark}

The next is a standard property of the $\Omega$-loop, which we prove for completeness.
\begin{lemma}\label{lem:kill}
Let $X\in\cC$.
An $\Omega$-loop cuts the strand it encircles:
  \[
  \begin{tikzpicture}[mb]
\draw[mstr] (0,0) -- (0,0.95); \draw[mstr] (0,1.13) -- (0,2.3);
\draw[mstr] (0,1.2) ++(0:0.42 and 0.16)
  arc[start angle=0, end angle=70, x radius=0.42, y radius=0.16];
\draw[mstr] (0,1.2) ++(110:0.42 and 0.16)
  arc[start angle=110, end angle=360, x radius=0.42, y radius=0.16];
\node[anchor=west,inner sep=1.5pt] at (0.42,1.2) {$\Omega$};
\blab{0}{X}\tlab{0}{2.3}{X}
\end{tikzpicture}
  =D^2\;
  \begin{tikzpicture}[mb]
\draw[mstr] (0,0) -- (0,0.45); \draw[mstr] (0,1.85) -- (0,2.3);
\node[mcpnc] at (0,0.7) {$e$}; \node[mcpnc] at (0,1.6) {$e$};
\blab{0}{X}\tlab{0}{2.3}{X}
\end{tikzpicture},
  \]
  where $e$ is the canonical element of $\Hom(X,\one)\otimes\Hom(\one,X)$.
The same holds for a strand passing horizontally through the loop:
  \[
  \begin{tikzpicture}[mb]
\draw[mstr] (0,1.2) -- (0.99,1.2); \draw[mstr] (1.17,1.2) -- (2.3,1.2);
\draw[mstr] (1.15,1.2) ++(20:0.16 and 0.42)
  arc[start angle=20, end angle=340, x radius=0.16, y radius=0.42];
\node[anchor=south,inner sep=1.5pt] at (1.15,1.62) {$\Omega$};
\node[anchor=east,inner sep=1.5pt] at (0,1.2) {$X$};
\node[anchor=west,inner sep=1.5pt] at (2.3,1.2) {$X$};
\end{tikzpicture}
  =D^2\;
  \begin{tikzpicture}[mb]
\draw[mstr] (0,1.2) -- (0.25,1.2) .. controls (0.55,1.2) .. (0.55,1.5) -- (0.55,1.6);
\node[mcpnc] at (0.55,1.85) {$e$};
\node[mcpnc] at (1.45,0.55) {$e$};
\draw[mstr] (1.45,0.8) -- (1.45,0.9) .. controls (1.45,1.2) .. (1.75,1.2) -- (2.0,1.2);
\node[anchor=east,inner sep=1.5pt] at (0,1.2) {$X$};
\node[anchor=west,inner sep=1.5pt] at (2.0,1.2) {$X$};
\end{tikzpicture}.
  \]
In particular
  $\ell_{\Omega,C}=D^2\,\delta_{C,\one}\,\id_C$ for simple $C$.
\end{lemma}

\begin{proof}
First let $X=C$ be simple.
As $\ell_{V,C}=\chi_C([V])\,\id_C$ and $d_V=D\,S_{\one,V}$,
\[
\ell_{\Omega,C}=\sum_{V\in\Irr(\cC)}d_V\,\frac{S_{V,C}}{S_{\one,C}}\,\id_C
=\frac{D^2}{d_C}\sum_{V\in\Irr(\cC)}S_{\one,V}S_{V,C}\,\id_C
=D^2\,\delta_{C,\one}\,\id_C ,
\]
since $S$ is symmetric and $S^2$ is the map $X\mapsto X^\vee$.
For general $X$,
\[
\begin{tikzpicture}[mb]
\draw[mstr] (0,0) -- (0,0.95); \draw[mstr] (0,1.13) -- (0,2.3);
\draw[mstr] (0,1.2) ++(0:0.42 and 0.16)
  arc[start angle=0, end angle=70, x radius=0.42, y radius=0.16];
\draw[mstr] (0,1.2) ++(110:0.42 and 0.16)
  arc[start angle=110, end angle=360, x radius=0.42, y radius=0.16];
\node[anchor=west,inner sep=1.5pt] at (0.42,1.2) {$\Omega$};
\blab{0}{X}\tlab{0}{2.3}{X}
\end{tikzpicture}
\overset{\text{\ref{item:E3}}}{=}\sum_{C\in\Irr(\cC)}d_C\;
\begin{tikzpicture}[mb]
\draw[mstr] (0,0) -- (0,0.45); \draw[mstr] (0,0.75) -- (0,1.61); \draw[mstr] (0,1.79) -- (0,2.35);
\draw[mstr] (0,2.65) -- (0,3.1);
\draw[mstr] (0,1.86) ++(0:0.42 and 0.16)
  arc[start angle=0, end angle=70, x radius=0.42, y radius=0.16];
\draw[mstr] (0,1.86) ++(110:0.42 and 0.16)
  arc[start angle=110, end angle=360, x radius=0.42, y radius=0.16];
\node[anchor=west,inner sep=1.5pt] at (0.42,1.86) {$\Omega$};
\node[mcpnc] at (0,0.6) {$e$}; \node[mcpnc] at (0,2.5) {$e$};
\node[anchor=east,inner sep=1.5pt] at (0,1.2) {$C$};
\blab{0}{X}\tlab{0}{3.1}{X}
\end{tikzpicture}
=D^2\;
\begin{tikzpicture}[mb]
\draw[mstr] (0,0) -- (0,0.45); \draw[mstr] (0,2.65) -- (0,3.1);
\node[mcpnc] at (0,0.6) {$e$}; \node[mcpnc] at (0,2.5) {$e$};
\blab{0}{X}\tlab{0}{3.1}{X}
\end{tikzpicture}
\]
by the simple case, as only $C=\one$ contributes and $d_\one=1$.
The second picture follows from the first by a planar isotopy which turns
the part of the picture around the loop by a right angle; this part
contains no coupon.
\end{proof}

We now need to recall the properties of $\ell_A$ for invertible objects $A$.

\begin{definition}\label{def:I}
Let $\Inv(\cC)\subset\Irr(\cC)$ be the set of invertible simple objects, that is,
of those $A$ with $A\otimes A^\vee\cong\one$.
A tensor product of invertible objects is invertible, hence simple, and
$A\otimes B\cong B\otimes A$ by the braiding; so the isomorphism classes of invertible objects form a finite abelian group
under the tensor product, with unit $\one$ and inverse $A\mapsto A^\vee$.
We identify this group with $\Inv(\cC)$, so that relations in $\Inv(\cC)$ are
isomorphisms in $\cC$.
We write $\operatorname{Aut}_\otimes(\id_\cC)$ for the group of monoidal
automorphisms of the identity functor, that is, of natural automorphisms
$\varphi_X$ of $X$ with $\varphi_{X\otimes Y}=\varphi_X\otimes\varphi_Y$.
\end{definition}

\begin{lemma}\label{lem:monodromy}
Let $A$ be an invertible object and $X\in\cC$.
Then $d_A=1$, and the monodromy of $A$ and $X$ is the $A$-loop on $X$:
\[
c_{X,A}\,c_{A,X}=\id_A\otimes\ell_{A,X},\qquad
c_{A,X}\,c_{X,A}=\ell_{A,X}\otimes\id_A,
\]
that is,
\[
\begin{tikzpicture}[mb]
\sstr{0.9}{0}{0}{1}
\sstr[mover]{0}{0}{0.9}{1}
\sstr{0.9}{1}{0}{2}
\sstr[mover]{0}{1}{0.9}{2}
\blab{0}{A}\blab{0.9}{X}\tlab{0}{2}{A}\tlab{0.9}{2}{X}
\end{tikzpicture}
\;=\;
\begin{tikzpicture}[mb]
\draw[mstr] (0,0) -- (0,2);
\begin{scope}[shift={(0.9,0.15)}]
\draw[mstr] (0,-0.15) -- (0,0.6); \draw[mstr] (0,0.78) -- (0,1.85);
\draw[mstr] (0,0.85) ++(0:0.42 and 0.16)
  arc[start angle=0, end angle=70, x radius=0.42, y radius=0.16];
\draw[mstr] (0,0.85) ++(110:0.42 and 0.16)
  arc[start angle=110, end angle=360, x radius=0.42, y radius=0.16];
\draw[mstr,-{Stealth[length=1.5mm,width=1.5mm]}] (0,0.85) ++(232:0.42 and 0.16)
  arc[start angle=232, end angle=256, x radius=0.42, y radius=0.16];
\node[anchor=west,inner sep=1.5pt] at (0.42,0.85) {$A$};
\end{scope}
\blab{0}{A}\blab{0.9}{X}\tlab{0}{2}{A}\tlab{0.9}{2}{X}
\end{tikzpicture}\,,
\qquad
\begin{tikzpicture}[mb]
\sstr{0.9}{0}{0}{1}
\sstr[mover]{0}{0}{0.9}{1}
\sstr{0.9}{1}{0}{2}
\sstr[mover]{0}{1}{0.9}{2}
\blab{0}{X}\blab{0.9}{A}\tlab{0}{2}{X}\tlab{0.9}{2}{A}
\end{tikzpicture}
\;=\;
\begin{tikzpicture}[mb]
\begin{scope}[shift={(0,0.15)}]
\draw[mstr] (0,-0.15) -- (0,0.6); \draw[mstr] (0,0.78) -- (0,1.85);
\draw[mstr] (0,0.85) ++(0:0.42 and 0.16)
  arc[start angle=0, end angle=70, x radius=0.42, y radius=0.16];
\draw[mstr] (0,0.85) ++(110:0.42 and 0.16)
  arc[start angle=110, end angle=360, x radius=0.42, y radius=0.16];
\draw[mstr,-{Stealth[length=1.5mm,width=1.5mm]}] (0,0.85) ++(232:0.42 and 0.16)
  arc[start angle=232, end angle=256, x radius=0.42, y radius=0.16];
\node[anchor=west,inner sep=1.5pt] at (0.42,0.85) {$A$};
\end{scope}
\draw[mstr] (1.2,0) -- (1.2,2);
\blab{0}{X}\blab{1.2}{A}\tlab{0}{2}{X}\tlab{1.2}{2}{A}
\end{tikzpicture}\,.
\]
\end{lemma}

\begin{proof}
From $A\otimes A^\vee\cong\one$ we get $d_A^2=1$, so $d_A=1$.
As $A\otimes-$ is an equivalence, $f\mapsto\id_A\otimes f$ is a bijection
$\End(X)\to\End(A\otimes X)$, so $c_{X,A}c_{A,X}=\id_A\otimes\varphi_X$ for a unique
$\varphi_X\in\End(X)$; the partial trace over $A$ gives $\ell_{A,X}=d_A\varphi_X=\varphi_X$.
Conjugating by $c_{X,A}$ and using the naturality of the braiding,
$c_{A,X}c_{X,A}=c_{X,A}^{-1}(\id_A\otimes\ell_{A,X})\,c_{X,A}=\ell_{A,X}\otimes\id_A$.
\end{proof}

The following is well known (see e.g.\ \cite{DGNO} and \cite[Lemma~8.22.9]{EGNO}); we include a proof for
completeness.

\begin{proposition}
  \label{prop:loop-automorphism}
For $A\in\Inv(\cC)$, the loop $\ell_A$ is a monoidal automorphism of $\id_\cC$, and
$A\mapsto\ell_A$ is a group isomorphism $\Inv(\cC)\to\operatorname{Aut}_\otimes(\id_\cC)$.
\end{proposition}

\begin{proof}
\noindent\emph{$\ell_A$ is monoidal.}
Let $A\in\Inv(\cC)$.
By the hexagon axioms and Lemma~\ref{lem:monodromy},
\[
c_{X\otimes Y,A}\,c_{A,X\otimes Y}
=(c_{X,A}\otimes\id_Y)(\id_X\otimes c_{Y,A}c_{A,Y})(c_{A,X}\otimes\id_Y)
=c_{X,A}c_{A,X}\otimes\ell_{A,Y}
=\id_A\otimes\ell_{A,X}\otimes\ell_{A,Y},
\]
so $\ell_{A,X\otimes Y}=\ell_{A,X}\otimes\ell_{A,Y}$ by Lemma~\ref{lem:monodromy}
for $X\otimes Y$, and $\ell_A$ is monoidal; it is natural because the braiding is.

\medskip\noindent\emph{Homomorphism.}
For $E,F,X\in\cC$,
\[
\ell_{E\otimes F,X}
=
\begin{tikzpicture}[mb]
\draw[mstr] (0,0) -- (0,0.6); \draw[mstr] (0,0.78) -- (0,1.7);
\draw[mstr] (0,0.85) ++(0:0.42 and 0.16)
  arc[start angle=0, end angle=70, x radius=0.42, y radius=0.16];
\draw[mstr] (0,0.85) ++(110:0.42 and 0.16)
  arc[start angle=110, end angle=360, x radius=0.42, y radius=0.16];
\draw[mstr,-{Stealth[length=1.5mm,width=1.5mm]}] (0,0.85) ++(232:0.42 and 0.16)
  arc[start angle=232, end angle=256, x radius=0.42, y radius=0.16];
\node[anchor=west,inner sep=1.5pt] at (0.42,0.85) {$E\otimes F$};
\blab{0}{X}\tlab{0}{1.7}{X}
\end{tikzpicture}
\overset{\mathrm{(1)}}{=}
\begin{tikzpicture}[mb]
\draw[mstr] (0,0) -- (0,0.42); \draw[mstr] (0,0.56) -- (0,0.64); \draw[mstr] (0,0.78) -- (0,1.7);
\draw[mstr] (0,0.85) ++(0:0.42 and 0.14)
  arc[start angle=0, end angle=70, x radius=0.42, y radius=0.14];
\draw[mstr] (0,0.85) ++(110:0.42 and 0.14)
  arc[start angle=110, end angle=360, x radius=0.42, y radius=0.14];
\draw[mstr,-{Stealth[length=1.5mm,width=1.5mm]}] (0,0.85) ++(232:0.42 and 0.14)
  arc[start angle=232, end angle=256, x radius=0.42, y radius=0.14];
\draw[mstr] (0,0.85) ++(0:1.0 and 0.36)
  arc[start angle=0, end angle=82, x radius=1.0, y radius=0.36];
\draw[mstr] (0,0.85) ++(98:1.0 and 0.36)
  arc[start angle=98, end angle=360, x radius=1.0, y radius=0.36];
\draw[mstr,-{Stealth[length=1.5mm,width=1.5mm]}] (0,0.85) ++(240:1.0 and 0.36)
  arc[start angle=240, end angle=252, x radius=1.0, y radius=0.36];
\node[anchor=west,inner sep=1pt] at (0.42,0.85) {$E$};
\node[anchor=west,inner sep=1.5pt] at (1.0,0.85) {$F$};
\blab{0}{X}\tlab{0}{1.7}{X}
\end{tikzpicture}
\overset{\text{isotopy}}{=}
\begin{tikzpicture}[mb]
\draw[mstr] (0,0) -- (0,0.5); \draw[mstr] (0,0.68) -- (0,1.25); \draw[mstr] (0,1.43) -- (0,2.3);
\draw[mstr] (0,0.75) ++(0:0.42 and 0.16)
  arc[start angle=0, end angle=70, x radius=0.42, y radius=0.16];
\draw[mstr] (0,0.75) ++(110:0.42 and 0.16)
  arc[start angle=110, end angle=360, x radius=0.42, y radius=0.16];
\draw[mstr,-{Stealth[length=1.5mm,width=1.5mm]}] (0,0.75) ++(232:0.42 and 0.16)
  arc[start angle=232, end angle=256, x radius=0.42, y radius=0.16];
\node[anchor=west,inner sep=1.5pt] at (0.42,0.75) {$F$};
\draw[mstr] (0,1.5) ++(0:0.42 and 0.16)
  arc[start angle=0, end angle=70, x radius=0.42, y radius=0.16];
\draw[mstr] (0,1.5) ++(110:0.42 and 0.16)
  arc[start angle=110, end angle=360, x radius=0.42, y radius=0.16];
\draw[mstr,-{Stealth[length=1.5mm,width=1.5mm]}] (0,1.5) ++(232:0.42 and 0.16)
  arc[start angle=232, end angle=256, x radius=0.42, y radius=0.16];
\node[anchor=west,inner sep=1.5pt] at (0.42,1.5) {$E$};
\blab{0}{X}\tlab{0}{2.3}{X}
\end{tikzpicture}
=\ell_{E,X}\,\ell_{F,X},
\]
where (1) draws the strand $E\otimes F$ as two parallel strands $E$, $F$ (the hexagon axioms for
the crossings, and the nested cap and cup of $E\otimes F$).
In particular $\ell_{A\otimes B}=\ell_A\ell_B$ for $A,B\in\Inv(\cC)$.
Since the loop does not change when its label is replaced by an isomorphic
object, $A\mapsto\ell_A$ is a homomorphism from $\Inv(\cC)$.

\medskip\noindent\emph{Injectivity.}
If $\ell_A=\ell_B$, then $\chi_Y([A])=\chi_Y([B])$ for all $Y\in\Irr(\cC)$, as
$\ell_{A,Y}=\chi_Y([A])\,\id_Y$.
The characters $\chi_Y$ separate the elements of $\fF$, so $[A]=[B]$.

\medskip\noindent\emph{Surjectivity.}
Let $\varphi\in\operatorname{Aut}_\otimes(\id_\cC)$.
For $Y\in\Irr(\cC)$, $\varphi_Y$ is a non-zero scalar.
If $Z\subset X\otimes Y$ for simple $X,Y,Z$, then $\varphi_Z=\varphi_X\varphi_Y$ by
the naturality of $\varphi$ and $\varphi_{X\otimes Y}=\varphi_X\otimes\varphi_Y$.
Hence $[Y]\mapsto\varphi_Yd_Y$ is a character of $\fF$, so there is
$A\in\Irr(\cC)$ with
\[
\varphi_Y\,d_Y=\chi_A([Y])=\frac{s_{Y,A}}{d_A}\qquad(Y\in\Irr(\cC)).
\]
We show $A\in\Inv(\cC)$.
Since $\one\subset Y\otimes Y^\vee$, we have $\varphi_{Y^\vee}=\varphi_Y^{-1}$,
while $s_{Y^\vee,A}=\overline{s_{Y,A}}$ and $d_{Y^\vee}=d_Y$.
Hence $\varphi_Y^{-1}=\overline{\varphi_Y}$, that is, $|\varphi_Y|=1$, and so
$|S_{Y,A}|=d_Ad_Y/D$.
Unitarity of $S$ gives
\[
1=\sum_{Y\in\Irr(\cC)}|S_{Y,A}|^2=\frac{d_A^2}{D^2}\sum_{Y\in\Irr(\cC)}d_Y^2=d_A^2 ,
\]
so $d_A=1$.
Then $A\otimes A^\vee$ contains $\one$ and has dimension $1$; as all
dimensions are positive, $A\otimes A^\vee\cong\one$, that is, $A\in\Inv(\cC)$.
Finally $\varphi_Y=s_{Y,A}/d_Y=\chi_Y([A])$, which is the scalar of
$\ell_{A,Y}$.
As $\cC$ is semisimple, $\varphi=\ell_A$.
\end{proof}

\subsection{Modular tensor categories and anti-monoidal involutions}

\begin{lemma}\label{lem:box-loop}
For $E,X\in\cC$, the box of the $E$-loop on $X$ is the $(\sP E)^\vee$-loop on
$\sP X$:
\[
\sP(\ell_{E,X})=\ell_{(\sP E)^\vee,\sP X}.
\]
In particular $\sP(\ell_{A,X})=\ell_{A,\sP X}^{-1}$ for $A\in\Inv(\cC)$ with
$\sP A\cong A$.
\end{lemma}

\begin{proof}
\[
\sP(\ell_{E,X})
=
\begin{tikzpicture}[mb]
\draw[mstr] (0,0) -- (0,0.6); \draw[mstr] (0,0.78) -- (0,1.7);
\draw[mstr] (0,0.85) ++(0:0.42 and 0.16)
  arc[start angle=0, end angle=70, x radius=0.42, y radius=0.16];
\draw[mstr] (0,0.85) ++(110:0.42 and 0.16)
  arc[start angle=110, end angle=360, x radius=0.42, y radius=0.16];
\draw[mstr,-{Stealth[length=1.5mm,width=1.5mm]}] (0,0.85) ++(308:0.42 and 0.16)
  arc[start angle=308, end angle=284, x radius=0.42, y radius=0.16];
\node[anchor=east,inner sep=1.5pt] at (-0.42,0.85) {$E$};
\ubox{-1.0}{0.3}{0.65}{1.4}
\blab{0}{\sP X}\tlab{0}{1.7}{\sP X}
\end{tikzpicture}
\overset{\text{\ref{item:B1}--\ref{item:B6}}}{=}
\begin{tikzpicture}[mb]
\draw[mstr] (0,0) -- (0,0.6); \draw[mstr] (0,0.78) -- (0,1.7);
\draw[mstr] (0,0.85) ++(0:0.42 and 0.16)
  arc[start angle=0, end angle=70, x radius=0.42, y radius=0.16];
\draw[mstr] (0,0.85) ++(110:0.42 and 0.16)
  arc[start angle=110, end angle=360, x radius=0.42, y radius=0.16];
\draw[mstr,-{Stealth[length=1.5mm,width=1.5mm]}] (0,0.85) ++(308:0.42 and 0.16)
  arc[start angle=308, end angle=284, x radius=0.42, y radius=0.16];
\node[anchor=east,inner sep=1.5pt] at (-0.42,0.85) {$\sP E$};
\blab{0}{\sP X}\tlab{0}{1.7}{\sP X}
\end{tikzpicture}
=
\begin{tikzpicture}[mb]
\draw[mstr] (0,0) -- (0,0.6); \draw[mstr] (0,0.78) -- (0,1.7);
\draw[mstr] (0,0.85) ++(0:0.42 and 0.16)
  arc[start angle=0, end angle=70, x radius=0.42, y radius=0.16];
\draw[mstr] (0,0.85) ++(110:0.42 and 0.16)
  arc[start angle=110, end angle=360, x radius=0.42, y radius=0.16];
\draw[mstr,-{Stealth[length=1.5mm,width=1.5mm]}] (0,0.85) ++(232:0.42 and 0.16)
  arc[start angle=232, end angle=256, x radius=0.42, y radius=0.16];
\node[anchor=west,inner sep=1.5pt] at (0.42,0.85) {$(\sP E)^\vee$};
\blab{0}{\sP X}\tlab{0}{1.7}{\sP X}
\end{tikzpicture}
=\ell_{(\sP E)^\vee,\sP X}.
\]
For $A\in\Inv(\cC)$ with $\sP A\cong A$, $\ell_{(\sP A)^\vee}=\ell_{A^\vee}=\ell_A^{-1}$ by
Proposition~\ref{prop:loop-automorphism}.
\end{proof}

\begin{proposition}\label{prop:kappa}
There is a unique $\obstruction\in\Inv(\cC)$, the \emph{coherence obstruction} of $\beta$, with
\begin{equation}\label{eq:obstruction}
\beta_{\sP X}=\sP(\beta_X)\circ\ell_{\obstruction,\sP X}\qquad(X\in\cC),
\end{equation}
that is,
\[
\begin{tikzpicture}[mb]
\draw[mstr] (0,0) -- (0,3.1);
\node[mcpn] at (0,1.55) {$\beta_{\sP X}$};
\blab{0}{\sP X}\tlab{0}{3.1}{\sP^3X}
\end{tikzpicture}
\;=\;
\begin{tikzpicture}[mb]
\draw[mstr] (0,0) -- (0,0.6); \draw[mstr] (0,0.78) -- (0,3.1);
\draw[mstr] (0,0.85) ++(0:0.42 and 0.16)
  arc[start angle=0, end angle=70, x radius=0.42, y radius=0.16];
\draw[mstr] (0,0.85) ++(110:0.42 and 0.16)
  arc[start angle=110, end angle=360, x radius=0.42, y radius=0.16];
\draw[mstr,-{Stealth[length=1.5mm,width=1.5mm]}] (0,0.85) ++(232:0.42 and 0.16)
  arc[start angle=232, end angle=256, x radius=0.42, y radius=0.16];
\node[anchor=west,inner sep=1.5pt] at (0.42,0.85) {$\obstruction$};
\node[mcpn] at (0,2.15) {$\beta_X$};
\ubox{-0.45}{1.6}{0.45}{2.7}
\blab{0}{\sP X}\tlab{0}{3.1}{\sP^3X}
\end{tikzpicture}\,.
\]
\end{proposition}

\begin{proof}
By \eqref{eq:beta-monoidal} and since $(\sP,U)$ is anti-monoidal,
$\beta_{\sP X}$ and $\sP(\beta_X)$ are monoidal natural
isomorphisms $\sP\Rightarrow\sP^3$, so, $\sP$ being an equivalence,
$\sP(\beta_X)^{-1}\beta_{\sP X}=\varphi_{\sP X}$ for a unique
$\varphi\in\operatorname{Aut}_\otimes(\id_\cC)$, which is $\ell_\obstruction$ for a unique
$\obstruction\in\Inv(\cC)$ by Proposition~\ref{prop:loop-automorphism}.
\end{proof}

\begin{definition}
  \label{def:coherent}
We call an anti-monoidal involution $(\sP,U,\beta)$ of $\cC$
\emph{coherent}
if it satisfies a further condition
\begin{equation}
\beta_{\sP X}=\sP(\beta_X),
 \tag{Coh}\label{eq:beta-coherence}
\end{equation}
that is, if its coherence obstruction is $\obstruction=\one$.
\end{definition}

The object $\obstruction$ is the counterpart of the
anyon $\mathfrak O(\mathbf g,\mathbf g,\mathbf g)$ of \cite[(105)]{Barkeshli_2019}
for the generator $\mathbf g$ of $\bZ/2$, whence the letter $O$.

\begin{proposition}
\label{prop:obstruction}
Let $\obstruction$ be the coherence obstruction of $\beta$ (Proposition~\ref{prop:kappa}).
\begin{enumerate}[label=(\alph*)]
\item $\sP\obstruction\cong\obstruction$.
\item The monoidal natural isomorphisms $\beta'\colon\id_\cC\Rightarrow\sP^2$
  are the $\beta\circ\ell_B$ with $B\in\Inv(\cC)$, and the coherence obstruction
  of $\beta\circ\ell_B$ is $\obstruction\otimes B\otimes\sP B$.
\item Regard $\Inv(\cC)$ as a $\bZ/2$-module, the generator acting by
  $A\mapsto(\sP A)^\vee$.
  The class of $\obstruction$ in
  \[
  \{A\in\Inv(\cC)\mid\sP A\cong A\}\big/\{B\otimes\sP B\mid B\in\Inv(\cC)\}\cong
  H^3(\bZ/2,\Inv(\cC))
  \]
  is the obstruction class of $(\sP,U)$: it depends only on $(\sP,U)$, and it
  vanishes if and only if $\beta$ can be chosen to be coherent.
\end{enumerate}
\end{proposition}

\begin{proof}
We use that $\sP(\ell_{A,X})=\ell_{(\sP A)^\vee,\sP X}$ (Lemma~\ref{lem:box-loop}).

\medskip\noindent (a)
The naturality of $\beta$ applied to the morphism $\beta_X\colon X\to\sP^2X$
gives $\beta_{\sP^2X}\beta_X=\sP^2(\beta_X)\beta_X$, so
$\beta_{\sP^2X}=\sP^2(\beta_X)$.
Applying $\sP$ to \eqref{eq:obstruction} and using
\eqref{eq:obstruction} for $\sP X$,
\[
\ell_{(\sP\obstruction)^\vee,\sP^2X}
=\sP^2(\beta_X)^{-1}\,\sP(\beta_{\sP X})
=\beta_{\sP^2X}^{-1}\,\sP(\beta_{\sP X})
=\ell_{\obstruction,\sP^2X}^{-1}
=\ell_{\obstruction^\vee,\sP^2X}.
\]
As $\sP^2$ is essentially surjective,
$\ell_{(\sP\obstruction)^\vee}=\ell_{\obstruction^\vee}$, and
$\sP\obstruction\cong\obstruction$ by Proposition~\ref{prop:loop-automorphism}.

\medskip\noindent (b)
If $\beta'$ is a monoidal natural isomorphism $\id_\cC\Rightarrow\sP^2$, then
$\beta^{-1}\circ\beta'\in\operatorname{Aut}_\otimes(\id_\cC)$, which is $\ell_B$
for a unique $B\in\Inv(\cC)$ by Proposition~\ref{prop:loop-automorphism};
conversely each $\beta\circ\ell_B$ is monoidal.
For $\beta'=\beta\circ\ell_B$, since the $\ell_A$ commute with each other,
\[
\sP(\beta'_X)^{-1}\beta'_{\sP X}
=\sP(\ell_{B,X})^{-1}\,\sP(\beta_X)^{-1}\beta_{\sP X}\,\ell_{B,\sP X}
=\ell_{\sP B,\sP X}\,\ell_{\obstruction,\sP X}\,\ell_{B,\sP X}
=\ell_{\obstruction\otimes B\otimes\sP B,\sP X}.
\]

\medskip\noindent (c)
For a cyclic group generated by $g$ with norm $N=1+g$,
$H^3(\bZ/2,\Inv(\cC))=\Ker N/\Ima(1-g)$; in the multiplicative notation of $\Inv(\cC)$,
$N(A)=A\otimes(\sP A)^\vee$ and $(1-g)(B)=B\otimes\sP B$, which gives the
displayed isomorphism.
By (a), $\obstruction$ lies in the numerator, and by (b) its class does not
depend on $\beta$.
If the class vanishes, $\obstruction\cong B\otimes\sP B$ for some $B\in\Inv(\cC)$, and then
the coherence obstruction of $\beta\circ\ell_{B^\vee}$ is
$\obstruction\otimes B^\vee\otimes\sP(B^\vee)\cong\one$, as $\sP$ preserves duals by
\ref{item:B5}; so $\beta\circ\ell_{B^\vee}$ is coherent by Definition~\ref{def:coherent}.
Conversely, if some $\beta'$ is coherent, its coherence obstruction is $\one$, and
the class vanishes.
\end{proof}

\begin{remark}\label{rem:finding-P}
In Appendix~\ref{app:skeletal} we give an explicit method to find an
anti-monoidal involution on a given $\cC$ with multiplicity-free fusion
rules, or to show that none exists, starting from the $F$- and $R$-symbols
of a skeletal model of $\cC$.
\end{remark}

\begin{lemma}\label{lem:fY}
Let $Y\in\Irr(\cC)$ and $X\in\Fix(\sP)$, put $Z:=Y\otimes\sP Y$, and let
$\zeta_Y:=\chi_{\sP Y}([\obstruction])$, so that
$\ell_{\obstruction,\sP Y}=\zeta_Y\,\id_{\sP Y}$.
The linear map
$R_{\xi_X,b_Y}\colon\Phi\mapsto b_Y^{-1}\,\sP(\Phi)\,\xi_X$ of $\Hom(X,Z)$,
introduced in Definition~\ref{def:R}, satisfies
\[
R_{\xi_X,b_Y}^2=\zeta_Y,\qquad
\operatorname{Tr}(R_{\xi_X,b_Y})=\varepsilon\bigl(\mob(X)\,\bdy(Y^\vee)\bigr).
\]
Here $\zeta_Y$ is a root of unity.
When $\beta$ is coherent, $\zeta_Y=1$, so $R_{\xi_X,b_Y}$ is an involution and its
trace is an integer.
\end{lemma}

\begin{proof}
Write $b:=b_Y$, $\xi:=\xi_X$ and $R:=R_{\xi_X,b_Y}$.

\noindent\emph{The square of $R$.}
We first show $\sP(b)\,b=\zeta_Y^{-1}\,\beta_Z$.
We have $b_Y=U_{Y,\sP Y}\,\hat b_Y$ with $\hat b_Y:=\beta_Y\otimes\id_{\sP Y}$.
By the argument of Remark~\ref{rem:P-structure-box}, it suffices to show
that $\hat b_Y$ followed by its box is $\zeta_Y^{-1}\,\beta_Y\otimes\beta_{\sP Y}$:
\[
\begin{tikzpicture}[mb]
\draw[mstr] (0,0) -- (0,3.7); \draw[mstr] (1.4,0) -- (1.4,3.7);
\node[mcpn] at (0,0.5) {$\beta_Y$};
\node[mcpn] at (1.4,2.3) {$\beta_Y$};
\ubox{-0.4}{1.4}{2.1}{3.2}
\node[anchor=east,inner sep=1pt] at (0,1.1) {$\sP^2Y$};
\node[anchor=west,inner sep=1pt] at (1.4,0.75) {$\sP Y$};
\node[anchor=west,inner sep=1pt] at (0,2.3) {$\sP Y$};
\node[anchor=east,inner sep=1pt] at (1.4,1.75) {$Y$};
\node[anchor=east,inner sep=1pt] at (1.4,2.9) {$\sP^2Y$};
\blab{0}{Y}\blab{1.4}{\sP Y}\tlab{0}{3.7}{\sP^2Y}\tlab{1.4}{3.7}{\sP^3Y}
\end{tikzpicture}
\overset{\text{\ref{item:B2}}}{\underset{\text{\ref{item:B3}}}{=}}
\begin{tikzpicture}[mb]
\draw[mstr] (0,0) -- (0,3.7); \draw[mstr] (1.4,0) -- (1.4,3.7);
\node[mcpn] at (0,0.5) {$\beta_Y$};
\node[mcpn] at (1.4,2.3) {$\sP\beta_Y$};
\blab{0}{Y}\blab{1.4}{\sP Y}\tlab{0}{3.7}{\sP^2Y}\tlab{1.4}{3.7}{\sP^3Y}
\end{tikzpicture}
\overset{\eqref{eq:obstruction}}{=}
\zeta_Y^{-1}\;
\begin{tikzpicture}[mb]
\draw[mstr] (0,0) -- (0,3.7); \draw[mstr] (1.4,0) -- (1.4,3.7);
\node[mcpn] at (0,0.5) {$\beta_Y$};
\node[mcpn] at (1.4,2.3) {$\beta_{\sP Y}$};
\blab{0}{Y}\blab{1.4}{\sP Y}\tlab{0}{3.7}{\sP^2Y}\tlab{1.4}{3.7}{\sP^3Y}
\end{tikzpicture}\,.
\]
Hence, by \eqref{eq:equivariant} for $\xi$ and the naturality of $\beta$,
\[
R^2(\Phi)=b^{-1}\sP(b)^{-1}\sP^2(\Phi)\,\sP(\xi)\,\xi
=(\sP(b)\,b)^{-1}\sP^2(\Phi)\,\beta_X
=(\sP(b)\,b)^{-1}\beta_Z\,\Phi
=\zeta_Y\,\Phi .
\]
As $\obstruction$ has finite order in $\Inv(\cC)$ and $A\mapsto\ell_A$ is a homomorphism
(Proposition~\ref{prop:loop-automorphism}), $\zeta_Y$ is a root of unity; and
$\zeta_Y=1$ when $\beta$ is coherent, as then $\obstruction=\one$
(Definition~\ref{def:coherent}).
So the eigenvalues of $R$, whose squares are $\zeta_Y$, have modulus one,
and $\operatorname{Tr}(R^{-1})=\overline{\operatorname{Tr}(R)}$.

\noindent\emph{The trace of $R$.}
By Proposition~\ref{prop:basis}, the coefficient of $\mob(X)$ in $\bdy(Y)=[Z,b]$ is the
closing $\operatorname{Tr}_X(Z,b)$ of Definition~\ref{def:TrX}.
Write the canonical element of $\Hom(X,Z)\otimes\Hom(Z,X)$ as
$\sum_if_i\otimes g^i$.
Then $\operatorname{Tr}_X(Z,b)=\sum_i\operatorname{tr}(\xi^{-1}\sP(g^i)\,b\,f_i)$.
As $b\,f=\sP(R^{-1}f)\,\xi$ and $g^iR^{-1}f_i\in\End(X)=\bC\id_X$, this is
$\sum_i\operatorname{tr}(g^iR^{-1}f_i)=\operatorname{Tr}(R^{-1})$ by \ref{item:E2}.
On the other hand, this coefficient is $\varepsilon(\mob(X)^\star\bdy(Y))$, and
$\varepsilon(u^\star)=\overline{\varepsilon(u)}$ (Proposition~\ref{prop:hilbert}(d)).
So
\[
\operatorname{Tr}(R)=\overline{\varepsilon\bigl(\mob(X)^\star\,\bdy(Y)\bigr)}
=\varepsilon\bigl(\bdy(Y)^\star\,\mob(X)\bigr)
=\varepsilon\bigl(\mob(X)\,\bdy(Y^\vee)\bigr),
\]
as $\bdy(Y)^\star=\bdy(Y^\vee)$ is central (Lemma~\ref{lem:L-star},
Proposition~\ref{prop:L}).
\end{proof}

\begin{definition}\label{def:InvT}
We set
\[
\Fix(\sT):=\{X\in\Irr(\cC)\mid X\cong(\sP X)^\vee\}.
\qedhere
\]
\end{definition}

\begin{remark}
The notation $\Fix(\sT)$ is motivated by the fact that $\Fix(\sT)$
is the subset of $\Irr(\cC)$ fixed, up to isomorphism, by
the contravariant monoidal functor $\sT:=\sP\circ(-)^\vee$.
For unitary $\cC$, composing with the dagger makes $\sT$ a covariant,
anti-linear monoidal functor, which is the time reversal of physics.
We do not need $\sT$ in the rest of the paper,
and the reader can think of $\Fix(\sT)$ as a single symbol. 
\end{remark}

\begin{definition}\label{def:eta}
For $Y\in\Irr(\cC)$ we put $\eta_Y:=\varepsilon(\bdy(Y^\vee))$.
\end{definition}

\begin{proposition}\label{prop:beta}
Let $Y\in\Irr(\cC)$.
\begin{enumerate}[label=(\alph*)]
\item If $Y\notin\Fix(\sT)$, then $\eta_Y=0$.
\item If $Y\in\Fix(\sT)$, then $\Hom(\one,Y\otimes\sP Y)$ is a line, and
  \[
  \eta_Y=\operatorname{Tr}(R_{U_0,b_Y})
  \]
  for the map $R_{U_0,b_Y}$ of Lemma~\ref{lem:fY} on this line; it is a
  root of unity.
  When $\beta$ is coherent, $R_{U_0,b_Y}$ is an involution, and
  $\eta_Y=\pm1$.
\item If $Y\in\Fix(\sT)$, then $\eta_{Y^\vee}=\overline{\eta_Y}=\eta_Y^{-1}$ and
  \[
  \eta_{Y^\vee}=\chi_Y([\obstruction])\,\eta_Y.
  \]
  When $\beta$ is coherent, $\eta_{Y^\vee}=\eta_Y$.
\end{enumerate}
\end{proposition}

\begin{proof}
Since $\mob(\one)$ is the unit of $\fX$ and $\xi_\one=U_0$, Lemma~\ref{lem:fY} for
$X=\one$ gives
$\eta_Y=\varepsilon(\bdy(Y^\vee))=\operatorname{Tr}(R_{U_0,b_Y})$,
where $R_{U_0,b_Y}$ acts on $\Hom(\one,Y\otimes\sP Y)\cong\Hom(Y^\vee,\sP Y)$.
As $\varepsilon(u^\star)=\overline{\varepsilon(u)}$ (Proposition~\ref{prop:hilbert}(d)),
Lemma~\ref{lem:L-star} gives
$\eta_{Y^\vee}=\varepsilon(\bdy(Y^\vee)^\star)=\overline{\eta_Y}$.

\medskip\noindent (a)
This space is zero unless $\sP Y\cong Y^\vee$, that is, unless $Y\in\Fix(\sT)$.

\medskip\noindent (b)
Otherwise the space is a line, on which $R_{U_0,b_Y}$ acts by a scalar $\rho$
with $\rho^2=\zeta_Y$, a root of unity (Lemma~\ref{lem:fY}); so
$\eta_Y=\rho$ is a root of unity.
When $\beta$ is coherent, $\zeta_Y=1$, so $\eta_Y=\rho=\pm1$.

\medskip\noindent (c)
We saw $\eta_{Y^\vee}=\overline{\eta_Y}=\eta_Y^{-1}$.
So $\eta_{Y^\vee}/\eta_Y=\eta_Y^{-2}=\overline{\zeta_Y}=\chi_Y([\obstruction])$,
the last since $\sP Y\cong Y^\vee$.
When $\beta$ is coherent, $\eta_Y=\pm1$ is real.
\end{proof}

\begin{definition}
  \label{def:MX}
We define $M_X$ for $X\in\Irr(\cC)$ by
\[M_X:=\sum_{Y\in \Fix(\sT)}S_{X,Y}\,\eta_Y.\qedhere\]
\end{definition}

\subsection{The structure of $\fX$}
\label{subsec:proof-cardy}

The aim of this section is to prove the following theorem,
which was stated in the introduction as Theorem~\ref{thm:intro}:
\begin{theorem}
  \label{thm:main}
\begin{enumerate}[label=(\alph*)]
\item The center of $\fX$ is $\bdy(\fF)$.
\item For each $X\in\Irr(\cC)$ with $p_X\neq 0$,
the algebra $p_X\fX$ is simple, 
\[p_X\fX\cong\End_\bC(V_X),\]
where $V_X$ is an irreducible $\fX$-module; we put $V_X:=0$ if $p_X=0$.
Furthermore,  every irreducible $\fX$-module arises in this manner for a unique $X$.
\item We have \[
\dim V_X=M_X.
\]
In
  particular the numbers $M_X$ in Definition~\ref{def:MX} are non-negative
  integers, and $\fX\cong\bigoplus_{X\in\Irr(\cC)}\Mat_{M_X}(\bC)$.
\end{enumerate}
\end{theorem}

The rest of the section is devoted to the proof of this theorem.
The essential idea is to construct
an averaging operator $\Cas:\fX\to\fX$
whose image is the center $Z(\fX)$, 
and show that $\Cas(\fX)$ is contained in $\bdy(\fF)$ 
using a manipulation of pictures.
As we already know $\bdy(\fF)\subset Z(\fX)$, we get Theorem~\ref{thm:main}(a), from which (b) follows.
For $u=\mob(\one)$ the same computation gives the relation
$\sum_X\mob(X)^\star \mob(X)=\sum_Y\eta_Y\bdy(Y)$, 
from which Theorem~\ref{thm:main}(c) follows by
computing the dimension of the blocks of $\fX$ in two ways.

Let us implement this idea.
We begin by constructing an averaging operator $\Cas$ on $\fX$.

\begin{definition}\label{def:Theta}
Let $\Cas\colon\fX\to\fX$ be the linear map
\[
\Cas(u):=\sum_{X\in\Fix(\sP)}\mob(X)^\star\,u\,\mob(X)\qquad(u\in\fX).
\qedhere
\]
\end{definition}

Let $V$ run over the
irreducible $\fX$-modules up to isomorphism, let $1_V$ be the unit of
the summand $\End_\bC(V)$ of $\fX\cong\bigoplus_V\End_\bC(V)$,
and let $\mu_V>0$ be the numbers with
$\varepsilon=\sum_V\mu_V\operatorname{Tr}_V$.
We then have:
\begin{lemma}\label{lem:cas}
$\Cas(u)=\sum_V\dfrac{\operatorname{Tr}_V(u)}{\mu_V}\,1_V$
for all $u\in\fX$.
Hence $\Cas(\fX)=Z(\fX)$.
\end{lemma}

\begin{proof}
By Proposition~\ref{prop:hilbert}(d),
$\varepsilon(\mob(X)^\star \mob(X'))=\delta_{X,X'}$; hence
$\Cas(u)=\sum_ib^i\,u\,b_i$ for any bases $\{b_i\}$,
$\{b^i\}$ of $\fX$ with $\varepsilon(b^ib_j)=\delta_{ij}$.
Choose matrix
units $E^V_{ij}$ of $\End_\bC(V)$.
Then $\{E^V_{ij}\}$ and
$\{E^V_{ji}/\mu_V\}$ are such bases, since
$\varepsilon(E^V_{ji}E^W_{kl})=\delta_{V,W}\delta_{ik}\delta_{jl}\mu_V$,
and
\[
\Cas(u)=\sum_V\frac1{\mu_V}\sum_{i,j}E^V_{ji}\,u\,E^V_{ij}
=\sum_V\frac{\operatorname{Tr}_V(u)}{\mu_V}\,1_V ,
\]
because $E^V_{ji}\,u\,E^V_{ij}=u^V_{ii}E^V_{jj}$, where $(u^V_{kl})$ is
the component of $u$ in $\End_\bC(V)$.
The $1_V$ form a basis of
$Z(\fX)$, and $\Cas(1_V)=(\dim V/\mu_V)\,1_V$, so the image
of $\Cas$ is $Z(\fX)$.
\end{proof}

\newcommand{\picI}{%
\begin{tikzpicture}[mb]
\mbars{-0.5}{2.6}{2.5}
\draw[mstr] (0,0) -- (0,0.9); \draw[mstr] (1.05,0) -- (1.05,0.9); \draw[mstr] (2.1,0) -- (2.1,0.9);
\draw[mstr] (2.1,0.9) -- (0,2.5);
\draw[mover] (1.05,0.9) .. controls (1.05,1.45) and (0.6,1.35) .. (0.65,1.7)
  .. controls (0.7,2.05) and (1.05,1.95) .. (1.05,2.5);
\draw[mover] (0,0.9) -- (2.1,2.5);
\node[mcpn] at (0,0.45) {$\xi_V^\natural$}; \node[mcpn] at (1.05,0.45) {$\xi$}; \node[mcpn] at (2.1,0.45) {$\xi_V$};
\blab{0}{V^\vee}\blab{1.05}{X}\blab{2.1}{V}
\tlab{0}{2.5}{\sP V}\tlab{1.05}{2.5}{\sP X}\tlab{2.1}{2.5}{\sP V^\vee}
\end{tikzpicture}}
\newcommand{\picIb}{%
\begin{tikzpicture}[mb,y=5.4mm]
\mbars{-0.5}{2.6}{5.3}
\draw[mstr] (0,0) -- (0,0.73) .. controls (0,1.12) and (0.15,1.3) .. (0.35,1.3)
  .. controls (0.5,1.3) and (0.5,2.1) .. (0.35,2.1) .. controls (0.15,2.1) and (0,2.28) .. (0,2.68) -- (0,3.8);
\draw[mstr] (2.1,0) -- (2.1,0.73) .. controls (2.1,1.12) and (1.8,1.3) .. (1.3,1.3) -- (0.85,1.3)
  .. controls (0.7,1.3) and (0.7,2.1) .. (0.85,2.1) -- (1.3,2.1)
  .. controls (1.8,2.1) and (2.1,2.28) .. (2.1,2.68) -- (2.1,3.8);
\draw[mover] (1.2,0) -- (1.2,3.8);
\draw[mstr] (2.1,3.8) -- (0,5.3);
\draw[mover] (1.2,3.8) .. controls (1.2,4.3) and (0.75,4.2) .. (0.8,4.55)
  .. controls (0.85,4.9) and (1.2,4.8) .. (1.2,5.3);
\draw[mover] (0,3.8) -- (2.1,5.3);
\node[mcpn] at (0,3.4) {$\xi_V^\natural$}; \node[mcpn] at (1.2,3.4) {$\xi$}; \node[mcpn] at (2.1,3.4) {$\xi_V$};
\blab{0}{V^\vee}\blab{1.2}{X}\blab{2.1}{V}
\tlab{0}{5.3}{\sP V}\tlab{1.2}{5.3}{\sP X}\tlab{2.1}{5.3}{\sP V^\vee}
\end{tikzpicture}}
\newcommand{\picII}{%
\begin{tikzpicture}[mb,y=5.4mm]
\mbars{-0.5}{2.6}{5.3}
\draw[mstr] (0.45,3.0) -- (0.45,0.5) .. controls (0.45,0.15) and (0.75,0.15) .. (0.75,0.5)
  -- (0.75,3.0) .. controls (0.75,3.35) and (0.45,3.35) .. (0.45,3.0);
\draw[mover] (2.1,0) -- (2.1,0.73) .. controls (2.1,1.12) and (1.65,1.3) .. (1.05,1.3)
          .. controls (0.45,1.3) and (0,1.12) .. (0,0.73) -- (0,0);
\draw[mover] (0,3.8) -- (0,2.68) .. controls (0,2.28) and (0.45,2.1) .. (1.05,2.1)
          .. controls (1.65,2.1) and (2.1,2.28) .. (2.1,2.68) -- (2.1,3.8);
\draw[mover] (0.75,0.9) -- (0.75,2.5);
\draw[mover] (1.2,0) -- (1.2,3.8);
\node[anchor=east,inner sep=1pt] at (0.45,1.7) {$\Omega$};
\draw[mstr] (2.1,3.8) -- (0,5.3);
\draw[mover] (1.2,3.8) .. controls (1.2,4.3) and (0.75,4.2) .. (0.8,4.55)
  .. controls (0.85,4.9) and (1.2,4.8) .. (1.2,5.3);
\draw[mover] (0,3.8) -- (2.1,5.3);
\node[mcpn] at (0,3.4) {$\xi_V^\natural$}; \node[mcpn] at (1.2,3.4) {$\xi$}; \node[mcpn] at (2.1,3.4) {$\xi_V$};
\blab{0}{V^\vee}\blab{1.2}{X}\blab{2.1}{V}
\tlab{0}{5.3}{\sP V}\tlab{1.2}{5.3}{\sP X}\tlab{2.1}{5.3}{\sP V^\vee}
\end{tikzpicture}}
\newcommand{\picIIIa}{%
\begin{tikzpicture}[mb,y=5.4mm]
\mbars{-0.6}{3.45}{5.3}
\draw[mstr] (0.45,2.45) -- (0.45,0.5) .. controls (0.45,0.15) and (0.75,0.15) .. (0.75,0.5)
  -- (0.75,2.45) .. controls (0.75,2.8) and (0.45,2.8) .. (0.45,2.45);
\draw[mover] (3,0) -- (3,0.73) .. controls (3,1.12) and (2.4,1.3) .. (1.5,1.3)
          .. controls (0.6,1.3) and (0,1.12) .. (0,0.73) -- (0,0);
\draw[mover] (0.75,1.05) -- (0.75,1.55);
\draw[mover] (0,2.6) .. controls (0,2.2) and (0.1,2.1) .. (0.5,2.1)
          .. controls (0.9,2.1) and (1.0,2.2) .. (1.0,2.6);
\draw[mover] (0.75,1.85) -- (0.75,2.35);
\draw[mstr] (0,2.6) -- (0,3.6); \draw[mstr] (1.0,2.6) -- (1.0,3.6);
\draw[mstr] (1.0,3.6) .. controls (1.0,4.0) and (0,3.9) .. (0,4.3) -- (0,5.3);
\draw[mover] (2.05,0) -- (2.05,5.3);
\draw[mover] (0,3.6) .. controls (0,4.0) and (1.0,3.9) .. (1.0,4.3)
  .. controls (1.0,4.8) and (3,4.8) .. (3,5.3);
\node[mcpn] at (0,3.2) {$\xi_V^\natural$}; \node[mcpn] at (1.0,3.2) {$\xi_V$}; \node[mcpn] at (2.05,3.2) {$\xi$};
\node[anchor=east,inner sep=1pt] at (0.42,1.75) {$\Omega$};
\blab{0}{V^\vee}\blab{2.05}{X}\blab{3.0}{V}
\tlab{0}{5.3}{\sP V}\tlab{2.05}{5.3}{\sP X}\tlab{3.0}{5.3}{\sP V^\vee}
\end{tikzpicture}}
\newcommand{\picIIIb}{%
\begin{tikzpicture}[mb,y=5.4mm]
\mbars{-0.6}{3.6}{5.5}
\draw[mstr] (1.45,5.0) -- (1.45,0.75) .. controls (1.45,0.42) and (1.7,0.42) .. (1.7,0.75)
  -- (1.7,5.0) .. controls (1.7,5.28) and (1.45,5.28) .. (1.45,5.0);
\draw[mover] (3.15,0) -- (3.15,0.73) .. controls (3.15,1.12) and (2.5,1.3) .. (1.55,1.3)
          .. controls (0.6,1.3) and (0,1.12) .. (0,0.73) -- (0,0);
\draw[mover] (1.7,1.05) -- (1.7,1.55);
\draw[mstr] (0,3.6) -- (0,2.6) .. controls (0,2.2) and (0.1,2.1) .. (0.5,2.1)
          .. controls (0.9,2.1) and (1.0,2.2) .. (1.0,2.6) -- (1.0,3.6);
\draw[mstr] (1.0,3.6) .. controls (1.0,4.0) and (0,3.9) .. (0,4.3) -- (0,5.5);
\draw[mover] (2.2,0) -- (2.2,5.5);
\draw[mover] (0,3.6) .. controls (0,4.0) and (1.0,3.9) .. (1.0,4.3)
  .. controls (1.0,4.85) and (3.15,4.85) .. (3.15,5.5);
\draw[mover] (1.45,4.42) -- (1.45,4.92);
\node[mcpn] at (0,3.2) {$\xi_V^\natural$}; \node[mcpn] at (1.0,3.2) {$\xi_V$}; \node[mcpn] at (2.2,3.2) {$\xi$};
\node[anchor=east,inner sep=1pt] at (1.42,1.85) {$\Omega$};
\blab{0}{V^\vee}\blab{2.2}{X}\blab{3.15}{V}
\tlab{0}{5.5}{\sP V}\tlab{2.2}{5.5}{\sP X}\tlab{3.15}{5.5}{\sP V^\vee}
\end{tikzpicture}}
\newcommand{\picIIIc}{%
\begin{tikzpicture}[mb,y=5.4mm]
\mbars{-0.6}{3.6}{5.5}
\draw[mstr] (1.45,5.0) -- (1.45,0.75) .. controls (1.45,0.42) and (1.7,0.42) .. (1.7,0.75)
  -- (1.7,5.0) .. controls (1.7,5.28) and (1.45,5.28) .. (1.45,5.0);
\draw[mover] (3.15,0) -- (3.15,0.73) .. controls (3.15,1.12) and (2.5,1.3) .. (1.55,1.3)
          .. controls (0.6,1.3) and (0,1.12) .. (0,0.73) -- (0,0);
\draw[mover] (1.7,1.05) -- (1.7,1.55);
\draw[mstr] (0,5.5) -- (0,2.6) .. controls (0,2.2) and (0.1,2.1) .. (0.5,2.1)
          .. controls (0.9,2.1) and (1.0,2.2) .. (1.0,2.6) -- (1.0,4.3);
\draw[mover] (2.2,0) -- (2.2,5.5);
\draw[mover] (1.0,4.3) .. controls (1.0,4.85) and (3.15,4.85) .. (3.15,5.5);
\draw[mover] (1.45,4.42) -- (1.45,4.92);
\node[mcpn] at (2.2,3.2) {$\xi$};
\node[anchor=east,inner sep=1pt] at (1.42,1.85) {$\Omega$};
\blab{0}{V^\vee}\blab{2.2}{X}\blab{3.15}{V}
\tlab{0}{5.5}{\sP V}\tlab{2.2}{5.5}{\sP X}\tlab{3.15}{5.5}{\sP V^\vee}
\end{tikzpicture}}
\newcommand{\picIII}{%
\begin{tikzpicture}[mb,x=6mm,y=5.4mm]
\mbars{-0.6}{3.45}{4.25}
\draw[mstr] (0.5,2.9) -- (0.5,0.6) .. controls (0.5,0.25) and (0.8,0.25) .. (0.8,0.6)
  -- (0.8,2.9) .. controls (0.8,3.25) and (0.5,3.25) .. (0.5,2.9);
\draw[mover] (3,0) -- (3,1.0) .. controls (3,1.42) and (2.4,1.6) .. (1.5,1.6)
          .. controls (0.6,1.6) and (0,1.42) .. (0,1.0) -- (0,0);
\draw[mover] (0.8,1.35) -- (0.8,1.85);
\draw[mover] (1.75,0) -- (1.75,4.25);
\draw[mover] (0,4.25) -- (0,3.1) .. controls (0,2.7) and (0.6,2.55) .. (1.5,2.55)
  .. controls (2.4,2.55) and (3,2.7) .. (3,3.1) -- (3,4.25);
\draw[mover] (0.5,2.3) -- (0.5,2.8);
\node[anchor=east,inner sep=1pt] at (0.45,2.08) {$\Omega$};
\node[mcpn] at (1.75,3.65) {$\xi$};
\blab{0}{V^\vee}\blab{1.75}{X}\blab{3.0}{V}
\tlab{0}{4.25}{\sP V}\tlab{1.75}{4.25}{\sP X}\tlab{3.0}{4.25}{\sP V^\vee}
\end{tikzpicture}}
\newcommand{\picIV}{%
\begin{tikzpicture}[mb,x=6mm,y=5.4mm]
\mbars{-0.6}{3.45}{4.25}
\draw[mstr] (0.5,2.9) -- (0.5,0.6) .. controls (0.5,0.25) and (0.8,0.25) .. (0.8,0.6)
  -- (0.8,2.9) .. controls (0.8,3.25) and (0.5,3.25) .. (0.5,2.9);
\draw[mover] (3,0) -- (3,1.0) .. controls (3,1.42) and (2.4,1.6) .. (1.5,1.6)
          .. controls (0.6,1.6) and (0,1.42) .. (0,1.0) -- (0,0);
\draw[mover] (0.8,1.35) -- (0.8,1.85);
\draw[mover] (1.75,0) -- (1.75,4.25);
\draw[mover] (0,4.25) -- (0,3.1) .. controls (0,2.7) and (0.6,2.55) .. (1.5,2.55)
  .. controls (2.4,2.55) and (3,2.7) .. (3,3.1) -- (3,4.25);
\draw[mover] (0.5,2.3) -- (0.5,2.8);
\node[anchor=east,inner sep=1pt] at (0.45,2.08) {$A$};
\node[anchor=west,inner sep=1pt] at (0.85,2.08) {$A^\vee$};
\node[mcpn] at (1.75,3.65) {$\xi$};
\blab{0}{V^\vee}\blab{1.75}{X}\blab{3.0}{V}
\tlab{0}{4.25}{\sP V}\tlab{1.75}{4.25}{\sP X}\tlab{3.0}{4.25}{\sP V^\vee}
\end{tikzpicture}}
\newcommand{\picIVa}{%
\begin{tikzpicture}[mb,x=6mm,y=5.4mm]
\mbars{-0.6}{3.6}{5.6}
\draw[mstr] (1.22,1.75) -- (1.22,0.45) .. controls (1.22,0.15) and (0.75,0.15) .. (0.75,0.45)
  -- (0.75,4.35) .. controls (0.75,4.65) and (1.22,4.65) .. (1.22,4.35) -- (1.22,3.15);
\draw[mover] (3,0) -- (3,0.55) .. controls (3,0.85) and (2.4,1.0) .. (1.5,1.0)
          .. controls (0.6,1.0) and (0,0.85) .. (0,0.55) -- (0,0);
\draw[mover] (1.22,0.75) -- (1.22,1.25);
\draw[mover] (1.77,0) -- (1.77,1.75);
\draw[mstr] (1.5,1.75) -- (1.5,3.15);
\draw[mstr] (1.77,3.15) -- (1.77,5.6);
\draw[mover] (0,5.6) -- (0,4.35) .. controls (0,4.1) and (0.6,3.95) .. (1.5,3.95)
  .. controls (2.4,3.95) and (3,4.1) .. (3,4.35) -- (3,5.6);
\draw[mover] (0.75,3.75) -- (0.75,4.25);
\node[mcpnc,minimum width=7mm] at (1.5,1.75) {$e'$};
\node[mcpnc,minimum width=7mm] at (1.5,3.15) {$e'$};
\node[mcpn] at (1.77,4.95) {$\xi$};
\node[anchor=east,inner sep=1pt] at (0.7,2.45) {$A$};
\node[anchor=west,inner sep=1pt] at (1.55,2.45) {$B$};
\blab{0}{V^\vee}\blab{1.77}{X}\blab{3.0}{V}
\tlab{0}{5.6}{\sP V}\tlab{1.77}{5.6}{\sP X}\tlab{3.0}{5.6}{\sP V^\vee}
\end{tikzpicture}}
\newcommand{\picIVb}{%
\begin{tikzpicture}[mb,x=6mm,y=5.4mm]
\mbars{-0.6}{3.6}{5.05}
\draw[mstr] (1.77,0) -- (1.77,0.6);
\draw[mstr] (1.22,0.6) -- (1.22,0.45) .. controls (1.22,0.15) and (0.75,0.15) .. (0.75,0.45)
  -- (0.75,4.0) .. controls (0.75,4.3) and (1.22,4.3) .. (1.22,4.0) -- (1.22,3.85);
\draw[mstr] (1.77,3.85) -- (1.77,5.05);
\draw[mover] (3,0) -- (3,1.0) .. controls (3,1.42) and (2.4,1.6) .. (1.5,1.6)
          .. controls (0.6,1.6) and (0,1.42) .. (0,1.0) -- (0,0);
\draw[mover] (1.5,1.0) -- (1.5,3.3);
\draw[mover] (0,5.05) -- (0,3.1) .. controls (0,2.7) and (0.6,2.55) .. (1.5,2.55)
  .. controls (2.4,2.55) and (3,2.7) .. (3,3.1) -- (3,5.05);
\draw[mover] (0.75,2.3) -- (0.75,2.8);
\node[mcpnc,minimum width=7mm] at (1.5,0.75) {$e'$};
\node[mcpnc,minimum width=7mm] at (1.5,3.55) {$e'$};
\node[mcpn] at (1.77,4.45) {$\xi$};
\node[anchor=east,inner sep=1pt] at (0.7,2.08) {$A$};
\node[anchor=west,inner sep=1pt] at (1.55,2.08) {$B$};
\blab{0}{V^\vee}\blab{1.77}{X}\blab{3.0}{V}
\tlab{0}{5.05}{\sP V}\tlab{1.77}{5.05}{\sP X}\tlab{3.0}{5.05}{\sP V^\vee}
\end{tikzpicture}}
\newcommand{\picV}{%
\begin{tikzpicture}[mb,x=6mm,y=5.4mm]
\mbars{-0.6}{3.6}{5.05}
\draw[mstr] (1.5,0) -- (1.5,0.45);
\draw[mstr] (1.25,0.75) .. controls (1.1,0.9) .. (1.1,1.1) -- (1.1,2.95) .. controls (1.1,3.1) .. (1.25,3.25);
\draw[mover] (3,0) -- (3,1.0) .. controls (3,1.42) and (2.4,1.6) .. (1.5,1.6)
          .. controls (0.6,1.6) and (0,1.42) .. (0,1.0) -- (0,0);
\draw[mover] (1.75,0.75) .. controls (1.9,0.9) .. (1.9,1.1) -- (1.9,2.95) .. controls (1.9,3.1) .. (1.75,3.25);
\draw[mover] (0,5.05) -- (0,3.1) .. controls (0,2.7) and (0.6,2.55) .. (1.5,2.55)
  .. controls (2.4,2.55) and (3,2.7) .. (3,3.1) -- (3,5.05);
\draw[mover] (1.1,2.3) -- (1.1,2.8);
\draw[mstr] (1.5,3.8) -- (1.5,5.05);
\node[mcpnc] at (1.5,0.6) {$e$}; \node[mcpnc] at (1.5,3.55) {$e$};
\node[mcpn] at (1.5,4.45) {$\xi$};
\node[anchor=east,inner sep=1pt] at (1.1,2.08) {$A$}; \node[anchor=west,inner sep=1pt] at (1.9,2.08) {$B$};
\blab{0}{V^\vee}\blab{1.5}{X}\blab{3.0}{V}
\tlab{0}{5.05}{\sP V}\tlab{1.5}{5.05}{\sP X}\tlab{3.0}{5.05}{\sP V^\vee}
\end{tikzpicture}}
\newcommand{\picVI}{%
\begin{tikzpicture}[mb,x=6mm,y=5.4mm]
\mbars{-0.6}{3.6}{6.45}
\draw[mstr] (2.0,0) -- (2.0,1.25) .. controls (2.0,1.4) .. (1.75,1.56);
\draw[mstr] (1.75,4.19) .. controls (2.0,4.35) .. (2.0,4.5) -- (2.0,6.45);
\draw[mover] (0,4.65) -- (0,1.3) .. controls (0,0.95) and (0.6,0.8) .. (1.5,0.8)
  .. controls (2.4,0.8) and (3,0.95) .. (3,1.3) -- (3,4.65);
\draw[mover] (0,4.65) .. controls (0,5.05) and (0.75,5.25) .. (1.5,5.25)
          .. controls (2.25,5.25) and (3.0,5.05) .. (3.0,4.65);
\draw[mover] (1.0,0) -- (1.0,1.25) .. controls (1.0,1.4) .. (1.25,1.56);
\draw[mover] (1.25,4.19) .. controls (1.0,4.35) .. (1.0,4.5) -- (1.0,6.45);
\draw[mstr] (1.5,2.34) -- (1.5,3.41);
\node[mcpnc] at (1.5,1.95) {$e$};
\node[mcpn] at (1.5,2.85) {$\xi$};
\node[mcpnc] at (1.5,3.8) {$e$};
\ubox{-0.4}{3.2}{3.4}{5.55}
\node[anchor=west,inner sep=1pt] at (0,3.85) {$V$}; \node[anchor=east,inner sep=1pt] at (3.0,3.85) {$V^\vee$};
\node[anchor=east,inner sep=1pt] at (1.0,4.65) {$B$}; \node[anchor=west,inner sep=1pt] at (2.0,4.65) {$A$};
\blab{1.0}{A}\blab{2.0}{B}\tlab{1.0}{6.45}{\sP B}\tlab{2.0}{6.45}{\sP A}
\end{tikzpicture}}
\newcommand{\picVIb}{%
\begin{tikzpicture}[mb,x=6mm,y=5.4mm]
\mbars{-0.6}{3.6}{6.45}
\draw[mstr] (2.0,0) -- (2.0,1.25) .. controls (2.0,1.4) .. (1.75,1.56);
\draw[mstr] (1.75,4.19) .. controls (2.0,4.35) .. (2.0,4.5) -- (2.0,6.45);
\draw[mover] (0,4.65) -- (0,1.3) .. controls (0,0.95) and (0.6,0.8) .. (1.5,0.8)
  .. controls (2.4,0.8) and (3,0.95) .. (3,1.3) -- (3,4.65)
  .. controls (3,5.05) and (2.4,5.25) .. (1.5,5.25) .. controls (0.6,5.25) and (0,5.05) .. (0,4.65);
\draw[mover] (1.0,0) -- (1.0,1.25) .. controls (1.0,1.4) .. (1.25,1.56);
\draw[mover] (1.25,4.19) .. controls (1.0,4.35) .. (1.0,4.5) -- (1.0,6.45);
\draw[mstr] (1.5,2.34) -- (1.5,3.41);
\node[mcpnc] at (1.5,1.95) {$e$};
\node[mcpn] at (1.5,2.85) {$\xi$};
\node[mcpnc] at (1.5,3.8) {$e$};
\ubox{0.75}{3.3}{2.25}{4.3}
\node[anchor=east,inner sep=1pt] at (0,2.85) {$\sP V$};
\blab{1.0}{A}\blab{2.0}{B}\tlab{1.0}{6.45}{\sP B}\tlab{2.0}{6.45}{\sP A}
\end{tikzpicture}}
\newcommand{\picVII}{%
\begin{tikzpicture}[mb,x=6.5mm]
\mbars{-1.35}{2.6}{3.5}
\draw[mstr] (1,0) .. controls (1,0.3) .. (0.7,0.45);
\draw[mstr] (0.7,2.4) .. controls (0.85,2.6) .. (1,2.85) -- (1,3.5);
\draw[mover] (0.5,1.75) ellipse (1.6 and 1.62);
\draw[mover] (0,0) .. controls (0,0.3) .. (0.3,0.45);
\draw[mover] (0.3,2.4) .. controls (0.15,2.6) .. (0,2.85) -- (0,3.5);
\draw[mstr] (0.5,0.9) -- (0.5,1.95);
\node[mcpnc] at (0.5,0.6) {$e$};
\node[mcpn] at (0.5,1.35) {$\xi$};
\node[mcpnc] at (0.5,2.25) {$e$};
\ubox{-0.3}{1.75}{1.3}{2.75}
\node at (2.35,1.75) {$\Omega$};
\blab{0}{A}\blab{1}{B}\tlab{0}{3.5}{\sP B}\tlab{1}{3.5}{\sP A}
\end{tikzpicture}}
\newcommand{\picVIIa}{%
\begin{tikzpicture}[mb,x=6.5mm]
\mbars{-1.0}{3.2}{3.5}
\draw[mstr] (0,0) .. controls (0,0.3) .. (0.3,0.45);
\draw[mstr] (0.3,2.4) .. controls (0.15,2.6) .. (0,2.85) -- (0,3.5);
\draw[mstr] (0.5,0.9) -- (0.5,1.95);
\draw[mstr] (1.85,2.95) arc[start angle=90, end angle=270, x radius=0.22, y radius=1.4];
\draw[mover] (0.7,2.4) .. controls (0.8,2.75) and (2.7,2.45) .. (2.7,3.5);
\draw[mover] (2.7,0) .. controls (2.7,0.95) and (0.8,0.1) .. (0.7,0.45);
\draw[mover] (1.85,0.15) arc[start angle=-90, end angle=90, x radius=0.22, y radius=1.4];
\node[mcpnc] at (0.5,0.6) {$e$};
\node[mcpn] at (0.5,1.35) {$\xi$};
\node[mcpnc] at (0.5,2.25) {$e$};
\ubox{-0.3}{1.75}{1.3}{2.75}
\node[anchor=west,inner sep=1pt] at (2.07,1.55) {$\Omega$};
\blab{0}{A}\blab{2.7}{B}\tlab{0}{3.5}{\sP B}\tlab{2.7}{3.5}{\sP A}
\end{tikzpicture}}
\newcommand{\picVIIb}{%
\begin{tikzpicture}[mb,x=6.5mm]
\mbars{-1.0}{3.6}{3.5}
\draw[mstr] (0,0) .. controls (0,0.3) .. (0.3,0.45);
\draw[mstr] (0.3,2.4) .. controls (0.15,2.6) .. (0,2.85) -- (0,3.5);
\draw[mstr] (0.5,0.9) -- (0.5,1.95);
\draw[mstr] (0.7,2.4) .. controls (0.85,2.65) and (1.45,2.75) .. (1.65,2.4) -- (1.65,0.6)
  .. controls (1.45,0.2) and (0.85,0.25) .. (0.7,0.45);
\draw[mstr] (3.1,0) -- (3.1,3.5);
\node[mcpnc] at (0.5,0.6) {$e$};
\node[mcpn] at (0.5,1.35) {$\xi$};
\node[mcpnc] at (0.5,2.25) {$e$};
\ubox{-0.3}{1.75}{1.3}{2.75}
\node[mcpnc] at (1.65,1.35) {$e''$};
\node[mcpnc] at (3.1,1.35) {$e''$};
\node[anchor=west,inner sep=1pt] at (1.72,2.05) {$\sP A^\vee$};
\node[anchor=west,inner sep=1pt] at (1.65,0.55) {$B^\vee$};
\blab{0}{A}\blab{3.1}{B}\tlab{0}{3.5}{\sP B}\tlab{3.1}{3.5}{\sP A}
\end{tikzpicture}}
\newcommand{\picVIIbY}{%
\begin{tikzpicture}[mb,x=6.5mm]
\mbars{-1.0}{3.6}{3.5}
\draw[mstr] (0,0) .. controls (0,0.3) .. (0.3,0.45);
\draw[mstr] (0.3,2.4) .. controls (0.15,2.6) .. (0,2.85) -- (0,3.5);
\draw[mstr] (0.5,0.9) -- (0.5,1.95);
\draw[mstr] (0.7,2.4) .. controls (0.85,2.65) and (1.45,2.75) .. (1.65,2.4) -- (1.65,0.6)
  .. controls (1.45,0.2) and (0.85,0.25) .. (0.7,0.45);
\draw[mstr] (3.1,0) -- (3.1,3.5);
\node[mcpnc] at (0.5,0.6) {$e$};
\node[mcpn] at (0.5,1.35) {$\xi$};
\node[mcpnc] at (0.5,2.25) {$e$};
\ubox{-0.3}{1.75}{1.3}{2.75}
\node[mcpnc] at (1.65,1.35) {$e''$};
\node[mcpnc] at (3.1,1.35) {$e''$};
\node[anchor=west,inner sep=1pt] at (1.72,2.05) {$\sP Y^\vee$};
\node[anchor=west,inner sep=1pt] at (1.65,0.55) {$\sP Y^\vee$};
\blab{0}{Y}\blab{3.1}{\sP Y}\tlab{0}{3.5}{\sP^2Y}\tlab{3.1}{3.5}{\sP Y}
\end{tikzpicture}}
\newcommand{\picVIIc}{%
\begin{tikzpicture}[mb,x=6.5mm]
\fill[rounded corners=3pt,fill=mshade] (-0.5,0.15) rectangle (1.85,2.95);
\begin{scope}
\clip (-0.5,0.15) rectangle (1.85,2.95);
\draw[white,line width=2.2pt] (0,0) .. controls (0,0.3) .. (0.3,0.45)
  (0.3,2.4) .. controls (0.15,2.6) .. (0,2.85) -- (0,3.5)
  (0.5,0.9) -- (0.5,1.95)
  (0.7,2.4) .. controls (0.85,2.65) and (1.45,2.75) .. (1.65,2.4) -- (1.65,0.6)
  .. controls (1.45,0.2) and (0.85,0.25) .. (0.7,0.45);
\end{scope}
\mbars{-1.0}{3.6}{3.5}
\draw[mstr] (0,0) .. controls (0,0.3) .. (0.3,0.45);
\draw[mstr] (0.3,2.4) .. controls (0.15,2.6) .. (0,2.85) -- (0,3.5);
\draw[mstr] (0.5,0.9) -- (0.5,1.95);
\draw[mstr] (0.7,2.4) .. controls (0.85,2.65) and (1.45,2.75) .. (1.65,2.4) -- (1.65,0.6)
  .. controls (1.45,0.2) and (0.85,0.25) .. (0.7,0.45);
\draw[mstr] (3.1,0) -- (3.1,3.5);
\node[mcpnc] at (0.5,0.6) {$e$};
\node[mcpn] at (0.5,1.35) {$\xi$};
\node[mcpnc] at (0.5,2.25) {$e$};
\ubox{-0.3}{1.75}{1.3}{2.75}
\node[anchor=west,inner sep=1pt] at (1.87,1.35) {$\sP Y^\vee$};
\node[anchor=south,inner sep=1pt] at (1.3,2.97) {$\varphi$};
\blab{0}{Y}\blab{3.1}{\sP Y}\tlab{0}{3.5}{\sP^2Y}\tlab{3.1}{3.5}{\sP Y}
\end{tikzpicture}}
\newcommand{\picIXa}{%
\begin{tikzpicture}[mb,x=6.5mm]
\draw[mstr] (0,0.2) -- (0,0.5) .. controls (0,0.8) .. (0.3,0.95);
\draw[mstr] (1,0.2) -- (1,0.5) .. controls (1,0.8) .. (0.7,0.95);
\draw[mstr] (0.5,1.4) -- (0.5,2.45);
\draw[mstr] (0.3,2.9) .. controls (0.15,3.1) .. (0,3.35) -- (0,4.4);
\draw[mstr] (0.7,2.9) .. controls (0.85,3.1) .. (1,3.35) -- (1,4.4);
\draw[mstr] (0,4.4) .. controls (0,4.9) and (-0.9,4.9) .. (-0.9,4.4) -- (-0.9,0.2)
  .. controls (-0.9,-0.3) and (0,-0.3) .. (0,0.2);
\draw[mstr] (1,4.4) .. controls (1,4.9) and (1.9,4.9) .. (1.9,4.4) -- (1.9,0.2)
  .. controls (1.9,-0.3) and (1,-0.3) .. (1,0.2);
\node[mcpnc] at (0.5,1.1) {$e$};
\node[mcpn] at (0.5,1.85) {$\xi$};
\node[mcpnc] at (0.5,2.75) {$e$};
\node[mcpn] at (0,4.0) {$\beta_Y^{-1}$};
\ubox{-0.3}{2.25}{1.3}{3.25}
\node[anchor=east,inner sep=1pt] at (0,0.35) {$Y$};
\node[anchor=west,inner sep=1pt] at (1,0.35) {$\sP Y$};
\end{tikzpicture}}
\newcommand{\picIXb}{%
\begin{tikzpicture}[mb,x=6.5mm]
\mbars{-0.5}{1.55}{2.2}
\draw[mstr] (0,0) -- (0,2.2); \draw[mstr] (1.05,0) -- (1.05,2.2);
\node[mcpn] at (0,1.1) {$\beta_Y$};
\blab{0}{Y}\blab{1.05}{\sP Y}\tlab{0}{2.2}{\sP^2Y}\tlab{1.05}{2.2}{\sP Y}
\end{tikzpicture}}

Our central statement is the following proposition,
which rewrites $\Cas(\mob(X))$ in terms of an explicit linear combination of $\bdy(Y)$:

\begin{proposition}\label{prop:two-steps}
For $X\in\Fix(\sP)$,
\[
\Cas(\mob(X))=\sum_{Y\in\Irr(\cC)}\varepsilon\bigl(\mob(X)\,\bdy(Y^\vee)\bigr)\,\bdy(Y).
\]
\end{proposition}

\begin{remark}
  The proof below is a long sequence of elementary picture manipulations,
  but the reader might wonder why such a statement is expected to hold in the first place.
  Here we give a brief explanation, which is more conceptual
  for readers familiar with the general theory of 
  modular tensor categories
  and their equivariantizations;
  those who are happy with the picture proof can skip this remark.
  Assume that $\beta$ is coherent.
  As mentioned in Remark~\ref{rem:equivariantization},
  our M\"obius fusion algebra $\fX$ can be obtained by
  considering the fusion algebra of $\cD:=\cC^{\bZ/2}$,
  the equivariantization of $\cC$ 
  by the monoidal (and not anti-monoidal) $\bZ/2$-action given by $(\sP,J)$
  of Remark~\ref{rem:equivariantization},
  and restricting to the subspace
  where the equivariant object $(\one,-U_0)$ acts by $-1$.
  Now, in general, the center of the fusion algebra of a fusion category $\cD$ 
  is known to be determined by the Drinfeld center $Z(\cD)$ of $\cD$,
  and the averaging map $\Cas$ in this context is known to be given by
  $F\circ I$,
  where $F$ is the forgetful functor from $Z(\cD)$ to $\cD$ and $I$ is its adjoint. 
Indeed,  $F(I(X))\simeq\bigoplus_{V\in\Irr(\cD)}V^\vee\otimes X\otimes V$.
  When $\cD=\cC^{\bZ/2}$, the Drinfeld center $Z(\cD)$ 
  can be explicitly computed, following \cite{GNN2009}.
  Combining all these and filling in the details, 
  this argument should provide an alternative proof, at least when $\beta$ is coherent.
  In this paper, the authors instead decided to use picture sequences in the proof below, 
  as this is more elementary and does not require the assumption that $\beta$ is coherent.
\end{remark}

\begin{proof}[Proof\/\protect\footnotemark{} of Proposition~\ref{prop:two-steps}]
\footnotetext{An animated version of this proof is available at \url{https://ippo0127.github.io/movies/mobius_prop4_23.html}.}%
Write $\xi:=\xi_X$, so that $\mob(X)=[X,\xi]$.

\medskip\noindent\emph{Step 1.}
We first make the following moves.
\begin{gather*}
\Cas([X,\xi])
\overset{\text{Def.~\ref{def:Theta}}}{\underset{\text{Prop.~\ref{prop:hilbert}(c)}}{=}}
\sum_{V\in\Fix(\sP)}
\picI
=\sum_{V\in\Fix(\sP)}
\picIb
\\
\overset{\text{Lem.~\ref{lem:kill}}}{=}\sum_{V\in\Fix(\sP)}\frac{d_V}{D^2}
\picII
=
\sum_{V\in\Fix(\sP)}\frac{d_V}{D^2}
\picIIIa
\\
\overset{\text{isotopy}}{=}
\sum_{V\in\Fix(\sP)}\frac{d_V}{D^2}
\picIIIb
\overset{\text{Lem.~\ref{lem:natural-dual}(b)}}{=}
\sum_{V\in\Fix(\sP)}\frac{d_V}{D^2}
\picIIIc
\\
\overset{\text{isotopy}}{=}
\sum_{V\in\Fix(\sP)}\frac{d_V}{D^2}
\picIII
=
\sum_{V\in\Irr(\cC)}\frac{d_V}{D^2}
\picIII
\end{gather*}
At the last step, 
the sum is extended from $\Fix(\sP)$ to $\Irr(\cC)$ because for
$V\notin\Fix(\sP)$ the picture vanishes:
the two strands passing through the $\Omega$-loop are the strand $V$ of the
cap and the leg $\sP V^\vee$ of the cup, and
the loop cuts them to zero by Lemma~\ref{lem:kill}, because
$\Hom(\one,V^\vee\otimes\sP V)\cong\Hom(V,\sP V)=0$.

\medskip\noindent\emph{Step 2.}
Below, we use the following canonical elements
for shaded coupons:
for $(A,B)\in\Irr(\cC)^2$ let $e$ be the canonical element of
$\Hom(X,A\otimes B)\otimes\Hom(A\otimes B,X)$, $e'$ that of
$\Hom(A^\vee\otimes X,B)\otimes\Hom(B,A^\vee\otimes X)$, and $e''$ that of
$\Hom(B,\sP A)\otimes\Hom(\sP A,B)$.
The end of $e''$ in $\Hom(\sP A,B)$ will sit on a strand running downward, where it is
drawn as its transpose $B^\vee\to\sP A^\vee$, see \eqref{eq:transpose}.
Then
\begin{gather*}
\sum_{V\in\Irr(\cC)}\frac{d_V}{D^2}
\picIII
\overset{\text{Def.~\ref{def:Omega}}}{=}\sum_{V\in\Irr(\cC)}\sum_{A}\frac{d_Ad_V}{D^2}
\picIV
\\
\overset{\text{\ref{item:E3}}}{=}\sum_{V\in\Irr(\cC)}\sum_{A,B}\frac{d_Ad_Bd_V}{D^2}
\picIVa
\overset{\text{nat.}}{=}\sum_{V\in\Irr(\cC)}\sum_{A,B}\frac{d_Ad_Bd_V}{D^2}
\picIVb
\\
\overset{\text{\ref{item:E4}}}{=}\sum_{V\in\Irr(\cC)}\sum_{A,B}\frac{d_Ad_Bd_V}{D^2}
\picV
\overset{\text{Lem.~\ref{lem:trace-property}}}{=}
\sum_{V\in\Irr(\cC)}\sum_{A,B}\frac{d_Ad_Bd_V}{D^2}
\picVI
\\
\overset{\text{\ref{item:B1},\ref{item:B2}}}{\underset{\text{\ref{item:B4},\ref{item:B6}}}{=}}
\sum_{V\in\Irr(\cC)}\sum_{A,B}\frac{d_Ad_Bd_V}{D^2}
\picVIb
\overset{\text{Rem.~\ref{rem:Omega-mirror}}}{=}
\sum_{A,B}\frac{d_Ad_B}{D^2}
\picVII
\\
\overset{\text{isotopy}}{=}\sum_{A,B}\frac{d_Ad_B}{D^2}
\picVIIa
\overset{\text{Lem.~\ref{lem:kill}}}{\underset{\text{\ref{item:E4}}}{=}}
\sum_{A,B}d_Ad_B
\picVIIb
\\
\overset{\mathrm{(i)}}{=}\sum_{Y\in\Irr(\cC)}d_Y^2
\picVIIbY
\overset{\mathrm{(ii)}}{=}\sum_{Y\in\Irr(\cC)}d_Y
\picVIIc
\overset{\mathrm{(iii)}}{=}\sum_{Y\in\Irr(\cC)}
\picIXa
\;
\picIXb\,.
\end{gather*}

\noindent The explanations of the last three equalities are as follows:

(i) The right strand of the left hand side of (i) is a morphism $B\to\sP A$, so the summand vanishes unless $B\cong\sP A$.
For such a pair put $Y:=A$,
and we may take $B=\sP Y$; then $d_Ad_B=d_Y^2$.

(ii) As $\End(\sP Y)$ is spanned by $\id_{\sP Y}$, and
$\operatorname{tr}(\id_{\sP Y})=d_{\sP Y}=d_Y$ (Corollary~\ref{cor:P-dim}),
we have
$e''=d_Y^{-1}\,\id_{\sP Y}\otimes\id_{\sP Y}$ (Definition~\ref{def:repeated}). Therefore the two coupons $e''$ disappear with a factor $d_Y^{-1}$.

(iii) The shaded part of the left hand side of (iii) 
is a morphism $\varphi\colon Y\to\sP^2Y$; as $Y$ is simple,
$\varphi=z\,\beta_Y$ for a scalar $z$.
Then the picture on the left hand side is $[Y\otimes\sP Y,U_{Y,\sP Y}(\varphi\otimes\id_{\sP Y})]=z\,\bdy(Y)$
(Definition~\ref{def:L}), so the summand is $z\,d_Y\,\bdy(Y)$.
On the right hand side, the coefficient of $\bdy(Y)$
is the result of composing $\varphi$ with $\beta_Y^{-1}$ and closing, 
which gives $z\,d_Y$.
This shows the equality (iii).

\medskip\noindent\emph{Step 3.}
It remains to evaluate the coefficient of $\bdy(Y)$.
Fix $Y\in\Irr(\cC)$.
As $b_Y^{-1}=(\beta_Y\otimes\id)^{-1}U^{-1}_{Y,\sP Y}$ (Definition~\ref{def:L}), the map
$R_{\xi,b_Y}$ of Lemma~\ref{lem:fY} is
$R_{\xi,b_Y}(\Phi)=(\beta_Y\otimes\id)^{-1}\Ubox{\Phi}\,\xi$ for $\Phi\in\Hom(X,Y\otimes\sP Y)$, so
by \ref{item:E2} the closed picture at the end of Step~2 is $\operatorname{Tr}(R_{\xi,b_Y})$. 
As $\operatorname{Tr}(R_{\xi,b_Y})=\varepsilon(\mob(X)\,\bdy(Y^\vee))$ by
Lemma~\ref{lem:fY}, this proves the proposition.
\end{proof}

\begin{corollary}\label{cor:Cas1}
\[
\Cas(\mob(\one))=\sum_{X\in\Fix(\sP)}\mob(X)^\star \mob(X)
=\sum_{Y\in\Fix(\sT)}\eta_Y\,\bdy(Y).
\]
\end{corollary}

\begin{proof}
Take $X=\one$ in Proposition~\ref{prop:two-steps}: $\mob(\one)$ is the unit, and
$\varepsilon(\bdy(Y^\vee))=\eta_Y$ (Definition~\ref{def:eta}), which is $0$
for $Y\notin\Fix(\sT)$ (Proposition~\ref{prop:beta}(a)).
\end{proof}


We can now prove Theorem~\ref{thm:main}.
We need one final lemma.
\begin{lemma}
\label{lem:eps-pX}
We have $\varepsilon(p_X)=S_{\one,X}\,M_X$.
\end{lemma}
\begin{proof}
Since $p_X=\bdy(e_X)$ and $\bdy$ is linear, \eqref{eq:fusion-idempotent} and
Proposition~\ref{prop:beta} give
\[
\varepsilon(p_X)=S_{\one,X}\sum_{A\in\Irr(\cC)}\overline{S_{A,X}}\,\varepsilon(\bdy(A))
=S_{\one,X}\sum_{A\in\Irr(\cC)}S_{A^\vee,X}\,\eta_{A^\vee}
=S_{\one,X}\,M_X,
\]
using $\overline{S_{A,X}}=S_{A^\vee,X}$ and Definition~\ref{def:eta}, and, in the
last step, the reindexing $A\mapsto A^\vee$, the symmetry of $S$, and
Proposition~\ref{prop:beta}(a).
\end{proof}
\begin{proof}[Proof of Theorem~\ref{thm:main}]
(a) By Proposition~\ref{prop:L}, $\bdy(\fF)\subset Z(\fX)$.
Conversely,
$Z(\fX)=\Cas(\fX)$ by Lemma~\ref{lem:cas}, and
$\Cas(\mob(X))\in \bdy(\fF)$ for every $X\in\Fix(\sP)$ by
Proposition~\ref{prop:two-steps}; since the $\mob(X)$ span $\fX$, we get
$\Cas(\fX)\subset \bdy(\fF)$.

\medskip\noindent (b) By (a) and Proposition~\ref{prop:L}, the non-zero $p_X$, which span
$\bdy(\fF)$, form a basis of $Z(\fX)=\bigoplus_V\bC1_V$ consisting of orthogonal idempotents.
Such a basis is $\{1_V\}$: each $p_X$ is a sum of some $1_V$, these sets are disjoint, and
there are as many $p_X$ as $1_V$.
We write $V_X$ for the $V$ with $p_X=1_V$; then $p_X\fX=1_{V_X}\fX\cong\End_\bC(V_X)$.

\medskip\noindent (c)
If $p_X=0$, then $M_X=0$ by Lemma~\ref{lem:eps-pX}.
Otherwise $p_X=1_{V_X}$ by (b); let $\mu:=\mu_{V_X}$.
Multiplying Lemma~\ref{lem:cas} for $u=\mob(\one)$ and Corollary~\ref{cor:Cas1} by $p_X$,
\[
\frac{\dim V_X}{\mu}\,p_X
\overset{\mathrm{(1)}}{=}\Cas(\mob(\one))\,p_X
\overset{\mathrm{(2)}}{=}\sum_{Y\in\Fix(\sT)}\eta_Y\,\bdy(Y)\,p_X
\overset{\mathrm{(3)}}{=}\frac{M_X}{S_{\one,X}}\,p_X .
\]

(1) Lemma~\ref{lem:cas} for $u=\mob(\one)$.

(2) Corollary~\ref{cor:Cas1}.

(3) $\bdy(Y)\,p_X=\chi_X([Y])\,p_X$ with $\chi_X([Y])=S_{Y,X}/S_{\one,X}$, by
Proposition~\ref{prop:L} and \eqref{eq:fusion-idempotent}; then $S_{Y,X}=S_{X,Y}$ and the
definition of $M_X$.

\noindent Also $\mu\dim V_X=\varepsilon(p_X)=S_{\one,X}M_X$ by
Proposition~\ref{prop:hilbert}(e) and Lemma~\ref{lem:eps-pX}.
Multiplying the two identities gives $(\dim V_X)^2=M_X^2$, and $M_X>0$ by the first, so
$\dim V_X=M_X$.
\end{proof}

\section{Bimodules over $\fX$ and the twist}
\label{sec:bimodules}

In this section
we attach to every object $E$ an $\fX$-bimodule $\cH_E$.
Let $\Sigma:=\{X\in\Irr(\cC)\mid V_X\neq0\}$, where $V_X$ is the irreducible
$\fX$-module of Theorem~\ref{thm:main}.
In
Sections~\ref{subsec:W} and~\ref{subsec:tube} we determine its structure,
$\cH_E\cong\bigoplus_{X,X'\in\Sigma}(V_X\otimes V_{X'}^*)^{\oplus N_{EX'}^X}$, and in
Section~\ref{subsec:theta} we use it to relate the twist $\theta_X$ on all $X$ with
$V_X\neq0$ to a central element $\mob(\obstruction)$; when $\beta$ is coherent, $\theta_X$
is a sign independent of such $X$.

\subsection{The bimodules $\cH_E$}
\label{subsec:W}

\begin{definition}\label{def:W}
For $E\in\cC$ let $\cH_E$ be the quotient of
$\bigoplus_{Z\in\Ob(\cC)}\Hom(Z,E\otimes\sP Z)$ by the relations (R1)
(direct sums) and (R2) (coinvariance) of Definition~\ref{def:X}, the latter in
the form \[
[Z,(\id_E\otimes\sP(a))g]=[Z,g\circ a],\quad a\in\End(Z).
\] 
Thus $\cH_\one=\fX$.
We draw
$[Z,g]$ as the picture of $g$ in which the strand $E$
leaves from the left side of the M\"obius band:
\[
[Z,g]=\picZg\,.
\qedpic
\]
\end{definition}

\begin{lemma}\label{lem:W-bimodule}
\begin{enumerate}[label=(\alph*)]
\item 
For $a\colon Z\to Z'$ and
  $h\colon Z'\to E\otimes\sP Z$, we have
\[
\begin{tikzpicture}[mb]
\mbars{-1}{0.6}{2.9}\medge{-1}{0}{2.9}
\draw[mstr] (0,0) -- (0,1.9);
\draw[mstr] (0.15,1.9) -- (0.15,2.9);
\draw[mstr] (-0.15,1.9) .. controls (-0.15,2.45) and (-0.6,2.45) .. (-1,2.45);
\node[mcpn] at (0,0.8) {$a$};
\node[mcpn] at (0,1.9) {$h$};
\legend{-1}{2.45}{E}
\blab{0}{Z}\tlab{0.15}{2.9}{\sP Z}
\end{tikzpicture}
=
\begin{tikzpicture}[mb]
\mbars{-1}{0.75}{2.9}\medge{-1}{0}{2.9}
\draw[mstr] (0,0) -- (0,0.8);
\draw[mstr] (0.15,0.8) -- (0.15,2.9);
\draw[mstr] (-0.15,0.8) .. controls (-0.15,1.3) and (-0.6,1.3) .. (-1,1.3);
\node[mcpn] at (0,0.8) {$h$};
\node[mcpn] at (0.15,1.95) {$a$};
\ubox{-0.3}{1.5}{0.6}{2.4}
\legend{-1}{1.3}{E}
\blab{0}{Z'}\tlab{0.15}{2.9}{\sP Z'}
\end{tikzpicture}\,.
\]
\item 
For $f\colon X_1\otimes\cdots\otimes X_n\to E\otimes\sP X_n\otimes\cdots\otimes\sP X_1$ and $1\le i<n$,
\[
\begin{tikzpicture}[mb]
\mbars{-1.3}{5.1}{2.4}\medge{-1.3}{0}{2.4}
\foreach \x in {0,1.7,2.9,4.6} \draw[mstr] (\x,0) -- (\x,2.4);
\draw[mstr] (-0.2,1.2) .. controls (-0.2,1.8) and (-0.6,1.8) .. (-1.3,1.8);
\node[mcpn, minimum width=45mm] at (2.3,1.2) {$f$};
\legend{-1.3}{1.8}{E}
\blab{0}{X_1}\blab{0.85}{\cdots}\blab{1.7}{X_i}\blab{2.9}{X_{i+1}}\blab{3.75}{\cdots}\blab{4.6}{X_n}
\tlab{0}{2.4}{\sP X_n}\tlab{0.85}{2.4}{\cdots}\tlab{1.7}{2.4}{\sP X_{i+1}}\tlab{2.9}{2.4}{\sP X_i}\tlab{3.75}{2.4}{\cdots}\tlab{4.6}{2.4}{\sP X_1}
\end{tikzpicture}
\;=\;
\begin{tikzpicture}[mb]
\mbars{-1.3}{5.1}{3.2}\medge{-1.3}{0}{3.2}
\draw[mstr] (0,0) -- (0,3.2); \draw[mstr] (4.6,0) -- (4.6,3.2);
\sstr{1.7}{0}{2.9}{1.1}
\sstr[mover]{2.9}{0}{1.7}{1.1}
\sstr{1.7}{2.1}{2.9}{3.2}
\sstr[mover]{2.9}{2.1}{1.7}{3.2}
\draw[mstr] (1.7,1.1) -- (1.7,2.1); \draw[mstr] (2.9,1.1) -- (2.9,2.1);
\draw[mstr] (-0.2,1.6) .. controls (-0.2,2.2) and (-0.6,2.2) .. (-1.3,2.2);
\node[mcpn, minimum width=45mm] at (2.3,1.6) {$f$};
\legend{-1.3}{2.2}{E}
\blab{0}{X_1}\blab{0.85}{\cdots}\blab{1.7}{X_{i+1}}\blab{2.9}{X_i}\blab{3.75}{\cdots}\blab{4.6}{X_n}
\tlab{0}{3.2}{\sP X_n}\tlab{0.85}{3.2}{\cdots}\tlab{1.7}{3.2}{\sP X_i}\tlab{2.9}{3.2}{\sP X_{i+1}}\tlab{3.75}{3.2}{\cdots}\tlab{4.6}{3.2}{\sP X_1}
\end{tikzpicture}
\;.
\]
\item The rules
\[
\mone{h}{Z}{\sP Z}\;\cdot\;\picXg
\;:=\;
\begin{tikzpicture}[mb]
\mbars{-1.2}{1.65}{2.8}\medge{-1.2}{0}{2.8}
\draw[mstr] (1.0,0) -- (1.0,0.8);
\draw[mstr] (1.15,0.8) -- (1.15,1.8);
\draw[mstr] (0.85,0.8) .. controls (0.85,1.35) and (0.4,1.35) .. (-1.2,1.35);
\draw[mover] (0,0) -- (0,1.8);
\sstr{1.15}{1.8}{0}{2.8}
\sstr[mover]{0}{1.8}{1.15}{2.8}
\node[mcpn] at (0,0.6) {$h$};
\node[mcpn] at (1.0,0.8) {$g$};
\legend{-1.2}{1.35}{E}
\blab{0}{Z}\blab{1.0}{X}\tlab{0}{2.8}{\sP X}\tlab{1.15}{2.8}{\sP Z}
\end{tikzpicture}\,,
\qquad
\picXg\;\cdot\;\mone{h}{Z}{\sP Z}
\;:=\;
\begin{tikzpicture}[mb]
\mbars{-1.2}{1.65}{2.8}\medge{-1.2}{0}{2.8}
\draw[mstr] (0,0) -- (0,0.8);
\draw[mstr] (0.15,0.8) -- (0.15,1.8);
\draw[mstr] (-0.15,0.8) .. controls (-0.15,1.35) and (-0.6,1.35) .. (-1.2,1.35);
\draw[mstr] (1.15,0) -- (1.15,1.8);
\sstr{1.15}{1.8}{0.15}{2.8}
\sstr[mover]{0.15}{1.8}{1.15}{2.8}
\node[mcpn] at (0,0.8) {$g$};
\node[mcpn] at (1.15,0.6) {$h$};
\legend{-1.2}{1.35}{E}
\blab{0}{X}\blab{1.15}{Z}\tlab{0.15}{2.8}{\sP Z}\tlab{1.15}{2.8}{\sP X}
\end{tikzpicture}
\]
  for $[Z,h]\in\fX$, $[X,g]\in\cH_E$ make $\cH_E$ an $\fX$-bimodule.
\end{enumerate}
\end{lemma}


\begin{proof}
  It is completely  parallel  to the proof for $\fX$ itself.
\end{proof}

To prove Proposition~\ref{prop:W-basis} below, we need the following linear forms on $\cH_E$:

\begin{definition}\label{def:kappa}
Let $E\in\cC$.
For simple $Y$ and $k\in\Hom(E\otimes\sP Y,Y)$ let
\[
\operatorname{Tr}_{Y,k}([Z,g]):=
\begin{tikzpicture}[mb]
\draw[mstr] (0,0.5) -- (0,1.5);
\draw[mstr] (-0.35,1.5) -- (-0.35,3.7);
\draw[mstr] (0.35,1.5) -- (0.35,3.7);
\draw[mstr] (0,3.7) -- (0,4.4) .. controls (0,4.85) and (-1.05,4.85) .. (-1.05,4.4)
  -- (-1.05,0.2) .. controls (-1.05,-0.25) and (0,-0.25) .. (0,0.2) -- (0,0.5);
\node[mcpnc] at (0,0.5) {$e$};
\node[mcpn,minimum width=10mm] at (0,1.5) {$g$};
\node[mcpnc] at (0.35,2.6) {$e$};
\ubox{0}{2.2}{0.7}{3.0}
\node[mcpn,minimum width=10mm] at (0,3.7) {$k$};
\node[anchor=west,inner sep=1pt] at (0,4.2) {$Y$};
\node[anchor=west,inner sep=1pt] at (0,1.0) {$Z$};
\node[anchor=east,inner sep=1pt] at (-0.35,2.6) {$E$};
\end{tikzpicture}\,,
\]
where $e$ is the canonical element of $\Hom(Y,Z)\otimes\Hom(Z,Y)$. 
\end{definition}

\begin{lemma}\label{lem:kappa}
$\operatorname{Tr}_{Y,k}$ is a well-defined linear form on $\cH_E$, and
$\operatorname{Tr}_{Y,k}([Y',h])=\delta_{Y,Y'}\operatorname{tr}(kh)/d_Y$ for simple $Y'$.
\end{lemma}

\begin{proof}
By \ref{item:E1} and \ref{item:B1}, $\operatorname{Tr}_{Y,k}$ vanishes on
(R2), as in the proof of Lemma~\ref{lem:TrX}, and it vanishes on (R1) since the
canonical element for $Z\oplus Z'$ is the sum of those for $Z$ and $Z'$.
For simple $Y'$ we have $\Hom(Y,Y')=0$ unless $Y'\cong Y$, and the canonical
element of $\End(Y)\otimes\End(Y)$ is $d_Y^{-1}\id_Y\otimes\id_Y$.
\end{proof}

\begin{proposition}\label{prop:W-basis}
Let $E\in\cC$.
The map
\[
\bigoplus_{Y\in\Irr(\cC)}\Hom(Y,E\otimes\sP Y)\to\cH_E,\qquad h\mapsto[Y,h],
\]
is an isomorphism.
\end{proposition}

\begin{proof}
Again this is parallel to the proof for $\fX$ itself,
if we use $\Tr_{Y,k}$ instead of $\Tr_X$.
\end{proof}

\begin{proposition}\label{prop:W-factor}
For $u,v\in\fX$,
\[
\sum_{E\in\Irr(\cC)}d_E\Tr_{\cH_E}(w\mapsto u\cdot w\cdot v)=D^2\varepsilon(u)\,\varepsilon(v).
\]
\end{proposition}


\begin{proof}[Proof\/\protect\footnotemark]
\footnotetext{An animated version of this proof is available at \url{https://ippo0127.github.io/movies/mobius_prop5_6.html}.}%
It suffices to
take $u=\mob(A)=[A,\xi_A]$ and $v=\mob(B)=[B,\xi_B]$ with $A,B\in\Fix(\sP)$.
By Definition~\ref{def:eps},  the right hand side
is $D^2\delta_{A,\one}\delta_{B,\one}$.
What remains then is to compute the left hand side. 

\medskip\noindent\emph{Step 1.}
We compute the trace over $\cH_E$ in the basis $[Y,h_i]$ of
Proposition~\ref{prop:W-basis}, where for each $Y\in\Irr(\cC)$ the $h_i$
form a basis of $\Hom(Y,E\otimes\sP Y)$, as the sum of the diagonal
coefficients with respect to the dual basis $d_Y\operatorname{Tr}_{Y,k_i}$
(Lemma~\ref{lem:kappa}, where $\sum_ih_i\otimes k_i=t$); thus we can write
\[
\sum_{E\in\Irr(\cC)}d_E\Tr_{\cH_E}(w\mapsto u\cdot w\cdot v)
=
\sum_{E,Y\in\Irr(\cC)}d_Ed_Y
\begin{tikzpicture}[mb,y=6mm]
\draw[mstr] (0,0.5) -- (0,1.6); \draw[mstr] (1.1,0.5) -- (1.1,1.6); \draw[mstr] (2.2,0.5) -- (2.2,1.6);
\draw[mstr] (0.95,1.6) .. controls (0.95,2.2) and (-0.8,2.0) .. (-0.8,2.7) -- (-0.8,5.9)
  .. controls (-0.8,6.2) and (0.1,6.0) .. (0.1,6.4);
\draw[mstr] (1.25,1.6) -- (1.25,2.5); \sstr{1.25}{2.5}{0}{3.1} \draw[mstr] (0,3.1) -- (0,3.7);
\draw[mstr] (2.2,1.6) -- (2.2,3.1); \sstr{2.2}{3.1}{1.1}{3.7} \sstr{1.1}{3.7}{0}{4.3}
\sstr[mover]{0}{3.7}{1.1}{4.3}
\draw[mover] (0,1.6) -- (0,2.5); \sstr[mover]{0}{2.5}{1.1}{3.1} \sstr[mover]{1.1}{3.1}{2.2}{3.7}
\draw[mstr] (2.2,3.7) -- (2.2,4.3);
\draw[mstr] (0,4.3) -- (0,4.9); \draw[mstr] (1.1,4.3) -- (1.1,4.9); \draw[mstr] (2.2,4.3) -- (2.2,4.9);
\draw[mstr] (1.1,4.9) -- (1.1,5.6) .. controls (1.1,6.0) and (0.7,5.95) .. (0.7,6.4);
\draw[mstr] (0.4,6.4) -- (0.4,7.25) .. controls (0.4,7.8) and (-1.5,7.8) .. (-1.5,7.15) -- (-1.5,0.4)
  .. controls (-1.5,-0.3) and (1.1,-0.3) .. (1.1,0.2) -- (1.1,0.5);
\node[mcpnc,minimum width=24mm] at (1.1,0.5) {$e$};
\node[mcpn] at (0,1.6) {$\xi_A$}; \node[mcpnc] at (1.1,1.6) {$t$}; \node[mcpn] at (2.2,1.6) {$\xi_B$};
\node[mcpnc,minimum width=24mm] at (1.1,4.9) {$e$};
\ubox{-0.55}{4.45}{2.75}{5.35}
\node[mcpnc,minimum width=10mm] at (0.4,6.4) {$t$};
\node[anchor=west,inner sep=1pt] at (0.4,7.0) {$Y$};
\node[anchor=east,inner sep=1pt] at (-0.8,4.0) {$E$};
\end{tikzpicture}
\]
where $t$ is the canonical element of $\Hom(Y,E\otimes\sP Y)\otimes\Hom(E\otimes\sP Y,Y)$
and $e$ that of $\Hom(Y,A\otimes Y\otimes B)\otimes\Hom(A\otimes Y\otimes B,Y)$.
This can be further simplified as follows, where $t'$ is the canonical element of
$\Hom(Y\otimes(\sP Y)^\vee,E)\otimes\Hom(E,Y\otimes(\sP Y)^\vee)$, obtained from $t$ by
bending $\sP Y$ \ref{item:E4}:
\begin{gather*}
=\sum_{E,Y\in\Irr(\cC)}d_Ed_Y
\begin{tikzpicture}[mb]
\draw[mstr] (0,0.5) -- (0,0.75) .. controls (0,1.15) and (-1.25,1.0) .. (-1.25,1.4) -- (-1.25,1.6);
\draw[mstr] (-0.7,5.35) -- (-0.7,6.0) .. controls (-0.7,6.45) and (-2.0,6.45) .. (-2.0,6.0)
  -- (-2.0,0.2) .. controls (-2.0,-0.25) and (0,-0.25) .. (0,0.2) -- (0,0.5);
\draw[mover] (-0.4,0.5) -- (-0.4,1.95);
\draw[mstr] (0.4,0.5) -- (0.4,1.95);
\draw[mstr] (-1.4,1.6) -- (-1.4,4.5) .. controls (-1.4,4.95) and (-1.0,4.9) .. (-1.0,5.35);
\sstr{0.4}{1.95}{-0.4}{2.5}
\sstr[mover]{-0.4}{1.95}{0.4}{2.5}
\draw[mstr] (-0.4,2.5) -- (-0.4,4.3); \draw[mstr] (0.4,2.5) -- (0.4,4.3);
\draw[mover] (-1.1,1.6) -- (-1.1,2.9) .. controls (-1.1,3.3) and (0,3.1) .. (0,3.5) -- (0,4.3);
\draw[mstr] (0,4.3) -- (0,4.75) .. controls (0,5.05) and (-0.4,4.95) .. (-0.4,5.35);
\node[mcpnc,minimum width=12mm] at (0,0.5) {$e$};
\node[mcpnc] at (-1.25,1.6) {$t$};
\node[mcpn] at (-0.4,1.6) {$\xi_A$};
\node[mcpn] at (0.4,1.6) {$\xi_B$};
\node[mcpnc,minimum width=12mm] at (0,4.3) {$e$};
\ubox{-0.88}{3.85}{0.88}{4.75}
\node[mcpnc,minimum width=9mm] at (-0.7,5.35) {$t$};
\node[anchor=west,inner sep=1pt] at (-0.7,5.83) {$Y$};
\node[anchor=east,inner sep=1pt] at (-1.4,2.8) {$E$};
\end{tikzpicture}
\overset{\text{\ref{item:E4}}}{=}
\sum_{E,Y\in\Irr(\cC)}d_Ed_Y
\begin{tikzpicture}[mb]
\draw[mstr] (0,0.5) -- (0,0.75) .. controls (0,1.1) and (-1.6,0.95) .. (-1.6,1.35) -- (-1.6,2.0);
\draw[mstr] (-1.15,5.35) -- (-1.15,6.3) .. controls (-1.15,6.75) and (-2.2,6.75) .. (-2.2,6.3)
  -- (-2.2,0.2) .. controls (-2.2,-0.25) and (0,-0.25) .. (0,0.2) -- (0,0.5);
\draw[mover] (-0.4,0.5) -- (-0.4,1.95);
\draw[mstr] (0.4,0.5) -- (0.4,1.95);
\draw[mstr] (-1.45,2.0) -- (-1.45,4.5) .. controls (-1.45,4.85) and (-1.0,4.8) .. (-1.0,5.35);
\sstr{0.4}{1.95}{-0.4}{2.5}
\sstr[mover]{-0.4}{1.95}{0.4}{2.5}
\draw[mstr] (-0.4,2.5) -- (-0.4,4.3); \draw[mstr] (0.4,2.5) -- (0.4,4.3);
\draw[mover] (-1.3,2.0) -- (-1.3,1.65) .. controls (-1.3,1.3) and (-1.0,1.3) .. (-1.0,1.65)
  -- (-1.0,2.9) .. controls (-1.0,3.3) and (0,3.1) .. (0,3.5) -- (0,4.3);
\draw[mstr] (0,4.3) -- (0,5.75) .. controls (0,6.2) and (-0.85,6.2) .. (-0.85,5.75) -- (-0.85,5.35);
\node[mcpnc,minimum width=12mm] at (0,0.5) {$e$};
\node[mcpnc,minimum width=6.5mm] at (-1.45,2.0) {$t'$};
\node[mcpn] at (-0.4,1.6) {$\xi_A$};
\node[mcpn] at (0.4,1.6) {$\xi_B$};
\node[mcpnc,minimum width=12mm] at (0,4.3) {$e$};
\ubox{-0.88}{3.85}{0.88}{4.75}
\node[mcpnc,minimum width=6.5mm] at (-1.0,5.35) {$t'$};
\node[anchor=east,inner sep=1pt] at (-1.15,6.1) {$Y$};
\node[anchor=east,inner sep=1pt] at (-1.45,3.2) {$E$};
\end{tikzpicture}
\\
\overset{\text{\ref{item:E3}}}{=}
\sum_{Y\in\Irr(\cC)}d_Y
\begin{tikzpicture}[mb]
\draw[mstr] (0,1.0) -- (0,1.5) .. controls (0,1.95) and (-1.0,1.95) .. (-1.0,1.5)
  -- (-1.0,0.3) .. controls (-1.0,-0.15) and (0,-0.15) .. (0,0.3) -- (0,1.0);
\draw[mover] (-0.4,1.0) -- (-0.4,2.65);
\draw[mstr] (0.4,1.0) -- (0.4,2.65);
\sstr{0.4}{2.65}{-0.4}{3.2}
\sstr[mover]{-0.4}{2.65}{0.4}{3.2}
\draw[mstr] (-0.4,3.2) -- (-0.4,4.6); \draw[mstr] (0.4,3.2) -- (0.4,4.6);
\draw[mover] (0,4.6) -- (0,5.3) .. controls (0,5.75) and (-1.0,5.75) .. (-1.0,5.3)
  -- (-1.0,4.1) .. controls (-1.0,3.65) and (0,3.65) .. (0,4.1) -- (0,4.6);
\node[mcpnc,minimum width=12mm] at (0,1.0) {$e$};
\node[mcpn] at (-0.4,2.3) {$\xi_A$};
\node[mcpn] at (0.4,2.3) {$\xi_B$};
\node[mcpnc,minimum width=12mm] at (0,4.6) {$e$};
\ubox{-0.88}{4.2}{0.88}{5.0}
\node[anchor=west,inner sep=1pt] at (0,0.5) {$Y$};
\node[anchor=west,inner sep=1pt] at (0,5.17) {$\sP Y$};
\end{tikzpicture}
\overset{\text{\ref{item:B1},\ref{item:B2}}}{\underset{\text{\ref{item:B4}--\ref{item:B6}}}{=}}
\sum_{Y\in\Irr(\cC)}d_Y
\begin{tikzpicture}[mb]
\draw[mstr] (0,1.0) -- (0,1.5) .. controls (0,1.95) and (-1.0,1.95) .. (-1.0,1.5)
  -- (-1.0,0.3) .. controls (-1.0,-0.15) and (0,-0.15) .. (0,0.3) -- (0,1.0);
\draw[mover] (-0.4,1.0) -- (-0.4,2.65);
\draw[mstr] (0.4,1.0) -- (0.4,2.65);
\sstr{0.4}{2.65}{-0.4}{3.2}
\sstr[mover]{-0.4}{2.65}{0.4}{3.2}
\draw[mstr] (-0.4,3.2) -- (-0.4,4.6); \draw[mstr] (0.4,3.2) -- (0.4,4.6);
\draw[mover] (0,4.6) -- (0,5.3) .. controls (0,5.75) and (-1.0,5.75) .. (-1.0,5.3)
  -- (-1.0,4.1) .. controls (-1.0,3.65) and (0,3.65) .. (0,4.1) -- (0,4.6);
\node[mcpnc,minimum width=12mm] at (0,1.0) {$e$};
\node[mcpn] at (-0.4,2.3) {$\xi_A$};
\node[mcpn] at (0.4,2.3) {$\xi_B$};
\node[mcpnc,minimum width=12mm] at (0,4.6) {$e$};
\ubox{-1.35}{3.45}{0.75}{5.8}
\node[anchor=west,inner sep=1pt] at (0,0.5) {$Y$};
\node[anchor=west,inner sep=1pt] at (0,5.1) {$Y$};
\end{tikzpicture}\;.
\end{gather*}
\medskip\noindent\emph{Step 2.}
The first picture below is the last picture of Step~1.
Let $\hat e$ be the canonical element of
$\Hom(\one,A\otimes Y\otimes B\otimes Y^\vee)\otimes\Hom(A\otimes Y\otimes B\otimes Y^\vee,\one)$ and
$e'$ that of $\Hom(\one,A\otimes B)\otimes\Hom(A\otimes B,\one)$, and put
\[
\iota:=(\id_A\otimes c^{-1}_{Y,B}\otimes\id_{Y^\vee})(\id_{A\otimes B}\otimes\coev_Y)\colon
A\otimes B\to A\otimes Y\otimes B\otimes Y^\vee,
\]
the cup on the right of $B$ whose strand $Y$ passes in front of $B$ to the middle.
For the upper end $q$ of $e$ and its bent form $\hat q:=\ev'_Y(q\otimes\id_{Y^\vee})$, the morphism
$\hat q\,\iota$ is $q$ with its output $Y$ joined to its middle input by a strand passing in front
of $B$: this is the upper end of $e$ together with its loop inside the box.
So, once both ends of $e$ are bent by \ref{item:E4}, the box contains $\hat q\,\iota$, and
\ref{item:E1} for $h:=\iota$, which also applies inside the box (see the remark after
Proposition~\ref{prop:box-rules}), moves $\iota$ to the lower end, which becomes $e'$
(the part $\iota$ is shaded):
\[
\sum_{Y\in\Irr(\cC)}d_Y
\begin{tikzpicture}[mb]
\draw[mstr] (0,1.0) -- (0,1.5) .. controls (0,1.95) and (-1.0,1.95) .. (-1.0,1.5)
  -- (-1.0,0.3) .. controls (-1.0,-0.15) and (0,-0.15) .. (0,0.3) -- (0,1.0);
\draw[mover] (-0.4,1.0) -- (-0.4,2.65);
\draw[mstr] (0.4,1.0) -- (0.4,2.65);
\sstr{0.4}{2.65}{-0.4}{3.20}
\sstr[mover]{-0.4}{2.65}{0.4}{3.20}
\draw[mstr] (-0.4,3.2) -- (-0.4,4.6); \draw[mstr] (0.4,3.2) -- (0.4,4.6);
\draw[mover] (0,4.6) -- (0,5.3) .. controls (0,5.75) and (-1.0,5.75) .. (-1.0,5.3)
  -- (-1.0,4.1) .. controls (-1.0,3.65) and (0,3.65) .. (0,4.1) -- (0,4.6);
\node[mcpnc,minimum width=12mm] at (0,1.0) {$e$};
\node[mcpn] at (-0.4,2.30) {$\xi_A$}; \node[mcpn] at (0.4,2.30) {$\xi_B$};
\node[mcpnc,minimum width=12mm] at (0,4.6) {$e$};
\ubox{-1.35}{3.45}{0.75}{5.8}
\node[anchor=west,inner sep=1pt] at (0,0.5) {$Y$};
\node[anchor=west,inner sep=1pt] at (0,5.1) {$Y$};
\end{tikzpicture}
\overset{\text{\ref{item:E4}}}{=}
\sum_{Y\in\Irr(\cC)}d_Y
\begin{tikzpicture}[mb]
\fill[rounded corners=3pt,fill=mshade] (-1.2,3.58) rectangle (0.6,4.42);
\begin{scope}
\clip (-1.2,3.58) rectangle (0.6,4.42);
\draw[white,line width=2.2pt] (-0.4,3.2) -- (-0.4,4.75) (0.4,3.2) -- (0.4,4.75)
  (0,4.75) -- (0,4.2) .. controls (0,3.75) and (-1.0,3.75) .. (-1.0,4.2) -- (-1.0,4.75);
\end{scope}
\draw[mstr] (0,1.0) -- (0,1.5) .. controls (0,1.95) and (-1.0,1.95) .. (-1.0,1.5)
  -- (-1.0,0.3) .. controls (-1.0,-0.3) and (1.35,-0.3) .. (1.35,0.3)
  -- (1.35,1.35) .. controls (1.35,1.85) and (0.8,1.85) .. (0.8,1.35) -- (0.8,1.0);
\draw[mover] (-0.4,1.0) -- (-0.4,2.65);
\draw[mstr] (0.4,1.0) -- (0.4,2.65);
\sstr{0.4}{2.65}{-0.4}{3.20}
\sstr[mover]{-0.4}{2.65}{0.4}{3.20}
\draw[mstr] (-0.4,3.2) -- (-0.4,4.75); \draw[mstr] (0.4,3.2) -- (0.4,4.75);
\draw[white,line width=3pt] (0,4.6) -- (0,4.2) .. controls (0,3.75) and (-1.0,3.75) .. (-1.0,4.2) -- (-1.0,4.6);
\draw[mstr] (0,4.75) -- (0,4.2) .. controls (0,3.75) and (-1.0,3.75) .. (-1.0,4.2) -- (-1.0,4.75);
\node[mcpnc,minimum width=15mm] at (0.2,1.0) {$\hat e$};
\node[mcpn] at (-0.4,2.30) {$\xi_A$}; \node[mcpn] at (0.4,2.30) {$\xi_B$};
\node[mcpnc,minimum width=15mm] at (-0.3,4.75) {$\hat e$};
\ubox{-1.35}{3.45}{0.75}{5.2}
\node[anchor=west,inner sep=1pt] at (1.35,0.6) {$Y$};
\node[anchor=west,inner sep=1pt] at (0,4.27) {$Y$};
\node[anchor=west,inner sep=1pt] at (0.8,4.0) {$\iota$};
\end{tikzpicture}
\overset{\text{\ref{item:E1}}}{=}
\sum_{Y\in\Irr(\cC)}d_Y
\begin{tikzpicture}[mb]
\fill[rounded corners=3pt,fill=mshade] (-0.6,0.9) rectangle (1.2,1.74);
\begin{scope}
\clip (-0.6,0.9) rectangle (1.2,1.74);
\draw[white,line width=2.2pt] (-0.4,0.5) -- (-0.4,2.3) (0.4,0.5) -- (0.4,2.3);
\end{scope}
\draw[mstr] (0.4,0.5) -- (0.4,3.65);
\draw[white,line width=3pt] (0,1.9) -- (0,1.52) .. controls (0,1.07) and (1.0,1.07) .. (1.0,1.52) -- (1.0,1.9);
\draw[mstr] (1.0,1.52) .. controls (1.0,1.07) and (0,1.07) .. (0,1.52)
  -- (0,2.15) .. controls (0,2.6) and (-1.0,2.6) .. (-1.0,2.15) -- (-1.0,0.3)
  .. controls (-1.0,-0.35) and (1.8,-0.35) .. (1.8,0.3) -- (1.8,2.15)
  .. controls (1.8,2.6) and (1.0,2.6) .. (1.0,2.15) -- (1.0,1.52);
\draw[mstr] (-0.4,0.5) -- (-0.4,2.05);
\draw[mover] (-0.4,2.05) -- (-0.4,3.65);
\sstr{0.4}{3.65}{-0.4}{4.2}
\sstr[mover]{-0.4}{3.65}{0.4}{4.2}
\draw[mstr] (-0.4,4.2) -- (-0.4,4.95); \draw[mstr] (0.4,4.2) -- (0.4,4.95);
\node[mcpnc,minimum width=9mm] at (0,0.5) {$e'$};
\node[mcpn] at (-0.4,3.3) {$\xi_A$}; \node[mcpn] at (0.4,3.3) {$\xi_B$};
\node[mcpnc,minimum width=10mm] at (0,4.95) {$e'$};
\ubox{-0.75}{4.45}{0.75}{5.4}
\node[anchor=west,inner sep=1pt] at (1.8,0.9) {$Y$};
\node[anchor=west,inner sep=1pt] at (1.23,1.3) {$\iota$};
\end{tikzpicture}
\]
\begin{gather*}
=\sum_{Y\in\Irr(\cC)}d_Y
\begin{tikzpicture}[mb]
\draw[mstr] (-0.4,0.5) -- (-0.4,2.7); \draw[mstr] (0.4,0.5) -- (0.4,2.7);
\draw[mover] (-0.4,1.35) ellipse (0.32 and 0.13);
\draw[mover] (-0.4,1.4) -- (-0.4,1.6);
\sstr{0.4}{2.70}{-0.4}{3.30}
\sstr[mover]{-0.4}{2.70}{0.4}{3.30}
\draw[mstr] (-0.4,3.3) -- (-0.4,3.9); \draw[mstr] (0.4,3.3) -- (0.4,3.9);
\node[mcpnc,minimum width=10mm] at (0,0.5) {$e'$};
\node[mcpn] at (-0.4,2.20) {$\xi_A$}; \node[mcpn] at (0.4,2.20) {$\xi_B$};
\node[mcpnc,minimum width=10mm] at (0,3.9) {$e'$};
\ubox{-0.75}{3.45}{0.75}{4.35}
\node[anchor=west,inner sep=1pt] at (-0.08,1.35) {$Y$};
\end{tikzpicture}
\overset{\text{Def.~\ref{def:Omega}}}{=}
\begin{tikzpicture}[mb]
\draw[mstr] (-0.4,0.5) -- (-0.4,2.7); \draw[mstr] (0.4,0.5) -- (0.4,2.7);
\draw[mover] (-0.4,1.35) ellipse (0.32 and 0.13);
\draw[mover] (-0.4,1.4) -- (-0.4,1.6);
\sstr{0.4}{2.7}{-0.4}{3.3}
\sstr[mover]{-0.4}{2.7}{0.4}{3.3}
\draw[mstr] (-0.4,3.3) -- (-0.4,3.9); \draw[mstr] (0.4,3.3) -- (0.4,3.9);
\node[mcpnc,minimum width=10mm] at (0,0.5) {$e'$};
\node[mcpn] at (-0.4,2.2) {$\xi_A$}; \node[mcpn] at (0.4,2.2) {$\xi_B$};
\node[mcpnc,minimum width=10mm] at (0,3.9) {$e'$};
\ubox{-0.75}{3.45}{0.75}{4.35}
\node[anchor=east,inner sep=1pt] at (-0.72,1.35) {$\Omega$};
\end{tikzpicture}
\overset{\text{Lem.~\ref{lem:kill}}}{=}
D^2\delta_{A,\one}
\begin{tikzpicture}[mb]
\draw[mstr] (-0.4,0.5) -- (-0.4,2.7); \draw[mstr] (0.4,0.5) -- (0.4,2.7);
\sstr{0.4}{2.7}{-0.4}{3.3}
\sstr[mover]{-0.4}{2.7}{0.4}{3.3}
\draw[mstr] (-0.4,3.3) -- (-0.4,3.9); \draw[mstr] (0.4,3.3) -- (0.4,3.9);
\node[mcpnc,minimum width=10mm] at (0,0.5) {$e'$};
\node[mcpn] at (-0.4,2.2) {$\xi_A$}; \node[mcpn] at (0.4,2.2) {$\xi_B$};
\node[mcpnc,minimum width=10mm] at (0,3.9) {$e'$};
\ubox{-0.75}{3.45}{0.75}{4.35}
\end{tikzpicture}
=D^2\delta_{A,\one}\delta_{B,\one}
\end{gather*}
\end{proof}

\begin{corollary}\label{cor:W-weighted}
For $X,X'\in\Sigma$,
$\sum_{E\in\Irr(\cC)}d_E\dim p_X\cH_Ep_{X'}=d_Xd_{X'}M_XM_{X'}$.
\end{corollary}

\begin{proof}
Take $u=p_X$ and $v=p_{X'}$ in Proposition~\ref{prop:W-factor}.
Since the left and right
actions commute, $w\mapsto p_Xwp_{X'}$ is an idempotent of $\cH_E$ with image
$p_X\cH_Ep_{X'}$, so the left side is $\sum_Ed_E\dim p_X\cH_Ep_{X'}$.
By
Lemma~\ref{lem:eps-pX}, $\varepsilon(p_X)=S_{\one,X}M_X$, so the right side is
$D^2S_{\one,X}S_{\one,X'}M_XM_{X'}=d_Xd_{X'}M_XM_{X'}$.
\end{proof}

\subsection{Tube algebra and the decomposition of the bimodules}
\label{subsec:tube}

We now show that $\bigoplus_{E\in\Irr(\cC)}\cH_E$ is a module over the tube algebra $\cT=\operatorname{Tube}(\cC)$ of $\cC$, commuting with the $\fX$-bimodule structure. 
We combine this information with Proposition~\ref{prop:W-factor} 
to deduce the decomposition of the bimodules $\cH_E$ into simple $\fX$-bimodules.
\subsubsection{Generalities on the tube algebra}
Here we collect the basic facts about the tube algebra \cite{Ocneanu1994,Izumi2000} of a modular tensor category $\cC$
that we need.
For completeness, we include the proofs, which are standard.

\begin{definition}\label{def:tube}
Let $\cT$ be the vector space
\[
\cT:=\bigoplus_{E,E',V\in\Irr(\cC)}\Hom(V\otimes E,\,E'\otimes V),
\]
where we introduce a special pictorial notation for the elements of $\cT$:
for $\tau\in\Hom(V\otimes E,E'\otimes V)$ we write
\[
\tubeelt{\tau}{V}{E}{E'}
\qquad\text{for}\qquad
\begin{tikzpicture}[mb]
\draw[mstr] (-0.25,0) -- (-0.25,3);
\draw[mstr] (0.25,0) -- (0.25,3);
\node[mcpn,minimum width=8mm] at (0,1.5) {$\tau$};
\blab{-0.25}{V}\blab{0.25}{E}\tlab{-0.25}{3}{E'}\tlab{0.25}{3}{V}
\end{tikzpicture}\,.
\]
The product on $\cT$ is defined by 
\[
\tubeelt{\tau'}{V'}{E'}{E''}
\;*\;
\tubeelt{\tau}{V}{E}{E'}
\;:=\;
\sum_{V''\in\Irr(\cC)}d_{V''}\;
\begin{tikzpicture}[mb]
\draw[mstr] (0,0) -- (0,0.9);
\draw[mstr] (-0.3,0.9) -- (-0.3,3.4);
\draw[mstr] (0.3,0.9) -- (0.3,4.6);
\draw[mstr,rounded corners=4pt] (2.0,1.3) -- (0.8,1.3) -- (0.8,2.2);
\draw[mstr] (0.8,2.2) -- (0.8,4.6);
\draw[mstr,rounded corners=4pt] (-0.3,3.4) -- (-0.3,4.0) -- (-1.5,4.0);
\draw[mstr] (0.55,4.6) -- (0.55,5.6);
\tubemark{0}{0.22}\tubemark{0.55}{5.38}
\node[mcpnc,minimum width=7mm] at (0,0.9) {$r$};
\node[mcpn,minimum width=7mm] at (0.55,2.2) {$\tau$};
\node[mcpn,minimum width=7mm] at (0,3.4) {$\tau'$};
\node[mcpnc,minimum width=7mm] at (0.55,4.6) {$r$};
\node[anchor=east,inner sep=1pt] at (-0.3,2.2) {$V'$};
\node[anchor=west,inner sep=1pt] at (0.3,1.58) {$V$};
\node[anchor=west,inner sep=1pt] at (0.3,2.8) {$E'$};
\node[anchor=west,inner sep=1pt] at (0.8,3.4) {$V$};
\node[anchor=west,inner sep=1pt] at (0.3,4.0) {$V'$};
\node[anchor=west,inner sep=1.5pt] at (2.0,1.3) {$E$};
\node[anchor=east,inner sep=1.5pt] at (-1.5,4.0) {$E''$};
\blab{0}{V''}\tlab{0.55}{5.6}{V''}
\end{tikzpicture}\,,
\]
where $r$ is the canonical element of $\Hom(V'',V'\otimes V)\otimes\Hom(V'\otimes V,V'')$.
We denote the summand for $V''$ on the right hand side by $(\tau'*\tau)_{V''}$.
If the middle objects differ, i.e.\ if the top object of $\tau$ is not the bottom object of $\tau'$, the product is declared to be zero, and the product is extended to all of $\cT$ by linearity.
The unit of $\cT$ is $\sum_E\id_E$, with $V=\one$. 
\end{definition}

\begin{remark}
Motivationally, in the special notation for the tube algebra elements,
the two ends of the strand $V$ marked by $/\!/$ are to be thought 
of as glued together around a tube.
Namely, the picture is drawn on the surface of a horizontal tube, along which $E$ runs from the right to the left, and the strand $V$ wraps once around it, passing in front of the tube through $\tau$ and coming back behind it (left picture):
\[
\begin{tikzpicture}[mb]
\draw[mstr,densely dashed,black!50] (1.8,2.6) arc[start angle=90,end angle=270,x radius=0.3,y radius=1.3];
\draw[mstr] (1.8,0) arc[start angle=-90,end angle=90,x radius=0.3,y radius=1.3];
\draw[mstr] (-1.8,1.3) ellipse[x radius=0.3,y radius=1.3];
\draw[mstr] (-1.8,0) -- (1.8,0);
\draw[mstr] (-1.8,2.6) -- (1.8,2.6);
\draw[mstr,densely dashed,black!60] (0.4,2.6) .. controls (-0.8,2.4) and (-0.9,0.3) .. (-0.3,0);
\draw[mover,rounded corners=4pt] (-0.25,1.3) -- (-0.25,1.85) -- (-1.53,1.85);
\draw[mstr,rounded corners=4pt] (2.07,0.75) -- (0.25,0.75) -- (0.25,1.3);
\draw[mstr] (0.25,1.3) -- (0.25,1.55) .. controls (0.25,1.9) and (0.4,2.2) .. (0.4,2.6);
\draw[mstr] (-0.3,0) .. controls (-0.3,0.4) and (-0.25,0.7) .. (-0.25,1.05) -- (-0.25,1.3);
\node[mcpn,minimum width=8mm] at (0,1.3) {$\tau$};
\node[anchor=west,inner sep=1.8pt] at (2.07,0.75) {$E$};
\node[anchor=east,inner sep=1.8pt] at (-1.53,1.85) {$E'$};
\node[anchor=west,inner sep=1.5pt] at (0.32,2.2) {$V$};
\end{tikzpicture}
\qquad\qquad
\begin{tikzpicture}[mb]
\draw[mstr,densely dashed,black!50] (2.7,2.6) arc[start angle=90,end angle=270,x radius=0.3,y radius=1.3];
\draw[mstr] (2.7,0) arc[start angle=-90,end angle=90,x radius=0.3,y radius=1.3];
\draw[mstr] (-2.7,1.3) ellipse[x radius=0.3,y radius=1.3];
\draw[mstr] (-2.7,0) -- (2.7,0);
\draw[mstr] (-2.7,2.6) -- (2.7,2.6);
\draw[mstr,densely dashed,black!60] (1.3,2.6) .. controls (0.3,2.4) and (0.1,0.35) .. (0.6,0);
\draw[mstr,densely dashed,black!60] (-0.55,2.6) .. controls (-1.55,2.4) and (-1.7,0.35) .. (-1.2,0);
\draw[mstr] (0.6,0) .. controls (0.6,0.3) and (0.65,0.5) .. (0.65,0.8);
\draw[mstr] (1.15,0.8) -- (1.15,1.05) .. controls (1.15,1.7) and (1.3,2.1) .. (1.3,2.6);
\draw[mstr] (-1.2,0) .. controls (-1.2,0.6) and (-1.15,1.1) .. (-1.15,1.7);
\draw[mstr] (-0.65,1.7) -- (-0.65,1.9) .. controls (-0.65,2.25) and (-0.55,2.45) .. (-0.55,2.6);
\draw[mstr,rounded corners=4pt] (2.905,0.35) -- (1.15,0.35) -- (1.15,0.8);
\draw[mover,rounded corners=4pt] (0.65,0.8) -- (0.65,1.25) -- (-0.65,1.25) -- (-0.65,1.7);
\draw[mover,rounded corners=4pt] (-1.15,1.7) -- (-1.15,2.1) -- (-2.46,2.1);
\node[mcpn,minimum width=8mm] at (0.9,0.8) {$\tau$};
\node[mcpn,minimum width=8mm] at (-0.9,1.7) {$\tau'$};
\node[anchor=west,inner sep=1.8pt] at (2.905,0.35) {$E$};
\node[anchor=south,inner sep=1pt] at (0,1.25) {$E'$};
\node[anchor=south,inner sep=1pt] at (-2.0,2.1) {$E''$};
\node[anchor=west,inner sep=1.5pt] at (1.22,1.95) {$V$};
\node[anchor=west,inner sep=1pt] at (-0.52,2.2) {$V'$};
\end{tikzpicture}
\]
Cutting the left tube open along a horizontal line on its back gives the picture on the left in Definition~\ref{def:tube}.
The product $\tau'*\tau$ is the right picture, in which the tube of $\tau'$ is glued to the left of the tube of $\tau$, so that $V$ and $V'$ both wrap around the tube.
Fusing $V'$ and $V$ into $V''$ by \ref{item:E3} on the back of the tube and cutting it open there gives the picture of $\tau'*\tau$ in Definition~\ref{def:tube}.
\end{remark}

\begin{definition}\label{def:pants}
For $A,B\in\Irr(\cC)$, let $W_{A,B}:=\bigoplus_{E\in\Irr(\cC)}\Hom(A,E\otimes B)$.
In pictures, we draw an element $\psi\in\Hom(A,E\otimes B)$ with the strand $E$ going out to the
left, and an element $\varphi\in\Hom(E\otimes B,A)$ with the strand $E$ coming in from the left:
\[
\begin{tikzpicture}[mb]
\draw[mstr] (0,-0.9) -- (0,0.9);
\draw[mstr] (0,0) -- (-0.9,0);
\node[mcpn] at (0,0) {$\psi$};
\node[anchor=east,inner sep=1.5pt] at (-0.9,0) {$E$};
\node[anchor=base,inner sep=0pt,yshift=-8.8pt] at (0,-0.9) {$A$};\tlab{0}{0.9}{B}
\end{tikzpicture}
\;:=\;
\begin{tikzpicture}[mb]
\draw[mstr] (0,-0.9) -- (0,0);
\draw[mstr] (-0.25,0) -- (-0.25,0.9); \draw[mstr] (0.25,0) -- (0.25,0.9);
\node[mcpn,minimum width=8mm] at (0,0) {$\psi$};
\node[anchor=base,inner sep=0pt,yshift=-8.8pt] at (0,-0.9) {$A$};\tlab{-0.25}{0.9}{E}\tlab{0.25}{0.9}{B}
\end{tikzpicture}\,,
\qquad
\begin{tikzpicture}[mb]
\draw[mstr] (0,-0.9) -- (0,0.9);
\draw[mstr] (0,0) -- (-0.9,0);
\node[mcpn] at (0,0) {$\varphi$};
\node[anchor=east,inner sep=1.5pt] at (-0.9,0) {$E$};
\node[anchor=base,inner sep=0pt,yshift=-8.8pt] at (0,-0.9) {$B$};\tlab{0}{0.9}{A}
\end{tikzpicture}
\;:=\;
\begin{tikzpicture}[mb]
\draw[mstr] (-0.25,-0.9) -- (-0.25,0); \draw[mstr] (0.25,-0.9) -- (0.25,0);
\draw[mstr] (0,0) -- (0,0.9);
\node[mcpn,minimum width=8mm] at (0,0) {$\varphi$};
\node[anchor=base,inner sep=0pt,yshift=-8.8pt] at (-0.25,-0.9) {$E$};\node[anchor=base,inner sep=0pt,yshift=-8.8pt] at (0.25,-0.9) {$B$};\tlab{0}{0.9}{A}
\end{tikzpicture}\,.
\]
When such pictures are composed, the end of $E$ is bent upwards (or downwards) on the far left,
crossing no other strand.
For
$\tau\in\Hom(V\otimes E,E'\otimes V)$ let $F_\tau$ act on $W_{A,B}$ by
\[
F_\tau\,\psi:=
\begin{tikzpicture}[mb]
\draw[mstr,rounded corners=6pt] (1.4,0.8) -- (1.4,-0.7) -- (-0.85,-0.7) -- (-0.85,0.75);
\draw[mover] (0.6,-1.3) -- (0.6,0.0);
\draw[mstr] (0.6,0.0) -- (0.6,2.2);
\draw[mstr,rounded corners=4pt] (0.6,0.0) -- (-0.35,0.0) -- (-0.35,0.75);
\draw[mover,rounded corners=6pt] (-0.35,0.75) -- (-0.35,1.6) -- (1.4,1.6) -- (1.4,0.8);
\draw[mstr,rounded corners=4pt] (-0.85,0.75) -- (-0.85,1.15) -- (-1.5,1.15);
\node[mcpn] at (0.6,0.0) {$\psi$};
\node[mcpn,minimum width=8mm] at (-0.6,0.75) {$\tau$};
\node[anchor=east,inner sep=1.5pt] at (-1.5,1.15) {$E'$};
\node[anchor=south,inner sep=1pt] at (0.05,0.0) {$E$};
\node[anchor=east,inner sep=1pt] at (-0.85,-0.1) {$V$};
\node[anchor=base,inner sep=0pt,yshift=-8.8pt] at (0.6,-1.3) {$A$};\tlab{0.6}{2.2}{B}
\end{tikzpicture}
\]
on the summand $E$, and by $0$ on the other summands.
\end{definition}

\begin{proposition}\label{prop:pants-rep}
For $A,B\in\Irr(\cC)$, $\tau\mapsto F_\tau$ is a representation of $\cT$ on $W_{A,B}$, that is,
$F_{\tau'*\tau}=F_{\tau'}\,F_\tau$ and $F_{\sum_E\id_E}=\id$.
\end{proposition}

\begin{proof}
By the definitions of $F$ and of the product,
\begin{gather*}
F_{\tau'*\tau}\,\psi=
\sum_{V''}d_{V''}
\begin{tikzpicture}[mb,y=6.1mm]
\draw[mstr,rounded corners=6pt] (2.4,1.6) -- (2.4,-0.9) -- (0,-0.9) -- (0,0.0);
\draw[mover] (1.6,-1.5) -- (1.6,0.6);
\draw[mover] (1.6,0.6) -- (1.6,4.8);
\draw[mstr,rounded corners=4pt] (1.6,0.6) -- (0.8,0.6) -- (0.8,1.2);
\draw[mstr] (-0.3,0.0) -- (-0.3,2.4);
\draw[mstr] (0.3,0.0) -- (0.3,1.2);
\draw[mstr] (0.3,1.2) -- (0.3,2.4);
\draw[mstr] (0.8,1.2) -- (0.8,3.5);
\draw[mstr] (0.3,2.4) -- (0.3,3.5);
\draw[mstr,rounded corners=4pt] (-0.3,2.4) -- (-0.3,2.9) -- (-0.9,2.9);
\draw[mover,rounded corners=6pt] (0.55,3.5) -- (0.55,4.2) -- (2.4,4.2) -- (2.4,1.6);
\node[mcpn] at (1.6,0.6) {$\psi$};
\node[mcpnc,minimum width=7mm] at (0,0.0) {$r$};
\node[mcpn,minimum width=7mm] at (0.55,1.2) {$\tau$};
\node[mcpn,minimum width=7mm] at (0,2.4) {$\tau'$};
\node[mcpnc,minimum width=7mm] at (0.55,3.5) {$r$};
\node[anchor=east,inner sep=1.5pt] at (-0.9,2.9) {$E''$};
\node[anchor=east,inner sep=1pt] at (-0.1,-0.75) {$V''$};
\node[anchor=east,inner sep=1pt] at (-0.3,1.2) {$V'$};
\node[anchor=west,inner sep=1pt] at (0.3,0.55) {$V$};
\node[anchor=west,inner sep=1pt] at (0.3,1.8) {$E'$};
\node[anchor=base,inner sep=0pt,yshift=-8.8pt] at (1.6,-1.5) {$A$};\tlab{1.6}{4.8}{B}
\end{tikzpicture}
\overset{\text{nat.}}{=}
\sum_{V''}d_{V''}
\begin{tikzpicture}[mb,y=6.1mm]
\draw[mstr,rounded corners=6pt] (2.7,1.6) -- (2.7,-0.35) -- (-0.3,-0.35) -- (-0.3,2.4);
\draw[mstr,rounded corners=6pt] (2.4,1.6) -- (2.4,-0.05) -- (0.3,-0.05) -- (0.3,1.2);
\draw[mover] (1.6,-0.9) -- (1.6,0.6);
\draw[mover] (1.6,0.6) -- (1.6,5.6);
\draw[mstr,rounded corners=4pt] (1.6,0.6) -- (0.8,0.6) -- (0.8,1.2);
\draw[mstr] (0.3,1.2) -- (0.3,2.4);
\draw[mstr] (0.8,1.2) -- (0.8,3.5);
\draw[mstr] (0.3,2.4) -- (0.3,3.5);
\draw[mstr] (0.55,3.5) -- (0.55,4.3);
\draw[mstr,rounded corners=4pt] (-0.3,2.4) -- (-0.3,2.9) -- (-0.9,2.9);
\draw[mover,rounded corners=6pt] (0.3,4.3) -- (0.3,5.1) -- (2.7,5.1) -- (2.7,1.6);
\draw[mover,rounded corners=6pt] (0.8,4.3) -- (0.8,4.8) -- (2.4,4.8) -- (2.4,1.6);
\node[mcpn] at (1.6,0.6) {$\psi$};
\node[mcpn,minimum width=7mm] at (0.55,1.2) {$\tau$};
\node[mcpn,minimum width=7mm] at (0,2.4) {$\tau'$};
\node[mcpnc,minimum width=7mm] at (0.55,3.5) {$r$};
\node[mcpnc,minimum width=7mm] at (0.55,4.3) {$r$};
\node[anchor=east,inner sep=1.5pt] at (-0.9,2.9) {$E''$};
\node[anchor=west,inner sep=1pt] at (0.3,0.35) {$V$};
\node[anchor=east,inner sep=1pt] at (-0.3,1.0) {$V'$};
\node[anchor=west,inner sep=1pt] at (0.3,1.8) {$E'$};
\node[anchor=base,inner sep=0pt,yshift=-8.8pt] at (1.6,-0.9) {$A$};\tlab{1.6}{5.6}{B}
\end{tikzpicture}
\\
\overset{\text{\ref{item:E3}}}{=}
\begin{tikzpicture}[mb,y=6.1mm]
\draw[mstr,rounded corners=6pt] (2.7,1.6) -- (2.7,-0.35) -- (-0.3,-0.35) -- (-0.3,2.4);
\draw[mstr,rounded corners=6pt] (2.4,1.6) -- (2.4,-0.05) -- (0.3,-0.05) -- (0.3,1.2);
\draw[mover] (1.6,-0.9) -- (1.6,0.6);
\draw[mover] (1.6,0.6) -- (1.6,4.6);
\draw[mstr,rounded corners=4pt] (1.6,0.6) -- (0.8,0.6) -- (0.8,1.2);
\draw[mstr] (0.3,1.2) -- (0.3,2.4);
\draw[mstr,rounded corners=4pt] (-0.3,2.4) -- (-0.3,2.9) -- (-0.9,2.9);
\draw[mover,rounded corners=6pt] (0.3,2.4) -- (0.3,4.0) -- (2.7,4.0) -- (2.7,1.6);
\draw[mover,rounded corners=6pt] (0.8,1.2) -- (0.8,3.7) -- (2.4,3.7) -- (2.4,1.6);
\node[mcpn] at (1.6,0.6) {$\psi$};
\node[mcpn,minimum width=7mm] at (0.55,1.2) {$\tau$};
\node[mcpn,minimum width=7mm] at (0,2.4) {$\tau'$};
\node[anchor=east,inner sep=1.5pt] at (-0.9,2.9) {$E''$};
\node[anchor=west,inner sep=1pt] at (0.3,0.35) {$V$};
\node[anchor=east,inner sep=1pt] at (-0.3,1.0) {$V'$};
\node[anchor=west,inner sep=1pt] at (0.3,1.8) {$E'$};
\node[anchor=base,inner sep=0pt,yshift=-8.8pt] at (1.6,-0.9) {$A$};\tlab{1.6}{4.6}{B}
\end{tikzpicture}
=F_{\tau'}F_\tau\,\psi .
\end{gather*}
For $\tau=\sum_E\id_E$ ($V=\one$) there is no loop, and $F_\tau\psi=\psi$.
\end{proof}

\begin{lemma}\label{lem:pants}
\begin{enumerate}[label=(\alph*)]
\item For $E,V\in\Irr(\cC)$, $F_{c_{V,E}}$ acts on $\Hom(A,E\otimes B)$ by $\chi_A([V])$, and
  $F_{c^{-1}_{E,V}}$ by $\chi_B([V])$.
\item For $\tau\in\Hom(V\otimes E,E'\otimes V)$,
\[
\sum_{A,B\in\Irr(\cC)}d_Ad_B\,\Tr_{W_{A,B}}F_\tau=D^2\,\delta_{V,\one}\,\delta_{E,E'}\operatorname{tr}(\tau).
\]
\end{enumerate}
\end{lemma}

\begin{proof}

(a) The loops are natural and slide off:
\[
F_{c_{V,E}}\,\psi=
\begin{tikzpicture}[mb,y=5.5mm]
\draw[mstr,rounded corners=5pt] (1.2,0.3) -- (1.2,-0.6) -- (0.0,-0.6) -- (0.0,-0.25);
\draw[mover] (0.6,-1.1) -- (0.6,0.0);
\draw[mstr] (0.6,0.0) -- (0.6,1.3);
\draw[mstr] (0.6,0.0) -- (-0.9,0.0);
\draw[mover,rounded corners=5pt] (0.0,-0.25) -- (0.0,0.65) -- (1.2,0.65) -- (1.2,0.3);
\node[mcpn] at (0.6,0.0) {$\psi$};
\node[anchor=east,inner sep=1.5pt] at (-0.9,0.0) {$E$};
\node[anchor=east,inner sep=1pt] at (0.0,0.4) {$V$};
\node[anchor=base,inner sep=0pt,yshift=-8.8pt] at (0.6,-1.1) {$A$};\tlab{0.6}{1.3}{B}
\end{tikzpicture}
=
\begin{tikzpicture}[mb,y=5.5mm]
\draw[mstr] (0.4,-0.75) arc[start angle=0,end angle=180,x radius=0.4,y radius=0.13];
\draw[mover] (0,-1.3) -- (0,0.0);
\draw[mstr] (0,0.0) -- (0,1.0);
\draw[mstr] (0,0.0) -- (-0.9,0.0);
\draw[mover] (-0.4,-0.75) arc[start angle=180,end angle=360,x radius=0.4,y radius=0.13];
\draw[mstr,-{Stealth[length=1.5mm,width=1.5mm]}] (0,-0.75) ++(232:0.4 and 0.13)
  arc[start angle=232, end angle=256, x radius=0.4, y radius=0.13];
\node[mcpn] at (0,0.0) {$\psi$};
\node[anchor=east,inner sep=1.5pt] at (-0.9,0.0) {$E$};
\node[anchor=west,inner sep=1pt] at (0.4,-0.75) {$V$};
\node[anchor=base,inner sep=0pt,yshift=-8.8pt] at (0,-1.3) {$A$};\tlab{0}{1.0}{B}
\end{tikzpicture}
=\chi_A([V])\,\psi,
\]
\[
F_{c^{-1}_{E,V}}\,\psi=
\begin{tikzpicture}[mb,y=5.5mm]
\draw[mstr,rounded corners=5pt] (1.2,0.3) -- (1.2,-0.6) -- (0.0,-0.6) -- (0.0,0.25);
\draw[mover] (0.6,-1.1) -- (0.6,0.0);
\draw[mstr] (0.6,0.0) -- (0.6,1.3);
\draw[mover] (0.6,0.0) -- (-0.9,0.0);
\draw[mover,rounded corners=5pt] (0.0,0.25) -- (0.0,0.65) -- (1.2,0.65) -- (1.2,0.3);
\node[mcpn] at (0.6,0.0) {$\psi$};
\node[anchor=east,inner sep=1.5pt] at (-0.9,0.0) {$E$};
\node[anchor=east,inner sep=1pt] at (0.0,0.4) {$V$};
\node[anchor=base,inner sep=0pt,yshift=-8.8pt] at (0.6,-1.1) {$A$};\tlab{0.6}{1.3}{B}
\end{tikzpicture}
=
\begin{tikzpicture}[mb,y=5.5mm]
\draw[mstr] (0.4,0.75) arc[start angle=0,end angle=180,x radius=0.4,y radius=0.13];
\draw[mstr] (0,-1.0) -- (0,0.0);
\draw[mover] (0,0.0) -- (0,1.3);
\draw[mstr] (0,0.0) -- (-0.9,0.0);
\draw[mover] (-0.4,0.75) arc[start angle=180,end angle=360,x radius=0.4,y radius=0.13];
\draw[mstr,-{Stealth[length=1.5mm,width=1.5mm]}] (0,0.75) ++(232:0.4 and 0.13)
  arc[start angle=232, end angle=256, x radius=0.4, y radius=0.13];
\node[mcpn] at (0,0.0) {$\psi$};
\node[anchor=east,inner sep=1.5pt] at (-0.9,0.0) {$E$};
\node[anchor=west,inner sep=1pt] at (0.4,0.75) {$V$};
\node[anchor=base,inner sep=0pt,yshift=-8.8pt] at (0,-1.0) {$A$};\tlab{0}{1.3}{B}
\end{tikzpicture}
=\chi_B([V])\,\psi .
\]

\medskip\noindent (b) Since $F_\tau$ maps $\Hom(A,E\otimes B)$ to $\Hom(A,E'\otimes B)$, only $E'=E$ contributes, and
\begin{align*}
\sum_{A,B}d_Ad_B\,\Tr_{W_{A,B}}F_\tau
&\overset{\text{\ref{item:E2}}}{=}\delta_{E,E'}\sum_{A,B}d_Ad_B
\begin{tikzpicture}[mb,y=6.1mm]
\draw[mstr,rounded corners=6pt] (1.4,0.8) -- (1.4,-0.7) -- (-0.85,-0.7) -- (-0.85,0.75);
\draw[mover] (0.6,-1.2) -- (0.6,0.0);
\draw[mstr] (0.6,0.0) -- (0.6,2.6);
\draw[mstr,rounded corners=4pt] (0.6,0.0) -- (-0.35,0.0) -- (-0.35,0.75);
\draw[mover,rounded corners=6pt] (-0.35,0.75) -- (-0.35,1.6) -- (1.4,1.6) -- (1.4,0.8);
\draw[mstr,rounded corners=4pt] (-0.85,0.75) -- (-0.85,1.15) -- (-1.5,1.15) -- (-1.5,2.6) -- (0.6,2.6);
\draw[mstr] (0.6,2.6) -- (0.6,3.0) .. controls (0.6,3.5) and (2.1,3.5) .. (2.1,3.0) -- (2.1,-1.2)
  .. controls (2.1,-1.7) and (0.6,-1.7) .. (0.6,-1.2);
\node[mcpnc] at (0.6,0.0) {$e$};
\node[mcpn,minimum width=8mm] at (-0.6,0.75) {$\tau$};
\node[mcpnc] at (0.6,2.6) {$e$};
\node[anchor=east,inner sep=1pt] at (-0.85,-0.1) {$V$};
\node[anchor=south,inner sep=1pt] at (0.05,0.0) {$E$};
\node[anchor=east,inner sep=1pt] at (-1.5,1.9) {$E$};
\node[anchor=west,inner sep=1pt] at (0.6,1.85) {$B$};
\node[anchor=west,inner sep=1pt] at (0.6,-1.0) {$A$};
\end{tikzpicture}
\\
&\overset{\text{\ref{item:E3}}}{=}\delta_{E,E'}\sum_Bd_B
\begin{tikzpicture}[mb]
\draw[mstr,rounded corners=6pt] (1.4,0.8) -- (1.4,-0.7) -- (-0.85,-0.7) -- (-0.85,0.75);
\draw[mover,rounded corners=6pt] (0.6,0.4) -- (0.6,2.2) -- (1.9,2.2) -- (1.9,-1.1) -- (0.6,-1.1) -- (0.6,0.4);
\draw[mover,rounded corners=6pt] (-0.35,0.75) -- (-0.35,1.6) -- (1.4,1.6) -- (1.4,0.8);
\draw[mover,rounded corners=6pt] (-0.85,0.75) -- (-0.85,1.15) -- (-1.5,1.15) -- (-1.5,2.7) -- (2.4,2.7)
  -- (2.4,0.0) -- (-0.35,0.0) -- (-0.35,0.75);
\node[mcpn,minimum width=8mm] at (-0.6,0.75) {$\tau$};
\node[anchor=east,inner sep=1pt] at (-0.85,-0.1) {$V$};
\node[anchor=west,inner sep=1pt] at (0.6,0.85) {$B$};
\node[anchor=east,inner sep=1pt] at (-1.5,1.9) {$E$};
\end{tikzpicture}
\overset{\text{Def.~\ref{def:Omega}}}{=}\delta_{E,E'}
\begin{tikzpicture}[mb]
\draw[mstr] (1.75,0.8) arc[start angle=0,end angle=180,x radius=0.35,y radius=0.12];
\draw[mover,rounded corners=6pt] (1.4,0.8) -- (1.4,-0.7) -- (-0.85,-0.7) -- (-0.85,0.75);
\draw[mover,rounded corners=6pt] (-0.35,0.75) -- (-0.35,1.6) -- (1.4,1.6) -- (1.4,0.8);
\draw[mover] (1.05,0.8) arc[start angle=180,end angle=360,x radius=0.35,y radius=0.12];
\draw[mover,rounded corners=6pt] (-0.85,0.75) -- (-0.85,1.15) -- (-1.5,1.15) -- (-1.5,2.2) -- (2.4,2.2)
  -- (2.4,0.0) -- (-0.35,0.0) -- (-0.35,0.75);
\node[mcpn,minimum width=8mm] at (-0.6,0.75) {$\tau$};
\node[anchor=west,inner sep=1pt] at (1.75,0.8) {$\Omega$};
\node[anchor=east,inner sep=1pt] at (1.4,-0.35) {$V^\vee$};
\node[anchor=east,inner sep=1pt] at (-1.5,1.65) {$E$};
\end{tikzpicture}
\\
&\overset{\text{Lem.~\ref{lem:kill}}}{=}\delta_{E,E'}\,D^2\;
\begin{tikzpicture}[mb]
\draw[mstr,rounded corners=6pt] (1.4,0.55) -- (1.4,-0.7) -- (-0.85,-0.7) -- (-0.85,0.75);
\draw[mstr,rounded corners=6pt] (-0.35,0.75) -- (-0.35,1.6) -- (1.4,1.6) -- (1.4,1.2);
\draw[mover,rounded corners=6pt] (-0.85,0.75) -- (-0.85,1.15) -- (-1.5,1.15) -- (-1.5,2.2) -- (2.4,2.2)
  -- (2.4,0.0) -- (-0.35,0.0) -- (-0.35,0.75);
\node[mcpn,minimum width=8mm] at (-0.6,0.75) {$\tau$};
\node[mcpnc] at (1.4,1.2) {$e$};
\node[mcpnc] at (1.4,0.55) {$e$};
\node[anchor=east,inner sep=1pt] at (1.4,-0.35) {$V^\vee$};
\node[anchor=east,inner sep=1pt] at (-1.5,1.65) {$E$};
\end{tikzpicture}
=D^2\,\delta_{V,\one}\,\delta_{E,E'}\operatorname{tr}(\tau).
\end{align*}
In the second picture $e$ is the canonical element of $\Hom(A,E\otimes B)\otimes\Hom(E\otimes B,A)$,
and in the last picture $e$ is the canonical element of $\Hom(V^\vee,\one)\otimes\Hom(\one,V^\vee)$, which vanishes unless $V=\one$. 
For $V=\one$ the picture is $\operatorname{tr}(\tau)$.
\end{proof}

\begin{proposition}\label{prop:tube}
The map $\tau\mapsto F_\tau$, $\cT\to\bigoplus_{A,B\in\Irr(\cC)}\End_\bC(W_{A,B})$, is an isomorphism of algebras.
In
particular $\cT$ is semisimple, and the $W_{A,B}$ are its simple modules.
\end{proposition}

\begin{proof}
For $\tau\in\cT$ let $\operatorname{tr}_0(\tau)$ be the sum of the traces of its summands with
$E'=E$ and $V=\one$.
By Proposition~\ref{prop:pants-rep} and Lemma~\ref{lem:pants}(b),
$D^2\operatorname{tr}_0(\tau'*\tau)=\sum_{A,B}d_Ad_B\Tr_{W_{A,B}}\bigl(F_{\tau'}F_\tau\bigr)$.
For $\tau\in\Hom(V\otimes E,E'\otimes V)$ and
$\tau'\in\Hom(V^\vee\otimes E',E\otimes V^\vee)$, only $V''=\one$ contributes to $\operatorname{tr}_0(\tau'*\tau)$, where
the canonical element $r$ is $d_V^{-1}\coev'_V\otimes\ev_V$, so
\[
\operatorname{tr}_0(\tau'*\tau)=\frac1{d_V}\;
\begin{tikzpicture}[mb]
\draw[mstr] (0,2.9) -- (0,0.9) .. controls (0,0.4) and (1.0,0.4) .. (1.0,0.9) -- (1.0,1.6);
\draw[mstr,rounded corners=4pt] (1.0,1.6) -- (1.0,2.25) -- (0.5,2.25) -- (0.5,2.9);
\draw[mstr] (1.5,1.6) -- (1.5,3.6) .. controls (1.5,4.1) and (0.5,4.1) .. (0.5,3.6) -- (0.5,2.9);
\draw[mstr,rounded corners=6pt] (0,2.9) -- (0,3.4) -- (-0.6,3.4) -- (-0.6,4.5) -- (2.2,4.5)
  -- (2.2,1.0) -- (1.5,1.0) -- (1.5,1.6);
\node[mcpn,minimum width=9mm] at (1.25,1.6) {$\tau$};
\node[mcpn,minimum width=9mm] at (0.25,2.9) {$\tau'$};
\node[anchor=east,inner sep=1pt] at (0,1.5) {$V^\vee$};
\node[anchor=east,inner sep=1pt] at (1.0,1.05) {$V$};
\node[anchor=south,inner sep=1pt] at (0.75,2.25) {$E'$};
\node[anchor=west,inner sep=1pt] at (1.5,2.8) {$V$};
\node[anchor=east,inner sep=1pt] at (-0.6,3.95) {$E$};
\end{tikzpicture}\,,
\]
which after bending $V$ is the trace pairing between $\Hom(V\otimes E\otimes V^\vee,E')$ and
$\Hom(E',V\otimes E\otimes V^\vee)$; it is non-degenerate.
So $F_\tau=0$ implies $\tau=0$, and $\tau\mapsto F_\tau$ is injective.
It is multiplicative by Proposition~\ref{prop:pants-rep}.
Finally, with
$N_{EB}^A=\dim\Hom(A,E\otimes B)$ and $N_{VE}^G=N_{EV}^G$,
\[
\dim\cT=\sum_{E,E',V,G}N_{VE}^G\,N_{E'V}^G=\sum_{A,B}\Bigl(\sum_EN_{EB}^A\Bigr)^2=\dim\bigoplus_{A,B}\End_\bC(W_{A,B}),
\]
so $\tau\mapsto F_\tau$ is an isomorphism.
\end{proof}

The simple modules of $\cT$ are thus labeled by pairs $(A,B)$ of simple objects of
$\cC$, in accordance with the description of the representations of the tube algebra
by the Drinfeld center \cite{Izumi2000} and with $\mathcal{Z}(\cC)\simeq\cC\boxtimes\overline{\cC}$
for modular $\cC$ \cite[Theorem~7.10]{Muger2003}, where $\overline{\cC}$ is $\cC$ with the inverse braiding.

\subsubsection{The tube algebra action on the bimodules}

\begin{definition}
\label{def:F-tube}
For $\tau\in\Hom(V\otimes E,E'\otimes V)$ and $w=[X,g]\in\cH_E$ let
\[
F_\tau(w):=
\begin{tikzpicture}[mb]
\mbars{-1.4}{2.0}{3.2}\medge{-1.4}{0}{3.2}
\draw[mstr] (-0.2,0) -- (-0.2,1.5);
\draw[mstr] (0.9,0) -- (0.9,0.55);
\draw[mstr,rounded corners=4pt] (0.75,0.55) -- (0.75,1.0) -- (0.2,1.0) -- (0.2,1.5);
\draw[mstr] (1.05,0.55) -- (1.05,3.2);
\draw[mstr,rounded corners=4pt] (-0.2,1.5) -- (-0.2,2.0) -- (-1.4,2.0);
\draw[mstr] (0.2,1.5) -- (0.2,3.2);
\draw[mstr] (1.7,0) -- (1.7,3.2);
\node[mcpn,minimum width=8mm] at (0,1.5) {$\tau$};
\node[mcpn] at (0.2,2.6) {$\beta_V$};
\node[mcpn] at (0.9,0.55) {$g$};
\legend{-1.4}{2.0}{E'}
\blab{-0.2}{V}\blab{0.9}{X}\blab{1.7}{\sP V}
\tlab{0.2}{3.2}{\sP^2V}\tlab{1.05}{3.2}{\sP X}\tlab{1.7}{3.2}{\sP V}
\end{tikzpicture}
\;\in\cH_{E'}.
\]
We extend $F_\tau$ by $0$ to the other summands of $\cH:=\bigoplus_{E\in\Irr(\cC)}\cH_E$, and linearly to $\cT$.
Here we have used the same symbol as for the action on $W_{A,B}$ (Definition~\ref{def:pants}); the module on which $F_\tau$ acts is always clear from its argument.
\end{definition}
\begin{lemma}\label{lem:tube-action}
\begin{enumerate}[label=(\alph*)]
\item $F_\tau$ is well defined, and $F_{\id_E}$ is the identity of $\cH_E$.
\item $F_{\tau'}\circ F_\tau=F_{\tau'*\tau}$.
\end{enumerate}
Hence $\cH$ is a unital $\cT$-module.
\end{lemma}

\begin{proof}
(a) For $a\in\End(X)$, $F_\tau[X,(\id_E\otimes\sP(a))g]$ and $F_\tau[X,ga]$ differ by moving
$\id_V\otimes a\otimes\id_{\sP V}$ once around the core (naturality of $U$), so they agree
by Lemma~\ref{lem:W-bimodule}(a); (R1) is respected by naturality.
For the unit, take $V=\one$ and $\tau=\id_E$: the strands $\one$ and $\sP\one$ are
removed by the trace property for $U_0\colon\one\to\sP\one$, \eqref{eq:U-unit} and
$\beta_\one=\sP(U_0)U_0$, as in the proof of $\bdy(\one)=\mob(\one)$
(Proposition~\ref{prop:L}(b)).

\medskip\noindent (b) We have
\begin{align*}
F_{\tau'}F_\tau(w)
&\overset{\text{\ref{item:E3}}}{=}\sum_{V''}d_{V''}
\begin{tikzpicture}[mb]
\mbars{-1.8}{3.0}{6.0}\medge{-1.8}{0}{6.0}
\draw[mstr] (-0.9,0) -- (-0.9,0.6); \draw[mstr] (0,0) -- (0,0.6);
\draw[mstr] (-0.45,0.6) -- (-0.45,1.6);
\draw[mstr] (-0.9,1.6) -- (-0.9,3.9);
\draw[mstr] (0,1.6) -- (0,2.7);
\draw[mstr] (1.0,0) -- (1.0,0.8);
\draw[mstr,rounded corners=4pt] (0.85,0.8) -- (0.85,2.1) -- (0.4,2.1) -- (0.4,2.7);
\draw[mstr] (1.15,0.8) -- (1.15,6.0);
\draw[mstr,rounded corners=4pt] (0,2.7) -- (0,3.3) -- (-0.5,3.3) -- (-0.5,3.9);
\draw[mstr] (0.4,2.7) -- (0.4,6.0);
\draw[mstr] (-0.5,3.9) -- (-0.5,6.0);
\draw[mstr,rounded corners=4pt] (-0.9,3.9) -- (-0.9,4.5) -- (-1.8,4.5);
\draw[mstr] (1.9,0) -- (1.9,6.0); \draw[mstr] (2.6,0) -- (2.6,6.0);
\node[mcpnc,minimum width=11mm] at (-0.45,0.6) {$r$};
\node[mcpnc,minimum width=11mm] at (-0.45,1.6) {$r$};
\node[mcpn,minimum width=8mm] at (0.2,2.7) {$\tau$};
\node[mcpn,minimum width=8mm] at (-0.7,3.9) {$\tau'$};
\node[mcpn] at (-0.5,5.1) {$\beta_{V'}$};
\node[mcpn] at (0.4,5.1) {$\beta_V$};
\node[mcpn] at (1.0,0.8) {$g$};
\legend{-1.8}{4.5}{E''}
\blab{-0.9}{V'}\blab{0}{V}\blab{1.0}{X}\blab{1.9}{\sP V}\blab{2.6}{\sP V'}
\tlab{-0.5}{6.0}{\sP^2V'}\tlab{0.4}{6.0}{\sP^2V}\tlab{1.15}{6.0}{\sP X}\tlab{1.9}{6.0}{\sP V}\tlab{2.6}{6.0}{\sP V'}
\end{tikzpicture}
\overset{\text{Lem.~\ref{lem:W-bimodule}(a)}}{=}\sum_{V''}d_{V''}
\begin{tikzpicture}[mb]
\mbars{-1.8}{3.0}{6.0}\medge{-1.8}{0}{6.0}
\draw[mstr] (-0.45,0) -- (-0.45,1.0);
\draw[mstr] (-0.9,1.0) -- (-0.9,3.9);
\draw[mstr] (0,1.0) -- (0,2.7);
\draw[mstr] (1.0,0) -- (1.0,0.8);
\draw[mstr,rounded corners=4pt] (0.85,0.8) -- (0.85,2.1) -- (0.4,2.1) -- (0.4,2.7);
\draw[mstr] (1.15,0.8) -- (1.15,6.0);
\draw[mstr,rounded corners=4pt] (0,2.7) -- (0,3.3) -- (-0.5,3.3) -- (-0.5,3.9);
\draw[mstr] (0.4,2.7) -- (0.4,6.0);
\draw[mstr] (-0.5,3.9) -- (-0.5,6.0);
\draw[mstr,rounded corners=4pt] (-0.9,3.9) -- (-0.9,4.5) -- (-1.8,4.5);
\draw[mstr] (1.9,0) -- (1.9,5.1); \draw[mstr] (2.6,0) -- (2.6,5.1);
\draw[mstr] (2.25,5.1) -- (2.25,6.0);
\node[mcpnc,minimum width=11mm] at (-0.45,1.0) {$r$};
\node[mcpn,minimum width=8mm] at (0.2,2.7) {$\tau$};
\node[mcpn,minimum width=8mm] at (-0.7,3.9) {$\tau'$};
\node[mcpn] at (-0.5,5.1) {$\beta_{V'}$};
\node[mcpn] at (0.4,5.1) {$\beta_V$};
\node[mcpn] at (1.0,0.8) {$g$};
\node[mcpnc,minimum width=9mm] at (2.25,5.1) {$r$};
\ubox{1.6}{4.7}{2.9}{5.5}
\legend{-1.8}{4.5}{E''}
\blab{-0.45}{V''}\blab{1.0}{X}\blab{1.9}{\sP V}\blab{2.6}{\sP V'}
\tlab{-0.5}{6.0}{\sP^2V'}\tlab{0.4}{6.0}{\sP^2V}\tlab{1.15}{6.0}{\sP X}\tlab{2.25}{6.0}{\sP V''}
\end{tikzpicture}
\\
&\overset{\text{Lem.~\ref{lem:W-bimodule}(a)}}{\underset{\text{\ref{item:B7}}}{=}}\sum_{V''}d_{V''}
\begin{tikzpicture}[mb]
\mbars{-1.8}{2.3}{6.0}\medge{-1.8}{0}{6.0}
\draw[mstr] (-0.45,0) -- (-0.45,1.0);
\draw[mstr] (-0.9,1.0) -- (-0.9,3.9);
\draw[mstr] (0,1.0) -- (0,2.7);
\draw[mstr] (1.0,0) -- (1.0,0.7);
\draw[mstr,rounded corners=4pt] (0.85,0.7) -- (0.85,2.1) -- (0.4,2.1) -- (0.4,2.7);
\draw[mstr] (1.15,0.7) -- (1.15,6.0);
\draw[mstr,rounded corners=4pt] (0,2.7) -- (0,3.3) -- (-0.5,3.3) -- (-0.5,3.9);
\draw[mstr] (0.4,2.7) -- (0.4,4.7);
\draw[mstr] (-0.5,3.9) -- (-0.5,4.7);
\draw[mstr] (-0.05,4.7) -- (-0.05,6.0);
\draw[mstr,rounded corners=4pt] (-0.9,3.9) -- (-0.9,4.5) -- (-1.8,4.5);
\draw[mstr] (1.9,0) -- (1.9,6.0);
\node[mcpnc,minimum width=11mm] at (-0.45,1.0) {$r$};
\node[mcpn,minimum width=8mm] at (0.2,2.7) {$\tau$};
\node[mcpn,minimum width=8mm] at (-0.7,3.9) {$\tau'$};
\node[mcpnc,minimum width=12mm] at (-0.05,4.7) {$r$};
\node[mcpn] at (-0.05,5.45) {$\beta_{V''}$};
\node[mcpn] at (1.0,0.7) {$g$};
\legend{-1.8}{4.5}{E''}
\blab{-0.45}{V''}\blab{1.0}{X}\blab{1.9}{\sP V''}
\tlab{-0.05}{6.0}{\sP^2V''}\tlab{1.15}{6.0}{\sP X}\tlab{1.9}{6.0}{\sP V''}
\end{tikzpicture}
=\sum_{V''}d_{V''}F_{(\tau'*\tau)_{V''}}(w)=F_{\tau'*\tau}(w).
\qedhere
\end{align*}
\end{proof}

\begin{lemma}\label{lem:tube-bimodule}
For $\tau\in\Hom(V\otimes E,E'\otimes V)$, $u,v\in\fX$ and $w\in\cH_E$,
$F_\tau(u\cdot w\cdot v)=u\cdot F_\tau(w)\cdot v$.
Hence $\cH$ is a module over $\cT\otimes\fX\otimes\fX^{\mathrm{op}}$.
\end{lemma}

\begin{proof}
Let $u=[Z,h]$ and $w=[X,g]$. For the left action,
\begin{align*}
u\cdot F_\tau(w)
&=
\begin{tikzpicture}[mb]
\mbars{-2.0}{2.0}{4.4}\medge{-2.0}{0}{4.4}
\draw[mstr] (-0.2,0) -- (-0.2,1.5);
\draw[mstr] (0.9,0) -- (0.9,0.55);
\draw[mstr,rounded corners=4pt] (0.75,0.55) -- (0.75,1.0) -- (0.2,1.0) -- (0.2,1.5);
\draw[mstr] (1.05,0.55) -- (1.05,3.2);
\draw[mstr,rounded corners=4pt] (-0.2,1.5) -- (-0.2,2.0) -- (-2.0,2.0);
\draw[mstr] (0.2,1.5) -- (0.2,3.2);
\draw[mstr] (1.7,0) -- (1.7,3.2);
\sstr{0.2}{3.2}{-1.1}{4.4}
\sstr{1.05}{3.2}{0.2}{4.4}
\sstr{1.7}{3.2}{1.05}{4.4}
\draw[mover] (-1.1,0) -- (-1.1,3.2);
\sstr[mover]{-1.1}{3.2}{1.7}{4.4}
\node[mcpn] at (-1.1,0.6) {$h$};
\node[mcpn,minimum width=8mm] at (0,1.5) {$\tau$};
\node[mcpn] at (0.2,2.6) {$\beta_V$};
\node[mcpn] at (0.9,0.55) {$g$};
\legend{-2.0}{2.0}{E'}
\blab{-1.1}{Z}\blab{-0.2}{V}\blab{0.9}{X}\blab{1.7}{\sP V}
\tlab{-1.1}{4.4}{\sP^2V}\tlab{0.2}{4.4}{\sP X}\tlab{1.05}{4.4}{\sP V}\tlab{1.7}{4.4}{\sP Z}
\end{tikzpicture}
\overset{\text{Lem.~\ref{lem:W-bimodule}(b)}}{=}
\begin{tikzpicture}[mb]
\mbars{-2.0}{2.0}{5.4}\medge{-2.0}{0}{5.4}
\sstr{-1.1}{0}{-0.2}{1.0}
\draw[mstr] (-0.2,1.0) -- (-0.2,2.5);
\draw[mstr] (0.9,0) -- (0.9,1.55);
\draw[mstr,rounded corners=4pt] (0.75,1.55) -- (0.75,2.0) -- (0.2,2.0) -- (0.2,2.5);
\draw[mstr] (1.05,1.55) -- (1.05,4.2);
\draw[mstr,rounded corners=4pt] (-0.2,2.5) -- (-0.2,3.0) -- (-2.0,3.0);
\draw[mstr] (0.2,2.5) -- (0.2,4.2);
\draw[mstr] (1.7,0) -- (1.7,5.4);
\sstr{0.2}{4.2}{-1.1}{5.4}
\sstr{1.05}{4.2}{0.2}{5.4}
\sstr[mover]{-0.2}{0}{-1.1}{1.0}
\draw[mover] (-1.1,1.0) -- (-1.1,4.2);
\sstr[mover]{-1.1}{4.2}{1.05}{5.4}
\node[mcpn] at (-1.1,1.7) {$h$};
\node[mcpn,minimum width=8mm] at (0,2.5) {$\tau$};
\node[mcpn] at (0.2,3.6) {$\beta_V$};
\node[mcpn] at (0.9,1.55) {$g$};
\legend{-2.0}{3.0}{E'}
\blab{-1.1}{V}\blab{-0.2}{Z}\blab{0.9}{X}\blab{1.7}{\sP V}
\tlab{-1.1}{5.4}{\sP^2V}\tlab{0.2}{5.4}{\sP X}\tlab{1.05}{5.4}{\sP Z}\tlab{1.7}{5.4}{\sP V}
\end{tikzpicture}
\\
&\overset{\text{isotopy}}{=}
\begin{tikzpicture}[mb]
\mbars{-1.4}{3.1}{3.8}\medge{-1.4}{0}{3.8}
\draw[mstr] (-0.2,0) -- (-0.2,1.9);
\draw[mstr] (1.9,0) -- (1.9,0.55);
\draw[mstr,rounded corners=4pt] (1.75,0.55) -- (1.75,1.25) -- (0.2,1.25) -- (0.2,1.9);
\draw[mstr] (2.05,0.55) -- (2.05,2.6);
\draw[mstr,rounded corners=4pt] (-0.2,1.9) -- (-0.2,2.4) -- (-1.4,2.4);
\draw[mstr] (0.2,1.9) -- (0.2,3.8);
\draw[mstr] (2.8,0) -- (2.8,3.8);
\sstr{2.05}{2.6}{1.0}{3.8}
\draw[mover] (1.0,0) -- (1.0,2.6);
\sstr[mover]{1.0}{2.6}{2.05}{3.8}
\node[mcpn] at (1.0,0.55) {$h$};
\node[mcpn,minimum width=8mm] at (0,1.9) {$\tau$};
\node[mcpn] at (0.2,2.95) {$\beta_V$};
\node[mcpn] at (1.9,0.55) {$g$};
\legend{-1.4}{2.4}{E'}
\blab{-0.2}{V}\blab{1.0}{Z}\blab{1.9}{X}\blab{2.8}{\sP V}
\tlab{0.2}{3.8}{\sP^2V}\tlab{1.0}{3.8}{\sP X}\tlab{2.05}{3.8}{\sP Z}\tlab{2.8}{3.8}{\sP V}
\end{tikzpicture}
=F_\tau(u\cdot w).
\end{align*}
For the right action the proof is basically the same: the crossing move for the strands $\sP V$
and $Z$ cancels the crossing of $\sP^2V$ with $\sP Z$, and the two remaining crossings of $Z$
with $\sP V$ cancel.
\end{proof}

\begin{lemma}\label{lem:tube-loops}
For $E,Y,V\in\Irr(\cC)$ and $w\in\cH_E$,
\[
F_{c_{Y,E}}(w)=\bdy(Y)\cdot w,\qquad F_{c^{-1}_{E,V}}(w)=w\cdot \bdy(V).
\]
\end{lemma}


\begin{proof}

The leg passes behind $Y$ and $\sP Y$, and
after the crossing move (Lemma~\ref{lem:W-bimodule}(b)) behind $Y$ only:
\[
\bdy(Y)\cdot w=
\begin{tikzpicture}[mb]
\mbars{-1.0}{2.4}{3.0}\medge{-1.0}{0}{3.0}
\draw[mstr] (1.8,0) -- (1.8,0.6);
\draw[mstr] (1.65,0.6) .. controls (1.65,1.1) and (0.5,1.05) .. (-1.0,1.05);
\draw[mstr] (1.95,0.6) -- (1.95,1.5);
\sstr{1.95}{1.5}{0}{2.5}
\draw[mstr] (0,2.5) -- (0,3.0);
\draw[mover] (0,0) -- (0,1.5);
\draw[mover] (0.8,0) -- (0.8,1.5);
\sstr[mover]{0}{1.5}{0.8}{2.5}
\sstr[mover]{0.8}{1.5}{1.6}{2.5}
\draw[mstr] (0.8,2.5) -- (0.8,3.0); \draw[mstr] (1.6,2.5) -- (1.6,3.0);
\node[mcpn] at (0,0.5) {$\beta_Y$};
\node[mcpn] at (1.8,0.6) {$g$};
\legend{-1.0}{1.05}{E}
\blab{0}{Y}\blab{0.8}{\sP Y}\blab{1.8}{X}
\tlab{0}{3.0}{\sP X}\tlab{0.8}{3.0}{\sP^2Y}\tlab{1.6}{3.0}{\sP Y}
\end{tikzpicture}
\;=\;
\begin{tikzpicture}[mb]
\mbars{-1.4}{1.8}{2.4}\medge{-1.4}{0}{2.4}
\draw[mstr] (0.5,0) -- (0.5,0.7);
\draw[mstr] (0.65,0.7) -- (0.65,2.4);
\draw[mstr] (0.35,0.7) .. controls (0.35,1.5) and (-0.2,1.5) .. (-1.4,1.5);
\draw[mover] (-0.6,0) -- (-0.6,2.4);
\node[mcpn] at (-0.6,0.6) {$\beta_Y$};
\draw[mstr] (1.4,0) -- (1.4,2.4);
\node[mcpn] at (0.5,0.7) {$g$};
\legend{-1.4}{1.5}{E}
\blab{-0.6}{Y}\blab{0.5}{X}\blab{1.4}{\sP Y}
\tlab{-0.6}{2.4}{\sP^2Y}\tlab{0.65}{2.4}{\sP X}\tlab{1.4}{2.4}{\sP Y}
\end{tikzpicture}
=F_{c_{Y,E}}(w).
\]
The second identity is entirely analogous.
\end{proof}

We can now determine the structure of $\cH$:
\begin{theorem}\label{thm:W-structure}
As modules over $\cT\otimes\fX\otimes\fX^{\mathrm{op}}$,
\[
\cH\;\cong\bigoplus_{X,X'\in\Sigma}W_{X,X'}\otimes V_X\otimes V_{X'}^*,
\]
where $\cT$ acts on $W_{X,X'}$ (Definition~\ref{def:pants}) and $\fX$ acts on
$V_X\otimes V_{X'}^*$ by $u\,(x\otimes\varphi)\,v=ux\otimes\varphi v$.
\end{theorem}


\begin{proof}
By Lemmas~\ref{lem:tube-action} and~\ref{lem:tube-bimodule} and Proposition~\ref{prop:tube},
\[
\cH\cong\bigoplus_{A,B\in\Irr(\cC)}W_{A,B}\otimes K_{A,B}
\]
as modules over $\cT\otimes\fX\otimes\fX^{\mathrm{op}}$, where the multiplicity spaces $K_{A,B}$
are $\fX$-bimodules.
By Lemma~\ref{lem:tube-loops}, left multiplication by $\bdy(Y)$ on $\cH$ is the action of the
element $\sum_Ec_{Y,E}$ of $\cT$, which acts on $W_{A,B}$ by $\chi_A([Y])$
(Lemma~\ref{lem:pants}(a)); so $\bdy(Y)$ acts on $K_{A,B}$ from the
left by $\chi_A([Y])$, and likewise from the right by $\chi_B([Y])$.
By \eqref{eq:fusion-idempotent}, $p_X=\bdy(e_X)$ therefore acts on $K_{A,B}$ from the left by
$\chi_A(e_X)=\delta_{A,X}$, and from the right by $\delta_{B,X}$.
So $K_{A,B}$ is a unital module over $p_A\fX\otimes(p_B\fX)^{\mathrm{op}}$, which is zero unless
$A,B\in\Sigma$, and is $\End_\bC(V_A)\otimes\End_\bC(V_B)^{\mathrm{op}}$ if $A,B\in\Sigma$
(Theorem~\ref{thm:main}(b)).
Hence $K_{A,B}=0$ unless $A,B\in\Sigma$, and
$K_{X,X'}\cong(V_X\otimes V_{X'}^*)^{\oplus k_{X,X'}}$ for $X,X'\in\Sigma$, with some $k_{X,X'}\ge0$.

It remains to show $k_{X,X'}=1$ for each fixed pair $X,X'\in\Sigma$; note that $k_{X,X'}$ does not depend on $E$.
For $E\in\Irr(\cC)$, $\id_E\in\cT$ acts on $\cH$ as the projection onto $\cH_E$
(Lemma~\ref{lem:tube-action}(a)), and on $W_{X,X'}$ as the projection onto $\Hom(X,E\otimes X')$.
So $p_X\cH_Ep_{X'}\cong\Hom(X,E\otimes X')\otimes K_{X,X'}$, of dimension $N_{EX'}^X\,k_{X,X'}\,M_XM_{X'}$
by Theorem~\ref{thm:main}(c), where $N_{EX'}^X=\dim\Hom(X,E\otimes X')$.
By Corollary~\ref{cor:W-weighted} and $\sum_Ed_EN_{EX'}^X=d_{X\otimes X'^\vee}=d_Xd_{X'}$,
\[
k_{X,X'}\,d_Xd_{X'}M_XM_{X'}=d_Xd_{X'}M_XM_{X'},
\]
and $M_X,M_{X'}>0$ as $X,X'\in\Sigma$. Hence $k_{X,X'}=1$.
\end{proof}

Applying the projection $\id_E\in\cT$, we get the form of the theorem stated in the introduction as
Theorem~\ref{thm:bimodule-intro}:
\begin{corollary}\label{cor:W-structure}
For $E\in\cC$ there is an isomorphism of $\fX$-bimodules
\[
\cH_E\;\cong\bigoplus_{X,X'\in\Sigma}\bigl(V_X\otimes V_{X'}^*\bigr)^{\oplus N_{EX'}^X},\qquad N_{EX'}^X=\dim\Hom(X,E\otimes X').
\]
In particular
$\dim p_X\cH_Ep_{X'}=N_{EX'}^XM_XM_{X'}$.
\end{corollary}

\begin{proof}
Both sides are additive in $E$, so let $E$ be simple.
Then $\id_E\in\cT$ projects $\cH$ onto $\cH_E$ and $W_{X,X'}$ onto $\Hom(X,E\otimes X')$.
\end{proof}

\subsection{Constancy of the twist}
\label{subsec:theta}

Recall that $\Sigma=\{X\in\Irr(\cC)\mid V_X\neq0\}$.
The aim of this subsection is to show that $\theta_X$ is, on $\Sigma$, a constant times
the scalar $\obs_X$ by which a central element $\mob(\obstruction)$ acts on $V_X$, and to give an
explicit formula for the constant.
When $\beta$ is coherent, $\obs_X=1$, so $\theta_X$ is constant on $\Sigma$.

We first record the square of $\mob(A)$ for an invertible $A$.

\begin{lemma}\label{lem:Ma-square}
Let $A\in\Inv(\cC)$ with $\sP A\cong A$, so that $A\in\Fix(\sP)$.
Then $\mob(A)^2=\theta_A^{-1}\bdy(A)$; here $\theta_A=\pm1$ by Lemma~\ref{lem:P-modular-data}.
\end{lemma}

\begin{proof}
Write $\xi:=\xi_A$.  Then
\begin{align*}
\mob(A)\cdot \mob(A)
&\overset{\text{Prop.~\ref{prop:product}}}{=}
\begin{tikzpicture}[mb]
\mbars{-0.5}{1.4}{2.4}
\draw[mstr] (0,0) -- (0,0.9); \draw[mstr] (0.9,0) -- (0.9,0.9);
\sstr{0.9}{0.9}{0}{2.4}
\sstr[mover]{0}{0.9}{0.9}{2.4}
\node[mcpn] at (0,0.6) {$\xi$}; \node[mcpn] at (0.9,0.6) {$\xi$};
\blab{0}{A}\blab{0.9}{A}\tlab{0}{2.4}{\sP A}\tlab{0.9}{2.4}{\sP A}
\end{tikzpicture}
\overset{\text{Lem.~\ref{lem:P-modular-data}}}{=}
\theta_A^{-1}\,
\begin{tikzpicture}[mb]
\mbars{-0.5}{1.4}{2.4}
\draw[mstr] (0,0) -- (0,2.4); \draw[mstr] (0.9,0) -- (0.9,2.4);
\node[mcpn] at (0,1.2) {$\xi$}; \node[mcpn] at (0.9,1.2) {$\xi$};
\blab{0}{A}\blab{0.9}{A}\tlab{0}{2.4}{\sP A}\tlab{0.9}{2.4}{\sP A}
\end{tikzpicture}
\\
&\overset{\text{Lem.~\ref{lem:trace-property}}}{=}
\theta_A^{-1}\,
\begin{tikzpicture}[mb]
\mbars{-0.5}{1.4}{2.4}
\draw[mstr] (0,0) -- (0,2.4); \draw[mstr] (0.9,0) -- (0.9,2.4);
\node[mcpn] at (0,0.6) {$\xi$};
\node[mcpn] at (0,1.6) {$\xi$};
\ubox{-0.45}{1.15}{0.45}{2.05}
\blab{0}{A}\blab{0.9}{\sP A}\tlab{0}{2.4}{\sP^2A}\tlab{0.9}{2.4}{\sP A}
\end{tikzpicture}
\overset{\eqref{eq:equivariant}}{=}
\theta_A^{-1}\,
\begin{tikzpicture}[mb]
\mbars{-0.5}{1.4}{2.4}
\draw[mstr] (0,0) -- (0,2.4); \draw[mstr] (0.9,0) -- (0.9,2.4);
\node[mcpn] at (0,1.2) {$\beta_A$};
\blab{0}{A}\blab{0.9}{\sP A}\tlab{0}{2.4}{\sP^2A}\tlab{0.9}{2.4}{\sP A}
\end{tikzpicture}
\overset{\text{Def.~\ref{def:L}}}{=}
\theta_A^{-1}\bdy(A).
\end{align*}
\end{proof}

We next describe the conjugation by $\mob(A)$ for an invertible $A$.

\begin{lemma}\label{lem:Ma-conj}
Let $A\in\Inv(\cC)$ with $\sP A\cong A$, so that $A\in\Fix(\sP)$.
Then $\mob(A)$ is invertible, with $\mob(A)^{-1}=\mob(A)^\star$, and for $E\in\cC$ and
$[Z,g]\in\cH_E$,
\[
\mob(A)\cdot[Z,g]\cdot \mob(A)^{-1}=[Z,\,g\circ\ell_{A,Z}].
\]
\end{lemma}

\begin{proof}
By Proposition~\ref{prop:hilbert}(c) and Definition~\ref{def:product},
$\mob(A)^\star \mob(A)$ is represented by a morphism on $A^\vee\otimes A\cong\one$, so by
Proposition~\ref{prop:basis} it is a multiple of the unit $\mob(\one)$, and the multiple
is $\varepsilon(\mob(A)^\star \mob(A))=1$ by Proposition~\ref{prop:hilbert}(d).
As $\fX$ is finite-dimensional, $\mob(A)^{-1}=\mob(A)^\star$.

For the conjugation, write $\xi:=\xi_A$.
Then
\begin{align*}
\mob(A)\cdot\picZg
&\overset{\text{Lem.~\ref{lem:W-bimodule}(c)}}{=}
\begin{tikzpicture}[mb]
\mbars{-1.2}{1.65}{2.8}\medge{-1.2}{0}{2.8}
\draw[mstr] (1.0,0) -- (1.0,0.8);
\draw[mstr] (1.15,0.8) -- (1.15,1.8);
\draw[mstr] (0.85,0.8) .. controls (0.85,1.35) and (0.4,1.35) .. (-1.2,1.35);
\draw[mover] (0,0) -- (0,1.8);
\sstr{1.15}{1.8}{0}{2.8}
\sstr[mover]{0}{1.8}{1.15}{2.8}
\node[mcpn] at (0,0.6) {$\xi$};
\node[mcpn] at (1.0,0.8) {$g$};
\legend{-1.2}{1.35}{E}
\blab{0}{A}\blab{1.0}{Z}\tlab{0}{2.8}{\sP Z}\tlab{1.15}{2.8}{\sP A}
\end{tikzpicture}
\overset{\text{Lem.~\ref{lem:W-bimodule}(b)}}{=}
\begin{tikzpicture}[mb]
\mbars{-1.2}{1.65}{3.0}\medge{-1.2}{0}{3.0}
\sstr{0}{0}{1.0}{0.9}
\sstr[mover]{1.0}{0}{0}{0.9}
\draw[mstr] (1.0,0.9) -- (1.0,1.6);
\draw[mstr] (1.15,1.6) -- (1.15,3.0);
\draw[mstr] (0.85,1.6) .. controls (0.85,2.15) and (0.4,2.15) .. (-1.2,2.15);
\draw[mover] (0,0.9) -- (0,3.0);
\node[mcpn] at (0,1.35) {$\xi$};
\node[mcpn] at (1.0,1.6) {$g$};
\legend{-1.2}{2.15}{E}
\blab{0}{Z}\blab{1.0}{A}\tlab{0}{3.0}{\sP A}\tlab{1.15}{3.0}{\sP Z}
\end{tikzpicture}
=
\begin{tikzpicture}[mb]
\mbars{-1.2}{1.65}{2.8}\medge{-1.2}{0}{2.8}
\draw[mstr] (0,0) -- (0,0.8);
\draw[mstr] (0.15,0.8) -- (0.15,1.8);
\draw[mstr] (-0.15,0.8) .. controls (-0.15,1.35) and (-0.6,1.35) .. (-1.2,1.35);
\draw[mstr] (1.15,0) -- (1.15,1.8);
\sstr{0.15}{1.8}{1.15}{2.8}
\sstr[mover]{1.15}{1.8}{0.15}{2.8}
\node[mcpn] at (0,0.8) {$g$};
\node[mcpn] at (1.15,0.6) {$\xi$};
\legend{-1.2}{1.35}{E}
\blab{0}{Z}\blab{1.15}{A}\tlab{0.15}{2.8}{\sP A}\tlab{1.15}{2.8}{\sP Z}
\end{tikzpicture}
\\
&\overset{\text{Lem.~\ref{lem:monodromy}}}{=}
\begin{tikzpicture}[mb]
\mbars{-1.2}{1.65}{3.0}\medge{-1.2}{0}{3.0}
\draw[mstr] (0,0) -- (0,0.7);
\draw[mstr] (0.15,0.7) -- (0.15,2.0);
\draw[mstr] (-0.15,0.7) .. controls (-0.15,1.2) and (-0.6,1.2) .. (-1.2,1.2);
\draw[mstr] (1.15,0) -- (1.15,2.0);
\sstr{1.15}{2.0}{0.15}{3.0}
\sstr[mover]{0.15}{2.0}{1.15}{3.0}
\node[mcpn] at (0,0.7) {$g$};
\node[mcpn] at (0.15,1.5) {$\ell_A^{-1}$};
\node[mcpn] at (1.15,0.6) {$\xi$};
\legend{-1.2}{1.2}{E}
\blab{0}{Z}\blab{1.15}{A}\tlab{0.15}{3.0}{\sP A}\tlab{1.15}{3.0}{\sP Z}
\end{tikzpicture}
\overset{\text{Lem.~\ref{lem:W-bimodule}(a)}}{\underset{\text{Lem.~\ref{lem:box-loop}}}{=}}
\begin{tikzpicture}[mb]
\mbars{-1.2}{1.65}{3.0}\medge{-1.2}{0}{3.0}
\draw[mstr] (0,0) -- (0,1.35);
\draw[mstr] (0.15,1.35) -- (0.15,2.0);
\draw[mstr] (-0.15,1.35) .. controls (-0.15,1.75) and (-0.6,1.75) .. (-1.2,1.75);
\draw[mstr] (1.15,0) -- (1.15,2.0);
\sstr{1.15}{2.0}{0.15}{3.0}
\sstr[mover]{0.15}{2.0}{1.15}{3.0}
\node[mcpn] at (0,0.55) {$\ell_A$};
\node[mcpn] at (0,1.35) {$g$};
\node[mcpn] at (1.15,0.6) {$\xi$};
\legend{-1.2}{1.75}{E}
\blab{0}{Z}\blab{1.15}{A}\tlab{0.15}{3.0}{\sP A}\tlab{1.15}{3.0}{\sP Z}
\end{tikzpicture}
\overset{\text{Lem.~\ref{lem:W-bimodule}(c)}}{=}
[Z,g\circ\ell_{A,Z}]\cdot \mob(A).
\end{align*}
Multiplying by $\mob(A)^{-1}$ on the right gives the claim.
\end{proof}

We apply this to the coherence obstruction $\obstruction$ of $\beta$ (Proposition~\ref{prop:kappa}).

\begin{proposition}\label{prop:Mkappa}
$\obstruction\in\Fix(\sP)$, so that $\mob(\obstruction)\in\fX$ is defined.
It is invertible, with $\mob(\obstruction)^{-1}=\mob(\obstruction)^\star$, and
$\mob(\obstruction)=\mob(\one)$ when $\beta$ is coherent.
\end{proposition}

\begin{proof}
$\sP\obstruction\cong\obstruction$ by Proposition~\ref{prop:obstruction}(a), and
Lemma~\ref{lem:Ma-conj} gives the inverse.
When $\beta$ is coherent, $\obstruction=\one$ (Definition~\ref{def:coherent}), so
$\mob(\obstruction)=\mob(\one)$.
\end{proof}

\begin{lemma}\label{lem:W-wind}
For $E\in\Irr(\cC)$, 
let $\omega_E\in\cT$ be $\id_{E\otimes E}$, regarded as an element of
$\Hom(V\otimes E,E\otimes V)$ with $V=E$:
\[
\omega_E=
\begin{tikzpicture}[mb]
\draw[mstr,rounded corners=4pt] (-0.25,0) -- (-0.25,2.25) -- (-1.3,2.25);
\draw[mstr,rounded corners=4pt] (1.3,0.75) -- (0.25,0.75) -- (0.25,3);
\tubemark{-0.25}{0.22}\tubemark{0.25}{2.78}
\node[anchor=west,inner sep=1.5pt] at (1.3,0.75) {$E$};
\node[anchor=east,inner sep=1.5pt] at (-1.3,2.25) {$E$};
\blab{-0.25}{E}\tlab{0.25}{3}{E}
\end{tikzpicture}\,.
\]
Then $F_{\omega_E}(w)=\mob(\obstruction)\,w\,\mob(\obstruction)^{-1}$ for $w\in\cH_E$;
in particular $F_{\omega_E}=\id$ when $\beta$ is coherent.
That is, for every object $Z$ and every $g\colon Z\to E\otimes\sP Z$,
\[
F_{\omega_E}[Z,g]=
\begin{tikzpicture}[mb]
\mbars{-1.3}{2.0}{2.6}\medge{-1.3}{0}{2.6}
\draw[mstr] (-0.7,0) .. controls (-0.7,0.45) and (-1,0.6) .. (-1.3,0.6);
\draw[mstr] (0.3,0) -- (0.3,0.8);
\draw[mstr] (0.15,0.8) -- (0.15,2.6);
\draw[mstr] (0.45,0.8) .. controls (0.45,1.3) and (0.95,1.2) .. (0.95,1.7) -- (0.95,2.6);
\draw[mstr] (1.6,0) -- (1.6,2.6);
\node[mcpn] at (0.3,0.8) {$g$};
\node[mcpn] at (0.15,1.8) {$\beta_E$};
\legend{-1.3}{0.6}{E}
\blab{-0.7}{E}\blab{0.3}{Z}\blab{1.6}{\sP E}
\tlab{0.15}{2.6}{\sP^2E}\tlab{0.95}{2.6}{\sP Z}\tlab{1.6}{2.6}{\sP E}
\end{tikzpicture}
\;=\;\mob(\obstruction)\cdot\picZg\cdot \mob(\obstruction)^{-1}\,.
\]
\end{lemma}
\begin{proof}
We move $g$ twice around the core.
With $\tilde g:=(\id_E\otimes\ell^{-1}_{\obstruction,\sP Z})\,g$ we have
$(\beta_E\otimes\id_{\sP Z})\,g=(\id_{\sP^2E}\otimes\sP(\beta_Z^{-1}))\,\Ubox{\Ubox{\tilde g}}\,\beta_Z$, and
\begin{multline*}
F_{\omega_E}[Z,g]=
\begin{tikzpicture}[mb]
\mbars{-1.3}{2.3}{3.4}\medge{-1.3}{0}{3.4}
\draw[mstr] (-0.7,0) .. controls (-0.7,0.45) and (-1,0.6) .. (-1.3,0.6);
\draw[mstr] (0.3,0) -- (0.3,0.9);
\draw[mstr] (0.15,0.9) -- (0.15,3.4);
\draw[mstr] (0.45,0.9) .. controls (0.45,1.4) and (1.05,1.3) .. (1.05,1.8) -- (1.05,3.4);
\draw[mstr] (1.9,0) -- (1.9,3.4);
\node[mcpn] at (0.3,0.9) {$g$};
\node[mcpn] at (0.15,2.2) {$\beta_E$};
\legend{-1.3}{0.6}{E}
\blab{-0.7}{E}\blab{0.3}{Z}\blab{1.9}{\sP E}
\tlab{0.15}{3.4}{\sP^2E}\tlab{1.05}{3.4}{\sP Z}\tlab{1.9}{3.4}{\sP E}
\end{tikzpicture}
\overset{\text{\ref{item:B7}}}{\underset{\eqref{eq:obstruction}}{=}}
\begin{tikzpicture}[mb]
\mbars{-1.3}{2.3}{3.4}\medge{-1.3}{0}{3.4}
\draw[mstr] (-0.7,0) .. controls (-0.7,0.45) and (-1,0.6) .. (-1.3,0.6);
\draw[mstr] (0.3,0) -- (0.3,1.5);
\draw[mstr] (0.15,1.5) -- (0.15,3.4);
\draw[mstr] (0.45,1.5) -- (0.45,2.1) .. controls (0.45,2.35) and (1.05,2.25) .. (1.05,2.5) -- (1.05,3.4);
\draw[mstr] (1.9,0) -- (1.9,3.4);
\node[mcpn] at (0.3,0.5) {$\beta_Z$};
\node[mcpn] at (0.3,1.5) {$\tilde g$};
\ubox{-0.15}{1.05}{0.75}{1.95}
\ubox{-0.3}{0.9}{0.9}{2.1}
\node[mcpn] at (1.05,2.85) {$\beta^{-1}_Z$};
\ubox{0.5}{2.42}{1.6}{3.28}
\legend{-1.3}{0.6}{E}
\blab{-0.7}{E}\blab{0.3}{Z}\blab{1.9}{\sP E}
\tlab{0.15}{3.4}{\sP^2E}\tlab{1.05}{3.4}{\sP Z}\tlab{1.9}{3.4}{\sP E}
\end{tikzpicture}
\overset{\text{\ref{item:B1}--\ref{item:B3}}}{=}
\begin{tikzpicture}[mb]
\mbars{-1.3}{2.3}{3.4}\medge{-1.3}{0}{3.4}
\draw[mstr] (-0.7,0) .. controls (-0.7,0.45) and (-1,0.6) .. (-1.3,0.6);
\draw[mstr] (0.3,0) -- (0.3,1.45);
\draw[mstr] (0.15,1.45) -- (0.15,3.4);
\draw[mstr] (0.45,1.45) -- (0.45,1.9) .. controls (0.45,2.1) and (1.05,2.05) .. (1.05,2.25) -- (1.05,3.4);
\draw[mstr] (1.75,0) -- (1.75,3.4);
\node[mcpn] at (0.3,0.5) {$\beta_Z$};
\node[mcpn] at (0.3,1.45) {$\tilde g$};
\ubox{-0.15}{1.0}{0.75}{1.9}
\node[mcpn] at (1.05,2.6) {$\beta^{-1}_Z$};
\ubox{-0.3}{0.85}{2.05}{3.0}
\legend{-1.3}{0.6}{E}
\blab{-0.7}{E}\blab{0.3}{Z}\blab{1.75}{\sP E}
\tlab{0.15}{3.4}{\sP^2E}\tlab{1.05}{3.4}{\sP Z}\tlab{1.75}{3.4}{\sP E}
\end{tikzpicture}
\\
\overset{\text{Lem.~\ref{lem:W-bimodule}(a)}}{=}
\begin{tikzpicture}[mb]
\mbars{-1.1}{1.4}{2.4}\medge{-1.1}{0}{2.4}
\draw[mstr] (-0.55,0) .. controls (-0.55,0.45) and (-0.85,0.6) .. (-1.1,0.6);
\draw[mstr] (0.3,0) -- (0.3,0.95);
\draw[mstr] (0.15,0.95) -- (0.15,2.4);
\draw[mstr] (0.45,0.95) -- (0.45,1.4) .. controls (0.45,1.7) and (1.05,1.6) .. (1.05,1.9) -- (1.05,2.4);
\node[mcpn] at (0.3,0.95) {$\tilde g$};
\ubox{-0.15}{0.5}{0.75}{1.4}
\legend{-1.1}{0.6}{E}
\blab{-0.55}{E}\blab{0.3}{\sP Z}
\tlab{0.15}{2.4}{\sP^2Z}\tlab{1.05}{2.4}{\sP E}
\end{tikzpicture}
\overset{\text{Lem.~\ref{lem:W-bimodule}(a)}}{=}
\begin{tikzpicture}[mb]
\mbars{-1}{0.6}{2.8}\medge{-1}{0}{2.8}
\draw[mstr] (0,0) -- (0,1);
\draw[mstr] (0.15,1) -- (0.15,2.8);
\draw[mstr] (-0.15,1) .. controls (-0.15,1.45) and (-0.6,1.45) .. (-1,1.45);
\node[mcpn] at (0,1) {$g$};
\node[mcpn] at (0.15,2.05) {$\ell^{-1}_\obstruction$};
\legend{-1}{1.45}{E}
\blab{0}{Z}\tlab{0.15}{2.8}{\sP Z}
\end{tikzpicture}
\overset{\text{Lem.~\ref{lem:box-loop}}}{\underset{\text{(R2)}}{=}}
\begin{tikzpicture}[mb]
\mbars{-1}{0.6}{2.8}\medge{-1}{0}{2.8}
\draw[mstr] (0,0) -- (0,1.75);
\draw[mstr] (0.15,1.75) -- (0.15,2.8);
\draw[mstr] (-0.15,1.75) .. controls (-0.15,2.2) and (-0.6,2.2) .. (-1,2.2);
\node[mcpn] at (0,0.75) {$\ell_\obstruction$};
\node[mcpn] at (0,1.75) {$g$};
\legend{-1}{2.2}{E}
\blab{0}{Z}\tlab{0.15}{2.8}{\sP Z}
\end{tikzpicture}
\overset{\text{Lem.~\ref{lem:Ma-conj}}}{=}
\mob(\obstruction)\cdot[Z,g]\cdot \mob(\obstruction)^{-1}.
\end{multline*}
\end{proof}

\begin{lemma}\label{lem:omega}
For $A,B,E\in\Irr(\cC)$, $F_{\omega_E}$ acts on $\Hom(A,E\otimes B)\subset W_{A,B}$ by
$\theta_A\theta_B^{-1}$.
\end{lemma}

\begin{proof}
We have
\[
F_{\omega_E}\,\psi\overset{\text{Def.~\ref{def:pants}}}{=}
\begin{tikzpicture}[mb]
\draw[mstr,rounded corners=6pt] (1.4,0.8) -- (1.4,-0.7) -- (-0.85,-0.7) -- (-0.85,1.15) -- (-1.5,1.15);
\draw[mover] (0.6,-1.3) -- (0.6,0.0);
\draw[mstr] (0.6,0.0) -- (0.6,2.2);
\draw[mover,rounded corners=6pt] (0.6,0.0) -- (-0.35,0.0) -- (-0.35,1.6) -- (1.4,1.6) -- (1.4,0.8);
\node[mcpn] at (0.6,0.0) {$\psi$};
\node[anchor=east,inner sep=1.5pt] at (-1.5,1.15) {$E$};
\node[anchor=base,inner sep=0pt,yshift=-8.8pt] at (0.6,-1.3) {$A$};
\tlab{0.6}{2.2}{B}
\end{tikzpicture}
=
\begin{tikzpicture}[mb]
\draw[mstr,rounded corners=6pt] (1.4,1.2) -- (1.4,0.8) -- (-1.2,0.8);
\draw[mstr] (0.6,-1.0) -- (0.6,0.0);
\draw[mover] (0.6,0.0) -- (0.6,2.2);
\draw[mover,rounded corners=6pt] (0.6,0.0) -- (-0.35,0.0) -- (-0.35,1.6) -- (1.4,1.6) -- (1.4,1.2);
\node[mcpn] at (0.6,0.0) {$\psi$};
\node[anchor=east,inner sep=1.5pt] at (-1.2,0.8) {$E$};
\node[anchor=base,inner sep=0pt,yshift=-8.8pt] at (0.6,-1.0) {$A$};
\tlab{0.6}{2.2}{B}
\end{tikzpicture}
=\theta_E\,
\begin{tikzpicture}[mb]
\draw[mstr] (0.5,-0.6) -- (0.5,0.2);
\draw[mstr] (0,0.2) -- (0,0.6); \draw[mstr] (1.0,0.2) -- (1.0,0.6);
\sstr{1.0}{0.6}{0}{1.4} \sstr[mover]{0}{0.6}{1.0}{1.4}
\sstr{1.0}{1.4}{0}{2.2} \sstr[mover]{0}{1.4}{1.0}{2.2}
\draw[mstr] (0,2.2) -- (0,2.6); \draw[mstr] (1.0,2.2) -- (1.0,2.6);
\node[mcpn,minimum width=12mm] at (0.5,0.2) {$\psi$};
\node[anchor=north,inner sep=1.8pt] at (0.5,-0.6) {$A$};
\tlab{0}{2.6}{E}\tlab{1.0}{2.6}{B}
\end{tikzpicture}
=\theta_A\theta_B^{-1}\,\psi .
\qedpic
\]
\end{proof}

\begin{lemma}\label{lem:mX}
$\mob(\obstruction)$ is central in $\fX$.
\end{lemma}

\begin{proof}
For $E=\one$, $\omega_\one=\id_\one$ with $V=\one$, so $F_{\omega_\one}$ is the identity of
$\cH_\one=\fX$ by Lemma~\ref{lem:tube-action}(a).
So $\mob(\obstruction)\,w\,\mob(\obstruction)^{-1}=F_{\omega_\one}(w)=w$ for $w\in\fX$ by
Lemma~\ref{lem:W-wind}.
\end{proof}

\begin{definition}\label{def:mX}
For $X\in\Sigma$ let $\obs_X$ be the scalar by which the central element
$\mob(\obstruction)$ acts on the irreducible module $V_X$, so that
$\mob(\obstruction)\,p_X=\obs_X\,p_X$.
By Lemma~\ref{lem:Ma-square}, $\mob(\obstruction)^2=\theta_\obstruction^{-1}\bdy(\obstruction)$, and
$\bdy(\obstruction)$ acts on $V_X$ by $\chi_X([\obstruction])$, so
\[
\obs_X^2=\theta_\obstruction^{-1}\,\chi_X([\obstruction]);
\]
in particular $|\obs_X|=1$, as $\theta_\obstruction=\pm1$ and $\chi_X([\obstruction])$ is a root of unity.
When $\beta$ is coherent, $\obs_X=1$.
\end{definition}

\begin{proposition}\label{prop:theta-const}
The ratio $\theta_X/\obs_X$ is constant on $\Sigma$.
In particular, $\theta_X$ is constant on $\Sigma$ when $\beta$ is coherent.
\end{proposition}

\begin{proof}
Let $X,X'\in\Sigma$, and choose $E\in\Irr(\cC)$ with
$N_{EX'}^X=\dim\Hom(E,X\otimes X'^\vee)\ge1$.
By Theorem~\ref{thm:W-structure}, $p_X\cH_Ep_{X'}$ is the non-zero summand
$\Hom(X,E\otimes X')\otimes V_X\otimes V_{X'}^*$ of $\cH$.
On it
$F_{\omega_E}$ acts by $\theta_X\theta_{X'}^{-1}$ (Lemma~\ref{lem:omega}) and by
$\obs_X\obs_{X'}^{-1}$ (Lemma~\ref{lem:W-wind} and Definition~\ref{def:mX}).
Hence $\theta_X/\obs_X=\theta_{X'}/\obs_{X'}$.
\end{proof}

We next determine the constant value of $\theta_X/\obs_X$.
Let $p_+:=\sum_{X\in\Irr(\cC)}d_X^2\theta_X$ be the Gauss sum,
and $\lambda:=p_+/D$.
It is well-known that $\lambda$ is a root of unity \cite[Corollary~8.18.2]{EGNO}.

\begin{lemma}\label{lem:gauss}
$\lambda=\pm1$ when $\cC$ has an anti-monoidal involution $\sP$.
\end{lemma}

\begin{proof}
Since $X\mapsto\sP X$ permutes
$\Irr(\cC)$ with $d_{\sP X}=d_X$ (Corollary~\ref{cor:P-dim}) and
$\theta_{\sP X}=\theta_X^{-1}=\overline{\theta_X}$
(Lemma~\ref{lem:P-modular-data}; the twists are roots of unity
\cite{Vafa1988,AndersonMoore1988}, \cite[Corollary~8.18.2]{EGNO}), $p_+=\overline{p_+}$ is real.
Hence $\lambda=\pm1$.
\end{proof}

We also define, for $Y\in\Irr(\cC)$,
\[
\eta^\obstruction_Y:=\varepsilon\bigl(\bdy(Y^\vee)\,\mob(\obstruction)\bigr),
\qquad
\sigma:=\frac1D\sum_{Y\in\Irr(\cC)}\overline{\eta^\obstruction_Y}\,d_Y\,\theta_Y^{-1} .
\]

We can now prove our second main result, which was stated in the introduction as Theorem~\ref{thm:constancy-intro} for coherent $\beta$:
\begin{theorem}\label{thm:theta}
We have $\theta_X=\lambda\sigma\,\obs_X$ for $X\in\Sigma$, and $|\sigma|=1$.
\end{theorem}

\begin{proof}
By Proposition~\ref{prop:theta-const}, $\theta_X=k\,\obs_X$ for $X\in\Sigma$, for some $k$
with $|k|=1$.
By the balancing axiom,
$s_{X,Y}=\theta_X^{-1}\theta_Y^{-1}\operatorname{tr}(\theta_{X\otimes Y})
=\theta_X^{-1}\theta_Y^{-1}\sum_ZN_{XY}^Z\,\theta_Zd_Z$,
where $N_{XY}^Z=\dim\Hom(X\otimes Y,Z)$.
As $\sum_XN_{XY}^Zd_X=d_{Z\otimes Y^\vee}=d_Zd_Y$, this gives
\[
\sum_Xd_X\theta_XS_{X,Y}
=\frac{\theta_Y^{-1}}{D}\sum_Z\theta_Zd_Z\sum_XN_{XY}^Zd_X
=\lambda\,d_Y\theta_Y^{-1}.
\]
Therefore
\[
k\,\mob(\obstruction)\overset{\mathrm{(1)}}{=}\sum_{X\in\Irr(\cC)}\theta_X\,p_X
\overset{\mathrm{(2)}}{=}\bdy\Bigl(\sum_{X\in\Irr(\cC)}\theta_X\,e_X\Bigr)
\overset{\mathrm{(3)}}{=}\frac{\lambda}{D}\sum_{A\in\Irr(\cC)}d_A\,\theta_A^{-1}\,\bdy(A).
\]

(1) $\mob(\obstruction)\,p_X=\obs_X\,p_X$ (Definition~\ref{def:mX}), $p_X=0$ for $X\notin\Sigma$,
and $\sum_Xp_X=\mob(\one)$.

(2) $p_X=\bdy(e_X)$.

(3) By \eqref{eq:fusion-idempotent} and $\overline{S_{A,X}}=S_{X,A^\vee}$,
$\sum_X\theta_Xe_X=\sum_A[A]\,\frac1D\sum_Xd_X\theta_XS_{X,A^\vee}$, and the relation above for
$Y=A^\vee$, with $d_{A^\vee}=d_A$ and $\theta_{A^\vee}=\theta_A$.

\noindent Multiply by $\mob(\obstruction)^{-1}=\mob(\obstruction)^\star$
(Proposition~\ref{prop:Mkappa}) and apply $\varepsilon$.
By Lemma~\ref{lem:L-star} and Proposition~\ref{prop:hilbert}(d),
$\varepsilon(\bdy(A)\mob(\obstruction)^\star)=\overline{\varepsilon(\mob(\obstruction)\bdy(A^\vee))}
=\overline{\eta^\obstruction_A}$, so
$k=\lambda\sigma$.
Hence $|\sigma|=|k|=1$ by Lemma~\ref{lem:gauss}.
\end{proof}

\begin{remark}\label{rem:theta-square}
Squaring, Definition~\ref{def:mX} gives
$\theta_X^2=\sigma^2\,\theta_\obstruction^{-1}\,\chi_X([\obstruction])$ for $X\in\Sigma$:
the ratio $\theta_X^2/\chi_X([\obstruction])$ is constant on $\Sigma$, with an explicit value.
For $\cC_1$ of Section~\ref{subsec:metaplectic}, with $\obstruction=\epsilon$, this gives
$\theta_{\psi_\pm}^2=-1$.
\end{remark}

In the coherent case this reduces to the following, which was
Theorem~\ref{thm:constancy-intro} in the introduction:
\begin{corollary}\label{cor:theta-coherent}
Suppose that $\beta$ is coherent.
Then $\theta_X=\lambda\sigma$ for all $X\in\Sigma$, where
\[
\sigma=\frac1D\sum_{Y\in\Fix(\sT)}\eta_Y\,d_Y\,\theta_Y ,
\]
and $\lambda\sigma=\pm1$.
\end{corollary}

\begin{proof}
By Proposition~\ref{prop:Mkappa} and Definition~\ref{def:mX}, $\mob(\obstruction)=\mob(\one)$ and
$\obs_X=1$, so $\eta^\obstruction=\eta$, which vanishes off $\Fix(\sT)$ and is $\pm1$ on it
(Proposition~\ref{prop:beta}).
For $Y\in\Fix(\sT)$, $\sP Y\cong Y^\vee$ gives
$\theta_Y^{-1}=\theta_{\sP Y}=\theta_{Y^\vee}=\theta_Y$ (Lemma~\ref{lem:P-modular-data}), so
$\theta_Y=\pm1$.
This gives the formula for $\sigma$, which is therefore real; as $|\sigma|=1$ and
$\lambda=\pm1$ (Theorem~\ref{thm:theta}, Lemma~\ref{lem:gauss}), $\lambda\sigma=\pm1$.
\end{proof}

\begin{remark}\label{rem:plain-indicator}
The sum in Corollary~\ref{cor:theta-coherent} is defined for any $\beta$, but it is
non-zero only in the coherent case:
\[
\frac1D\sum_{Y\in\Fix(\sT)}\eta_Y\,d_Y\,\theta_Y=
\begin{cases}\sigma=\pm1, & \text{if $\beta$ is coherent},\\
0, & \text{otherwise}.\end{cases}
\]
Indeed, apply $\varepsilon$ to the identity
$k\,\mob(\obstruction)=\frac{\lambda}{D}\sum_Ad_A\theta_A^{-1}\bdy(A)$ in the proof of
Theorem~\ref{thm:theta}, without multiplying by $\mob(\obstruction)^{-1}$.
As $k=\lambda\sigma$, $\lambda^2=1$, and
$\varepsilon(\bdy(A))=\eta_{A^\vee}=\overline{\eta_A}$
(Definition~\ref{def:eta} and Proposition~\ref{prop:beta}(a), (c)), this gives
$\frac1D\sum_A\overline{\eta_A}\,d_A\,\theta_A^{-1}=\sigma\,\varepsilon(\mob(\obstruction))$, and
$\varepsilon(\mob(\obstruction))$, the coefficient of $\mob(\one)$ in $\mob(\obstruction)$
(Definition~\ref{def:eps}), is $1$ if $\obstruction\cong\one$ and $0$ otherwise.
Taking complex conjugates gives the claim, as $\eta$ vanishes
off $\Fix(\sT)$ (Proposition~\ref{prop:beta}(a)) and $\sigma$ is real in the coherent case.
So the expression of the anomaly indicator of \cite{Wang_2017}, which is defined
for coherent $\beta$, vanishes when it is evaluated for non-coherent $\beta$, and
Theorem~\ref{thm:theta} uses in its place the sum $\sigma$, twisted by $\mob(\obstruction)$.
\end{remark}

\section{Examples}
\label{sec:examples}

We illustrate the results on three examples.
For pointed
categories everything reduces to a finite abelian group with a
non-degenerate quadratic form and an involution, and we recover the
results of \cite{Lee:2018eqa,Orii_2025,Orii2025H3}.
The second example is the category 
$\cC=\cE\boxtimes\overline{\cE}$ where $\cE$ itself is modular,
with the involution exchanging the two
factors.
There $\fX$ turns out to be the fusion algebra of $\cE$.
The
third concerns the two modular categories with the fusion
rules of $\mathfrak{so}_5$ at level $2$ studied by Barkeshli and Cheng
\cite{Barkeshli:2017rzd}: using Corollary~\ref{cor:theta-coherent}, we show from the modular
data alone that the modular category $\mathfrak{so}_5$ at level $2$ itself admits no coherent anti-monoidal involution lifting
the only possible permutation of its simple objects.
For the other case
we construct $(\sP,U,\beta)$ from the $F$- and $R$-symbols by the
method of Appendix~\ref{app:skeletal}, obtaining $\fX\cong\Mat_2(\bC)$.
We also show that $\mathfrak{so}_5$ at level $2$ does admit anti-monoidal involutions
without coherence, and work out its $\fX\cong\bC\oplus\bC$.

\subsection{Pointed categories}
\label{subsec:pointed}\label{rem:pointed}

Let $\cC$ be pointed.
Then $\Irr(\cC)=A$ is a finite abelian group,
written additively, and the monodromy is a non-degenerate bicharacter
$B(a,b)=s_{a,b}$ with $\theta_{a+b}=\theta_a\theta_bB(a,b)$.
In this subsection, $P$ denotes an automorphism of $A$ with $P^2=1$ and
$\theta_{Pa}=\theta_a^{-1}$ for all $a\in A$, so that
$B(Pa,Pb)=B(a,b)^{-1}$.
An anti-monoidal involution $\sP$ induces such a $P$ by
Lemma~\ref{lem:P-modular-data}, and conversely every such $P$ is induced by
an anti-monoidal involution with coherence (Proposition~\ref{prop:pointed-lift}).
For pointed $\cC$ the numbers $M_a$, and in particular their integrality,
were obtained in \cite{Orii_2025} (see also \cite[(8)--(9)]{Wang_2017}); we rederive them from
Theorem~\ref{thm:main} and explain the value of the non-zero $M_a$ found there
(Proposition~\ref{prop:pointed-M}).

\begin{lemma}\label{lem:pointed-groups}
Put
\[
K_\pm:=\Ker(1\mp P),\qquad R_\pm:=\Ima(1\pm P).
\] Then we have:
\begin{enumerate}[label=(\alph*)]
\item If $P$ is induced by an anti-monoidal involution $\sP$, then
  $\Fix(\sP)=K_+$ and $\Fix(\sT)=K_-$.
\item $R_\pm\subset K_\pm$, and $2x\in R_\pm$ for $x\in K_\pm$, so
  $K_\pm/R_\pm$ are $\FF_2$-vector spaces.
\item With respect to $B$, $K_-^\perp=R_-$ and $K_+^\perp=R_+$.
\item $|K_+|/|R_+|=|K_-|/|R_-|$.
\end{enumerate}
\end{lemma}

\begin{proof}
(a) $\sP a=Pa$ and $(\sP a)^\vee=-Pa$.

\medskip\noindent (b) $(1\mp P)(1\pm P)=1-P^2=0$, and $2x=x\pm Px$ for $x\in K_\pm$.

\medskip\noindent (c) $B(x,y\mp Py)=B(x,y)\,B(Px,y)^{\pm1}=B(x\pm Px,y)$, since
$B(x,Py)=B(Px,y)^{-1}$.
Hence $R_\mp^\perp=K_\mp$, that is,
$K_\mp^\perp=R_\mp$.

\medskip\noindent (d) $|R_\pm|=|A|/|K_\mp|$, as $R_\pm\cong A/\ker(1\pm P)$.
\end{proof}

\begin{proposition}\label{prop:pointed-lift}
For every $P$ as above there is an anti-monoidal involution $(\sP,U,\beta)$ with
coherence on $\cC$ with $\sP a\cong Pa$ for all $a\in A$.
\end{proposition}

The only non-trivial point is \eqref{eq:beta-coherence}, that is, the
vanishing of the obstruction class (Proposition~\ref{prop:obstruction}).
This was shown in \cite{Orii2025H3} by a computation with explicit
cocycles; the proof below is a coordinate-free version of it.

\begin{proof}
(1) Let $\cC^{\otimes\op}$ be $\cC$ with the opposite tensor product
$X\otimes^{\op}Y:=Y\otimes X$ and the braiding
$c^{-1}_{X,Y}\colon X\otimes^{\op}Y\to Y\otimes^{\op}X$.
Then \eqref{eq:U-coherence}--\eqref{eq:U-braiding} say exactly that
$(\sP,U)$ is a braided monoidal functor $\cC\to\cC^{\otimes\op}$.
Both $\cC$ and $\cC^{\otimes\op}$ are pointed braided fusion categories with group of
simple objects $A$, and their quadratic forms, the self-braidings of
the simple objects, are $c_{a,a}=\theta_a$ and $c_{a,a}^{-1}=\theta_a^{-1}$.
By the theorem of Joyal and Street \cite[Theorem~8.4.12]{EGNO}, a
homomorphism of the groups of simple objects which preserves the
quadratic forms is induced by a braided monoidal functor, unique up to
monoidal natural isomorphism.
Applied to $P$, this gives $(\sP,U)$ satisfying
\eqref{eq:U-coherence}--\eqref{eq:U-braiding} with $\sP a\cong Pa$; it
is an equivalence since $P$ is bijective.

(2) The functor $\sP^2$, with the monoidal structure
$\sP(U_{X,Y})\circ U_{\sP Y,\sP X}$, is a braided autoequivalence of
$\cC$ inducing $P^2=1$ on $A$.
By the uniqueness in (1) it is monoidally isomorphic to $\id_\cC$, so
there is $\beta$ satisfying \eqref{eq:beta-monoidal}.

(3) By Proposition~\ref{prop:kappa},
$\sP(\beta_X)^{-1}\beta_{\sP X}=\ell_{a,\sP X}$, where $a\in A$ is the coherence
obstruction of $\beta$, and
by Proposition~\ref{prop:obstruction}(b) it suffices to show $a\in R_+=\Ima(1+P)$: if $a=b+Pb$,
then $\beta\circ\ell_{-b}$ is coherent.
Let $x\in K_+=\Ker(1-P)$, and let $\xi$ be a $\sP$-structure on $x$ (Lemma~\ref{lem:xiT}).
The naturality of $\beta$ and $\beta_x=\sP(\xi)\,\xi$ give
\[
\beta_{\sP x}=\sP^2(\xi)\,\beta_x\,\xi^{-1}=\sP^2(\xi)\,\sP(\xi)=\sP(\beta_x).
\]
Hence $B(a,x)=\ell_{a,\sP x}=1$ for all $x\in K_+$, that is,
$a\in K_+^\perp=R_+$ by Lemma~\ref{lem:pointed-groups}(c).
\end{proof}

From now on we fix an anti-monoidal involution $(\sP,U,\beta)$ with coherence
inducing $P$; it exists by Proposition~\ref{prop:pointed-lift}.

\begin{lemma}\label{lem:pointed-eta}
$\eta|_{K_-}$ is a character of $K_-$.
\end{lemma}

\begin{proof}
Let $b\in K_-$.
Then $b\otimes\sP b\cong b+Pb=0$, so by
Proposition~\ref{prop:basis} $\bdy(b)$ is a multiple of the unit $\mob(0)$, and the
multiple is $\varepsilon(\bdy(b))=\eta_{-b}=\eta_b$ (Definition~\ref{def:eta}, and
Proposition~\ref{prop:beta}(c) as $\beta$ is coherent).
Thus $\bdy(b)=\eta_b\,\mob(0)$, and $\eta|_{K_-}$ is a character because $\bdy$ is
multiplicative (Proposition~\ref{prop:L}).
\end{proof}

Since $B$ is non-degenerate, the character $\eta|_{K_-}$ extends to a
character of $A$, so there is $w\in A$ with $\eta_b=B(w,b)$ for
$b\in K_-$.
Define $n$ by $2^{n}:=|K_-|/|R_-|$.

\begin{proposition}\label{prop:pointed-M}
\[
M_a=\begin{cases}2^{n/2},&a\in-w+R_-,\\0,&\text{otherwise}.\end{cases}
\]
\end{proposition}

\begin{proof}
By Theorem~\ref{thm:main}(c) and Lemma~\ref{lem:pointed-groups}(a),
\[
M_a=|A|^{-1/2}\sum_{b\in K_-}B(a,b)\,\eta_b=|A|^{-1/2}\sum_{b\in K_-}B(a+w,b).
\]
The sum is $|K_-|$ if $a+w\in K_-^\perp=R_-$ and $0$ otherwise, and
$|K_-|\,|R_-|=|A|$ by Lemma~\ref{lem:pointed-groups}(c), so
$|A|^{-1/2}|K_-|=(|K_-|/|R_-|)^{1/2}=2^{n/2}$.
\end{proof}
As $M_a$ should be an integer, $n$ should be even.
We will explain this, together with the reason why $M_a$ is a power of $2$,
in the next proposition.
\begin{proposition}\label{prop:pointed-clifford}
\begin{enumerate}[label=(\alph*)]
\item For $x,y\in K_+$, $\mob(x)\mob(y)$ is a non-zero multiple of $\mob(x+y)$,
  and $\mob(y)\mob(x)=B(x,y)\,\mob(x)\mob(y)$ with $B(x,y)=\pm1$.
\item The center of $\fX$ is spanned by the $\mob(r)$, $r\in R_+$, and $B$
  induces a non-degenerate alternating form on $K_+/R_+$, whose
  dimension is $n$, which is necessarily even.
\item For every $a\in A$ with $p_a\neq0$,
  $p_a\fX\cong\mathrm{Cl}_{n}(\bC)\cong\Mat_{2^{n/2}}(\bC)$.
\end{enumerate}
\end{proposition}

\begin{proof}
(a) $\xi_x\cdot\xi_y$ is an isomorphism $x+y\to\sP(x+y)$, so
$\mob(x)\mob(y)=[x+y,\xi_x\cdot\xi_y]$ is a non-zero multiple of $\mob(x+y)$
(Proposition~\ref{prop:product}).
By Lemma~\ref{lem:Ma-conj} for $E=\one$,
$\mob(x)\cdot \mob(y)\cdot \mob(x)^{-1}=[y,\xi_y\circ\ell_{x,y}]=\chi_y([x])\,\mob(y)=B(x,y)\,\mob(y)$.
Finally $B(x,y)=B(Px,Py)^{-1}=B(x,y)^{-1}$, so $\mob(y)\mob(x)=B(x,y)\,\mob(x)\mob(y)$.

\medskip\noindent (b) By (a), $\fX$ is a twisted group algebra of $K_+$ with commutator
pairing $B|_{K_+}$, so its center is spanned by the $\mob(x)$ with $x$ in
the radical $K_+\cap K_+^\perp=R_+$ (Lemma~\ref{lem:pointed-groups}(b),(c)).
The induced form on $K_+/R_+$ is non-degenerate, and alternating since
$\mob(x)$ commutes with itself.
Its dimension is $n$ by
Lemma~\ref{lem:pointed-groups}(d).

\medskip\noindent (c) Choose $x_1,\dots,x_{n}\in K_+$ whose images form a basis of
$K_+/R_+$ with $B(x_i,x_j)=-1$ for all $i\neq j$; this is possible
because all non-degenerate alternating forms on $\FF_2^{n}$ are
equivalent, and the matrix $J-1$, with $J$ the all-ones matrix, is one
of them, as $(J-1)^2=1$ over $\FF_2$.
By (a) the $\mob(x_i)$ pairwise
anticommute, and $\mob(x_i)^2=\theta_{x_i}^{-1}\bdy(x_i)$ by Lemma~\ref{lem:Ma-square}, which acts
on $p_a\fX$ by the scalar $\theta_{x_i}^{-1}\chi_a([x_i])=\theta_{x_i}^{-1}B(a,x_i)\neq0$.
So after dividing each $\mob(x_i)$ by a square root of this scalar, the $p_a\mob(x_i)$ satisfy the
relations of $\mathrm{Cl}_{n}(\bC)$.
They generate $p_a\fX$, since
every $\mob(x)$ is a multiple of a product of the $\mob(x_i)$ and some $\mob(r)$
with $r\in R_+$.
As $\mathrm{Cl}_{n}(\bC)\cong\Mat_{2^{n/2}}(\bC)$ is
simple, $p_a\fX\cong\mathrm{Cl}_{n}(\bC)$.
\end{proof}

So the irreducible module $V_a$ is the spinor module of
$\mathrm{Cl}_{n}(\bC)$, of dimension $2^{n/2}=M_a$, in agreement with
Proposition~\ref{prop:pointed-M}.

Recall that a \emph{quadratic refinement} of a non-degenerate
alternating form $\langle-,-\rangle$ on $V=\FF_2^{n}$ is a map
$q\colon V\to\{\pm1\}$ with $q(x+y)=q(x)\,q(y)\,\langle x,y\rangle$, and
that its \emph{Arf invariant}
$\operatorname{Arf}(q):=2^{-n/2}\sum_{x\in V}q(x)$ is $\pm1$.

\begin{proposition}\label{prop:pointed-theta}
The twist is constant on $\Sigma=-w+R_-$, with value
\[
\theta_w=\lambda\,\operatorname{Arf}(q)=\pm1,
\qquad
\lambda=|A|^{-1/2}\sum_{a\in A}\theta_a ,
\]
where $q$ is the quadratic refinement $q(b)=\theta_b\,\eta_b$ of the
form induced by $B$ on $K_-/R_-$.
This is Corollary~\ref{cor:theta-coherent}, with
$\sigma=\operatorname{Arf}(q)$, which is the formula of
\cite{Lee:2018eqa}.
\end{proposition}

\begin{proof}
By Proposition~\ref{prop:pointed-M}, $\Sigma=-w+R_-$, and by
Corollary~\ref{cor:theta-coherent}, $\theta$ is constant on $\Sigma$ with value
$\lambda\sigma=\pm1$; as $\theta_{-x}=\theta_x$, this value is $\theta_w$.
For $r\in R_-$, constancy gives $\theta_{-w}=\theta_{-w+r}=\theta_{-w}\,\theta_r\,B(-w,r)$, so
$\theta_r=B(w,r)=\eta_r$.
It remains to show $\sigma=\operatorname{Arf}(q)$.

Let $b\in K_-$.
Then $\theta_b=\theta_{-b}=\theta_{Pb}=\theta_b^{-1}$, so $\theta_b=\pm1$,
and $\eta_b^2=\eta_{2b}=\theta_{2b}=\theta_b^4=1$, as $2b\in R_-$.
For $r\in R_-$ we have $B(b,r)=1$ and $\eta_r=\theta_r$, so
$q(b+r)=q(b)\,\theta_r^2=q(b)$, and $q$ descends to $K_-/R_-$.
There $B$ induces a non-degenerate form
(Lemma~\ref{lem:pointed-groups}(c)), which is alternating since
$B(b,b)=\theta_{2b}\theta_b^{-2}=1$, and
$q(b+b')=q(b)\,q(b')\,B(b,b')$.
As $|K_-/R_-|=2^{n}$ and $|A|=|K_-|\,|R_-|$, we get
\[
\sigma=|A|^{-1/2}\sum_{b\in K_-}\eta_b\,\theta_b
=\frac{|R_-|}{|A|^{1/2}}\sum_{x\in K_-/R_-}q(x)=\operatorname{Arf}(q).
\]
\end{proof}

\subsection{Categories of the form $\cE\boxtimes\overline{\cE}$}
\label{subsec:swap}

Let $\cE$ be a pseudo-unitary modular tensor category, which we take strict,
let $n:=|\Irr(\cE)|$, and let $w,x,y$ run over $\Irr(\cE)$.
Let
$\overline{\cE}$ be $\cE$ with the inverse braiding
$\bar c_{v,v'}:=c^{-1}_{v',v}$; its twist is $\theta^{-1}$, and its
$S$-matrix is $\overline S$, as its monodromies are the inverses of
those of $\cE$ (Section~\ref{subsec:statement}).
Let $\cC:=\cE\boxtimes\overline{\cE}$: its objects are
direct sums of $u\boxtimes v$ with $u,v\in\cE$, with
$\Hom(u\boxtimes v,u'\boxtimes v')=\Hom(u,u')\otimes\Hom(v,v')$,
the componentwise tensor product, and the braiding
\[
c_{u\boxtimes v,\,u'\boxtimes v'}=c_{u,u'}\boxtimes c^{-1}_{v',v}\colon
(u\otimes u')\boxtimes(v\otimes v')\to(u'\otimes u)\boxtimes(v'\otimes v).
\]
It is a strict pseudo-unitary modular tensor category with
$\Irr(\cC)=\{x\boxtimes y\}$, $(x\boxtimes y)^\vee=x^\vee\boxtimes y^\vee$,
$d_{x\boxtimes y}=d_xd_y$, $\theta_{x\boxtimes y}=\theta_x\theta_y^{-1}$,
$S_{x\boxtimes y,\,x'\boxtimes y'}=S_{x,x'}\,\overline{S_{y,y'}}$ and
$D_\cC=D_\cE^2$.

Let $\sP$ exchange the two factors:
\[
\sP(u\boxtimes v):=v\boxtimes u,\qquad \sP(f\boxtimes g):=g\boxtimes f .
\]
Then $\sP^2=\id_\cC$, $\sP\one=\one$, and $\sP$ is strictly monoidal,
$\sP(X\otimes Y)=\sP X\otimes\sP Y$, but it inverts the braiding:
\begin{equation}\label{eq:swap-braiding}
\sP(c_{X,Y})=c^{-1}_{\sP Y,\sP X},
\end{equation}
since $\sP(c_{u,u'}\boxtimes c^{-1}_{v',v})=c^{-1}_{v',v}\boxtimes c_{u,u'}$
while $c_{\sP Y,\sP X}=c_{v'\boxtimes u',\,v\boxtimes u}=c_{v',v}\boxtimes c^{-1}_{u,u'}$.
The braiding of $\cC$ makes $\sP$ anti-monoidal:

\begin{lemma}\label{lem:swap-P}
\[
U_{X,Y}:=c_{\sP Y,\sP X}\colon\sP Y\otimes\sP X\to\sP X\otimes\sP Y=\sP(X\otimes Y),
\qquad U_0:=\id_\one,\qquad\beta:=\id
\]
form an anti-monoidal involution in the sense of
Definition~\ref{def:anti-monoidal}.
\end{lemma}

\begin{proof}
\eqref{eq:U-coherence} holds because both sides are the positive half
twist $\sP Z\otimes\sP Y\otimes\sP X\to\sP X\otimes\sP Y\otimes\sP Z$, by
the hexagon axioms, and \eqref{eq:U-unit} because $c_{\one,-}=c_{-,\one}=\id$.
By \eqref{eq:swap-braiding}, both sides of \eqref{eq:U-braiding} are the
identity of $\sP Y\otimes\sP X$:
$\sP(c_{X,Y})\,U_{X,Y}=c^{-1}_{\sP Y,\sP X}c_{\sP Y,\sP X}$ and
$U_{Y,X}\,c^{-1}_{\sP X,\sP Y}=c_{\sP X,\sP Y}c^{-1}_{\sP X,\sP Y}$.
Likewise \eqref{eq:beta-monoidal} with $\beta=\id$ reads
$\sP(U_{X,Y})\,U_{\sP Y,\sP X}=\sP(c_{\sP Y,\sP X})\,c_{X,Y}=c^{-1}_{X,Y}c_{X,Y}=\id$
and $\sP(U_0)U_0=\id$, and \eqref{eq:beta-coherence} is trivial.
\end{proof}

From now on $u$ and $v$ also run over $\Irr(\cE)$.
Since $(x\boxtimes y)^\vee=x^\vee\boxtimes y^\vee$, we have
$(\sP(x\boxtimes y))^\vee=y^\vee\boxtimes x^\vee$, so
\[
\Fix(\sP)=\{u\boxtimes u\},\qquad \Fix(\sT)=\{u\boxtimes u^\vee\},
\qquad \dim\fX=n .
\]
Since $\beta=\id$, the identity of $u\boxtimes u$ is a $\sP$-structure,
and we take $\mob(u\boxtimes u):=[u\boxtimes u,\id]$.

\begin{proposition}\label{prop:swap-algebra}
\[
\mob(u\boxtimes u)\,\mob(v\boxtimes v)=\bdy(u\boxtimes v)=\sum_wN_{uv}^w\,\mob(w\boxtimes w),
\]
so that $\mob(u\boxtimes u)\mapsto[u]$ is an isomorphism of algebras
$\fX\cong\fF(\cE)$, under which $\bdy(u\boxtimes v)\mapsto[u][v]$,
$\varepsilon$ becomes the coefficient of $[\one]$, and $\star$ becomes
$[u]\mapsto[u^\vee]$.
Moreover $\eta\equiv1$ on $\Fix(\sT)$.
\end{proposition}

\begin{proof}
For $Y=u\boxtimes\one$ we have $\sP Y=\one\boxtimes u$, $Y\otimes\sP Y=u\boxtimes u$, and
$b_Y=U_{Y,\sP Y}=c_{u\boxtimes\one,\,\one\boxtimes u}=c_{u,\one}\boxtimes c^{-1}_{u,\one}=\id$, as $\beta=\id$;
so $\bdy(u\boxtimes\one)=[u\boxtimes u,\id]=\mob(u\boxtimes u)$, and in the same way
$\bdy(\one\boxtimes u)=\mob(u\boxtimes u)$.
By Proposition~\ref{prop:L},
\begin{align*}
\mob(u\boxtimes u)\,\mob(v\boxtimes v)&=\bdy(u\boxtimes\one)\,\bdy(\one\boxtimes v)=\bdy(u\boxtimes v),\\
\qquad
\bdy(u\boxtimes\one)\,\bdy(v\boxtimes\one)&=\bdy\bigl((u\otimes v)\boxtimes\one\bigr)=\sum_wN_{uv}^w\,\mob(w\boxtimes w).
\end{align*}
This gives the displayed formula and the isomorphism, and $\varepsilon$ is the
coefficient of $\mob(\one\boxtimes\one)\mapsto[\one]$.
By Lemma~\ref{lem:L-star},
$\mob(u\boxtimes u)^\star=\bdy(u\boxtimes\one)^\star=\bdy(u^\vee\boxtimes\one)=\mob(u^\vee\boxtimes u^\vee)$.
Finally, by Definition~\ref{def:eta},
$\eta_{u\boxtimes u^\vee}=\varepsilon(\bdy(u^\vee\boxtimes u))$ is the
coefficient of $[\one]$ in $[u^\vee][u]$, which is $N_{u^\vee u}^\one=1$.
\end{proof}

\begin{proposition}\label{prop:swap-M}
$M_{x\boxtimes y}=\delta_{y,x^\vee}$, so $\Sigma=\{x\boxtimes x^\vee\}$
and $\fX\cong\bC^{\oplus n}$.
Under the isomorphism of
Proposition~\ref{prop:swap-algebra}, $V_{x\boxtimes x^\vee}$ is the
one-dimensional module on which $\fF(\cE)$ acts by the character
$\chi_x$.
\end{proposition}

\begin{proof}
Under the isomorphism of
Proposition~\ref{prop:swap-algebra},
$\bdy\colon\fF(\cC)=\fF(\cE)\otimes\fF(\overline{\cE})\to\fF(\cE)$ is the
multiplication; as the $S$-matrix of $\overline{\cE}$ is $\overline S$,
its primitive idempotent labeled by $y$ is $e_{y^\vee}$, so
$p_{x\boxtimes y}=\bdy(e_x\otimes e_{y^\vee})\mapsto e_xe_{y^\vee}=\delta_{y,x^\vee}e_x$.
So $p_{x\boxtimes y}\fX$ is $\bC e_x$ for $y=x^\vee$ and $0$ otherwise, and
$M_{x\boxtimes y}=\dim V_{x\boxtimes y}=\delta_{y,x^\vee}$ by Theorem~\ref{thm:main}(b), (c);
on $V_{x\boxtimes x^\vee}=\bC e_x$, $\fF(\cE)$ acts by $\chi_x$.
\end{proof}

\begin{proposition}\label{prop:swap-theta}
$\theta_{x\boxtimes x^\vee}=1$ for every $x\boxtimes x^\vee\in\Sigma$,
and $\lambda=\sigma=1$, in accordance with Corollary~\ref{cor:theta-coherent}.
\end{proposition}

\begin{proof}
$\theta_{x\boxtimes x^\vee}=\theta_x\theta_{x^\vee}^{-1}=1$, and
\begin{gather*}
p_+=\sum_{x,y}d_x^2d_y^2\,\theta_x\theta_y^{-1}=\Bigl|\sum_xd_x^2\theta_x\Bigr|^2=D_\cE^2=D_\cC,\\
\sigma=\frac1{D_\cC}\sum_u\eta_{u\boxtimes u^\vee}\,d_u^2\,\theta_{u\boxtimes u^\vee}
=\frac1{D_\cE^2}\sum_ud_u^2=1 .\qedhere
\end{gather*}
\end{proof}

\begin{remark}
These results reflect the fact that the orientation double cover of the
M\"obius band is an annulus.
A diagram on the band colored by
$\cE\boxtimes\overline{\cE}$, glued at the cut through $\sP$, is the same as
a diagram on the annulus colored by $\cE$: the core loop $\mob(u\boxtimes u)$
lifts to a loop colored by $u$ around the annulus, and the boundary loop
$\bdy(u\boxtimes v)$ to the two boundary circles colored by $u$ and $v$; the
algebra of such loops on the annulus is $\fF(\cE)$.
Likewise, the
dimension $\delta_{y,x^\vee}=M_{x\boxtimes y}$ is that of
$\Hom_\cE(\one,x\otimes y)$, the space which the Reshetikhin--Turaev
theory of $\cE$ assigns to the sphere with two marked points colored by
$x$ and $y$, the double cover of $\mathbb{RP}^2$ with one marked point
colored by $x\boxtimes y$.
\end{remark}

\subsection{Certain non-pointed categories of rank six}
\label{subsec:metaplectic}

Consider the pseudo-unitary modular categories with six simple objects
$\one,\epsilon,\phi_1,\phi_2,\psi_+,\psi_-$ of dimensions
$1,1,2,2,\sqrt5,\sqrt5$ and the fusion rules
\begin{align*}
\epsilon\otimes\epsilon&=\one,&
\epsilon\otimes\phi_i&=\phi_i,&
\epsilon\otimes\psi_\pm&=\psi_\mp,\\
\phi_i\otimes\phi_i&=\one\oplus\epsilon\oplus\phi_{3-i},&
\phi_1\otimes\phi_2&=\phi_1\oplus\phi_2,&
\phi_i\otimes\psi_\pm&=\psi_+\oplus\psi_-,\\
\psi_+\otimes\psi_+&=\one\oplus\phi_1\oplus\phi_2,&
\psi_-\otimes\psi_-&=\one\oplus\phi_1\oplus\phi_2,&
\psi_+\otimes\psi_-&=\epsilon\oplus\phi_1\oplus\phi_2,
\end{align*}
those of $\mathfrak{so}_5$ at level $2$.
Four fusion categories with these fusion rules, each with a braiding,
were constructed in \cite[\S6.2]{AFT2016}, and by the census of
multiplicity-free fusion categories of low rank in
\cite{VercleyenSlingerland2026} there are no others.
We need two of the resulting
modular categories, $\cC_1$ and $\cC_2$.
(The other two are their Galois conjugates under $\sqrt{5}\mapsto-\sqrt5$,
which are still pseudo-unitary.)
Both have
$\theta_\one=\theta_\epsilon=1$ and
$\theta_{\phi_1}=\overline{\theta_{\phi_2}}=e^{4\pi i/5}$; they differ in
\[
\cC_1:\ \theta_{\psi_\pm}=\pm i,\qquad\qquad \cC_2:\ \theta_{\psi_\pm}=\pm1 .
\]
$\cC_1$ is the category of integrable representations of
$\widehat{\mathfrak{so}}_5$ at level $2$, and $\cC_2$ is the category
denoted $\mathrm{USp}(4)_2^\vee$ in \cite{Barkeshli:2017rzd}.
For both,
$D=\sqrt{20}$ and $\lambda=p_+/D=-1$.
The entries of $S$ we
need are $\sqrt{20}\,S_{X,\one}=d_X$; $\sqrt{20}\,S_{X,\epsilon}=d_X$ for
$X\in\{\one,\epsilon,\phi_i\}$ and $\sqrt{20}\,S_{\psi_\pm,\epsilon}=-\sqrt5$;
$S_{\phi_i,\psi_\pm}=0$; and $\sqrt{20}\,S_{\psi_s,\psi_t}=st\sqrt5$ for $\cC_1$ but
$-st\sqrt5$ for $\cC_2$ ($s,t=\pm$), as the balancing axiom requires.

\subsubsection{Anti-monoidal involutions with coherence}

Let us study whether these categories admit anti-monoidal involutions with
coherence.
We first show, using only the modular data, that $\cC_1$ admits no
anti-monoidal involution with coherence.

\begin{proposition}\label{prop:meta-C1}
$\cC_1$ admits no anti-monoidal involution with coherence in the sense of
Definition~\ref{def:anti-monoidal} and Definition~\ref{def:coherent}.
\end{proposition}

\begin{proof}
Suppose $(\sP,U,\beta)$ is one.
The only permutation of the simple
objects of $\cC_1$ preserving the fusion rules with
$\theta_{\sP X}=\overline{\theta_X}$ (Lemma~\ref{lem:P-modular-data}) is
$\phi_1\leftrightarrow\phi_2$, $\psi_+\leftrightarrow\psi_-$, so
$\Fix(\sP)=\Fix(\sT)=\{\one,\epsilon\}$.
Then
$M_{\phi_i}=2(1+\eta_{\epsilon})/\sqrt{20}$ is an integer only if
$\eta_{\epsilon}=-1$, and then $M_\one=M_\epsilon=M_{\phi_i}=0$ and
$M_{\psi_\pm}=1$.
So $\Sigma=\{\psi_+,\psi_-\}$, on which
$\theta=\pm i$ is neither constant nor $\pm1$, contradicting
Corollary~\ref{cor:theta-coherent}.
\end{proof}

\begin{remark}\label{rem:C1-BC}
This was first found by Barkeshli and Cheng \cite{Barkeshli:2017rzd}.
\end{remark}

For $\cC_2$ the modular data determine everything except one sign.

\begin{proposition}\label{prop:meta-C2-data}
Let $(\sP,U,\beta)$ be an anti-monoidal involution with coherence on $\cC_2$.
Then
$\Fix(\sP)=\Fix(\sT)=\{\one,\epsilon,\psi_+,\psi_-\}$, $\dim\fX=4$, and
there is a sign $\eta_0$ with
\[
\eta_{\epsilon}=-1,\qquad \eta_{\psi_\pm}=\pm\eta_0,\qquad
\Sigma=\{\psi_{-\eta_0}\},\qquad M_{\psi_{-\eta_0}}=2,\qquad
\fX\cong\Mat_2(\bC).
\]
The twist on $\Sigma$ is $-\eta_0=\lambda\sigma$ with $\sigma=\eta_0$,
as required by Corollary~\ref{cor:theta-coherent}.
Replacing $\beta$ by
$\beta\circ\ell_\epsilon$ gives another anti-monoidal involution, with
the opposite sign of $\eta_0$.
\end{proposition}

\begin{proof}
The permutation is now $\phi_1\leftrightarrow\phi_2$ with $\psi_\pm$
fixed, which gives $\Fix(\sP)$ and $\Fix(\sT)$, and $\dim\fX=4$ by
Proposition~\ref{prop:basis}.
Integrality of $M_{\phi_i}$ again gives
$\eta_{\epsilon}=-1$, as $S_{\phi_i,\psi_\pm}=0$.
Then
$M_\one=-M_\epsilon=\frac12(\eta_{\psi_+}+\eta_{\psi_-})$, and both are
non-negative only if $\eta_{\psi_-}=-\eta_{\psi_+}=:-\eta_0$; then
$M_{\psi_\pm}=1\mp\eta_0$.
By Theorem~\ref{thm:main}(c) we then have
$\fX\cong\Mat_2(\bC)$.
The twist of $\psi_{-\eta_0}$ is $-\eta_0$, and
\[
\sigma=\frac1{\sqrt{20}}\bigl(1-1+\sqrt5\,\eta_0+\sqrt5\,(-\eta_0)(-1)\bigr)=\eta_0 .
\]
Finally, $\beta\circ\ell_\epsilon$ is again monoidal, and satisfies
\eqref{eq:beta-coherence} since $\epsilon\otimes\sP\epsilon\cong\one$
(Definition~\ref{def:coherent} and Proposition~\ref{prop:obstruction}(b)); it multiplies $\eta_Y$ by
$\chi_Y([\epsilon])=S_{\epsilon,Y}/S_{\one,Y}$, which is $-1$ for
$Y=\psi_\pm$.
\end{proof}

It remains to show that $\cC_2$ does admit an anti-monoidal involution with coherence.

\begin{proposition}\label{prop:meta-C2-exists}
$\cC_2$ admits an anti-monoidal involution with coherence, and both signs $\eta_0$ of
Proposition~\ref{prop:meta-C2-data} occur.
\end{proposition}

\begin{proof}
See Appendix~\ref{app:so5}.
\end{proof}

The algebra $\fX$ can then be written down explicitly by
Corollary~\ref{cor:skeletal-product}.

\begin{proposition}\label{prop:meta-C2-algebra}
For $\eta_0=1$, $\Fix(\sP)\cong\bZ/2\times\bZ/2$ under the fusion rules,
and
\begin{gather*}
\mob(\epsilon)^2=-\mob(\one),\qquad \mob(\psi_\pm)^2=\mob(\one),\\
\mob(\epsilon)\mob(\psi_+)=\mob(\psi_-)=-\mob(\psi_+)\mob(\epsilon),\qquad
\mob(\psi_+)\mob(\psi_-)=-\mob(\epsilon)=-\mob(\psi_-)\mob(\psi_+).
\end{gather*}
Thus $\fX$ is the Clifford algebra $\mathrm{Cl}_2(\bC)\cong\Mat_2(\bC)$
generated by the anticommuting $\mob(\psi_\pm)$, with $\bdy(\epsilon)=-\mob(\one)$,
$\bdy(\phi_i)=0$, $\bdy(\psi_\pm)=\pm \mob(\one)$ and $p_{\psi_-}=\mob(\one)$.
\end{proposition}

\begin{proof}
In $\cC_2$ the only $\sP$-invariant summand of $X\otimes Y$, for
$X,Y\in\Fix(\sP)$, is the product in $\bZ/2\times\bZ/2$, so each
product in Corollary~\ref{cor:skeletal-product} has a single term, whose
coefficient $\xi_X\xi_YR^{XY}_Z/(u^{XY}_Z\xi_Z)$, in the notation of
Appendix~\ref{app:skeletal}, we evaluated with the data of
Appendix~\ref{app:so5}.
The values of $\bdy$ follow from
Proposition~\ref{prop:beta}, since each $\bdy(Y)$ is central, hence a
scalar, and $\varepsilon(\bdy(Y))=\eta_{Y^\vee}=\eta_Y$ as $Y^\vee\cong Y$; and $p_{\psi_-}=\mob(\one)$
since $\Sigma=\{\psi_-\}$.
The relation $\mob(\epsilon)^2=\theta_\epsilon^{-1}\bdy(\epsilon)=-\mob(\one)$ also follows from
Lemma~\ref{lem:Ma-square}, without the data.
\end{proof}

Thus $\fX$ has a single irreducible module $V_{\psi_-}$, of dimension
two, although $\psi_-$ is not invertible.

\subsubsection{Anti-monoidal involutions without coherence}

Although $\cC_1$ admits no anti-monoidal involution with coherence, it does
admit one without coherence, and Theorems~\ref{thm:main} and~\ref{thm:theta}
apply to it.

\begin{proposition}\label{prop:meta-C1-noncoh}
$\cC_1$ admits an anti-monoidal involution $(\sP,U,\beta)$, unique up to
isomorphism and up to replacing $\beta$ by $\beta\circ\ell_\epsilon$.
None of them is coherent: $\obstruction=\epsilon$.
For each of them:
\begin{enumerate}[label=(\alph*)]
\item $\Fix(\sP)=\Fix(\sT)=\{\one,\epsilon\}$, $\eta_{\epsilon}=-1$,
  $M_\one=M_\epsilon=M_{\phi_i}=0$ and $M_{\psi_\pm}=1$, so
  $\Sigma=\{\psi_+,\psi_-\}$.
\item $\fX$ has the basis $\mob(\one)$, $\mob(\epsilon)$, with $\mob(\epsilon)^2=-\mob(\one)$, so
  $\fX\cong\bC\oplus\bC$.
  There is a sign $s$ with
  \[
  \bdy(\epsilon)=-\mob(\one),\qquad \bdy(\phi_i)=0,\qquad \bdy(\psi_\pm)=\pm is\,\mob(\epsilon),\qquad
  p_{\psi_\pm}=\tfrac12\bigl(\mob(\one)\pm is\,\mob(\epsilon)\bigr).
  \]
  The sign $s$ changes when $\xi_\epsilon$ is replaced by $-\xi_\epsilon$, and also when
  $\beta$ is replaced by $\beta\circ\ell_\epsilon$, which has the same
  $\sP$-structures on $\epsilon$.
\item $\mob(\obstruction)=\mob(\epsilon)$ acts on $V_{\psi_\pm}$ by $\obs_{\psi_\pm}=\mp is$, and
  $\sigma=s$.
  So $\theta_{\psi_\pm}=\pm i=\lambda\sigma\,\obs_{\psi_\pm}$, as required by
  Theorem~\ref{thm:theta}, although $\theta$ is not constant on $\Sigma$.
\end{enumerate}
\end{proposition}

\begin{proof}
The existence and the uniqueness are the computation of
Appendix~\ref{app:so5}.
As $\obstruction\in\Inv(\cC_1)=\{\one,\epsilon\}$ and $\cC_1$ has no coherent
anti-monoidal involution (Proposition~\ref{prop:meta-C1}), $\obstruction=\epsilon$.

\medskip\noindent (a)
$\Fix(\sP)$ is read off from the permutation of
Proposition~\ref{prop:meta-C1}, and $\Fix(\sT)=\Fix(\sP)$ as all simple objects
are self-dual.
By Proposition~\ref{prop:beta}(c), $\eta_{\epsilon}^{-2}=\chi_\epsilon([\epsilon])=1$, and
$\eta_{\epsilon^\vee}=\eta_{\epsilon}$.
The rest is the computation in the proof of Proposition~\ref{prop:meta-C1}, which
uses only Theorem~\ref{thm:main}(c) and the modular data.

\medskip\noindent (b)
$\dim\fX=|\Fix(\sP)|=2$ by Proposition~\ref{prop:basis}.
Each $\bdy(Y)$ is central, and acts on $V_{\psi_t}$ by
$\chi_{\psi_t}([Y])=S_{Y,\psi_t}/S_{\one,\psi_t}$, which is $-1$ for $Y=\epsilon$, $0$ for
$Y=\phi_i$ and $\pm t$ for $Y=\psi_\pm$.
So $\bdy(\epsilon)=-\mob(\one)$, $\bdy(\phi_i)=0$, $\bdy(\psi_\pm)=\pm(p_{\psi_+}-p_{\psi_-})$, and
$\bdy(\psi_+)^2=\mob(\one)$; and $\mob(\epsilon)^2=\theta_\epsilon^{-1}\bdy(\epsilon)=-\mob(\one)$ by
Lemma~\ref{lem:Ma-square}.
As $\psi_+\otimes\sP\psi_+=\epsilon\oplus\phi_1\oplus\phi_2$ contains $\epsilon$ but not
$\one$, $\bdy(\psi_+)$ is a multiple of $\mob(\epsilon)$ by Proposition~\ref{prop:basis}; as
$\mob(\epsilon)^2=-\mob(\one)$, the multiple is $is$ for a sign $s$.
Replacing $\xi_\epsilon$ by $-\xi_\epsilon$ changes the sign of $\mob(\epsilon)$, hence of $s$.
For $\beta'=\beta\circ\ell_\epsilon$ we have
$\beta'_\epsilon=\chi_\epsilon([\epsilon])\,\beta_\epsilon=\beta_\epsilon$, so $\xi_\epsilon$ and
$\mob(\epsilon)$ can be kept, while
$\beta'_{\psi_+}=\chi_{\psi_+}([\epsilon])\,\beta_{\psi_+}=-\beta_{\psi_+}$ changes the sign of
$b_{\psi_+}$ and of $\bdy(\psi_+)$ (Definition~\ref{def:L}); so $s$ changes sign.
Finally $p_{\psi_\pm}=\frac12(\mob(\one)\pm \bdy(\psi_+))$.

\medskip\noindent (c)
By (b), $\mob(\epsilon)=-is\,\bdy(\psi_+)$, which acts on $V_{\psi_\pm}$ by $\mp is$.
Then $\eta^\obstruction_Y=\varepsilon(\bdy(Y)\mob(\epsilon))$ (as $Y^\vee\cong Y$) is $\mp is$ for $Y=\psi_\pm$ and $0$
otherwise, and
\[
\sigma=\frac1{\sqrt{20}}\Bigl(\overline{(-is)}\,\sqrt5\,(i)^{-1}+\overline{(is)}\,\sqrt5\,(-i)^{-1}\Bigr)
=\frac12(s+s)=s .
\]
With $\lambda=-1$, $\lambda\sigma\,\obs_{\psi_\pm}=(-s)(\mp is)=\pm i=\theta_{\psi_\pm}$.
\end{proof}

\appendix

\section{Finding an anti-monoidal involution in a skeletal model}
\label{app:skeletal}
\addtocontents{toc}{\protect\setcounter{tocdepth}{1}}

We describe how to decide, from the $F$- and $R$-symbols of a skeletal
model of $\cC$, whether a given permutation of $\Irr(\cC)$ is induced by an
anti-monoidal involution, with or without coherence, how to find all of them,
and how to read off $\obstruction$ and $\eta$.
We do this in two steps, first $(\sP,U)$ and then $\beta$.
For simplicity we assume that the fusion rules are multiplicity-free;
some of the arguments below do not generalize to the case with multiplicities.

\subsection{Setting}
\label{app:setting}

As in Corollary~\ref{cor:skeletal-product}, we drop the strictness
assumption and let $\cC$ be skeletal, now with $\dim\Hom(c,a\otimes b)\le1$
for simple $a,b,c$.
Choose a non-zero $v^{ab}_c\in\Hom(c,a\otimes b)$ whenever this space is
non-zero, with $v^{\one a}_a=v^{a\one}_a=\id_a$.
The data of $\cC$ are then the $F$- and $R$-symbols, defined by
\begin{equation}
(v^{ab}_e\otimes\id_c)\,v^{ec}_d=\sum_fF^{abc}_{d;ef}\,(\id_a\otimes v^{bc}_f)\,v^{af}_d,
\qquad
c_{a,b}\,v^{ab}_c=R^{ab}_c\,v^{ba}_c .
\label{eq:FR}
\end{equation}

Let $P$ be a permutation of $\Irr(\cC)$ with $P^2=1$ which preserves the
fusion rules and satisfies $\theta_{Pa}=\theta_a^{-1}$.
These conditions are necessary: an anti-monoidal involution permutes
$\Irr(\cC)$ preserving the fusion rules, $\beta$ gives $\sP^2X\cong X$, and
$\theta_{\sP X}=\theta_X^{-1}$ by Lemma~\ref{lem:P-modular-data}.
They imply $S_{Pa,Pb}=\overline{S_{a,b}}$ by the balancing axiom.

\subsection{Finding the functor $(\sP,U)$}
\label{app:step1}

In the skeleton, $\sP$ acts on simple objects by $P$.
The linear map $u^{ab}_c$ of Corollary~\ref{cor:skeletal-product} is a
scalar, which we denote by the same letter:
\[
U_{a,b}^{-1}\sP(v^{ab}_c)=\Ubox{v^{ab}_c}=u^{ab}_c\,v^{Pb,Pa}_{Pc}.
\]
Taking the boxes of the two relations \eqref{eq:FR} by
\ref{item:B1}--\ref{item:B4}, the conditions
\eqref{eq:U-coherence}--\eqref{eq:U-braiding} become
\begin{align*}
\text{\eqref{eq:U-coherence}}\quad&
F^{abc}_{d;ef}\,u^{bc}_f\,u^{af}_d
=\bigl[(F^{Pc\,Pb\,Pa}_{Pd})^{-1}\bigr]_{Pe,Pf}\,u^{ab}_e\,u^{ec}_d ,\\
\text{\eqref{eq:U-unit}}\quad& u^{\one a}_a=u^{a\one}_a=1,\\
\text{\eqref{eq:U-braiding}}\quad& u^{ab}_c=R^{ab}_c\,R^{Pa\,Pb}_{Pc}\,u^{ba}_c .
\end{align*}
Here $u^{\one a}_a=u^{a\one}_a=U_0$ by \eqref{eq:U-unit}, which we
normalize to $1$ by the natural isomorphism at $\one$.
Replacing $(\sP,U)$ by an isomorphic anti-monoidal functor, via scalars
$\varphi_a\in\bC^\times$ with $\varphi_\one=1$, changes $u^{ab}_c$ to
$u^{ab}_c\,\varphi_a\varphi_b/\varphi_c$.
We call this the \emph{gauge}.

Write $u^{ab}_c=|u^{ab}_c|\,e^{2\pi ix^{ab}_c}$.
The system becomes
\[
A\,x\equiv b\pmod{\bZ},\qquad A\,\log|u|=\log|\rho|,
\]
with an integer matrix $A$, with entries in $\{0,\pm1,\pm2\}$.
Here $\rho$ is the ratio of the two coefficients of an equation, and
$b=\arg\rho/2\pi$.
The second system is ordinary linear algebra over $\bR$; in a unitary gauge
it is solved by $|u|=1$.
For the first, choose integer matrices $V$ and $W$, invertible over $\bZ$,
such that $WAV$ is diagonal, with non-zero entries $d_1,\dots,d_r$ (Smith
normal form).
With $y=V^{-1}x$ and $c=Wb$, the system reads
\[
d_i\,y_i\equiv c_i\quad(i\le r),\qquad 0\equiv c_i\quad(i>r)\pmod{\bZ}.
\]

\begin{proposition}\label{prop:app-step1}
\begin{enumerate}[label=(\alph*)]
\item There is $(\sP,U)$ inducing $P$ if and only if $c_i\in\bZ$ for all $i>r$.
\item If so, the solutions $x$ form a torsor over
$(\bR/\bZ)^{n-r}\times\prod_{i\le r}\bZ/|d_i|$, where $n$ is the number of
unknowns.
\item The solutions up to isomorphism are in bijection with
$\prod_{i\le r}\bZ/|d_i|$.
\end{enumerate}
\end{proposition}

\begin{proof}
(a) and (b) are read off from the diagonal system.

(c) Isomorphisms between the $(\sP,U)$ inducing $P$ are exactly the gauge
transformations, since a natural isomorphism of $\sP$ is a family of scalars
$\varphi_a\in\End(Pa)=\bC$.
On the phases $x$ the gauge acts through the compact torus
$\{(\varphi_a)\mid\varphi_\one=1,\ |\varphi_a|=1\}$, which is connected and of
dimension $|\Irr(\cC)|-1$, and its stabilizer at any solution is finite: it
consists of the $\varphi$ with $\varphi_a\varphi_b=\varphi_c$ whenever
$c\subset a\otimes b$, that is, of the characters of the universal grading
group of $\cC$, which is finite.
So every gauge orbit has dimension $|\Irr(\cC)|-1$.
On the other hand, the gauge orbits lie in the solution set, whose components
have dimension $n-r$ by (b).
An anti-monoidal functor is a tensor functor $\cC\to\cC^{\otimes\op}$,
and tensor functors between fusion categories admit no non-trivial
deformations, by Ocneanu rigidity \cite[Theorem~2.31]{ENO2005}.
The condition \eqref{eq:U-braiding} only adds equations, so every
deformation of a solution is a gauge transformation, and
$|\Irr(\cC)|-1=n-r$.
Hence each gauge orbit, being connected and of full dimension in a compact
component, is the whole component; and the connected gauge group does not
permute the components.
So the isomorphism classes correspond to the components, that is, to
$\prod_{i\le r}\bZ/|d_i|$ by (b).
\end{proof}

The conditions in (a) have a direct meaning.
A row $w$ of $W$ with $i>r$ satisfies $wA=0$.
So it is an integer combination of the equations in which all $u$'s cancel,
that is, a gauge-invariant monomial in the $F$- and $R$-symbols.
The condition says that this invariant takes the same value on the data of
$\cC$ as on the data of the reversed category, relabeled by $P$.
The simplest such invariants are $R^{ab}_cR^{ba}_c=\theta_c/\theta_a\theta_b$,
and for these the condition is $\theta_{Pa}=\theta_a^{-1}$.
The remaining invariants involve $F$.
For pointed $\cC$ there are no further conditions, by
Proposition~\ref{prop:pointed-lift}.
We do not know whether the conditions on $P$ suffice in general.

\subsection{Finding $\beta$}
\label{app:step2}

Fix a solution $u$ of Step~1, and put
$\gamma^{ab}_c:=u^{ab}_c\,u^{Pb\,Pa}_{Pc}$.
Since $P^2=1$, each $\beta_a\in\End(a)=\bC$ is a scalar.
By \ref{item:B7}, $\Ubox{\Ubox{v^{ab}_c}}=\gamma^{ab}_c\,v^{ab}_c$ must equal
$(\beta_a\otimes\beta_b)\,v^{ab}_c\,\beta_c^{-1}$.
So \eqref{eq:beta-monoidal}, and \eqref{eq:beta-coherence} if we require coherence,
become
\[
\beta_a\beta_b/\beta_c=\gamma^{ab}_c,\quad \beta_\one=1,
\qquad\text{and}\qquad
\beta_{Pa}=\beta_a .
\]
We solve them
by the same method as above.
The solutions of \eqref{eq:beta-monoidal}, if any, form a torsor over the
group $\Inv(\cC)$ of Definition~\ref{def:I}.
Adding \eqref{eq:beta-coherence} as further rows, the condition of
Proposition~\ref{prop:app-step1}(a) for a solution with coherence is the
vanishing of the obstruction class in
$H^3(\bZ/2,\Inv(\cC))$ of Proposition~\ref{prop:obstruction}.
This class depends on $(\sP,U)$, not only on $P$.
So Step~2 has to be run for each of the finitely many $(\sP,U)$ of
Proposition~\ref{prop:app-step1}(c).

Without coherence, every solution $\beta$ of \eqref{eq:beta-monoidal} gives an
anti-monoidal involution, and $\obstruction$ is read off from it.
The box of the scalar $\beta_a$ is the same scalar, so \eqref{eq:obstruction} becomes
$\beta_{Pa}=\beta_a\,\chi_{Pa}([\obstruction])$, and $\obstruction$ is the element of $\Inv(\cC)$ with
$\chi_{a}([\obstruction])=\beta_a/\beta_{Pa}$ for all $a$.

Finally, for $Y\in\Fix(\sT)$ we have $PY=Y^\vee$, and $\Hom(\one,Y\otimes PY)$ is
spanned by $v^{Y,PY}_\one$.
With $U_0=1$, the map $R_{U_0,b_Y}$ of Lemma~\ref{lem:fY} is
$\Phi\mapsto(\beta_Y^{-1}\otimes\id)\,\Ubox{\Phi}$, which multiplies $v^{Y,PY}_\one$ by
$u^{Y,PY}_\one/\beta_Y$.
So by Proposition~\ref{prop:beta},
\[
\eta_Y=\frac{u^{Y,PY}_\one}{\beta_Y}\quad(Y\in\Fix(\sT)),\qquad
\eta_Y=0\quad(Y\notin\Fix(\sT)).
\]
The structure constants of $\fX$ are given by Corollary~\ref{cor:skeletal-product}.

\subsection{The categories $\cC_1$ and $\cC_2$}
\label{app:so5}

We applied this to the categories $\cC_1$ and $\cC_2$ of
Section~\ref{subsec:metaplectic}.
Their $F$- and $R$-symbols were determined in \cite{AFT2016}.
We used them in the numerical form of the AnyonWiki \cite{AnyonWiki},
whose data for this fusion ring come from \cite{VercleyenSlingerland2026}:
$\cC_2$ and $\cC_1$ are the categories with the codes $[6,1,0,9,1,1,1]$
and $[6,1,0,9,3,1,1]$ there, i.e.\ the first and the third categorification of the
fusion ring $\mathrm{FR}^{6,1,0}_9$, with the first pivotal (spherical) structure.
We checked the pentagon and hexagon equations independently.
In both cases $P$ is the permutation of Proposition~\ref{prop:meta-C1},
respectively Proposition~\ref{prop:meta-C2-data}.

Step~1 has $58$ unknowns and $946$ equations, of rank $53$.
The $5$ free directions are the gauge orbits, as they must be by
Proposition~\ref{prop:app-step1}(c), and there is no torsion.
So in both cases $(\sP,U)$ exists and is unique up to isomorphism.
In the consistent rows, the $c_i$ are integers to within $10^{-16}$, the
rounding error of the data.
In Step~2 we find two solutions $\beta$ of \eqref{eq:beta-monoidal} in both
cases, as $\Inv(\cC)=\{\one,\epsilon\}$; they differ by $\ell_\epsilon$.

For $\cC_2$ both satisfy \eqref{eq:beta-coherence}, so both are anti-monoidal
involutions with coherence.
This proves Proposition~\ref{prop:meta-C2-exists}.

For $\cC_1$ neither satisfies \eqref{eq:beta-coherence}: the added rows
deviate from $\bZ$ by exactly $1/2$, the non-zero class in
$H^3(\bZ/2,\Inv(\cC))\cong\bZ/2$.
As $(\sP,U)$ is unique, $\cC_1$ has no anti-monoidal involution with coherence, in
agreement with Proposition~\ref{prop:meta-C1}.
Both $\beta$ are anti-monoidal involutions without coherence.
For both, $\beta_a/\beta_{Pa}$ is $-1$ for $a=\psi_\pm$ and $1$ otherwise, which is
$\chi_a([\epsilon])$; so $\obstruction=\epsilon$.
Corollary~\ref{cor:skeletal-product} gives $\mob(\epsilon)^2=-\mob(\one)$, and
$\bdy(\psi_+)=is\,\mob(\epsilon)$ with opposite signs $s$ for the two $\beta$,
in agreement with Proposition~\ref{prop:meta-C1-noncoh}.

The values $\eta_Y=u^{Y,PY}_\one/\beta_Y$ on $\Fix(\sT)$ are as follows; $\eta$ vanishes
on the other simple objects.
\begin{center}
\begin{tabular}{l|c|cccc}
 & $\obstruction$ & $\eta_{\one}$ & $\eta_{\epsilon}$ & $\eta_{\psi_+}$ & $\eta_{\psi_-}$\\
\hline
$\cC_2$, first $\beta$ & $\one$ & $1$ & $-1$ & $1$ & $-1$\\
$\cC_2$, second $\beta$ & $\one$ & $1$ & $-1$ & $-1$ & $1$\\
$\cC_1$, both $\beta$ & $\epsilon$ & $1$ & $-1$ & -- & --
\end{tabular}
\end{center}
For $\cC_1$, $\psi_\pm\notin\Fix(\sT)$.
The rows for $\cC_2$ are the two cases of Proposition~\ref{prop:meta-C2-data}.

The other two categorifications behave like $\cC_2$ or $\cC_1$, according
as $\theta_{\psi_\pm}=\pm1$ or $\pm i$.
The final $u^{ab}_c$ and $\beta_a$ can be checked directly against the
equations above.

\def\arxivfont{\rm}
\bibliographystyle{ytamsalpha}
\bibliography{ref}

\end{document}